%% file: main.tex
\documentclass[11pt]{article}
\usepackage[numbers]{natbib}
\usepackage{amsmath,amssymb,amsthm}
\usepackage{array,booktabs,longtable}
\usepackage{needspace}
\usepackage{xcolor}
\usepackage{hyperref}
\usepackage[margin=1in]{geometry}
\hypersetup{colorlinks=true,citecolor=blue,linkcolor=blue}
\allowdisplaybreaks
\makeatletter
\renewcommand{\l@subsection}{\@dottedtocline{2}{1.5em}{3em}}
\makeatother

\input{_0_vars}
\input{_1_maths}
\theoremstyle{plain}
\newtheorem{theorem}{Theorem}[section]
\newtheorem{lemma}[theorem]{Lemma}
\newtheorem{proposition}[theorem]{Proposition}
\begin{document}
\date{October 2, 2026}
\title{\paperTitle}
\author{\paperAuthor}
\maketitle
\begin{abstract}
\input{00_abstract}
\end{abstract}

\input{_2_body}

\section*{Acknowlegements and AI Disclosure}
The authors thank Josh Alman, Haotian Jiang, Yin Tat Lee, Victor Reis, Santosh Vempala, Omri Weinstein, and Lichen Zhang for helpful discussions. The authors also thank Eleanor M. Song and Yaping Zhang for their general support throughout this project. The second author dedicates this work to her Ph.D. advisors, Anna Karlin and Yin Tat Lee.

The AI tools used in this work were GPT-6 Astra, GPT-5.6 Sol, Claude Fable 5, and Fable 5.1. The first author learned of this problem from Haotian Jiang while visiting the University of Washington nearly a decade ago. Subsequent attempts to solve the problem, first as a postdoc at the Institute for Advanced Study around 2019–2020 and again from roughly 2022 to 2024, were unsuccessful. After seeing Letwin’s result on arXiv on July 27, 2026, and carefully reading its proof, the authors noticed that it might be possible to improve the bound from $\log^{1/4}$ to $\log^{1/6}$. Beginning on July 28, 2026, this effort involved exploring more than 100 approaches in collaboration with AI tools. Candidate approaches were generated collaboratively, with the authors deciding which ones to prioritize. The authors and the AI tools ultimately obtained the $\log^{1/6}$ bound. 

After fully understanding the $\log^{1/6}$ proof, the authors found, on September 30, 2026, a way to reach an $f(\log^* n)$ bound for some fixed function $f$.  The authors thank Lijie Chen for the inspiring discussions in 2018 on the recursive use of the Chinese Remainder Theorem, a technique from Chen’s work in fine-grained complexity. Shortly after releasing the $f(\log^* (n))$ version, the authors realized that the log-star refinement could be applied repeatedly to further reduce the remaining height.  This height-reduction idea was inspired by a series of works by the first author on solving linear programs, which repeatedly reduce the matrix size in order to maintain the matrix inverse dynamically. Finally, the authors completed the proof by sharpening the analysis so that the costs decay geometrically and their total is bounded by a universal constant.

All of the proofs have been carefully verified by the authors and several AI tools. The authors take full responsibility for the correctness of the proofs.

\bibliographystyle{alpha}
\bibliography{ref}

\end{document}

%% file: _0_vars.tex
\newcommand{\paperTitle}{An $O(1)$ Bound for the KLS Constant}
\newcommand{\paperAuthor}{Zhao Song\\
{\normalsize Independent Researcher}\thanks{\texttt{magic.linuxkde@gmail.com}.}
\and Xinzhi Zhang\\
{\normalsize Microsoft Research}\thanks{\texttt{xinzhizhang@microsoft.com}.}
}

%% file: _1_maths.tex
\DeclareMathOperator{\Var}{Var}
\DeclareMathOperator{\E}{\mathbb{E}}
\DeclareMathOperator{\tr}{tr}
\renewcommand{\d}{\mathrm{d}}

%% file: 00_abstract.tex
The Kannan--Lov\'asz--Simonovits (KLS) conjecture asks whether every isotropic log-concave probability measure on $\mathbb R^n$ has a Cheeger constant bounded below by a universal positive constant.
We prove that $\psi_n\le C$ for a universal constant $C>0$, resolving the KLS conjecture.
We also prove that $C_P(\mu)\le C'$ for every isotropic log-concave probability measure $\mu$ on $\mathbb R^n$, with a universal constant $C'>0$.

%% file: _2_body.tex
\input{01_intro}
\input{10_preli}
\input{20_upper}
\input{30_upper_cont}

%% file: 01_intro.tex
\section{Introduction}
\label{sec:introduction}

The Kannan--Lov\'asz--Simonovits (KLS) conjecture~\cite{kls_1995} predicts that the Poincar\'e constant of an isotropic log-concave probability measure is bounded by a universal constant. In probabilistic terms, if $X\in\mathbb R^n$ is log-concave, $\E[X]=0$, and $\operatorname{Cov}(X)=I_n$, then
\begin{equation*}
 \Var[f(X)]\le C\cdot\E[|\nabla f(X)|^2]
\end{equation*}
for every locally Lipschitz function $f$ with finite Dirichlet energy, with $C$ independent of the dimension and the distribution. For log-concave measures, the Poincar\'e constant is comparable to the square of the reciprocal Cheeger constant~\cite{cheeger_1970,milman_2009,klartag_lehec_2025_survey}. Thus the conjecture is equivalently a dimension-free isoperimetric statement.

Write $C_P(\mu)$ for the Poincar\'e constant of $\mu$, and let $\psi_n$ be the supremum of the reciprocal Cheeger constant over all isotropic log-concave probability measures on $\mathbb R^n$.

The conjecture arose in the study of random walks on convex bodies. Bounds on the KLS constant control conductance and mixing times, and thereby the complexity of sampling, volume computation, and log-concave integration. These connections also underlie algorithms for convex optimization and for putting convex bodies in approximately isotropic position~\cite{lee_vempala_2019_survey,jia_laddha_lee_vempala_2020}.

\subsection{Previous bounds}
\label{sec:previous_bounds}

Kannan, Lov\'asz, and Simonovits established $\psi_n\lesssim\sqrt n$ by localization~\cite{kls_1995}. Subsequent progress used thin-shell concentration~\cite{bobkov_2007_isoperimetric,guedon_milman_2011} and Eldan's stochastic localization method~\cite{eldan_2013}. Lee and Vempala obtained $\psi_n\lesssim n^{1/4}$~\cite{lee_vempala_2024}, Chen proved the first subpolynomial bound~\cite{chen_2021}, and Klartag and Lehec obtained the first polylogarithmic bound~\cite{klartag_lehec_2022}. Further improvements were given by Jambulapati, Lee, and Vempala~\cite{jambulapati_lee_vempala_2022} and Klartag~\cite{klartag_2023}. Letwin's quadratic-form inequality~\cite{letwin_2026} yielded the strongest previous bound,
\begin{equation}\label{eq:introduction_letwin_bound}
 \psi_n\lesssim(\log n)^{1/4},\qquad n\ge2.
\end{equation}
Table~\ref{tab:introduction_kls_progress} summarizes these bounds together with the result of this work.

\begin{table}[htbp]
  \centering 
  \begin{tabular}{|l|l|l|l|}
    \hline
    Year & Author(s) & References & Upper bound for $\psi_n$ \\
    \hline
    1995 & Kannan--Lov\'asz--Simonovits
      & \cite{kls_1995} & $\lesssim \sqrt n$ \\
    2011 & Bobkov; Gu\'edon--Milman
      & \cite{bobkov_2007_isoperimetric,guedon_milman_2011} & $\lesssim n^{5/12}$ \\
    2013 & Eldan
      & \cite{eldan_2013} & $\lesssim n^{1/3}\cdot\sqrt{\log n}$ \\
    2017/2024 & Lee--Vempala
      & \cite{lee_vempala_2024} & $\lesssim n^{1/4}$ \\
    2021 & Chen
      & \cite{chen_2021} & $\leq \exp(C\cdot\sqrt{\log n\cdot\log\log n})$ \\
    2022 & Klartag--Lehec
      & \cite{klartag_lehec_2022} & $\lesssim \log^5 n$ \\
    2022 & Jambulapati--Lee--Vempala
      & \cite{jambulapati_lee_vempala_2022} & $\lesssim \log^{3.2226} n$ \\
    2023 & Klartag
      & \cite{klartag_2023} & $\lesssim \sqrt{\log n}$ \\
    2026 & Letwin
      & \cite{letwin_2026} & $\lesssim \log^{1/4} n$ \\
    2026 & This work
      & Theorem~\ref{thm:main_informal} & $O(1)$ \\
    \hline
  \end{tabular}
  \caption{General upper bounds for the KLS constant $\psi_n$. The notation $f(n)\lesssim g(n)$ means that $f(n)\le C\cdot g(n)$ for a universal constant $C>0$ independent of $n$.}
  \label{tab:introduction_kls_progress}
\end{table}

\subsection{Main result}
\label{sec:our_results}

We prove the KLS conjecture in its Poincar\'e formulation.

\begin{theorem}[Universal Poincar\'e and KLS bounds, informal version of Theorem~\ref{tg:thm:main}]
\label{thm:main_informal}
There is a universal constant $C>0$ such that, for every $n\ge1$ and every isotropic log-concave probability measure $\mu$ on $\mathbb R^n$, each locally Lipschitz function $f$ with $\int|\nabla f|^2\d\mu<\infty$ belongs to $L^2(\mu)$ and satisfies
\begin{equation*}
 \Var_\mu[f]\le C\cdot\int|\nabla f|^2\d\mu.
\end{equation*}
In particular, after increasing the universal constant if necessary,
\begin{equation*}
 C_P(\mu)\le C,\qquad \psi_n\le C.
\end{equation*}
\end{theorem}

The constant is independent of the dimension, the measure, and the test function. The theorem applies to arbitrary isotropic log-concave measures, including measures whose densities are not smooth or have bounded support. 

The proof begins with smooth, uniformly log-concave measures of covariance at most $I$. Our method of iterative refinement improves the dependence of the Poincar\'e bound on the lower curvature bound while keeping the total cost bounded. Once this dependence is removed, approximation gives the theorem in full generality.

\subsection{Technical overview}
\label{sec:technical_overview}

Our iterative refinement connects variance estimates for polynomials of
arbitrary degree with curvature estimates for the full Poincar\'e
constant. By alternately improving the polynomial and curvature estimates,
we introduce successive ordinary logarithms and obtain a first log-star bound. We then repeat the entire refinement,
using the improved polynomial estimates from the preceding round.
This outer iteration introduces successive log-stars. Finally, we
control the accumulated losses to obtain a universal bound.

The curvature arguments are carried out for centered regular
log-concave measures $\nu$ with covariance at most $I$ and density
$e^{-W}$ satisfying $D^2W\succeq aI$, where $0<a\le1$ is a lower
bound on the curvature of $W$.
The estimates hold uniformly over this class, so they also apply to
the affinely normalized posteriors arising in stochastic localization.
We use the analytic and approximation tools in
Section~\ref{sec:analytic_foundations} to pass at the end to
arbitrary isotropic log-concave measures.

\subsubsection{Polynomial variance estimates}
\label{sec:overview_polynomial}

Fix an isotropic log-concave random vector $X\sim\mu$ on
$\mathbb R^n$. Recall that isotropy already controls every linear fluctuation, making $\Var\langle v,X\rangle=|v|^2$ where $v\in\mathbb R^n$. 
We seek to extend this control to general functions, ultimately proving
$\Var_\mu f\le C\cdot\E_\mu|\nabla f|^2$ with a universal constant.

We first bound polynomial variances, with constants that depend explicitly on the degree but not on the dimension.
For an integer $d\ge1$ and a fully symmetric $d$-tensor $T$, subtract
lower-degree corrections from the homogeneous polynomial
$\langle T,x^{\otimes d}\rangle$ so that the resulting polynomial
$P_d^\mu[T]$ and its lower-order derivatives have zero mean.
Equivalently,
\begin{equation*}
 \E_\mu[D^kP_d^\mu[T]]=0\quad\forall\,0\le k<d,\qquad
 D^dP_d^\mu[T]=d!\cdot T.
\end{equation*}

A key estimate controlling fluctuations of the posterior covariance is
Letwin's quadratic-form inequality~\cite[Theorem~1.2]{letwin_2026}:
for every symmetric matrix $M\in\mathbb R^{n\times n}$, $ \Var[X^\top MX]\le2\cdot\E[|\nabla(X^\top MX)|^2] =8\cdot\tr[M^2]$. 
We use this inequality to bound the third-moment tensors that govern
covariance fluctuations during localization, within an induction on
polynomial degree whose initial case $d=1$ follows from isotropy.
Using Lemma~\ref{lem:construction_localization_process}, we localize the
measure and apply the estimates at smaller degrees to each affinely
normalized posterior. To close the induction, we use the derivative
hierarchy in Lemma~\ref{lem:hierarchy} and obtain
\begin{equation*}
 \Var_\mu[P_d^\mu[T]]\le C^d\cdot(d!)^4\cdot\|T\|_{\mathrm{HS}}^2,
\end{equation*}
where $\|\cdot\|_{\mathrm{HS}}$ is the Hilbert--Schmidt norm.
This is the polynomial variance bound of Theorem~\ref{thm:polynomial},
proved in Section~\ref{sec:polynomial_estimates}.

\subsubsection{From polynomial estimates to curvature}
\label{sec:overview_curvature}

We now combine polynomial variance estimates with the curvature
lower bound $D^2W\succeq aI$ to obtain an upper bound on $C_P(\nu)$.
For a centered log-concave measure $\nu$, define the normalized
polynomial fluctuation coefficients by
\begin{equation*}
 K_k(\nu):=\sup_{\|T\|_{\mathrm{HS}}=1}
                 \Var_\nu[P_k^\nu[T]],\qquad
 c_k(\nu):=\frac{\sqrt{K_k(\nu)}}{k!},
 \qquad \forall k\ge1.
\end{equation*}
The supremum is over fully symmetric $k$-tensors.
In Theorem~\ref{thm:dyadic_curvature}, we prove a comparison with
explicit dependence on the polynomial degree. If, for a nondecreasing function
$\ell:\mathbb N\to[1,\infty)$ describing the degree growth and a
sufficiently large coefficient scale $R>0$, we have
\begin{equation*}
 c_k(\nu)\le R^k\cdot\frac{\ell(k)^k}{(k+1)^2}\qquad\forall k\ge1,
\end{equation*}
then, for every dyadic $d\ge2$,
\begin{equation*}
 C_P(\nu)\lesssim R^2\cdot\ell(d)^2\cdot a^{-1/(d+1)}.
\end{equation*}
Increasing $d$ weakens the dependence on curvature, at the expense of
the coefficient factor $\ell(d)^2$. This balance explains why estimates
at all degrees are useful.

We consider the nonnegative diffusion operator
$H=-\Delta+\nabla W\cdot\nabla$ associated with $\nu$ and choose an
eigenfunction $f$ corresponding to its first positive eigenvalue
(Lemma~\ref{lem:first_eigenfunction}):
\begin{equation*}
 Hf=\lambda f,\qquad \lambda=C_P(\nu)^{-1},\qquad
 \E_\nu[f]=0,\qquad \|f\|_{L^2(\nu)}=1.
\end{equation*}
Starting from $F_0=f$, we form $\nabla H^{-1/2}F_r$ at each step and
subtract the mean of each component. We then apply $H^{-1/2}$ to
the centered family and rescale all components by a common factor
that preserves its squared $L^2(\nu)$ norm
(Section~\ref{tg:sec:families}). The squared norms of the discarded
means are the centering losses, whose sum we bound using the
polynomial variance estimates and Lemmas~\ref{dt_lem_inverse},
\ref{dt_lem_block}, and~\ref{lem:dyadic_loss}. Using Bochner's formula
in Lemma~\ref{lem:operator_core_bochner}, we show that the energy
decreases by at least $a\|F_r\|_2^2$ at each step. Since the initial
energy is $\lambda$, we can rule out a spectral gap that is too
small and obtain the stated bound on $C_P(\nu)$.

For the initial polynomial estimate, take $\ell(k)=k$ and enlarge
the universal constant $R$ if necessary. Applying
Theorem~\ref{thm:dyadic_curvature} then gives
\begin{equation*}
 C_P(\nu)\lesssim d^2\cdot a^{-1/(d+1)}
 \qquad\forall d\ge2\text{ dyadic}.
\end{equation*}
Choosing $d\asymp\log(e+a^{-1})$ keeps the curvature factor bounded
and gives
\begin{equation*}
 C_P(\nu)\lesssim[\log(e+a^{-1})]^2.
\end{equation*}
We prove the comparison and this first logarithmic improvement in
Section~\ref{sec:polynomial_curvature}, then apply the same comparison
again with better coefficient bounds.

\subsubsection{Iterative logarithmic curvature refinement}
\label{sec:overview_iteration}
\label{sec:overview_main}

We now use the logarithmic bound on $C_P(\nu)$ from the preceding
step to improve the polynomial variance estimates. In the localization
argument for these estimates, we apply this bound to the posterior
measures and obtain improved bounds on the same coefficients $c_d$
(Theorem~\ref{thm:profile_coefficients}). Applying
Theorem~\ref{thm:dyadic_curvature} with these coefficients then gives
a new bound on $C_P(\nu)$ with one more logarithm in $a^{-1}$.
Thus, we can repeat these two steps to introduce successive logarithms.
Obtaining a first log-star bound also requires control of the
accumulated cost of these repetitions.

For $x\ge0$, define the iterated logarithms recursively by
\begin{equation*}
 \ell_0(x)=x,\qquad
 \ell_{r+1}(x)=\log(e+\ell_r(x))\qquad\forall r\ge0.
\end{equation*}
At inner depth $r\ge1$, the estimate has the form
\begin{equation*}
 C_P(\nu)\le\Gamma_r^2\cdot\ell_r(a^{-1})^2.
\end{equation*}
Here $\Gamma_r$ is a positive amplitude depending only on $r$.
Using Theorem~\ref{thm:profile_coefficients}, we transfer the factor
$\ell_r(a^{-1})$ to a factor $\ell_r(d)$ in the degree-$d$ polynomial
estimate. Choosing $d$ on a logarithmic scale
in $a^{-1}$ in the curvature comparison then produces
$\ell_{r+1}(a^{-1})$. Thus, one refinement replaces the bound
$\Gamma_r^2\cdot\ell_r(a^{-1})^2$ by
$\Gamma_{r+1}^2\cdot[\log(e+\ell_r(a^{-1}))]^2$.

Controlling the amplitude $\Gamma_r$ requires a separate argument. A fixed
centering loss at every step would accumulate exponentially in $r$.
We instead begin with an averaged block of inverse-gradient iterates,
which makes both the initial mass loss and the excess energy small
(Proposition~\ref{prop:block_averaged_start}). For the later losses,
we use joint estimates for the derivative families and a
comparison through the subsequent blocks
(Lemmas~\ref{lem:block_joint_local_loss} and~\ref{lem:delayed_block_comparison}).
Combining these estimates with the coefficient transfer in
Proposition~\ref{prop:block_profile_coefficients}, we keep the total
amplitude polynomial in the refinement depth.
The averaged block has initial mass loss and excess energy of order
$z^{-3}$ at block radius $z$, which measures the squared operator
growth per iterate.
Balancing these against a loss budget of order
$(r+1)^{-1}$ gives the cube-root depth cost.
At the end of Section~\ref{sec:iterated_curvature}, we obtain the full
family of bounds in Proposition~\ref{prop:sharp_iterated_curvature}:
\begin{equation*}
 C_P(\nu)\le\Gamma_r^2\cdot\ell_r(a^{-1})^2,\qquad
 \Gamma_r^2\le C\cdot(r+1)^{1/3}\qquad\forall r\ge1.
\end{equation*}
The exponent here measures the cost of controlling the losses through
the iteration; it has a different role from the curvature exponent
$1/(d+1)$ in the preceding comparison.

A depth $r$ with universally bounded $\ell_r(a^{-1})$ suffices. Writing
$\log^*x$ for the least number of successive natural logarithms
needed to bring $x$ to at most one, we can take $r$ of order
$1+\log^*(a^{-1}+2)$. The logarithmic factor in the preceding bound
is then bounded, leaving only the depth cost $(r+1)^{1/3}$. We obtain
\begin{equation*}
 C_P(\nu)\le C\cdot[1+\log^*(a^{-1}+2)]^{1/3}.
\end{equation*}

For the dimensional estimate, we return to Gaussian localization.
Fix a regular isotropic measure $\mu$ and choose a bounded test
function whose variance-to-energy ratio approximates $C_P(\mu)$.
At a localization time of order $1/\log(en)$, a posterior with
controlled covariance preserves a fixed fraction of this variance
(Lemmas~\ref{lem:posterior_identities} and~\ref{tg:lem:stop}). After
affine normalization, the posterior has curvature at least a constant
multiple of $1/\log(en)$. Applying the inner curvature estimate and
transferring the energy back gives
\begin{equation*}
 C_P(\mu)\le C\cdot(r+1)^{1/3}\cdot\ell_r(\log(en))^2
 \qquad\forall r\ge1.
\end{equation*}
Optimizing $r$, passing by regular approximation, and using the
Cheeger--Poincar\'e comparison yield the first log-star bound of
Theorem~\ref{thm:initial_dimension_bound} in
Section~\ref{sec:initial_dimension_bound}:
\begin{equation*}
 C_P(\mu)\le C\cdot[1+\log^*(n+2)]^{1/3},\qquad
 \psi_n\le C'\cdot[1+\log^*(n+2)]^{1/6}.
\end{equation*}
The exponent $1/6$ comes from taking the square root in the
Cheeger--Poincar\'e comparison.

For the next stage, we use more than this optimized dimensional
conclusion. The curvature estimate remains available at every inner
depth, uniformly over the class of measures. This allows us to
choose the depth separately for each polynomial degree, extract
better coefficients, and run the entire refinement again.

\subsubsection{Repeated height reduction}
\label{sec:overview_repeated_height}

We now repeat iterative logarithmic curvature refinement using the improved
polynomial estimates from the preceding round. At outer stage $j$,
we choose an inner depth $r$ separately for each polynomial degree $d$
and use higher-order block estimates to control the cost of the
inner iteration. Write $\mathfrak t(x):=1+\log^*(x+2)$ for the shifted
logarithmic height function, and set $\mathfrak t_0(x):=x$ for all $x\ge1$.
After optimizing the depth and block order, the outer iteration gives
the following bound for any fixed odd block order $q\ge3$:
\begin{equation*}
 C_P(\nu)\le C_q\cdot C_*^{j-1}\cdot\mathfrak t_j(1+a^{-1})^{1/q},
 \qquad \mathfrak t_{j}=\mathfrak t\circ\mathfrak t_{j-1},
 \qquad\forall j\ge1.
\end{equation*}

Here $C_*$ is universal, and $C_q$ includes the initialization cost.
Each additional outer round contributes a factor $C_*$, giving
the prefactor $C_q C_*^{j-1}$ at stage $j$.
In Proposition~\ref{prop:repeated_height_outer_profiles}, we obtain
another log-star with each complete repetition, at the cost of
the fixed multiplicative factor $C_*$.

We first control how the cost grows with the inner depth $r$, using
the higher-order block estimates in Section~\ref{sec:repeated_height_reduction}.
Let $Q$ denote the variable block order, with $q$ reserved for the
fixed order above. For each odd $Q>1$, we bound the initial centering
loss by a multiple of $z^{-Q}$, where $z$ is the block radius.
Requiring this bound to fit a loss budget of order $(r+1)^{-1}$
gives the depth factor $(r+1)^{1/Q}$.
The order $Q$ thus controls the cost of depth, whereas the polynomial
degree $d$ controls the curvature exponent in the comparison theorem.

For the next outer round, we need improved bounds on the polynomial
coefficients. Fix an odd integer $q\ge3$ as the block order. For each
degree $d$, we choose the inner depth $r$ of order $\mathfrak t(d)$,
so that $\ell_r(d)$ is bounded. We apply the higher-order curvature
estimate in Proposition~\ref{prop:repeated_height_original_seed_profile}
with $Q=q$ and use the coefficient-extraction argument in
Lemma~\ref{lem:iterated_coefficient_seed} to obtain
\begin{equation}\label{eq:overview_initial_coefficients}
 c_d\le[C_q\cdot\mathfrak t(d)^{1/(2q)}]^{d-1},
\end{equation}
for any $d\ge2$ and centered log-concave measure of covariance at most $I$.
Here the Poincar\'e profile at fixed order $q$ has a depth factor
$(r+1)^{1/q}$; the square root in the definition of $c_d$ gives
the exponent $1/(2q)$.
We can now repeat the entire refinement, with each outer round
incurring a multiplicative cost.

For fixed odd $q\ge3$, let $A_{j,q}$ denote the accumulated
multiplicative factor in the curvature bound at outer stage $j$.
Here $j$ counts the outer stages, each repeating the entire inner
refinement; $r$ is the depth within one such stage. Each additional
outer round contributes a universal factor $C_*$, so we define
$A_{j,q}:=A_{1,q}\cdot C_*^{j-1}$ for $j\ge1$, where $A_{1,q}>0$
is the initialization cost and depends only on $q$.
This factor is independent of $r$, $Q$, and $a$; we keep their
dependence explicit in the bound
\begin{equation*}
 C_P(\nu)\le A_{j,q}\cdot\mathfrak t_j(Q)^{1/q}
                   \cdot(r+1)^{1/Q}\cdot\ell_r(a^{-1})^2.
\end{equation*}
We obtain this bound from Proposition~\ref{prop:repeated_height_outer_profiles}
after including the fixed initialization costs in $A_{1,q}$.
It holds for every $j\ge1$, odd $Q>1$, and integer $r\ge r_Q$,
where the lower threshold $r_Q$ depends only on $Q$.

For fixed $Q$, the inner-depth factor $(r+1)^{1/Q}$ grows
polynomially in $r$. We choose a suitable
$Q\asymp\log(r+1)$ to keep this factor bounded, leaving the
dependence $\mathfrak t_j(Q)^{1/q}$. Choosing the depth and order
separately for each polynomial degree $d$, we obtain
\begin{equation}\label{eq:overview_outer_coefficients}
 c_d\le[C\cdot\sqrt{A_{j,q}}\cdot
                \mathfrak t_{j+1}(d)^{1/(2q)}]^{d-1}.
\end{equation}
We apply the block comparison in Proposition~\ref{prop:second_logstar_profile}
to these improved coefficient bounds to complete the next outer round.
This replaces $\mathfrak t_j$ by $\mathfrak t_{j+1}$ and gives
$A_{j+1,q}=C_*\cdot A_{j,q}$.

To keep $C_*$ independent of $j$, we use the $q=3$ case of
Eq.~\eqref{eq:overview_initial_coefficients} as a fixed coefficient
bound, along with the improved bound in
Eq.~\eqref{eq:overview_outer_coefficients}, taking the smaller of the
two at each degree. With the fixed bound, we keep the constants in the
block estimates uniform as $A_{j,q}$ grows; with the improved bound,
we reduce the height dependence to $\mathfrak t_{j+1}$.

For each fixed $x$, repeated applications of $\mathfrak t$ reduce
the height to a universal constant
(Lemma~\ref{lem:repeated_height_termination}). Define the required
number of rounds by $\kappa(x)=\min\{j\ge0:\mathfrak t_j(x)\le3\}$.
Although finite for each $x$, this number grows without bound with $x$.
The height is then bounded, but the accumulated outer factor can
still contain $C_*^{\kappa(x)}$. In the next step, we reorganize
the repetitions so that their total multiplicative cost stays bounded.

\subsubsection{Controlling repetition costs through geometric decay}
\label{sec:overview_repetition_costs}

We now keep the accumulated cost bounded by replacing the fixed factor
$C_*$ with losses that approach one and assigning geometrically
decreasing loss budgets to successive groups of outer rounds. Throughout this
argument, we work with the common coefficient radius
\begin{equation*}
 \mathcal A(\mu)=\max\{1,\sup_{d\ge2}c_d(\mu)^{2/(d-1)}\},
 \qquad C_P(\mu)\le C\cdot\mathcal A(\mu),
\end{equation*}
and use the comparison in Lemma~\ref{lem:coefficient_radius_comparison}
only at the end, so that its fixed factor does not accumulate through
the repetitions.

We first strengthen the block and coefficient-transfer estimates so
that one outer round costs only $e^{C\epsilon}$
(Proposition~\ref{prop:slow_offset_outer_profiles}). A smaller
$\epsilon>0$ requires a larger degree cutoff
$K_\epsilon=\lceil C_{\mathrm{cut}}\cdot\epsilon^{-2}\rceil$
and a larger starting depth, where $C_{\mathrm{cut}}$ is a fixed
universal constant. We control these costs using earlier coefficient
bounds below the cutoffs and obtain curvature estimates at prescribed
depths with a loss close to one
(Propositions~\ref{prop:nearunit_secondary_blocks}
and~\ref{prop:nearunit_secondary_prescribed_profile}).

The block constants must remain uniform as the number of retained
coefficient bounds grows. For this purpose, we assign these bounds to
successive, disjoint degree ranges and bound their errors together
by a single convergent tail of the form
\begin{equation*}
 \sum_{k>K_0}k^M\cdot e^{-c\sqrt{k}}<\infty.
\end{equation*}
Here $K_0$, $M\ge0$, and $c > 0$ are fixed. In
Proposition~\ref{prop:constant_chain_blocks}, we use this estimate to
obtain block constants independent of the number of retained bounds.

Next, we control the remaining height dependence introduced by the
larger cutoffs and depths. Recall that $\mathfrak t$ is the shifted log-star
function, and $\chi(x)=\max\{3,1+\kappa(x)\}$ counts, up to fixed
shifts, the outer rounds needed to reduce the height to a bounded
range. We smooth both functions by setting
\begin{equation*}
 \overline f(x)=\max\left\{4,\frac14\int_0^4 f(x+t)\,\d t\right\},
 \qquad f=\mathfrak t,\chi,\quad\forall x\ge1.
\end{equation*}
In Lemma~\ref{lem:smoothed_secondary_maps}, we verify that these
averages are nondecreasing Lipschitz majorants. For integers $i\ge0$, set
\begin{equation*}
 W_i(x)=\overline{\mathfrak t}(\overline\chi^{\circ i}(x)),
 \qquad W_{i+1}(x)=W_i(\overline\chi(x)).
\end{equation*}
After $i$ grouped refinements, $W_i$ describes the remaining dependence
on inner depth. We use its slow growth and contraction to control
changes in the cutoff and depth parameters
(Lemma~\ref{lem:smoothed_secondary_costs}).

Choose geometrically decreasing loss budgets
$\alpha_i=\alpha_0\cdot2^{-i}$, with $\alpha_0=1/16$.
We normalize the logarithms by defining
$\mathcal L_r(x)=\ell_r(x)/\varrho$, where $\varrho$ is the unique
positive solution of $\varrho=\log(e+\varrho)$.
This normalization lets us make the logarithmic factor arbitrarily
close to one by increasing the depth.
After $i$ grouped refinements, let $A_i$ denote the accumulated
multiplicative factor and $S_i$ the additive offset to $W_i$.
Starting from a fixed universal $A_0>0$ and $S_0=1$, define
\begin{equation*}
 A_{i+1}:=e^{C_A\cdot\alpha_i}\cdot A_i,\qquad
 S_{i+1}:=S_i+B\cdot2^{-i},
\end{equation*}
where $C_A,B>0$ are suitable universal constants.
Let $R_i$ denote the lower threshold for the inner depth at this stage.
In Proposition~\ref{prop:constant_chain_profiles}, we prove the inductive estimate
\begin{equation*}
 \mathcal A\le A_i\cdot[W_i(r)+S_i]^{1/3}
                 \cdot\mathcal L_r(a^{-1})^2\qquad\forall r\ge R_i.
\end{equation*}

To pass from stage $i$ to stage $i+1$ at depth $r\ge R_{i+1}$,
take $j=\kappa(r)$ outer rounds and divide the loss budget by setting
$\epsilon=\frac{\alpha_i}{j+1}$. Then
$\epsilon\cdot j\le\alpha_i$, so the loss for the entire group
is at most $e^{C\alpha_i}$. A suitable block order also makes the
inner-depth factor at most $e^{\alpha_i}$.
The smaller $\epsilon$ increases the cutoffs and starting depths;
we control the resulting height dependence by
\begin{equation*}
 W_i(\epsilon^{-1})\le W_{i+1}(r)+O(2^{-i}).
\end{equation*}
Thus we improve the remaining height function from $W_i$ to $W_{i+1}$
with only the multiplicative and additive increments specified above.

We bound the cumulative multiplicative cost by summing its logarithms,
and the additive cost by summing its increments. Both sums are
geometric, giving
\begin{equation*}
 A_i\le A_0\cdot e^{2C_A\cdot\alpha_0},\qquad
 S_i\le S_0+2B
\end{equation*}
for every finite $i$.

It remains to control the growing starting depths $R_i$, which are
fixed powers of $\alpha_i^{-1}$. For a fixed sufficiently large $C$,
we show in Lemma~\ref{lem:constant_chain_arguments} that
\begin{equation*}
 W_i(R_i+C\cdot\mathfrak t(x))\le W_i(x)+b\cdot2^{-i}
 \qquad\forall x\ge1,
\end{equation*}
with a fixed universal constant $b>0$. Using the slow growth and
contraction of $W_i$, we therefore keep the additional error from
these thresholds geometrically decreasing as well.

Finally, we choose the number of groups and the inner depth to obtain
a bound independent of curvature. For $x=\max\{1,a^{-1}\}$, a finite
choice of $i$ makes $W_i(x)$ universally bounded
(Lemma~\ref{lem:smoothed_secondary_termination}). At the allowed
depth $r=R_i+\lceil C\cdot\mathfrak t(x)\rceil$,
both $W_i(r)$ and $\mathcal L_r(a^{-1})$ are universally bounded. Applying the inductive profile, we obtain
$\mathcal A\le C$, hence $C_P\le C$
(Proposition~\ref{prop:constant_chain_endpoint}). We then use regular
isotropic approximation (Lemmas~\ref{lem:poincare_limit}
and~\ref{lem:regular_isotropic_approximation}) and the
Cheeger--Poincar\'e comparison to deduce Theorem~\ref{tg:thm:main} for
arbitrary isotropic log-concave measures.

\paragraph{Organization of the paper.}
Sections~\ref{sec:preliminaries} and~\ref{sec:analytic_foundations}
provide the notation, standard inequalities, and analytic foundations.
Sections~\ref{sec:polynomial_estimates} and~\ref{sec:polynomial_curvature}
establish the polynomial estimates and the curvature comparison.
Section~\ref{sec:iterated_curvature} establishes iterative logarithmic curvature refinement,
and Section~\ref{sec:initial_dimension_bound} derives its dimensional
consequence. Section~\ref{sec:repeated_height_reduction} proves
repeated height reduction. Section~\ref{sec:controlling_repetition_costs}
controls the repetition costs and completes the iterative refinement
to a universal bound.

%% file: 10_preli.tex
\section{Preliminaries}
\label{sec:preliminaries}

We fix our notation and define the functional inequalities studied in this paper in Sections~\ref{sec:preliminaries_notation}--\ref{sec:preliminaries_functional}. Section~\ref{sec:preliminaries_appell} introduces the Appell polynomials and their basic identities. Section~\ref{pre:sec:inputs} collects the geometric and functional inequalities, including the weighted Bochner identity. Sections~\ref{sec:analytic_tools}--\ref{sec:stochastic_tools} record the analytic inequalities and results from stochastic calculus used in the proofs.

\subsection{Notation}
\label{sec:preliminaries_notation}

\paragraph{Vectors and matrices.}
We work over $\mathbb R$. For $x,y\in\mathbb R^n$, write $\langle x,y\rangle:=x^\top y$ and $|x|:=\sqrt{\langle x,x\rangle}$, and let $S^{n-1}:=\{x\in\mathbb R^n:|x|=1\}$.
We denote the identity matrix by $I_n$, or by $I$ when the dimension is clear, and the transpose of $M$ by $M^\top$.
For $M,N\in\mathbb R^{n\times n}$, define
\begin{equation*}
\begin{aligned}
  \tr[M]&:=\sum_{i=1}^n M_{ii},
  &\qquad \|M\|_{\mathrm{op}}&:=\sup_{|x|=1}|Mx|,\\
  \|M\|_{\mathrm{HS}}&:=\sqrt{\tr[M^\top M]},
  &\qquad \langle M,N\rangle_{\mathrm{HS}}&:=\tr[M^\top N].
\end{aligned}
\end{equation*}
For a finite tensor, the Hilbert--Schmidt norm is the square root of the sum of the squares of its entries.
For symmetric matrices $A,B$, we write $A\preceq B$ if $x^\top(B-A)x\ge0$ for every $x\in\mathbb R^n$.

\paragraph{Expectation and covariance.}
For a probability measure $\mu$, write
\begin{equation*}
  \E_\mu[f]:=\int f\d\mu,
  \qquad
  \Var_\mu[f]:=\E_\mu\bigl[(f-\E_\mu[f])^2\bigr],
\end{equation*}
whenever these quantities are defined.
For square-integrable random vectors $U,V$ on this probability space, define
\begin{equation*}
  \operatorname{Cov}_\mu(U,V)
  :=\E_\mu\bigl[(U-\E_\mu[U])(V-\E_\mu[V])^\top\bigr],
  \qquad
  \operatorname{Cov}_\mu(U):=\operatorname{Cov}_\mu(U,U).
\end{equation*}
Scalar arguments are viewed as one-dimensional vectors.
We omit the subscript $\mu$ when the underlying probability measure is clear, and write $\|f\|_2:=\sqrt{\E_\mu[f^2]}$ for the real $L^2(\mu)$ norm.
A function is centered if $\E_\mu[f]=0$.

\paragraph{Derivatives.}
For a scalar function $f$ on $\mathbb R^n$, $\partial_i f$ denotes its derivative in the $i$th coordinate, $\nabla f$ its gradient, and $D^2f$ its Hessian.
We write
\begin{equation*}
  \Delta f:=\sum_{i=1}^n\partial_i^2f,
  \qquad
  \operatorname{Lip}(f):=\sup_{x\ne y}\frac{|f(x)-f(y)|}{|x-y|}.
\end{equation*}
Derivatives are interpreted weakly when the function is not smooth and the relevant Sobolev regularity is specified.

\paragraph{Asymptotic notation.}
Throughout, $\log$ denotes the natural logarithm. A positive integer is dyadic if it is a power of two, including $1=2^0$.
The symbols $C,c>0$ denote universal constants, whose values may change from one occurrence to the next, unless a dependence is stated explicitly.
For nonnegative quantities $a,b$, the notation $a\lesssim b$ means $a\le Cb$ for a universal constant $C$, and $a\asymp b$ means both $a\lesssim b$ and $b\lesssim a$.

\subsection{Log-concavity and functional inequalities}
\label{sec:preliminaries_functional}

A probability measure $\mu$ on $\mathbb R^n$ is log-concave if
\begin{equation*}
  \mu(sA+(1-s)B)
  \geq \mu(A)^s\cdot\mu(B)^{1-s},
  \qquad 0\leq s\leq 1,
\end{equation*}
for compact sets $A,B\subseteq\mathbb R^n$.  Equivalently, $\mu$ has a density $e^{-V}$ with respect to Lebesgue measure on the affine hull of its support, where $V$ is an extended-real-valued convex function and may equal $+\infty$ outside the support.  We call $\mu$ isotropic when its mean is zero and its covariance is the identity:
\begin{equation*}
  \int x \d\mu(x)=0,
  \qquad
  \int x\otimes x \d\mu(x)=I_n.
\end{equation*}
Every full-dimensional log-concave measure with nonsingular covariance has an isotropic image, unique up to an orthogonal transformation.

We use the standard closure properties of log-concavity: affine images and convolutions preserve log-concavity by Pr\'ekopa's theorem~\cite[Theorems~6--8]{prekopa_1973}. Conditioning on a convex set of positive probability and multiplying a density by $e^{-q}$ for convex $q$ also preserve log-concavity, provided the normalizing integral is finite and positive. A full-dimensional log-concave probability density has exponential decay and hence finite moments of every order; see~\cite[Section~1]{klartag_lehec_2025_survey}.

The Poincar\'e constant $C_P(\mu)$ is the least $C$ for which
\begin{equation}
\label{eq:introduction_poincare}
  \Var_{\mu}[f]
  \leq C\cdot\int |\nabla f|^2 \d\mu
\end{equation}
holds for every locally Lipschitz function $f$ with finite Dirichlet energy $\int|\nabla f|^2\d\mu$, with the additional requirement that every such $f$ belongs to $L^2(\mu)$. If no finite $C$ has this property, set $C_P(\mu):=+\infty$. For full-dimensional log-concave measures, Lemma~\ref{lem:poincare_limit} justifies extending the inequality from compactly supported smooth functions to all such test functions.

For a Borel set $A\subseteq\mathbb R^n$, define its outer boundary measure by
\begin{equation*}
  \mu^+(A)
  :=\liminf_{\varepsilon\downarrow 0}
    \frac{\mu(A+\varepsilon B_2^n)-\mu(A)}{\varepsilon},
\end{equation*}
where $B_2^n:=\{x\in\mathbb R^n:|x|\le1\}$. The Cheeger constant of $\mu$ and its reciprocal are
\begin{equation*}
  h_\mu
  :=\inf_{0<\mu(A)<1}
    \frac{\mu^+(A)}{\min\{\mu(A),1-\mu(A)\}},
  \qquad
  \psi_\mu:=h_\mu^{-1}.
\end{equation*}
For full-dimensional log-concave measures, Cheeger's inequality and the reverse Cheeger--Buser inequality show that $C_P(\mu)$ and $\psi_\mu^2$ are comparable up to universal factors~\cite{cheeger_1970,buser_1982,ledoux_2004,milman_2009}. Using the explicit reverse comparison of De Ponti and Mondino~\cite{deponti_mondino_2021}, we obtain
\begin{equation*}
  \frac14
  \leq \frac{\psi_\mu^2}{C_P(\mu)}
  \leq \pi.
\end{equation*}
Thus, if $\psi_n:=\sup\{\psi_\mu:\mu\text{ is isotropic and log-concave on }\mathbb R^n\}$, then the KLS conjecture is precisely $\sup_{n\geq 1}\psi_n<\infty$.
Equivalently, every isotropic log-concave measure should satisfy $C_P(\mu)\leq C$. In affine-covariant form, if $\Sigma$ is the covariance of an arbitrary log-concave measure, the conjecture predicts $C_P(\mu)\leq C\cdot\|\Sigma\|_{\mathrm{op}}$.
The covariance dependence is necessary: testing Eq.~\eqref{eq:introduction_poincare} against linear functions gives $C_P(\mu)\geq\|\Sigma\|_{\mathrm{op}}$.

\subsection{Appell polynomials}
\label{sec:preliminaries_appell}

Let $\mu$ be isotropic, and let $\eta(x):=x$ denote the coordinate map. We write the multivariate Appell polynomials associated with $\mu$~\cite{a04} as symmetric tensor-valued polynomials $\mathcal A_k^\mu(x)$, defined by the formal generating identity
\begin{equation}\label{eq:appell}
 \frac{\exp(\langle z,x\rangle)}{\E_\mu[\exp(\langle z,\eta\rangle)]}
 =\sum_{k\ge0}\frac1{k!}\cdot\langle \mathcal A_k^\mu(x),z^{\otimes k}\rangle.
\end{equation}
Only finitely many coefficients are used at any fixed degree; no convergence of the entire series is asserted or needed. For a fully symmetric coefficient tensor $T$, write
\begin{equation*}
 P_k^\mu[T](x):=\langle T,\mathcal A_k^\mu(x)\rangle,\qquad
 K_k(\mu):=\sup_{\|T\|_{\mathrm{HS}}=1}\Var_\mu P_k^\mu[T].
\end{equation*}
These definitions also apply to a probability measure with finite moments. The norm sums over all ordered indices. For an isotropic measure, the first polynomials are
$\mathcal A_1=x$, $\mathcal A_2=x^{\otimes2}-I$, and
\begin{equation*}
 (\mathcal A_3)_{ijk}=x_ix_jx_k-\delta_{ij}x_k-\delta_{ik}x_j-\delta_{jk}x_i
 -\E[\eta_i\eta_j\eta_k].
\end{equation*}
The correction terms of lower degree are determined by the moments of $\mu$ and need not coincide with the Gaussian Hermite corrections.

Formal differentiation of Eq.~\eqref{eq:appell} gives
\begin{equation}\label{eq:appell_centering}
 \E[P_k[T]]=0,\quad
 \E[D^jP_k[T]]=0\ (1\le j<k),\quad
 D^kP_k[T]=k!\cdot T.
\end{equation}
Differentiating once contracts $k\cdot\mathcal A_{k-1}$ with $T$. Thus every polynomial $q$ of degree at most $s$ has the exact expansion
\begin{equation*}
 q-\E[q]=\sum_{k=1}^s\frac1{k!}\cdot P_k^\mu[\E[D^kq]].
\end{equation*}
Indeed, both sides have the same expected derivatives of every positive order and the same mean; descend from the highest coefficient. Consequently, if $K_k\le\overline K_k$ up to degree $s$, then
\begin{equation}\label{eq:polynomial_minkowski}
 \sqrt{\Var_\mu q}\le\sum_{k=1}^s c_k\cdot\|\E[D^kq]\|_{\mathrm{HS}},
 \qquad c_k:=\sqrt{\overline K_k}/k!.
\end{equation}
This is an $L^2$ triangle inequality, not an orthogonality assertion for Appell polynomials.

\subsection{Standard inequalities}
\label{pre:sec:inputs}
\label{tg:sec:normalization}

We collect the quadratic-form estimate of Letwin~\cite{letwin_2026}, the functional inequalities recorded by Klartag and Lehec~\cite{klartag_lehec_2025_survey}, a matrix form of the Brascamp--Lieb inequality~\cite{brascamp_lieb_1976}, and the weighted Bochner identity~\cite{bakry_gentil_ledoux_2014}.
Letwin's quadratic-form inequality controls the moment terms in the polynomial localization argument of Section~\ref{sec:polynomial_estimates} and the covariance fluctuations in Section~\ref{sec:gaussian_localization}.

\begin{lemma}[Quadratic forms, {\cite[Theorem~1.2]{letwin_2026}}]
\label{lem:quadratic_forms}
Let $X$ have an isotropic log-concave probability distribution $\mu$ on $\mathbb R^n$.
For every real symmetric matrix $B\in\mathbb R^{n\times n}$,
\begin{equation}\label{tg:eq:quadratic}
  \Var_\mu[X^\top B X]\le8\cdot\tr[B^2].
\end{equation}
\end{lemma}

To convert our Poincar\'e bounds into bounds on $\psi_n$, we use the Cheeger--Poincar\'e comparison with the explicit constants recorded after Corollary~21 of Klartag and Lehec~\cite{klartag_lehec_2025_survey}.

\begin{lemma}[Cheeger--Poincar\'e comparison, {\cite[Corollary~21 and the following remark]{klartag_lehec_2025_survey}}]
\label{lem:cheeger_poincare_comparison}
For a log-concave probability measure $\mu$ on $\mathbb R^n$, let $C_P(\mu)$ be its Poincar\'e constant and let $\psi_\mu$ be its reciprocal Cheeger constant.
Then
\begin{equation}\label{tg:eq:cheeger}
  C_P(\mu)\le4\cdot\psi_\mu^2,\qquad \psi_\mu^2\le\pi\cdot C_P(\mu).
\end{equation}
\end{lemma}

The localization argument also requires a test function whose variance and Lipschitz constant detect the size of $C_P(\mu)$ and whose absolute value has a universal bound. This last bound controls the contribution of paths along which the covariance leaves the prescribed interval. We state the result after centering the function and normalizing its variance to one. The short proof below derives this normalization from the cited result.

\begin{lemma}[Bounded test function for the Poincar\'e constant, {\cite[Theorem~20 and Corollary~21]{klartag_lehec_2025_survey}}]
\label{lem:bounded_witness}
There is a universal constant $B_0\ge1$ with the following property.
For every log-concave probability measure $\mu$ on $\mathbb R^n$ with $0<C_P(\mu)<\infty$, set $k:=C_P(\mu)$. Then there exists a bounded Lipschitz function $f:\mathbb R^n\to\mathbb R$ satisfying
\begin{equation}\label{tg:eq:witness}
  \E_\mu[f]=0,\qquad \E_\mu[f^2]=1,\qquad \sup_{x\in\mathbb R^n}|f(x)|\le B_0,\qquad \operatorname{Lip}(f)\le B_0/\sqrt{k}.
\end{equation}
\end{lemma}
\begin{proof}[Derivation of the normalization]
Take the bounded $1$-Lipschitz witness $g$ from~\cite[Theorem~20]{klartag_lehec_2025_survey} and set $\sigma^2:=\Var_\mu[g]$.
There are universal constants $c_1,C_1>0$ such that $\sigma^2\ge c_1\cdot\psi_\mu^2$ and $\|g\|_\infty^2\le C_1\cdot\sigma^2$.
If the supremum bound is initially essential, clip $g$ to its essential range; this preserves its values almost everywhere and its Lipschitz constant and makes the bound global.
The Cheeger inequality stated in Lemma~\ref{lem:cheeger_poincare_comparison} gives $\psi_\mu^2\ge k/4$, hence $\sigma^2\ge c_1\cdot k/4>0$.
Define $f:=(g-\E_\mu[g])/\sigma$.
Then $\E_\mu[f]=0$ and $\E_\mu[f^2]=1$, while
\begin{equation*} \|f\|_\infty\le2\cdot\|g\|_\infty/\sigma\le2\cdot\sqrt{C_1},\qquad \operatorname{Lip}(f)\le1/\sigma\le2/\sqrt{c_1\cdot k}. \end{equation*}
In the first chain, centering gives the first inequality and the supremum bound for $g$ gives the second.
In the second chain, the $1$-Lipschitz property gives the first inequality and the variance lower bound gives the second.
Thus $B_0:=\max\{1,2\cdot\sqrt{C_1},2/\sqrt{c_1}\}$ gives Eq.~\eqref{tg:eq:witness}.
\end{proof}

The Gaussian tilts used in localization add a positive quadratic term to the negative logarithm of the density. The following form of the Brascamp--Lieb inequality converts the matrix of this quadratic term into a weighted Poincar\'e inequality, including when the measure is supported on a convex set.

\begin{lemma}[Matrix Brascamp--Lieb inequality, {\cite{brascamp_lieb_1976}}]
\label{lem:matrix_strong_poincare}
Let $B$ be positive definite and let $\nu$ have density proportional to
$\exp(-V(x)-x^\top Bx/2)$, where $V$ is an extended-valued convex function.
For every locally Lipschitz $f\in L^2(\nu)$ with finite right side,
\begin{equation*}
 \Var_\nu[f]\le\E_\nu[(\nabla f)^\top B^{-1}\nabla f].
\end{equation*}
\end{lemma}

The weighted Bochner identity expresses the squared norm of the weighted Laplacian applied to a function as the sum of a Hessian term and a curvature term. We record its integrated form for smooth compactly supported functions; see~\cite[Section~C.5, Eqs.~(C.5.2)--(C.5.3)]{bakry_gentil_ledoux_2014} for the pointwise identity.

\begin{lemma}[Integrated weighted Bochner identity]\label{lem:weighted_bochner}
Let $\d\mu=e^{-W}\d x$ be a probability measure on $\mathbb R^n$ with $W\in C^\infty(\mathbb R^n)$, and define $H:=-\Delta+\nabla W\cdot\nabla$.
For $f,g\in C_c^\infty(\mathbb R^n)$, integration by parts gives
\begin{equation*}
 \langle f,Hg\rangle_{L^2(\mu)}=\E_\mu[\nabla f\cdot\nabla g].
\end{equation*}
Moreover,
\begin{equation}\label{tg:eq:bochner}
 \|Hg\|_2^2=\E[\|D^2g\|_\mathrm{HS}^2]
       +\E[(\nabla g)^\top D^2W\nabla g].
\end{equation}
\end{lemma}

\subsection{Analytic tools}
\label{sec:analytic_tools}

For the operator arguments in Section~\ref{sec:analytic_foundations}, we use two standard facts about Sobolev spaces~\cite[Sections~5.7 and~6.3.1]{evans_2010}. First, if $W$ is smooth, $g\in H^1_{\mathrm{loc}}$, and $-\Delta g+\nabla W\cdot\nabla g\in L^2_{\mathrm{loc}}$ in the distributional sense, then $g\in H^2_{\mathrm{loc}}$. Second, the embedding $H^1(B_R)\hookrightarrow L^2(B_R)$ is compact. For a nonnegative self-adjoint operator $H$ with compact resolvent, applying the spectral theorem for compact self-adjoint operators~\cite[Appendix~D.6, Theorem~7]{evans_2010} to $(I+H)^{-1}$ gives an orthonormal basis of eigenvectors of $H$. Expansion in this basis gives the Rayleigh characterization of its least positive eigenvalue on $(\ker H)^\perp$. On a spectral subspace where $H\ge\lambda I$ with $\lambda>0$, it also gives $\|H^{-1/2}\|_{\mathrm{op}}\le\lambda^{-1/2}$. Section~\ref{sec:analytic_foundations} verifies the hypotheses for the weighted Laplacians considered here.

The polynomial estimates use the following forms of Minkowski's and Bessel's inequalities. They allow us to estimate all tensor components together without introducing a factor depending on the dimension.

\begin{lemma}[Norm inequalities in Hilbert spaces]
\label{lem:hilbert_norm_tools}
Let $\mathcal H$ be a finite-dimensional real Hilbert space. For square-integrable $\mathcal H$-valued random variables $U_1,\ldots,U_r$,
\begin{equation*}
 \left(\E\left\|\sum_{k=1}^r U_k\right\|^2\right)^{1/2}
 \le\sum_{k=1}^r(\E\|U_k\|^2)^{1/2}.
\end{equation*}
If $Z_1,\ldots,Z_m$ are orthonormal real random variables in $L^2$, then every $U\in L^2(\mathcal H)$ satisfies
\begin{equation*}
 \sum_{i=1}^m\|\E[Z_i\cdot U]\|^2\le\E\|U\|^2.
\end{equation*}
Both statements remain valid for the seminorm $\|u\|_K:=\langle u,Ku\rangle^{1/2}$ associated with a fixed positive semidefinite operator $K$ on $\mathcal H$.
\end{lemma}
\begin{proof}
The first inequality is the triangle inequality in $L^2(\mathcal H)$. The second follows by applying scalar Bessel's inequality in an orthonormal basis of $\mathcal H$ and summing. Applying these statements to $K^{1/2}U$ proves their weighted versions.
\end{proof}

We repeatedly apply Gronwall's inequality to the derivative energies along localization. The version with a forcing term accounts for the zero initial conditions in Section~\ref{sec:iterated_curvature}.

\begin{lemma}[Gronwall's inequality]
\label{lem:gronwall}
Let $u:[0,T]\to[0,\infty)$ be continuous. If
\begin{equation*}
 u(t)\le u_0+\alpha\cdot\int_0^t u(s)\,\d s
 \qquad(0\le t\le T)
\end{equation*}
for $u_0,\alpha\ge0$, then $u(t)\le u_0\cdot e^{\alpha t}$. If $u$ is absolutely continuous and $u'(t)\le\alpha\cdot u(t)+f(t)$ almost everywhere, where $f\in L^1([0,T])$, then
\begin{equation*}
 u(t)\le e^{\alpha t}\cdot u(0)+\int_0^t e^{\alpha(t-s)}\cdot f(s)\,\d s.
\end{equation*}
\end{lemma}
\begin{proof}
For the first statement, set $v(t):=u_0+\alpha\cdot\int_0^t u(s)\,\d s$. Then $u\le v$ and $v'\le\alpha v$, so $e^{-\alpha t}\cdot v(t)$ is nonincreasing. For the second, multiply the differential inequality by $e^{-\alpha t}$ and integrate.
\end{proof}

\subsection{Stochastic calculus}
\label{sec:stochastic_tools}

We record the results from stochastic calculus used to construct localization and estimate its fluctuations. All filtrations satisfy the usual conditions. Matrix and tensor processes are treated through their coordinates. For It\^o's formula, the product rule, and the isometry below, see~\cite[Sections~B.1--B.2]{bakry_gentil_ledoux_2014}.

\begin{lemma}[It\^o's formula and isometry]
\label{lem:ito_calculus}
Let $B_t$ be an $m$-dimensional Brownian motion and let $X_t$ be an $\mathbb R^k$-valued process of the form
\begin{equation*}
 X_t=X_0+\int_0^t b_s\,\d s+\int_0^t\sigma_s\,\d B_s,
\end{equation*}
where $X_0$ is $\mathcal F_0$-measurable, $b$ and $\sigma$ are progressively measurable, and
$\int_0^T(|b_s|+\|\sigma_s\|_{\mathrm{HS}}^2)\,\d s<\infty$ almost surely for every finite $T$. If $F(t,x)$ is once continuously differentiable in $t$ and twice continuously differentiable in $x$, then
\begin{equation*}
 \d F(t,X_t)
 =\left(\partial_tF+\nabla F\cdot b_t
   +\frac12\cdot\tr[\sigma_t\sigma_t^\top D^2F]\right)\d t
   +\nabla F^\top\sigma_t\,\d B_t,
\end{equation*}
where the derivatives are evaluated at $(t,X_t)$. The formula also holds when $F$ is defined on an open domain, up to the exit time from that domain. In particular,
$\d[X^i,X^j]_t=(\sigma_t\sigma_t^\top)_{ij}\,\d t$, where $[\cdot,\cdot]$ denotes quadratic covariation. For continuous real semimartingales $U,V$,
\begin{equation*}
 \d(UV)=U\,\d V+V\,\d U+\d[U,V].
\end{equation*}
If $G$ is a progressively measurable matrix-valued process with
$\E\int_0^T\|G_s\|_{\mathrm{HS}}^2\,\d s<\infty$, then $\int_0^tG_s\,\d B_s$, $0\le t\le T$, is a square-integrable martingale and
\begin{equation*}
 \E\left|\int_0^T G_s\,\d B_s\right|^2
 =\E\int_0^T\|G_s\|_{\mathrm{HS}}^2\,\d s.
\end{equation*}
\end{lemma}

The next facts justify stopping the localization process, removing the stopping times, and identifying the Brownian motion in the Gaussian observation construction; see~\cite[Sections~3.4, 4.2--4.3, and 5.3]{le_gall_2016}.

\begin{lemma}[Martingale convergence, stopping, and L\'evy's characterization]
\label{lem:martingale_tools}
A bounded continuous local martingale is a martingale. A nonnegative continuous local supermartingale $Z$ with $\E Z_0<\infty$ is a supermartingale, and $\E Z_\tau\le\E Z_0$ for every bounded stopping time $\tau$. If $M$ is a continuous martingale and $\sigma\le\tau$ are bounded stopping times, then
\begin{equation*}
 \E[M_\tau\mid\mathcal F_\sigma]=M_\sigma.
\end{equation*}
A real continuous local martingale $M$ on $[0,\zeta)$ has almost surely a finite limit as $t\uparrow\zeta$ on the event $\{[M]_\zeta<\infty\}$, where $[M]_\zeta:=\lim_{t\uparrow\zeta}[M]_t$, provided $\zeta$ is an increasing limit of stopping times strictly below it.
Finally, an $\mathbb R^m$-valued continuous local martingale $M$ with $M_0=0$ and $[M^i,M^j]_t=\delta_{ij}t$ is an $m$-dimensional Brownian motion.
\end{lemma}

We also need to interchange a spatial integral with a stochastic integral. The following square-integrable form suffices after stopping on compact parameter sets.

\begin{lemma}[Stochastic Fubini]
\label{lem:stochastic_fubini}
Let $B_t$ be an $m$-dimensional Brownian motion and $\rho$ a probability measure on a measurable space $E$. Suppose $G(s,\omega,x)\in\mathbb R^m$ is jointly measurable with respect to the predictable $\sigma$-field and the $\sigma$-field on $E$, and
\begin{equation*}
 \int_E\E\int_0^T|G(s,x)|^2\,\d s\,\d\rho(x)<\infty.
\end{equation*}
Then, for $0\le t\le T$,
\begin{equation*}
 \int_E\left(\int_0^t G(s,x)\cdot\d B_s\right)\d\rho(x)
 =\int_0^t\left(\int_EG(s,x)\,\d\rho(x)\right)\cdot\d B_s,
\end{equation*}
where equality is understood in $L^2$; the processes admit indistinguishable continuous versions.
\end{lemma}
\begin{proof}
The identity holds for functions that are simple in $x$. Approximation in $L^2(\rho\otimes\mathbb P\otimes[0,T])$, Jensen's inequality, and the It\^o isometry in Lemma~\ref{lem:ito_calculus} extend it to the stated class. The right side supplies the continuous version.
\end{proof}

The localization parameters solve a finite-dimensional stochastic differential equation. We use the following local existence and continuation theorem; see~\cite[Sections~B.3--B.4]{bakry_gentil_ledoux_2014}.

\begin{lemma}[Local existence and continuation of SDE solutions]
\label{lem:local_sde}
Let $D\subseteq\mathbb R^k$ be open, and let $b:D\to\mathbb R^k$ and $\sigma:D\to\mathbb R^{k\times m}$ be locally Lipschitz. Let $B_t$ be an $m$-dimensional Brownian motion. For each initial point $X_0\in D$, the equation
\begin{equation*}
 \d X_t=b(X_t)\,\d t+\sigma(X_t)\,\d B_t
\end{equation*}
has a pathwise unique strong solution up to a maximal lifetime $\zeta$. On $\{\zeta<\infty\}$, the solution leaves every compact subset of $D$ before $\zeta$. In particular, a solution with a limit in $D$ at a finite lifetime can be continued.
\end{lemma}

%% file: 20_upper.tex
\section{Analytic Preparation}
\label{sec:analytic_foundations}
\label{sec:four_generation}
\label{tg:sec:domains}

We prepare the operator identities, approximation results, and spectral tools used in the polynomial variance estimates, the curvature comparison, its iteration, and the Gaussian localization transfer. Section~\ref{sec:operator_foundations} establishes the operator domains and Bochner identity needed to control the energies of the inverse-operator families. Section~\ref{sec:stability_approximation} provides regular approximation and stability of the Poincar\'e inequality, so that uniform bounds proved for regular measures extend to general isotropic log-concave measures. Section~\ref{sec:spectral_foundations} constructs a nonconstant eigenfunction attaining the Poincar\'e constant and justifies the inverse square root used in the curvature comparison.

\paragraph{Regular measures and the associated operator.}
We call a probability measure $\mu$ \emph{regular} if it has the form $\d\mu=e^{-W}\d x$, where $W$ is smooth
and there are constants $0<a\le b<\infty$, depending on $\mu$, such that
\begin{equation}\label{tg:eq:regular}
 aI\preceq D^2W\preceq bI.
\end{equation}
Neither $a$ nor $b$ is assumed uniform.
The curvature argument also allows centered regular measures with covariance at most $I$.
We define
\begin{equation*}
 H:=-\Delta+\nabla W\cdot\nabla,\qquad
 \lambda:=C_P(\mu)^{-1},\qquad P_+h:=h-\E_\mu[h].
\end{equation*}
$H$ is the nonnegative self-adjoint operator associated to the Dirichlet
form $\E[\nabla f\cdot\nabla g]$ on $L^2(\mu)$; inverse powers below act
only on the centered subspace. 

\subsection{Operator domains and the Bochner identity}
\label{sec:operator_foundations}

We extend the weighted Bochner identity in Lemma~\ref{lem:weighted_bochner} from smooth compactly supported functions to the full operator domain. The following lemma also establishes the domain properties needed to apply this identity throughout the proof.

\begin{lemma}[Operator core and Bochner identity]\label{lem:operator_core_bochner}
Let $\d\mu=e^{-W}\d x$ be a regular measure satisfying
Eq.~\eqref{tg:eq:regular}, and let $H$ be its nonnegative self-adjoint
operator associated with the Dirichlet form on $L^2(\mu)$.
Its form domain is the weighted Sobolev space
$\{g\in L^2(\mu):\nabla g\in L^2(\mu)\}$, with weak derivatives,
and its kernel consists exactly of the constants.
The space $C_c^\infty(\mathbb R^n)$ is a graph core for $H$, and
$P_+C_c^\infty(\mathbb R^n)$ is a graph core on the centered subspace.
For every $g\in\operatorname{Dom}(H)$, all second weak derivatives belong
to $L^2(\mu)$, mixed derivatives commute, and Eq.~\eqref{tg:eq:bochner} holds.
In particular, each $\partial_i g$ belongs to the form domain.
Every polynomial belongs to the form domain, and
$\langle \phi,Hg\rangle_{L^2(\mu)}=\E_\mu[\nabla\phi\cdot\nabla g]$
for every test function $\phi$ in the form domain.
\end{lemma}

We first identify the form and operator domains by cutoff and mollification. Integration by parts then proves the Bochner identity on smooth functions, and approximation extends it to the full operator domain.

\begin{proof}
On each compact set, the smooth positive weight is bounded above and below
by positive constants. Cutoff and mollification therefore identify the
closure of $C_c^\infty$ in the form norm with the stated weighted Sobolev
space. More explicitly, choose smooth cutoffs $\chi_R$ with
$0\le\chi_R\le1$, equal to one on the ball of radius $R$ and zero outside
the ball of radius $2R$, with $|\nabla\chi_R|\le C/R$ and
$|D^2\chi_R|\le C/R^2$. For a weighted Sobolev function $g$, the
functions $\chi_Rg$ converge to $g$ in the form norm: the terms involving
$1-\chi_R$ vanish by integrability, and
$\|g\nabla\chi_R\|_2\le(C/R)\cdot\|g\|_2$.
Each compactly supported member is then approximated by ordinary
mollification on a fixed compact neighborhood of its support.
The same cutoffs approximate the constant one, which has zero energy.
Conversely, a function of zero energy has weak gradient zero and hence is
constant on the connected space $\mathbb R^n$.
This also identifies $\ker H$.

Let $g\in\operatorname{Dom}(H)$. The weak equation
$-\Delta g+\nabla W\cdot\nabla g=Hg$ and the local elliptic regularity stated in Section~\ref{sec:analytic_tools} give
$g\in H^2_{\mathrm{loc}}$. The product rule gives
\begin{equation*} H(\chi_Rg)=\chi_RHg-2\nabla\chi_R\cdot\nabla g+gH\chi_R. \end{equation*}
This identity initially holds distributionally. Its right side is in
$L^2(\mu)$, so the defining weak characterization of the operator shows
that $\chi_Rg\in\operatorname{Dom}(H)$.
Since $|\nabla W(x)|\le|\nabla W(0)|+b|x|$, the function $H\chi_R$ is
bounded by a constant depending on $\mu$ but independent of $R\ge1$ and is
supported in the annulus $R\le|x|\le2R$.
Consequently $\chi_RHg\to Hg$ in $L^2(\mu)$,
$\|\nabla\chi_R\cdot\nabla g\|_2\to0$, and
$\|gH\chi_R\|_2\to0$. Thus $\chi_Rg\to g$ in the graph norm
$\|g\|_2+\|Hg\|_2$.
Mollifying each compactly supported $H^2$ function gives smooth compactly
supported functions converging in $H^2$ on a common compact set.
The bounded coefficients there imply convergence of their images under
$H$ in $L^2(\mu)$. A diagonal choice proves the assertion about the graph core.
Subtracting the means gives the centered graph core, because $H1=0$ and
convergence in $L^2(\mu)$ implies convergence of the means.

For a smooth compactly supported $g$, differentiation gives
$\partial_i(Hg)=H(\partial_i g)+\sum_j(\partial_i\partial_jW)\partial_jg$.
Integration by parts in each coordinate proves Eq.~\eqref{tg:eq:bochner}.
Let $g_j$ be a sequence in this core that converges in the graph norm to a general $g\in\operatorname{Dom}(H)$.
For a difference $u:=g_j-g_k$, the form identity gives
\begin{equation*} \|\nabla u\|_2^2=\langle u,Hu\rangle_{L^2(\mu)}\le\|u\|_2\cdot\|Hu\|_2. \end{equation*}
The equality is the definition of the operator through its form, and the
inequality is Cauchy--Schwarz. Hence the gradients converge in $L^2(\mu)$.
Applying Eq.~\eqref{tg:eq:bochner} to $u$ and using $D^2W\succeq0$ shows
that the Hessians are Cauchy in $L^2(\mu)$ as well.
Their limits are the distributional second derivatives of $g$.
The bound $D^2W\preceq bI$ passes the curvature term to the limit and
proves Eq.~\eqref{tg:eq:bochner} on all of $\operatorname{Dom}(H)$.
Distributional mixed derivatives commute, and the characterization of the form domain now applies to each $\partial_i g$.

Finally, strong convexity gives
$W(x)\ge W(0)+\nabla W(0)\cdot x+a|x|^2/2$.
Thus every polynomial and its gradient are square integrable, which places
all polynomials in the form domain.
The asserted testing identity is the defining operator identity on that
domain. In the isotropic case, the pairings
$\E[X_i\partial_j\partial_k g]$ are also finite by Cauchy--Schwarz.
\end{proof}

\subsection{Stability of the Poincar\'e inequality}
\label{sec:stability_approximation}

We first record the stability statement used to pass uniform bounds to weak limits, including every locally Lipschitz test function of finite Dirichlet energy.

\begin{lemma}[Stability of the Poincar\'e inequality]\label{lem:poincare_limit}
Let $\mu_j\Rightarrow\mu$ be probability measures on $\mathbb R^n$,
where $\mu$ is absolutely continuous with respect to Lebesgue measure.
Suppose that a finite $K\ge0$ satisfies
$\Var_{\mu_j}[q]\le K\cdot\E_{\mu_j}[|\nabla q|^2]$
for every $q\in C_c^\infty(\mathbb R^n)$ and every $j$.
Then every real locally Lipschitz $q$ with
$\E_\mu[|\nabla q|^2]<\infty$ belongs to $L^2(\mu)$ and satisfies
\begin{equation*} \Var_\mu[q]\le K\cdot\E_\mu[|\nabla q|^2]. \end{equation*}
In particular, $C_P(\mu)\le K$.
\end{lemma}
\begin{proof}
For a fixed $q\in C_c^\infty$, the functions $q$, $q^2$, and
$|\nabla q|^2$ are bounded and continuous. Weak convergence therefore
passes the assumed inequality to $\mu$ with the same constant $K$.

Let now $q$ be locally Lipschitz with finite Dirichlet energy, and set
$q_M:=\max\{-M,\min\{q,M\}\}$ for $M>0$.
The weak chain rule gives $|\nabla q_M|\le|\nabla q|$ almost everywhere.
With the cutoffs used in Lemma~\ref{lem:operator_core_bochner}, the
function $q_{M,R}:=\chi_Rq_M$ is compactly supported and globally
Lipschitz. Convolve it with a smooth compactly supported approximate
identity. The mollifications and their gradients converge almost everywhere with respect to Lebesgue measure to $q_{M,R}$ and its gradient, respectively, and are uniformly
bounded for fixed $M,R$, with supports in one compact set.
Absolute continuity of $\mu$ and dominated convergence pass the
inequality for compactly supported smooth functions to $q_{M,R}$.

As $R\to\infty$, bounded convergence gives $q_{M,R}\to q_M$ in
$L^2(\mu)$. The identity
$\nabla q_{M,R}=\chi_R\nabla q_M+q_M\nabla\chi_R$ shows convergence
of the gradients in $L^2(\mu)$: dominated convergence handles the first
term, and the second has norm at most $CM/R$.
Hence $\Var_\mu[q_M]\le K\cdot\E_\mu[|\nabla q_M|^2]$.
For independent $X,X'\sim\mu$, Fatou's lemma now gives
\begin{equation*} \frac12\cdot\E[(q(X)-q(X'))^2]\le\liminf_{M\to\infty}\Var_\mu[q_M]\le K\cdot\E_\mu[|\nabla q|^2]. \end{equation*}
The first inequality uses pointwise convergence of the truncations and
the variance formula for two independent copies. The second uses the inequality
for $q_M$ and $|\nabla q_M|\le|\nabla q|$.
The right side is finite, so Fubini's theorem gives a point $y$ for which
$\int(q(x)-q(y))^2\d\mu(x)<\infty$.
Since $q(y)$ is finite, this implies $q\in L^2(\mu)$; the expression with two independent copies is then exactly $\Var_\mu[q]$.
\end{proof}

\subsection{Regular isotropic approximation}

Gaussian convolution smooths the density, and a positive quadratic tilt makes the potential uniformly convex. An affine normalization restores isotropy while preserving convergence to the original measure.

\begin{lemma}[Regular isotropic approximation]\label{lem:regular_isotropic_approximation}
Every isotropic log-concave probability measure $\mu$ on $\mathbb R^n$
is the weak limit of regular isotropic measures
$\d\mu_j=e^{-W_j}\d x$ satisfying
$a_jI\preceq D^2W_j\preceq b_jI$ for constants $0<a_j\le b_j<\infty$.
The measures may be obtained by Gaussian convolution, a positive quadratic
tilt, and an affine normalization whose translation tends to zero and
whose linear part tends to $I$.
\end{lemma}

\begin{proof}
Let $X\sim\mu$ and let $Z$ be an independent standard Gaussian.
For $\delta>0$, write $\nu_\delta$ for the probability distribution of
$X+\sqrt\delta Z$. Its density $e^{-V_\delta}$ is positive and smooth.
It is log-concave by preservation of log-concavity under convolution.
Differentiating the Gaussian convolution under the integral gives
\begin{equation*} D^2V_\delta(y)=\delta^{-1}I-\delta^{-2}\cdot\operatorname{Cov}(X\mid X+\sqrt\delta Z=y). \end{equation*}
Indeed, the gradient of the logarithm of the convolution density is
$\delta^{-1}\cdot(\E[X\mid X+\sqrt\delta Z=y]-y)$, and the derivative of the
conditional mean is the conditional covariance divided by $\delta$.
Gaussian damping justifies these differentiations locally uniformly in
$y$. Convexity of $V_\delta$ and positivity of the conditional covariance
therefore give $0\preceq D^2V_\delta\preceq\delta^{-1}I$.

For $\epsilon>0$, define
\begin{equation*} \d\nu_{\delta,\epsilon}(y):=Z_{\delta,\epsilon}^{-1}\cdot e^{-\epsilon|y|^2/2}\d\nu_\delta(y), \end{equation*}
where $Z_{\delta,\epsilon}:=\int e^{-\epsilon|y|^2/2}\d\nu_\delta(y)$.
The new potential $W_{\delta,\epsilon}:=V_\delta+\epsilon|y|^2/2+\log Z_{\delta,\epsilon}$
has Hessian between $\epsilon I$ and $(\delta^{-1}+\epsilon)I$.
For fixed $\delta$, dominated convergence as $\epsilon\downarrow0$
gives convergence in total variation to $\nu_\delta$ and convergence of its
first two moments: apply it with weights $1$, $y_i$, and $y_iy_k$.
These weights are integrable since $X+\sqrt\delta Z$ has finite second moment.
Choose $\delta_j\downarrow0$ and then $\epsilon_j>0$ so that the
error in total variation and the errors in the first two moments are at most
$1/j$. Write $\nu_j:=\nu_{\delta_j,\epsilon_j}$, and denote its mean by
$m_j$ and covariance by $A_j$.
The coupling $X+\sqrt{\delta_j}Z\to X$ in $L^2$ implies
$\nu_j\Rightarrow\mu$, $m_j\to0$, and $A_j\to I$.

Each $A_j$ is positive definite because $\nu_j$ has a positive density.
Let $\mu_j$ be the probability distribution of $A_j^{-1/2}(Y-m_j)$ for $Y\sim\nu_j$.
It is isotropic, and its smooth potential has Hessian
$A_j^{1/2}D^2W_{\delta_j,\epsilon_j}(m_j+A_j^{1/2}x)A_j^{1/2}$.
Thus one may take
$a_j:=\epsilon_j\cdot\lambda_{\min}(A_j)>0$ and
$b_j:=(\delta_j^{-1}+\epsilon_j)\cdot\lambda_{\max}(A_j)<\infty$.
Since $A_j^{-1/2}\to I$ and $m_j\to0$, these affine maps converge to the
identity uniformly on compact sets. Tightness of $\nu_j$ then implies
$\mu_j\Rightarrow\mu$, as asserted.
\end{proof}

\subsection{Spectral gap and inverse square root}
\label{sec:spectral_foundations}

The following lemma identifies the reciprocal Poincar\'e constant with the first positive eigenvalue of $H$ and justifies the inverse square root used in the curvature argument.

\begin{lemma}[Spectral gap and inverse square root]\label{lem:first_eigenfunction}
For every regular measure satisfying Eq.~\eqref{tg:eq:regular}, the operator
$H$ has compact resolvent and a positive first nonzero eigenvalue
$\lambda=C_P(\mu)^{-1}$.
There is a real $f\in\operatorname{Dom}(H)$ with
$\E_\mu[f]=0$, $\|f\|_2=1$, and $Hf=\lambda f$.
The operator $H^{-1/2}$ is bounded on the centered subspace.
For every centered $h\in\operatorname{Dom}(H^{1/2})$, the function
$g:=H^{-1/2}h$ belongs to $\operatorname{Dom}(H)$ and satisfies
$\|Hg\|_2^2=\|H^{1/2}h\|_2^2$ and
$\E_\mu[|\nabla g|^2]=\|h\|_2^2$.
\end{lemma}

We conjugate $H$ to a Schr\"odinger operator with a confining potential and prove compactness of the embedding of its form domain. The variational characterization of the first positive eigenvalue then identifies the spectral gap, and spectral calculus gives the identities for the inverse square root.

\begin{proof}
The unitary map $Ug:=e^{-W/2}\cdot g$ from $L^2(\mu)$ to $L^2(\d x)$
transforms $H$ into the Friedrichs realization of $-\Delta+V$, where
$V:=|\nabla W|^2/4-\Delta W/2$.
This follows by direct differentiation on the core from
Lemma~\ref{lem:operator_core_bochner}, followed by closure of the forms.
Strong convexity gives a unique minimizer $x_*$ of $W$ and
$|\nabla W(x)|\ge a|x-x_*|$, while $\Delta W\le nb$.
It follows that $V(x)\ge a^2\cdot|x-x_*|^2/4-nb/2$.

After adding $nb/2+1$ to this potential, its form norm controls
$\|u\|_2$, $\|\nabla u\|_2$, and $\||x-x_*|u\|_2$.
A sequence bounded in this norm has a subsequence converging in $L^2$
on each fixed ball by the Rellich compactness theorem recalled in Section~\ref{sec:analytic_tools}.
Outside the ball of radius $R$ about $x_*$, its squared $L^2$ norm is at most
$R^{-2}\cdot\||x-x_*|u\|_2^2$, uniformly along the sequence.
A diagonal subsequence and this tail estimate give convergence in all of
$L^2(\d x)$. The embedding of the form domain is therefore compact, which gives
compact resolvent for $H$ by the unitary transformation.

Lemma~\ref{lem:operator_core_bochner} identifies the kernel as the
one-dimensional space of constants. The discrete spectrum consequently
has a positive first nonzero eigenvalue, attained on the centered subspace.
Denote this eigenvalue temporarily by $\lambda_1$.
The Rayleigh variational principle gives $\Var_\mu[q]\le\lambda_1^{-1}\cdot\E_\mu[|\nabla q|^2]$ for every compactly supported smooth $q$.
Apply Lemma~\ref{lem:poincare_limit} to the constant sequence $\mu_j:=\mu$ with $K:=\lambda_1^{-1}$.
Its truncation and cutoff argument extends the inequality to every locally Lipschitz function of finite energy and proves $C_P(\mu)\le\lambda_1^{-1}$.
Conversely, approximate a unit first nonconstant eigenfunction in the form norm by compactly supported smooth functions using Lemma~\ref{lem:operator_core_bochner}.
Their variances tend to one and their energies tend to $\lambda_1$, so the definition of $C_P(\mu)$ gives $C_P(\mu)\ge\lambda_1^{-1}$.
Thus $\lambda_1=C_P(\mu)^{-1}=\lambda$, with the full class of test functions in the definition of $C_P$ accounted for.
Since the operator has real coefficients, a nonzero real or imaginary
part of a corresponding eigenfunction is again an eigenfunction.
Normalize it in $L^2(\mu)$; orthogonality to the constants gives mean zero.
On the centered subspace the spectrum lies in $[\lambda,\infty)$, so
$\|H^{-1/2}\|\le\lambda^{-1/2}$.
Spectral calculus then gives $Hg=H^{1/2}h$ and $H^{1/2}g=h$ for centered
$h\in\operatorname{Dom}(H^{1/2})$.
The first identity proves $g\in\operatorname{Dom}(H)$ and the asserted
norm equality; the second, together with the identity for the Dirichlet form,
proves the gradient equality.
\end{proof}

\section{Polynomial Variance Estimates}\label{sec:polynomial_estimates}

We extend the quadratic-form estimate to polynomials of every degree whose lower-order derivatives have zero mean. The resulting variance bound depends only on the degree and the leading coefficient tensor, uniformly over all isotropic log-concave measures and dimensions. This supplies the initial polynomial estimates for the curvature comparison in the next section.

\begin{theorem}[Polynomial variance bound]\label{thm:polynomial}
Let $\mu$ be an isotropic log-concave probability measure on $\mathbb R^n$, let $d\ge1$ be an integer, and let $T$ be a fully symmetric $d$-tensor. Write $P_d^\mu[T]$ for the unique polynomial of degree at most $d$ such that
\begin{equation*}
 \E_\mu[D^jP_d^\mu[T]]=0,\ \forall 0\le j<d,\qquad
 D^dP_d^\mu[T]=d!\cdot T,
\end{equation*}
where $D^0P_d^\mu[T]:=P_d^\mu[T]$. Under Eq.~\eqref{tg:eq:quadratic},
\begin{equation*}
 \Var_\mu P_d^\mu[T]\le 1024^d\cdot(d!)^4\cdot\|T\|_{\mathrm{HS}}^2.
\end{equation*}
In particular, every polynomial $q$ of degree at most $s$ satisfies
\begin{equation*}
 \sqrt{\Var_\mu q}\le\sum_{k=1}^s32^k\cdot k!\cdot\|\E_\mu[D^kq]\|_{\mathrm{HS}}.
\end{equation*}
\end{theorem}

Using Lemma~\ref{lem:construction_localization_process}, we combine Letwin's quadratic inequality with bounds on the expected derivatives of the polynomial along localization. The posterior becomes strongly log-concave, so its Poincar\'e inequality converts these derivative estimates into a variance bound. Sections~\ref{sec:construction_localization_process} and~\ref{sec:polynomial_variance} establish the stochastic identities and derivative estimates needed for this argument.

\subsection{Construction of the localization process}\label{sec:construction_localization_process}

We construct a localization process for tracking polynomial derivative energies. Its evolution identities control these energies, while the Gaussian tilt creates the curvature needed to turn them into a variance bound.

\begin{lemma}[Construction of the localization process]
\label{lem:construction_localization_process}
Let $\mu_0$ be a compactly supported isotropic log-concave measure on $\mathbb R^n$. Write $\eta(x):=x$ for the coordinate map.
There is an adapted process of probability measures
\begin{equation*}
 \d\mu_t(x):=Z_t^{-1}\cdot\exp(\theta_t^\top x-x^\top \Lambda_tx/2)\d\mu_0(x)
\end{equation*}
with positive-definite covariance $\Gamma_t:=\operatorname{Cov}_{\mu_t}(\eta)$ at every finite time, mean $m_t:=\E_{\mu_t}[\eta]$, and
\begin{equation*}
 \theta_0=0,\qquad
 \d \theta_t=\Gamma_t^{-1}m_t\d t+\Gamma_t^{-1/2}\d \beta_t,\qquad
 \Lambda_t:=\int_0^t \Gamma_s^{-1}\d s,
\end{equation*}
where $\beta_t$ is a standard Brownian motion and $Z_t$ normalizes the measure.
For every bounded test function $q$, the posterior expectation is a martingale satisfying
\begin{equation*}
 \d\E_{\mu_t}[q]
 =\operatorname{Cov}_{\mu_t}(q,\Gamma_t^{-1/2}(\eta-m_t))\d \beta_t.
\end{equation*}
In particular,
\begin{equation*}
 \d m_t=\Gamma_t^{1/2}\d \beta_t,\qquad
 \d \Gamma_t=\sum_i U_{i,t}\d \beta_{i,t}-\Gamma_t\d t,
\end{equation*}
where
\begin{equation*}
 U_{i,t}:=\E_{\mu_t}[(\eta-m_t)(\eta-m_t)^\top(\Gamma_t^{-1/2}(\eta-m_t))_i].
\end{equation*}
For $t>0$, the posterior is $\Lambda_t$-strongly log-concave: its potential is the sum of $x^\top \Lambda_tx/2$ and an extended-valued convex function.
\end{lemma}

\begin{proof}
Choose $\rho\ge1$ such that the support of $\mu_0$ is contained in $\{x:|x|\le \rho\}$.
For $\theta\in\mathbb R^n$ and a real symmetric matrix $\Lambda$, define
\begin{equation*}
 r_{\theta,\Lambda}(x):=\exp(\theta^\top x-x^\top \Lambda x/2),\qquad
 Z(\theta,\Lambda):=\int r_{\theta,\Lambda}\d\mu_0,\text{ and} \qquad 
 \d\mu_{\theta,\Lambda}:=Z(\theta,\Lambda)^{-1}\cdot r_{\theta,\Lambda}\d\mu_0.
\end{equation*}
Let $m(\theta,\Lambda)$ and $\Gamma(\theta,\Lambda)$ denote the mean and covariance of this measure.
Compact support implies that $Z(\theta,\Lambda)>0$ and that these functions are smooth in all their parameters, by differentiation under the integral on compact parameter sets.
Moreover, $\mu_{\theta,\Lambda}$ and $\mu_0$ have the same null sets.
Since isotropy excludes support in an affine hyperplane, $\Gamma(\theta,\Lambda)$ is positive definite for every finite $(\theta,\Lambda)$.
Consequently the coefficients in
\begin{equation*}
 \d \theta_t=\Gamma(\theta_t,\Lambda_t)^{-1}m(\theta_t,\Lambda_t)\d t
       +\Gamma(\theta_t,\Lambda_t)^{-1/2}\d \beta_t,\qquad
 \d \Lambda_t=\Gamma(\theta_t,\Lambda_t)^{-1}\d t,\qquad (\theta_0,\Lambda_0)=(0,0)
\end{equation*}
are locally Lipschitz.
Lemma~\ref{lem:local_sde} gives a unique solution up to its maximal lifetime $\zeta$.
We first derive the identities before $\zeta$, and then prove that $\zeta=\infty$ almost surely.
All calculations before this conclusion are stopped on compact parameter sets.

Write $r_t:=r_{\theta_t,\Lambda_t}$, $Z_t:=Z(\theta_t,\Lambda_t)$, and $\mu_t:=\mu_{\theta_t,\Lambda_t}$.
It\^o's formula (Lemma~\ref{lem:ito_calculus}) gives
\begin{equation*}
 \frac{\d r_t(x)}{r_t(x)}
 =x^\top \Gamma_t^{-1}m_t\d t+x^\top \Gamma_t^{-1/2}\d \beta_t,
 \qquad
 \frac{\d Z_t}{Z_t}
 =m_t^\top \Gamma_t^{-1}m_t\d t+m_t^\top \Gamma_t^{-1/2}\d \beta_t.
\end{equation*}
Indeed, the quadratic variation of $\theta_t^\top x$ cancels the term $-x^\top\d \Lambda_tx/2$ in the first identity.
Applying It\^o's formula to the quotient cancels its drift and yields
\begin{equation*}
 \d(r_t(x)/Z_t)
 =(r_t(x)/Z_t)(x-m_t)^\top \Gamma_t^{-1/2}\d \beta_t.
\end{equation*}
Integration against $q\d\mu_0$ proves the asserted stochastic identity for every bounded measurable $q$.
On each compact parameter set, the integrands are uniformly bounded on the support, which justifies applying stochastic Fubini (Lemma~\ref{lem:stochastic_fubini}).
Taking $q$ to be the coordinates of $x$ gives $\d m_t=\Gamma_t^{1/2}\d \beta_t$.
Taking the coordinate products and subtracting $\d(m_tm_t^\top)$ gives
\begin{equation*}
 \d \Gamma_t=\sum_i U_{i,t}\d \beta_{i,t}-\Gamma_t\d t,
\end{equation*}
because $\d m_t\d m_t^\top=\Gamma_t\d t$ and the remaining noise is precisely the centered third moment in the statement.

We next exclude a finite lifetime, including covariance degeneration.
For $t<\zeta$, the matrix $\Lambda_t$ is positive semidefinite, so $\mu_t$ is log-concave.
Under this distribution set $\xi_t:=\Gamma_t^{-1/2}(\eta-m_t)$ and
\begin{equation*}
 S_{i,t}:=\E_{\mu_t}[\xi_{i,t}\xi_t\xi_t^\top],\qquad
 S_{z,t}:=\sum_i z_iS_{i,t},\ \forall z\in\mathbb R^n.
\end{equation*}
Thus $\xi_t$ is isotropic and $U_{i,t}=\Gamma_t^{1/2}S_{i,t}\Gamma_t^{1/2}$.
For every symmetric matrix $D$, Cauchy--Schwarz and Lemma~\ref{lem:quadratic_forms} give
\begin{equation*}
 |\langle D,S_{z,t}\rangle_{\mathrm{HS}}|
 =|\E_{\mu_t}[(z^\top \xi_t)(\xi_t^\top D\xi_t-\tr[D])]|\le\sqrt8\cdot|z|\cdot\|D\|_{\mathrm{HS}}.
\end{equation*}
Duality implies $\|S_{z,t}\|_{\mathrm{HS}}^2\le8|z|^2$.
The complete symmetry of the tensor $\E_{\mu_t}[\xi_{i,t}\xi_{j,t}\xi_{k,t}]$ also gives
$\sum_i|S_{i,t}z|^2=\|S_{z,t}\|_{\mathrm{HS}}^2$.
It follows that
\begin{equation*}
 \sum_i S_{i,t}^2\preceq8I,\qquad
 \sum_i\tr[S_{i,t}^2]\le8n,\qquad
 \sum_i(\tr[S_{i,t}])^2\le8n^2.
\end{equation*}
Applying It\^o's formula to $\log\det \Gamma_t$ gives
\begin{equation*}
 \log\det \Gamma_t
 =M_t-nt-\frac12\int_0^t\sum_i\tr[S_{i,s}^2]\d s
 \ge M_t-5nt,
 \qquad
 M_t:=\sum_i\int_0^t\tr[S_{i,s}]\d \beta_{i,s},
\end{equation*}
where $\Gamma_0=I$ and the quadratic variation of $M$ is at most $8n^2\cdot t$.
Extend the integrand of $M$ by zero after $\zeta$.
The resulting stochastic integral is a continuous square-integrable martingale on every finite time interval by Lemma~\ref{lem:ito_calculus}.
In particular, for each deterministic $T<\infty$, $\det \Gamma_t$ has a strictly positive lower bound for $0\le t<T\wedge\zeta$, possibly depending on the path.
The support bound gives $\Gamma_t\preceq \rho^2I$, and hence
\begin{equation*}
 \lambda_{\min}(\Gamma_t)\ge\frac{\det \Gamma_t}{(\rho^2)^{n-1}}.
\end{equation*}
Therefore $\Gamma_t^{-1}$ is bounded on the same interval.
Since $|m_t|\le \rho$, the drift of $\theta_t$ and the derivative of $\Lambda_t$ are bounded there, while the martingale part of $\theta_t$ has finite quadratic variation.
Localizing further on this bound and applying the martingale convergence statement in Lemma~\ref{lem:martingale_tools} show that $\theta_t$ and $\Lambda_t$ have finite limits as $t\uparrow\zeta$ on $\{\zeta\le T\}$.
At every such finite limiting parameter pair the covariance is positive definite and the coefficients are locally Lipschitz, so the local solution can be continued.
This contradicts maximality of $\zeta$.
Since $T$ was arbitrary, the process is defined for all $t\ge0$, and its covariance matrix remains positive definite at every finite time.

For bounded $q$, the process $\E_{\mu_t}[q]$ is a bounded local martingale, so it is a true martingale.
Finally, write $\d\mu_0=e^{-V_0}\d x$ with $V_0$ extended-valued and convex.
The posterior potential is
\begin{equation*}
 V_0(x)-\theta_t^\top x+\log Z_t+x^\top \Lambda_tx/2.
\end{equation*}
The first three terms form a convex function, and the matrix $\Lambda_t=\int_0^t \Gamma_s^{-1}\d s$ is positive definite for $t>0$.
This proves the strong log-concavity assertion.
\end{proof}

\subsection{Derivative estimates along localization}\label{sec:polynomial_variance}

We bound the growth of the expected derivatives of a polynomial in the evolving covariance metric. This estimate will let the induction on degree control the derivative energies over a short localization interval.

We use the notation of Lemma~\ref{lem:construction_localization_process} for a compactly supported isotropic initial measure.
Write $\E_t[\cdot]$ for posterior expectation and $\E_{\mathrm{loc}}[\cdot]$ for expectation over the process.
We suppress the time index in $S_i:=S_{i,t}$ and $U_i:=U_{i,t}$. With $\xi_t:=\Gamma_t^{-1/2}(\eta-m_t)$, we have $S_i=\E_t[\xi_{i,t}\xi_t\xi_t^\top]$ and $U_i=\Gamma_t^{1/2}S_i\Gamma_t^{1/2}$.
The proof of that lemma gives $\sum_iS_i^2\preceq8I$.
Since $2S_i\otimes S_i\preceq S_i^2\otimes I+I\otimes S_i^2$, congruence by $\Gamma_t^{1/2}\otimes \Gamma_t^{1/2}$ gives
\begin{equation}\label{eq:inverse_covariance_noise}
 \sum_iS_i^2\preceq8I,\qquad \sum_iU_i\otimes U_i\preceq8\Gamma_t\otimes \Gamma_t.
\end{equation}
These are inequalities on the full tensor spaces.

Fix a degree-$d$ polynomial $f$. For $1\le j\le d$, set
\begin{align*}
 h_j(t)&:=\E_t[D^jf],\\
 N_j(t)&:=\E_{\mathrm{loc}}[\langle h_j(t),\Gamma_t^{\otimes j}h_j(t)\rangle],\\
 L_j(t)&:=\E_{\mathrm{loc}}[\E_t[\langle D^jf-h_j(t),\Gamma_t^{\otimes j}(D^jf-h_j(t))\rangle]].
\end{align*}
Here $\E_{\mathrm{loc}}[\cdot]$ averages localization paths. All norms and products are on full finite tensor spaces. The energy of $D^jf$ in the covariance metric is $N_j+L_j$.

The next lemma controls the coupled evolution of the covariance metric and the conditional derivative mean.
The estimate includes their cross variation and applies to all tensor components simultaneously.

\begin{lemma}[Drift bound for the derivative moments]\label{lem:hierarchy}
For almost every $t$,
\begin{equation*}
 N_j'(t)\le4j^2\cdot N_j(t)+9j^2\cdot L_j(t).
\end{equation*}
\end{lemma}
\begin{proof}
Write $K:=\Gamma_t^{\otimes j}$. Its drift is at most
\begin{equation*}
 (-j+8\binom j2)K=(4j^2-5j)K
\end{equation*}
by Eq.~\eqref{eq:inverse_covariance_noise}. After congruence by $K^{-1/2}$, the $i$th noise coefficient is
$V_i:=\sum_{r=1}^j S_i^{(r)}$, where the superscript selects a tensor slot. The operator Cauchy--Schwarz inequality and Eq.~\eqref{eq:inverse_covariance_noise} give
\begin{equation*}
 \sum_iV_i^2\preceq8j^2I.
\end{equation*}
Indeed, the product rule for $\Gamma_t^{\otimes j}$ contributes $-K$ once from each slot and one term $U_i^{(r)}U_i^{(s)}$ for each pair $r<s$ and each noise coordinate $i$.
After congruence, the latter term is $S_i^{(r)}S_i^{(s)}$, whose sum over $i$ is at most $8I$.
Also $(\sum_rS_i^{(r)})^2\preceq j\sum_r(S_i^{(r)})^2$, which proves the noise bound without requiring different $S_i$ to commute.
Let $\d h_j=\sum_i c_i\d \beta_{i,t}$, where $c_i:=\E_t[(D^jf-h_j)\xi_{i,t}]$.
The weighted Bessel inequality in Lemma~\ref{lem:hilbert_norm_tools}, applied in the coordinate directions with weight $K$, gives
\begin{equation*}
 Q:=\sum_i\langle c_i,Kc_i\rangle\le
 \E_t[\langle D^jf-h_j,K(D^jf-h_j)\rangle].
\end{equation*}
The full drift of $\langle h,Kh\rangle$ is bounded by
\begin{equation*}
 (4j^2-5j)\cdot\langle h,Kh\rangle+Q+2\sum_i\langle K^{1/2}h,V_iK^{1/2}c_i\rangle.
\end{equation*}
The last term is at most
$\langle h,Kh\rangle+8j^2\cdot Q$.
Since $4j^2-5j+1\le4j^2$ and $1+8j^2\le9j^2$, taking expectations proves the lemma.
To justify taking expectations, we first stop the parameter process on compact sets.
The fixed spatial support bounds every derivative test and $\Gamma_t$; Bessel's inequality and Eq.~\eqref{eq:inverse_covariance_noise} bound the stochastic integrands.
The global existence in Lemma~\ref{lem:construction_localization_process} therefore permits removal of the stops.
In particular, the calculation includes the cross variation between the covariance and the conditional derivative mean.
\end{proof}

\subsection{Proof of the polynomial variance bound}

We prove the variance bound by induction on the degree, using the lower-degree bounds to close the derivative estimates. The Poincar\'e inequality supplied by the terminal measure's curvature then converts the resulting energy bound into the required variance estimate for the initial measure.

\begin{proof}[Proof of Theorem~\ref{thm:polynomial}]
Degree one follows from isotropy. Fix $d\ge2$, assume the assertion through degree $d-1$, and first take a compactly supported initial measure. Use Lemma~\ref{lem:construction_localization_process} for this measure. Apply Eq.~\eqref{eq:polynomial_minkowski} to $D^jf$ under each whitened posterior, and then apply Minkowski's inequality (Lemma~\ref{lem:hilbert_norm_tools}) in the localization probability space and the direct sum over tensor components. The chain rule inserts the needed covariance factor for each new derivative slot. Thus
\begin{equation}\label{eq:derivative_fluctuations}
 \sqrt{L_j(t)}\le\sum_{k=1}^{d-j}32^k\cdot k!\cdot\sqrt{N_{j+k}(t)},
 \forall 1\le j<d.
\end{equation}
The componentwise estimates introduce no factor from the number of tensor components and require no posterior centering of $D^jf$.

Now take $f:=P_d^{\mu_0}[T]$, and suppose $T\ne0$. Let
\begin{equation*}
 Q:=(d!)^2\cdot\|T\|_{\mathrm{HS}}^2,\qquad
 w_j:=1024^{d-j}\cdot((d-j)!)^2,\qquad
 M(t):=\max_{1\le j\le d}\frac{N_j(t)}{Q\cdot w_j}.
\end{equation*}
Initially $N_j(0)=0$ for $j<d$, and $N_d(0)=Q$. If $s:=d-j\ge1$, then Eq.~\eqref{eq:derivative_fluctuations} implies
\begin{equation*}
 \frac{\sqrt{L_j(t)}}{\sqrt{Q\cdot w_j}}
 \le\sqrt{M(t)}\cdot\sum_{k=1}^s\frac{k!\cdot(s-k)!}{s!}
 \le2\cdot\sqrt{M(t)}.
\end{equation*}
For $s=1$ the sum is one. For $s\ge2$, use $\binom sk\ge s$ on $1\le k<s$, and include the final term $k=s$. Hence
$L_j\le4Q\cdot w_j\cdot M$. Also $L_d=0$, because $D^df=d!\cdot T$ is constant.
Lemma~\ref{lem:hierarchy} now implies, for every $j$,
\begin{equation*}
 N_j(t)/(Q\cdot w_j)\le N_j(0)/(Q\cdot w_j)+40d^2\cdot\int_0^tM(s)\d s.
\end{equation*}
Taking the maximum and applying Gronwall's inequality (Lemma~\ref{lem:gronwall}) gives
\begin{equation*}
 M(t)\le e^{40d^2\cdot t}.
\end{equation*}
In particular, with $r:=\nabla f$,
\begin{equation*}
 G_1(t):=\E_{\mathrm{loc}}[\E_t[r^\top \Gamma_tr]]=N_1(t)+L_1(t)
 \le5\cdot(d!)^2\cdot\|T\|^2\cdot 1024^{d-1}\cdot((d-1)!)^2\cdot e^{40d^2\cdot t}.
\end{equation*}

The variance process $V(t):=\E_{\mathrm{loc}}[\Var_t f]$ obeys
\begin{equation*}
 V'(t)=-\E_{\mathrm{loc}}[|\operatorname{Cov}_t(f,\xi_t)|^2]\ge-V(t)
\end{equation*}
by Bessel's inequality under the whitened posterior. Hence $V(\tau)\ge e^{-\tau}\cdot\Var_{\mu_0}f$.
The terminal measure is $\Lambda_\tau$-strongly log-concave by Lemma~\ref{lem:construction_localization_process}, so Lemma~\ref{lem:matrix_strong_poincare}, applied with curvature matrix $\Lambda_\tau$, gives
\begin{equation*}
 \Var_\tau f\le\E_\tau[r^\top \Lambda_\tau^{-1}r],\qquad
 \Lambda_\tau^{-1}\preceq\tau^{-2}\cdot\int_0^\tau \Gamma_s\d s.
\end{equation*}
The matrix inequality follows, for each vector $v$, from
\begin{equation*}
 \int_0^\tau v^\top \Gamma_sv\d s-\tau^2v^\top \Lambda_\tau^{-1}v
 =\int_0^\tau|\Gamma_s^{1/2}v-\tau \Gamma_s^{-1/2}\Lambda_\tau^{-1}v|^2\d s.
\end{equation*}
For the fixed functions $r_ir_j$, the tower property gives
\begin{equation*}
 \E_{\mathrm{loc}}[\E_\tau[r^\top \Gamma_sr]]
 =\E_{\mathrm{loc}}[\E_s[r^\top \Gamma_sr]]=G_1(s).
\end{equation*}
Therefore
\begin{equation*}
 \Var_{\mu_0}f\le\frac{e^\tau}{\tau^2}\cdot\int_0^\tau G_1(s)\d s.
\end{equation*}
The conditional expectation identity applies with the random matrix $\Gamma_s$ held fixed at time $s$, so it also accounts for its dependence on the terminal function.

Choose $\tau:=1/(40d^2)$. Since $e^{1+\tau}<3$,
\begin{align*}
 \Var f&\le600d^2\cdot1024^{d-1}\cdot((d-1)!)^2\cdot(d!)^2\cdot\|T\|^2\\
 &=600\cdot1024^{d-1}\cdot(d!)^4\cdot\|T\|^2
 \le1024^d\cdot(d!)^4\cdot\|T\|^2.
\end{align*}
Degree one is simply isotropy. This closes the induction for compactly supported measures.

We remove the assumption of compact support by conditioning on growing balls, centering, and whitening. Moments through degree $2d$ converge. The coefficients of $\mathcal A_d^\mu$ are polynomials in the moments through degree $d$, so the varying Appell polynomials converge coefficientwise as well. Passing their variances to the limit proves Theorem~\ref{thm:polynomial} for all isotropic log-concave measures. The matrix inequality for compact nonsmooth posteriors is precisely Lemma~\ref{lem:matrix_strong_poincare}.

\end{proof}

\section{From Polynomial Estimates to Curvature}\label{sec:polynomial_curvature}

We turn polynomial variance bounds into an estimate for the full Poincar\'e constant in terms of a positive curvature lower bound. The comparison makes the dependence on the polynomial degree explicit, so that improved polynomial estimates give weaker dependence on inverse curvature. It supplies both the initial curvature estimate and the comparison used in the successive refinements of Section~\ref{sec:iterated_curvature}.

\begin{theorem}[Coefficient profiles and curvature]\label{thm:dyadic_curvature}
Let $0<\epsilon\le1$, let $\ell:\mathbb N\to[1,\infty)$ be nondecreasing, and let $R\ge2^{40}\cdot\epsilon^{-2}$.
Use the polynomial coefficients $c_k(\nu):=\sqrt{K_k(\nu)}/k!$ defined in Section~\ref{sec:preliminaries_appell}.
Suppose a centered regular log-concave probability measure $\d\nu=e^{-W}\d x$ with covariance at most $I$ satisfies, for every $k\ge1$,
\begin{equation}\label{eq:coefficient_profile}
 c_k(\nu)\le R^k\cdot\frac{\ell(k)^k}{(k+1)^2}.
\end{equation}
If $D^2W\succeq aI$ for $a>0$, then for every dyadic integer $d\ge2$,
\begin{equation}\label{eq:dyadic_curvature}
 C_P(\nu)\le16\cdot(1+\epsilon)\cdot R^2\cdot\ell(d)^2\cdot\max\{1,a^{-1/(d+1)}\}.
\end{equation}
\end{theorem}

The main difficulty is to control the mass lost when the families obtained by differentiation and the inverse operator are centered. Their energy decreases according to the Bochner identity, and polynomial testing bounds the discarded mass. We first establish the normalization identities, then estimate the errors in partial symmetrization. The proof in Section~\ref{tg:sec:iteration} combines these estimates across the polynomial degrees while controlling the total loss.

\subsection{Normalized families from the inverse operator}\label{tg:sec:families}

We construct families of derivatives whose mass and energy can be tracked throughout the curvature comparison. The normalization identifies the energy loss at each step, which will also control the symmetry errors in polynomial testing.

We allow centered regular measures with covariance $\Sigma\preceq I$ throughout the argument.
The coefficient bounds also hold for this class: whitening and Eq.~\eqref{eq:appell} give
$P_k^\mu[T](\Sigma^{1/2}y)=P_k^\rho[(\Sigma^{1/2})^{\otimes k}T](y)$ for the isotropic measure $\rho$.
The tensor map is a contraction, so every isotropic coefficient bound applies with the same constant.
For a finite family of functions, $\|\cdot\|_2^2$ denotes the sum of the squared $L^2(\mu)$ norms of its components; the operator and derivatives act componentwise.
Define $D_i:=P_+\partial_iH^{-1/2}$, $D:=(D_i)_i$, and $Lf:=\E[Xf]$.
The Dirichlet identity makes $D$ a contraction into the family direct sum.
Also $L$ is a contraction, since $\E[\langle z,X\rangle^2]=z^\top\Sigma z\le|z|^2$.
Integration by parts gives $\E[\nabla H^{-1/2}h]=LH^{1/2}h$ for centered $h$ in the form domain.
We normalize the inverse square root by one scalar shared by the entire family.
This preserves the family's squared $L^2$ norm and makes the polynomial testing identity exact.
Let $u$ be a nonzero centered finite family in the form domain. Define
\begin{equation}\label{dt_eq_inverse}
 v:=\|u\|_2^2,\quad b:=\|H^{-1/2}u\|_2^2,\quad
 \beta:=v/b,\quad h:=\sqrt\beta\cdot H^{-1/2}u,
\end{equation}
We further define
\begin{equation*}
 \chi:=\|H^{1/2}u\|_2^2-\beta v.
\end{equation*}
The scalar $\beta$ is shared by the whole family. The following lemma identifies the normalization error exactly and relates it to the energy decrease, which will also control the asymmetry between successive derivative indices.

\begin{lemma}[Exact normalization and defect]\label{dt_lem_inverse}
One has
\begin{equation}\label{dt_eq_inverseidentities}
 \|h\|_2^2=v,\qquad H^{1/2}h=\sqrt\beta\cdot u,\qquad
 \|H^{1/2}h\|_2^2=\beta v,\qquad \chi\ge0.
\end{equation}
For
\begin{equation*}
 z:=H^{-1/2}h-\beta^{-1/2}\cdot u
   =\sqrt\beta\cdot H^{-1}u-\beta^{-1/2}\cdot u,
\end{equation*}
there is the exact identity
\begin{equation}\label{dt_eq_defect}
 \|\nabla z\|_2^2=\frac\chi\beta.
\end{equation}
Moreover $\lambda\le\beta\le\|H^{1/2}u\|_2^2/v$.
The linear multiplier $\sqrt\beta\cdot H^{-1/2}$ has norm at most $\sqrt{\beta/\lambda}$.
\end{lemma}
\begin{proof}
The norm identity and the square-root identity are immediate. Cauchy--Schwarz gives
\begin{equation*}
 v^2\le\|H^{1/2}u\|_2^2\cdot\|H^{-1/2}u\|_2^2,
\end{equation*}
which proves $\chi\ge0$ and the upper bound on $\beta$.
The spectral inequality $b\le v/\lambda$ gives its lower bound. Spectral calculus gives
\begin{equation*}
 \|\nabla z\|_2^2
 =\beta b-2v+\beta^{-1}\cdot\|H^{1/2}u\|_2^2
 =\chi/\beta.
\end{equation*}
This calculation follows from the Dirichlet identity and spectral calculus. The multiplier bound follows from $H\succeq\lambda I$.
\end{proof}

Start $F_0$ at the centered unit first nonconstant eigenfunction from Lemma~\ref{lem:first_eigenfunction}. At generation $r$, form $u^r:=DF_r$ and apply Eq.~\eqref{dt_eq_inverse} to obtain $F_{r+1}$, with mean energy $\beta_{r+1}$ and energy drop $\chi_r$. Define
\begin{equation*}
 v_r:=\|F_r\|_2^2,\quad e_r:=\|H^{1/2}F_r\|_2^2,\quad
 \epsilon_r:=e_r-\lambda\cdot v_r,\quad
 p_r:=\|LH^{1/2}F_r\|^2.
\end{equation*}
Using the Bochner identity from Lemma~\ref{lem:operator_core_bochner} and the normalization identities from Lemma~\ref{dt_lem_inverse}, we obtain
\begin{align}
 v_{r+1}&=v_r-p_r,\notag\\
 e_{r+1}&\le e_r-a\cdot v_r-\chi_r,\notag\\
 0\le\epsilon_{r+1}&\le\epsilon_r+\lambda\cdot p_r-a\cdot v_r-\chi_r.\label{dt_eq_excess}
\end{align}
Moreover,
\begin{equation*}
 0\le p_r\le\min(v_r,e_r),\quad \chi_r\ge0,\quad
 v_r\le1,\quad e_r\le\lambda,\quad p_r\le\lambda,\quad a\le\lambda^2.
\end{equation*}
The last inequality follows from the first step of Eq.~\eqref{dt_eq_excess}. The uncentered gradient family has squared norm $v_r$ and mean squared norm $p_r$, which gives the identity for the squared norm.
Using the Bochner identity~\eqref{tg:eq:bochner} with $g:=H^{-1/2}F_r$, we bound the successor energy before normalization by
$\E[\|D^2g\|_{\mathrm{HS}}^2]\le e_r-a\cdot v_r$; normalization subtracts $\chi_r$.
Subtracting $\lambda v_{r+1}$ gives the excess recurrence.
The bounds $p_r\le v_r$ and $p_r\le e_r$ follow respectively from centering and contraction of $L$.
Before normalization at the first step, $\lambda\cdot(1-p_0)\le\|H^{1/2}u^0\|_2^2\le\lambda-a$, so $a\le\lambda\cdot p_0\le\lambda^2$.

These operations have the required domains.
If $F_r$ is in the form domain, then $H^{-1/2}F_r\in\operatorname{Dom}(H)$, and the domain statement in Lemma~\ref{lem:operator_core_bochner} implies that $u^r$ belongs to the form domain.
The next normalized family belongs to $\operatorname{Dom}(H)$, and the defect is in the form domain.
Thus the second weak derivatives used in the symmetry cancellations are well defined.
If $u^r=0$, we define $F_{r+1}$ and every subsequent family to be zero, and set $\beta_{r+1}:=\lambda$ and $\chi_r:=0$ at every zero step. The identities above remain valid. The finite stopping argument shows that this case cannot occur before the required horizon.

\subsection{A spectral estimate for the incidence operator}\label{tg:sec:tensor}

We prove a lower bound for the smallest singular value of the incidence operator on subsets. This is the combinatorial estimate needed to recover a tensor from a symmetrization of only some of its slots in the next subsection.

\begin{lemma}[Incidence-spectrum estimate]\label{lem:incidence_spectrum}
Let $N,s,q$ be integers with $0\le s\le q\le N-s$.
Let $\mathcal M_k$ be the Euclidean space of functions on the $k$-subsets
of $\{1,\ldots,N\}$, and define
\begin{equation*}
 (I_{q,s}f)(B):=\sum_{\substack{A\subset B\\|A|=s}}f(A).
\end{equation*}
Then
\begin{equation}\label{dt_eq_incidence}
 I_{q,s}^*I_{q,s}\succeq\binom{N-2s}{q-s}I.
\end{equation}
\end{lemma}
\begin{proof}
Let $U_k:\mathcal M_k\to\mathcal M_{k+1}$ sum over contained $k$-subsets, and let $D_{k+1}:=U_k^*$. Direct counting gives
\begin{equation*}
 D_{k+1}U_k-U_{k-1}D_k=(N-2k)I.
\end{equation*}
For a harmonic vector $v\in\ker D_j$, repeated use of this identity gives
\begin{equation*}
 D U^r v=r(N-2j-r+1)U^{r-1}v,\qquad
 \|U^rv\|^2=r!\cdot\frac{(N-2j)!}{(N-2j-r)!}\cdot\|v\|^2.
\end{equation*}
For $k<N/2$, the same commutator shows $U_k$ is injective; hence $D_{k+1}$ is surjective. Inductively, $\mathcal M_s$ is the orthogonal sum of the lifted harmonic spaces $U^{s-j}\ker D_j$, $0\le j\le s$. Orthogonality follows by pulling lowering operators to the other factor until a harmonic vector is paired with the range of an upper operator.

Since $I_{q,s}=U^{q-s}/(q-s)!$, the squared singular value on $U^{s-j}\ker D_j$ is
\begin{equation*}
 \binom{q-j}{s-j}\cdot\binom{N-s-j}{q-s},\ \forall 0\le j\le s,
\end{equation*}
which is a nonincreasing function of $j$, since the ratio of the $(j+1)$st value to the $j$th is
$\frac{s-j}{q-j}\cdot\frac{N-q-j}{N-s-j}\le1$. The last value is $\binom{N-2s}{q-s}$, proving Eq.~\eqref{dt_eq_incidence}.
\end{proof}

\subsection{Block recovery}

We bound a tensor with two symmetric blocks by its partial symmetrization, and include an error term when the second block is only approximately symmetric. This estimate lets polynomial tests at successively doubled degrees control the mass lost by centering the derivative families.
For an integer $q$, let $\mathsf P_q$ average permutations of the first $q$ slots. All permutations act orthogonally. Falling factorials are denoted by $(x)_{\underline s}$.
\begin{lemma}[Block recovery]\label{dt_lem_block}
Let $S\in\operatorname{Sym}^s(\mathbb R^n)\otimes\operatorname{Sym}^\ell(\mathbb R^n)$, where $s\le q\le\ell$. Then
\begin{equation}\label{dt_eq_block}
 \|\mathsf P_qS\|^2\ge
 \frac1{\binom qs}\cdot\frac{(s+\ell-q)_{\underline s}}{\ell_{\underline s}}\cdot\|S\|^2.
\end{equation}
If $T$ is symmetric in its first $s$ slots and $S$ is its symmetrization in the last $\ell$ slots, and $c$ is the reciprocal square root of the coefficient in Eq.~\eqref{dt_eq_block}, then
\begin{equation}\label{dt_eq_robustblock}
 \|T\|\le c\cdot\|\mathsf P_qT\|+(c+1)\cdot\|T-S\|.
\end{equation}
All assertions hold after adding finite direct sums over indices from earlier steps.
\end{lemma}
\begin{proof}
Define $N:=s+\ell$. The assumptions $s\le q\le\ell$ give
$0\le s\le q\le N-s$, so Lemma~\ref{lem:incidence_spectrum} applies.
For every $s$-subset $A$ of the slots, let $S_A$ be $S$ with its first symmetric block moved to $A$. This is independent of the internal ordering choices. Every $S_A$ has norm $\|S\|$. For a $q$-subset $B$, the sum $\sum_{A\subset B}S_A$ is $\binom qs$ times an orthogonal permutation of $\mathsf P_qS$. Apply Eq.~\eqref{dt_eq_incidence}, tensored with the identity on the tensor space, to this vector-valued list. It gives
\begin{equation*}
 \binom Nq\cdot\binom qs^2\cdot\|\mathsf P_qS\|^2
 \ge\binom{N-2s}{q-s}\cdot\binom Ns\cdot\|S\|^2.
\end{equation*}
Simplifying the factorials proves Eq.~\eqref{dt_eq_block}. The triangle inequality and contraction of $\mathsf P_q$ give Eq.~\eqref{dt_eq_robustblock}.
\end{proof}

The error in full symmetrization can be bounded by the errors in adjacent swaps. If $s_j$ swaps two adjacent designated slots of an $\ell$-slot array $Y$, then
\begin{equation}\label{dt_eq_sort}
 \|Y-\operatorname{Sym}_\ell Y\|\le\ell\cdot\sum_{j=1}^{\ell-1}\|Y-s_jY\|.
\end{equation}
Insertion sort expresses each permutation with each adjacent generator used at most $\ell$ times. Telescope the resulting product of isometries and average. We apply this estimate within each dyadic stage and count the iterations through which each error from an adjacent swap is propagated.

\subsection{Bounds on centering errors}\label{sec:centering_loss}

We bound the mass removed by centering in terms of polynomial coefficients and normalization errors. Block recovery at successively doubled degrees provides the estimate while keeping track of how each symmetry error propagates.

At degree $k$, we use a block of $L_0\cdot k$ derivative slots, where the integer $L_0\ge3$ will be chosen later.
All measures in Section~\ref{sec:centering_loss} are regular, centered, and have covariance at most $I$.
Write $c_k:=\sqrt{K_k(\mu)}/k!$ and $Q_kh:=\E[\mathcal A_k^\mu h]$.
The operator norm of $Q_k$ is at most $\sqrt{K_k(\mu)}$, by duality on the symmetric coefficient tensor space.
These bounds also hold componentwise in every finite direct sum over indices from earlier steps.

We first record the exact polynomial testing identity. Appell differentiation and integration by parts give, for $j\ge1$,
\begin{equation}\label{eq:dyadic_testing}
 (k+1)\mathsf P_{k+1}(Q_kDF_j)
 =Q_{k+1}H^{1/2}F_j
 =\sqrt{\beta_j}\cdot Q_{k+1}u^{j-1}.
\end{equation}
Indeed, pair with an arbitrary symmetric $(k+1)$-tensor and use
$\nabla P_{k+1}[T]=(k+1)\cdot P_k[T]$, with the derivative index left free.
The Dirichlet form identity in Lemma~\ref{lem:operator_core_bochner} gives the first equality.
Centering in $D$ contributes zero because $\E[\mathcal A_k^\mu]=0$.
The second equality is Eq.~\eqref{dt_eq_inverseidentities}.

Suppose locally that the relevant normalizers obey $\beta_j\le b\cdot\lambda$, where $b\ge1$.
For the defect $z^j:=H^{-1/2}F_{j+1}-\beta_{j+1}^{-1/2}\cdot u^j$, Eq.~\eqref{dt_eq_defect} gives
$\|\nabla z^j\|_2^2=\chi_j/\beta_{j+1}$.
Before normalization, the successor decomposes as
\begin{equation*}
 u^{j+1}=\beta_{j+1}^{-1/2}\cdot P_+\nabla u^j+P_+\nabla z^j.
\end{equation*}
Since $u^j=P_+\nabla H^{-1/2}F_j$, the first term is an ordinary weak Hessian and is symmetric in its two newest derivative slots.
Consequently the newest adjacent swap has norm at most $2\cdot\sqrt{\chi_j/\lambda}$.
Each subsequent map $D\sqrt{\beta}H^{-1/2}$ has norm at most $\sqrt b$.
It commutes with permutations of slots introduced at earlier steps, using the fixed scalar of the original complete family.
Thus a swap of derivative slots $a,a+1$, newest first, in $u^J$ satisfies
\begin{equation}\label{eq:dyadic_swap}
 \|u^J-s_au^J\|_2\le2\cdot b^{(a-1)/2}\cdot\sqrt{\chi_{J-a}/\lambda}.
\end{equation}
There are exactly $a-1$ later maps after this swap is created.
No old error is differentiated, and no derivative is commuted with an inverse power of $H$.

The next lemma combines polynomial testing with partial symmetrization to bound the norm removed by centering at each step. The bound consists of a contribution from the degree-$d$ polynomial estimate and contributions from the normalization errors at earlier steps.

\begin{lemma}[Dyadic local loss]\label{lem:dyadic_loss}
Let $L_0\ge3$ be an integer, let $b\ge1$, and define
\begin{equation}\label{eq:dyadic_block_constants}
 B:=2\cdot\sqrt{\frac{L_0}{L_0-2}},
 A:=b\cdot B^2, J:=B^4\cdot b^{L_0+1}.
\end{equation}
Let $d\ge2$ be a power of two, and suppose all relevant normalizers satisfy $\beta_j\le b\cdot\lambda$.
For $j\ge m_d:=(L_0+1)\cdot d/2-1$,
\begin{equation}\label{eq:dyadic_local_loss}
 \begin{aligned}
 \sqrt{p_j}&\le(A\cdot\lambda)^{d/2}\cdot c_d
 +\sum_{\substack{1\le k<d\\ k\text{ a power of two}}}
 t_k\cdot\sum_{a=1}^{L_0\cdot k-1}\sqrt{\chi_{j-k-a}},\\
 t_k&:=\frac{4L_0\cdot k}{B}\cdot J^{k/2}\cdot c_k\cdot\lambda^{(k-1)/2}.
 \end{aligned}
\end{equation}
The coefficients $c_k$ may be replaced throughout by any upper bounds.
\end{lemma}
\begin{proof}
For $k=1,2,4,\ldots,d$, define $T_k:=Q_ku^{j-k}$.
It has $k$ symmetric coordinate slots and $j-k+1$ derivative slots.
For $k<d$, the hypothesis on $j$ gives $j\ge(L_0+1)\cdot k-1$, so at least $L_0\cdot k$ derivative slots are available.
Let $S_k$ symmetrize the newest $L_0\cdot k$ derivative slots of $T_k$.
Lemma~\ref{dt_lem_block}, with coordinate block size $k$, derivative block size $L_0\cdot k$, and $q:=2k$, gives
\begin{equation*}
 \|T_k\|\le C_k\cdot\|\mathsf P_{2k}T_k\|+(C_k+1)\cdot\|T_k-S_k\|,
\end{equation*}
where
\begin{equation*}
 C_k^2:=\binom{2k}{k}\cdot\frac{(L_0\cdot k)_{\underline k}}{((L_0-1)\cdot k)_{\underline k}}
 \le4^k\cdot(\frac{L_0}{L_0-2})^k=B^{2k}.
\end{equation*}
Indeed, every factor in the falling-factorial ratio is at most $L_0/(L_0-2)$.
Repeated use of Eq.~\eqref{eq:dyadic_testing} gives
\begin{equation*}
 \|\mathsf P_{2k}T_k\|
 \le(b\cdot\lambda)^{k/2}\cdot\frac{k!}{(2k)!}\cdot\|T_{2k}\|.
\end{equation*}
The outer symmetrization absorbs every smaller symmetrization in this iteration.
Also, Eq.~\eqref{dt_eq_sort} and Eq.~\eqref{eq:dyadic_swap} imply
\begin{equation*}
 E_k:=\|T_k-S_k\|
 \le2L_0\cdot k\cdot b^{L_0\cdot k/2}\cdot\sqrt{K_k/\lambda}\cdot
 \sum_{a=1}^{L_0\cdot k-1}\sqrt{\chi_{j-k-a}}.
\end{equation*}
We used $b^{(a-1)/2}\le b^{L_0\cdot k/2}$ throughout the designated block.
Therefore
\begin{equation*}
 \|T_k\|\le B^k\cdot(b\cdot\lambda)^{k/2}\cdot\frac{k!}{(2k)!}\cdot\|T_{2k}\|
 +2B^k\cdot E_k.
\end{equation*}

The starting factor is $\sqrt{p_j}=\sqrt{\beta_j}\cdot\|T_1\|\le\sqrt{b\cdot\lambda}\cdot\|T_1\|$.
The preceding dyadic degrees sum to $k-1$, while the factorial ratios telescope.
Including the starting factor, the coefficient before stage $k$ is consequently at most
\begin{equation}\label{eq:dyadic_prefix}
 \frac{(b\cdot\lambda)^{k/2}\cdot B^{k-1}}{k!}.
\end{equation}
At degree $d$, $\|T_d\|\le\sqrt{K_d}$ because $\|u^{j-d}\|_2\le1$.
Eq.~\eqref{eq:dyadic_prefix} gives the terminal term $(A\cdot\lambda)^{d/2}\cdot c_d/B$, which is at most the one in Eq.~\eqref{eq:dyadic_local_loss}.
For an error at stage $k$, multiplying Eq.~\eqref{eq:dyadic_prefix} by $2B^k\cdot E_k$ gives the coefficient
\begin{equation*}
 \frac{(b\cdot\lambda)^{k/2}\cdot B^{k-1}}{k!}\cdot
 (2B^k)\cdot(2L_0\cdot k\cdot b^{L_0\cdot k/2})\cdot\sqrt{K_k/\lambda}
 =\frac{4L_0\cdot k}{B}\cdot J^{k/2}\cdot c_k\cdot\lambda^{(k-1)/2}.
\end{equation*}
All defect indices are nonnegative since $j\ge(L_0+1)\cdot k-1$.
This proves the lemma.
\end{proof}

\subsection{Proof of the curvature comparison}\label{tg:sec:iteration}

We prove the curvature comparison by showing that a smaller spectral gap would force more energy loss than the initial family contains. Summing the centering and normalization errors across dyadic degrees ensures that enough mass remains for this contradiction.

\begin{proof}[Proof of Theorem~\ref{thm:dyadic_curvature}]
Define $C:=16\cdot(1+\epsilon)$, and choose
\begin{equation*}
 L_0:=\lceil64/\epsilon\rceil+2,
 p_*:=\frac{\epsilon}{64\cdot(L_0+1)},
 b:=(1-p_*)^{-1}.
\end{equation*}
Use $B,A,J$ from Eq.~\eqref{eq:dyadic_block_constants}.
These choices give $L_0\le67/\epsilon$ and $p_*\ge\epsilon^2/4352$.
Moreover,
\begin{equation*}
 \log(J/16)
 =2\log\frac{L_0}{L_0-2}+(L_0+1)\cdot\log\frac1{1-p_*}
 \le\frac{\epsilon}{16}+\frac{\epsilon}{32}
 =\frac{3\epsilon}{32}.
\end{equation*}
Here $\log(1+x)\le x$ and $-\log(1-p_*)\le2p_*$ were used.
Since $e^{3\epsilon/32}\le1+\epsilon/4$, it follows that
\begin{equation*}
 q:=J/C\le1-3\epsilon/8,
 \sqrt q\le1-3\epsilon/16, A/C\le1/3.
\end{equation*}
The bound on $A/C$ also uses $L_0\ge66$ and $p_*\le1/4096$.
In particular $A<J<C$.

Fix $d$ and first suppose $\lambda\le1/(C\cdot R^2\cdot\ell(d)^2)$.
Define
\begin{equation*}
 P_N:=\sum_{j<N}p_j, X_N:=\sum_{j<N}\chi_j, V_N:=\sum_{j<N}v_j.
\end{equation*}
The exact identity for the squared norm and the recurrence for the energy excess give
\begin{equation}\label{eq:dyadic_budget}
 v_N=1-P_N, a\cdot V_N+X_N\le\lambda\cdot P_N.
\end{equation}
For each $j\ge L_0=m_2$, choose the largest dyadic $D_j\le d$ with $m_{D_j}\le j$.
Before degree $d$ becomes available, degree $k<d$ is used for exactly $(L_0+1)\cdot k/2$ generations.
The first $L_0$ losses are at most $\lambda$ each.
On a prefix where $\beta_j\le b\cdot\lambda$, apply Lemma~\ref{lem:dyadic_loss} at degree $D_j$.

One nonnegative lag kernel dominates the errors at all these degrees:
\begin{equation*}
 w_s:=\sum_{\substack{k<d\\k\text{ dyadic}}}t_k\cdot\mathbf1_{\{k+1\le s\le(L_0+1)\cdot k-1\}},
  W:=\sum_s w_s.
\end{equation*}
The number of lags contributed by degree $k$ is at most $L_0\cdot k$, so
\begin{equation}\label{eq:dyadic_kernel}
 \sqrt\lambda\cdot W
 \le\frac{4L_0^2}{B}\cdot\sum_{\substack{k<d\\k\text{ dyadic}}}k^2\cdot c_k\cdot(J\cdot\lambda)^{k/2}.
\end{equation}
Define the energy drops with negative indices to be zero.
Weighted Cauchy--Schwarz and summation give
\begin{equation*}
 \sum_{j<N}(\sum_s w_s\cdot\sqrt{\chi_{j-s}})^2\le W^2\cdot X_N.
\end{equation*}
Thus overlapping error intervals introduce no factor depending on the degree or the number of stages.
Squaring Eq.~\eqref{eq:dyadic_local_loss} and summing gives, on every justified prefix,
\begin{equation}\label{eq:dyadic_prefix_sum}
 P_N\le\delta_d+N\cdot B_d\cdot\lambda^d+2W^2\cdot X_N,
\end{equation}
where
\begin{equation*}
 B_d:=2A^d\cdot c_d^2,
 \delta_d:=L_0\cdot\lambda+(L_0+1)\cdot
 \sum_{\substack{2\le k<d\\k\text{ dyadic}}}k\cdot(A\cdot\lambda)^k\cdot c_k^2.
\end{equation*}
This also holds for $N<m_d$, since the term at the terminal degree then only adds a nonnegative quantity.

We next bound the kernel and initial terms uniformly in $d$.
Theorem~\ref{thm:polynomial} and the affine inheritance in Section~\ref{tg:sec:families} give, for every $k\ge1$,
\begin{equation}\label{eq:dyadic_initial_coefficients}
 c_k\le32^k\cdot k!.
\end{equation}
Set
\begin{equation*}
 k_0:=\lceil2^{10}\cdot\epsilon^{-1}\cdot\log(2/\epsilon)\rceil,
  C_0:=4L_0^2/B.
\end{equation*}
Then $k_0\le2^{12}\cdot\epsilon^{-2}$, $C_0\le2^{14}\cdot\epsilon^{-2}$, and $R\ge64k_0$.
For $k\le k_0$, Eq.~\eqref{eq:dyadic_initial_coefficients}, $k!\le k^k$, and $32k/R\le1/2$ give
\begin{equation*}
 C_0\cdot\sum_{k\le k_0}k^2\cdot c_k\cdot(J\cdot\lambda)^{k/2}
 \le\frac{32C_0}{R}\cdot\sum_{k\ge1}k^3\cdot 2^{-k+1}
 =\frac{1664C_0}{R}\le2^{-15}.
\end{equation*}
We used $J\cdot\lambda\le R^{-2}$.
For $k>k_0$, Eq.~\eqref{eq:coefficient_profile} and $\ell(k)\le\ell(d)$ for the degrees that occur in the kernel give the bound
\begin{equation*}
 C_0\cdot\sum_{k>k_0}q^{k/2}
 \le\frac{16C_0}{3\epsilon}\cdot e^{-3\epsilon k_0/16}
 \le2^{17}\cdot\epsilon^{-3}\cdot(\epsilon/2)^{192}<\frac18.
\end{equation*}
The sums over all positive integers bound the corresponding dyadic sums.
Together with Eq.~\eqref{eq:dyadic_kernel}, these estimates show
\begin{equation}\label{eq:dyadic_absorption}
 \theta:=2\lambda\cdot W^2\le\frac18.
\end{equation}

The same split bounds the initial terms.
For $2\le k\le k_0$, use Eq.~\eqref{eq:dyadic_initial_coefficients} and $A\cdot\lambda\le R^{-2}$ to obtain
\begin{equation*}
 \sum_{2\le k\le k_0}k\cdot(A\cdot\lambda)^k\cdot c_k^2
 \le(32/R)^4\cdot\sum_{k\ge2}k^5\cdot 4^{-k+2}
 \le2^{36}/R^4.
\end{equation*}
The last numerical series is less than $2^{16}$, as follows by differentiating the geometric series five times.
The choice $R\ge2^{40}\cdot\epsilon^{-2}$ therefore gives
\begin{equation*}
 \frac{L_0}{C\cdot R^2}+\frac{(L_0+1)\cdot 2^{36}}{R^4}\le\frac{p_*}{32}.
\end{equation*}
For the remaining degrees, Eq.~\eqref{eq:coefficient_profile} and $A/C\le1/3$ give
\begin{equation*}
 (L_0+1)\cdot\sum_{k>k_0}\frac{k\cdot(A/C)^k}{(k+1)^4}
 \le(L_0+1)\cdot 2^{-k_0}
 \le68\epsilon^{-1}\cdot(\epsilon/2)^{100}\le\frac{p_*}{32}.
\end{equation*}
Consequently
\begin{equation}\label{eq:dyadic_startup}
 \delta_d\le p_*/16.
\end{equation}
The same parameter bounds imply $\lambda\le p_*/2$ and $R^2\ge64/(p_*\cdot C)$.
All these estimates are independent of the terminal degree.

Suppose for contradiction that $a\ge32B_d\cdot\lambda^{d+1}/p_*$, and define $M:=\lceil\lambda/a\rceil$.
Stop at the first $N\le M$ with $P_N>p_*/2$, if one exists.
Since $p_j\le\lambda\le p_*/2$, through a possible exit $P_j\le p_*$.
Hence $v_j\ge1-p_*$ and $\beta_j=e_j/v_j\le b\cdot\lambda$ for $j\ge1$.
The successor before normalization is nonzero because its squared $L^2$ norm is $v_j-p_j\ge1-p_*-\lambda>0$.
These inequalities justify every local assumption on the entire stopped prefix before it is used.

Combining Eq.~\eqref{eq:dyadic_budget}, Eq.~\eqref{eq:dyadic_prefix_sum}, and Eq.~\eqref{eq:dyadic_absorption} gives
\begin{equation*}
 (1-\theta)\cdot P_N\le\delta_d+N\cdot B_d\cdot\lambda^d.
\end{equation*}
For $N\le M$, the contradictory assumption and $a\le\lambda^2$ imply
\begin{equation*}
 N\cdot B_d\cdot\lambda^d\le\frac{p_*}{32}\cdot(1+a/\lambda)
 \le\frac{p_*}{32}\cdot(1+\lambda)\le\frac{p_*}{16}.
\end{equation*}
Eq.~\eqref{eq:dyadic_startup} now yields $P_N<p_*/6$, contradicting a first exit.
No exit occurs through $M$, and Eq.~\eqref{eq:dyadic_budget} gives the contradiction
\begin{equation*}
 a\cdot V_M\ge(1-p_*/6)\cdot aM\ge(1-p_*/6)\cdot\lambda,
  a\cdot V_M\le\lambda\cdot P_M<p_*\cdot\lambda/6.
\end{equation*}
Thus, under the assumed bound on the spectral gap,
\begin{equation*}
 a\cdot C_P(\nu)^{d+1}<\frac{32}{p_*}\cdot B_d
 =\frac{64}{p_*}\cdot A^d\cdot c_d^2.
\end{equation*}
Using Eq.~\eqref{eq:coefficient_profile}, $A<C$, $\ell(d)\ge1$, and $R^2\ge64/(p_*\cdot C)$, the right side is at most
\begin{equation*}
 \frac{64}{p_*}\cdot\frac{A^d\cdot R^{2d}\cdot\ell(d)^{2d}}{(d+1)^4}
 \le(C\cdot R^2\cdot\ell(d)^2)^{d+1}.
\end{equation*}
Taking a $(d+1)$st root proves Eq.~\eqref{eq:dyadic_curvature} in this regime.
In the complementary regime, $C_P(\nu)<C\cdot R^2\cdot\ell(d)^2$ directly.
This completes the proof.
\end{proof}

\section{Iterative Logarithmic Curvature Refinement}\label{sec:iterated_curvature}

We repeatedly convert a curvature bound into improved polynomial coefficients and use these coefficients to add another logarithm to the curvature bound. The main task is to control the accumulated cost as the number of iterations grows. We prove that the squared amplitude grows only as $(r+1)^{1/3}$, uniformly in the measure and the curvature lower bound. This control will allow the iteration depth to be chosen in terms of the dimension in the next section.

\begin{proposition}[Iterative logarithmic curvature refinement]\label{prop:sharp_iterated_curvature}
For $x\ge0$ and every integer $r\ge1$, let $\ell_r(x)$ denote the $r$-fold iterate of $x\mapsto\log(e+x)$.
There is a universal constant $C_0$ and a sequence of universal constants $(\Gamma_r)_{r\ge1}$ such that every centered regular measure $\nu$ with covariance at most $I$ and curvature at least $aI$, where $a>0$, satisfies
\begin{equation}\label{eq:sharp_iterated_curvature}
 C_P(\nu)\le\Gamma_r^2 \cdot \ell_r(a^{-1})^2,
 \qquad \Gamma_r^2\le C_0 \cdot (r+1)^{1/3}
 ,\ \forall r\ge1.
\end{equation}
\end{proposition}

The recurrence is proved for auxiliary amplitudes $\widetilde\Gamma_r$.
For a fixed universal depth $r_0$ and universal constants $G,C_T>0$, it has the form
\begin{equation*}
 \begin{aligned}
 B_{r_0}&:=1,\qquad
 B_{r+1}:=B_r \cdot (1+r^{-2})^2 \cdot \zeta_r^2,\\
 \zeta_r&:=\exp \bigl(2^{-r} \cdot \log(2 \cdot C_T \cdot (r+1)^2)\bigr),\\
 \widetilde\Gamma_r^2&:=G^2 \cdot (r+1)^{1/3} \cdot B_r
 ,\ \forall r\ge r_0.
 \end{aligned}
\end{equation*}
The factors defining $B_r$ have a bounded product. The comparison with the Poincar\'e constant then gives $\Gamma_r^2:=\widetilde\Gamma_r^2+2$; the finitely many depths below $r_0$ are covered by fixed initial estimates. Thus $r$ counts the successive logarithmic improvements, with one sequence of amplitudes valid for every measure and every curvature lower bound. The coefficient estimates and the control of centering losses establish this recurrence in the proof of Proposition~\ref{prop:sharp_iterated_curvature} at the end of the section.

We use the notation
\begin{equation*}
 g(x):=\log(e+x),\qquad \ell_0(x):=x,\qquad
 \ell_r(x):=g^{\circ r}(x),\ \forall x\ge0,\ r\ge1.
\end{equation*}
A curvature profile $F$ will mean an upper bound $C_P(\nu)\le F(a)$ valid for every centered regular probability measure $\nu$ with covariance at most $I$ and potential Hessian at least $aI$.

\subsection{From curvature bounds to polynomial coefficients}\label{sec:curvature_feedback}

We first transfer a curvature bound to polynomial coefficient estimates whose growth is measured by the degree. These estimates provide the input for the next curvature comparison. The affine extension below allows the given curvature bound to be applied to each normalized posterior along localization.

\begin{lemma}[Affine extension of a curvature bound]\label{lem:profile_inflation}
Suppose a continuous function $F:(0,\infty)\to(0,\infty)$ bounds $C_P(\nu)$ by $F(a)$ for every centered regular measure with covariance at most $I$ and curvature at least $aI$.
Let $\nu$ have covariance $A\succ0$ and density proportional to $\exp(-x^\top Bx/2-V(x))$, where $B\succ0$ and $V$ is extended-valued and convex.
For every $\delta>0$ and every locally Lipschitz test function $f$ of finite energy,
\begin{equation}\label{eq:profile_inflation}
 \Var_\nu[f]\le F(\delta) \cdot \E_\nu[(\nabla f)^\top(A+\delta B^{-1})\nabla f].
\end{equation}
\end{lemma}
\begin{proof}
First extend the profile to a possibly nonsmooth centered measure with covariance at most $I$ and curvature at least $aI$.
Convolution with $N(0,\epsilon I)$ has covariance at most $(1+\epsilon)I$ and smooth potential Hessian between $a/(1+a\epsilon)I$ and $\epsilon^{-1}I$.
For the lower bound, the Hessian formula for Gaussian convolution in Lemma~\ref{lem:regular_isotropic_approximation} expresses the Hessian as
$\epsilon^{-1}I-\epsilon^{-2}\operatorname{Cov}(X\mid X+\sqrt\epsilon Z=y)$.
The conditional distribution has curvature at least $(a+\epsilon^{-1})I$.
Lemma~\ref{lem:matrix_strong_poincare}, applied to linear tests, bounds its covariance by $(a+\epsilon^{-1})^{-1}I$ and proves the asserted Hessian bound.
Divide the convolved vector by $\sqrt{1+\epsilon}$.
The resulting regular measure has covariance at most $I$ and curvature at least
$a_\epsilon:=a(1+\epsilon)/(1+a\epsilon)$, which tends to $a$.
Continuity of $F$ passes its Poincar\'e inequality to compactly supported smooth tests in the limit.
The truncation and approximation argument in Lemma~\ref{lem:poincare_limit} extends it to every locally Lipschitz test of finite energy.
The same limit works for every fixed polynomial degree, since strong log-concavity gives moment convergence of all fixed orders under this convolution and rescaling.

We define $M:=A+\delta B^{-1}$ and $Y:=M^{-1/2}(X-\E[X])$.
Then $\operatorname{Cov}(Y)=M^{-1/2}AM^{-1/2}\preceq I$.
Since $M\succeq\delta B^{-1}$, inversion gives $B\succeq\delta M^{-1}$, hence $M^{1/2}BM^{1/2}\succeq\delta I$.
Apply the extended profile to the distribution of $Y$ and the test $y\mapsto f(\E[X]+M^{1/2}y)$.
Its gradient energy is the right side of Eq.~\eqref{eq:profile_inflation} divided by $F(\delta)$.
This argument applies to arbitrary positive definite matrices, including noncommuting ones.
\end{proof}

We apply this extension to profiles defined by iterated logarithms.
The coefficient induction below uses the localization time
$\tau=(2308\cdot d^2)^{-1}$, where $2308=4+9\cdot16^2$ comes from
Lemma~\ref{lem:hierarchy}, the convolution bound in
Eq.~\eqref{eq:profile_convolution}, and Gronwall's inequality
(Lemma~\ref{lem:gronwall}).
For every $r\ge1$, we use monotonicity and concavity of the iterated
logarithms to control the resulting change of scale from $d$ to $\tau^{-1}$:
\begin{equation}\label{eq:iterated_log_scaling}
 \ell_r(x)\ge1,\qquad \ell_r(cx)\le c \cdot \ell_r(x),\ \forall c\ge1,\qquad
 \ell_r(2308 \cdot d^2)\le10 \cdot \ell_r(d),\ \forall d\ge1.
\end{equation}
The scaling inequality holds for $g$ since it is concave and $g(0)=1$, and then follows by composition.
For the last inequality first use $g(2308 \cdot d^2)\le\log(2309)+2 \cdot g(d)\le10 \cdot g(d)$, and apply the scaling inequality to the remaining compositions.
In particular, the factor ten is independent of $r$.

We also need a contraction that improves with the number of compositions.
For $x>0$,
\begin{equation*}
 \frac{x \cdot g'(x)}{g(x)}=\frac{x}{(e+x) \cdot \log(e+x)}\le\frac12.
\end{equation*}
Indeed, $t \cdot \log t-2t+2e$ is positive and increasing for $t\ge e$.
Integrating the logarithmic derivative and composing gives
\begin{equation}\label{eq:iterated_log_contraction}
 \log\frac{\ell_r(u)}{\ell_r(v)}\le2^{-r} \cdot \log\frac uv
 ,\ \forall u\ge v>0,\ r\ge0.
\end{equation}

The coefficient estimate below replaces the curvature variable in an $r$-fold logarithmic bound by the polynomial degree. Theorem~\ref{thm:dyadic_curvature} then adds one more logarithm to the curvature bound.

\begin{theorem}[Polynomial coefficients from curvature bounds]\label{thm:profile_coefficients}
Suppose $r\ge1$, $\Gamma\ge1$, and every centered regular measure with covariance at most $I$ and curvature at least $aI$ satisfies
$C_P\le\Gamma^2 \cdot \ell_r(a^{-1})^2$.
Define $R:=2^{12} \cdot \Gamma$.
Then for every isotropic log-concave measure, every integer $d\ge1$, and every symmetric coefficient tensor $T$,
\begin{equation}\label{eq:profile_coefficients}
 \Var_\mu[P_d^\mu[T]]\le R^{2d} \cdot (d!)^2 \cdot \frac{\ell_r(d)^{2d}}{(d+1)^4} \cdot \|T\|_{\mathrm{HS}}^2.
\end{equation}
The same assertion holds for centered measures with covariance at most $I$.
\end{theorem}
\begin{proof}
Write $\ell:=\ell_r$, $b_0:=1$, and $b_s:=\ell(s)^s/(s+1)^2$ for $s\ge1$.
Monotonicity gives
\begin{equation}\label{eq:profile_convolution}
 \sum_{k=1}^s\frac{b_k \cdot b_{s-k}}{b_s}
 \le\sum_{k=1}^s\frac{(s+1)^2}{(k+1)^2 \cdot (s-k+1)^2}\le16,
 \qquad \frac{b_d}{b_{d-1}}\ge\frac{\ell(d)}4.
\end{equation}
For the sum, split at $s/2$; the factor involving the larger index is at most four, and the reciprocal squares involving the smaller index are summable.
For the ratio, use $\ell(d)^{d-1}\ge\ell(d-1)^{d-1}$ and $d^2/(d+1)^2\ge1/4$.

We prove $c_d\le R^d \cdot b_d$ by induction on $d$.
Degree one follows from isotropy and $R\ge4$.
Assume the assertion for every isotropic measure at smaller degrees and first take a compactly supported isotropic initial measure $\mu_0$.
Use Lemma~\ref{lem:construction_localization_process} and the hierarchy of Section~\ref{sec:polynomial_variance}.
Applying Eq.~\eqref{eq:polynomial_minkowski} to $D^jf$ under each whitened posterior gives
\begin{equation*}
 \sqrt{L_j(t)}\le\sum_{k=1}^{d-j}R^k \cdot b_k \cdot \sqrt{N_{j+k}(t)}.
\end{equation*}
As in Eq.~\eqref{eq:derivative_fluctuations}, this is an inequality on the full output tensor direct sum, with the covariance factor inserted in every new derivative slot.
Take $f:=P_d^{\mu_0}[T]$, with $T\ne0$, and define
\begin{equation*}
 Q:=(d!)^2 \cdot \|T\|_{\mathrm{HS}}^2,\qquad
 w_s:=R^{2s} \cdot b_s^2,\qquad
 M(t):=\max_{1\le j\le d}\frac{N_j(t)}{Q \cdot w_{d-j}}.
\end{equation*}
Eq.~\eqref{eq:appell_centering} gives $M(0)=1$.
Eq.~\eqref{eq:profile_convolution} implies $L_j\le256 \cdot Q \cdot w_{d-j} \cdot M$; also $L_d=0$.
Lemma~\ref{lem:hierarchy}, followed by integration and an application of Lemma~\ref{lem:gronwall}, gives $M(t)\le e^{2308 \cdot d^2 \cdot t}$, where $2308=4+9\cdot256$.
Thus, at $\tau:=1/(2308 \cdot d^2)$, with $w:=\nabla f$,
\begin{equation}\label{eq:profile_gradient}
 G_1(s):=\E_{\mathrm{loc}}[\E_s[w^\top \Gamma_sw]]
 =N_1(s)+L_1(s)\le257 \cdot e \cdot Q \cdot  R^{2d-2} \cdot b_{d-1}^2,\ \forall 0\le s\le\tau.
\end{equation}

The proof of Theorem~\ref{thm:polynomial} gives the lower bound on the remaining variance
$\E_{\mathrm{loc}}[\Var_{\mu_\tau}[f]]\ge e^{-\tau} \cdot \Var_{\mu_0}[f]$.
Apply Lemma~\ref{lem:profile_inflation} at time $\tau$, with $A:=\Gamma_\tau$, $B:=\Lambda_\tau$, $\delta:=\tau$, and $F(\delta):=\Gamma^2 \cdot \ell_r(\delta^{-1})^2$.
The terminal measure has curvature $\Lambda_\tau$, and
$\Lambda_\tau^{-1}\preceq\tau^{-2}\int_0^\tau \Gamma_s\d s$.
For the fixed derivative products $w_i \cdot w_j$, the posterior martingale and tower property give
$\E_{\mathrm{loc}}[\E_\tau[w^\top \Gamma_sw]]=G_1(s)$ for $s\le\tau$.
This keeps the correlation with the earlier covariance matrix.
Consequently
\begin{equation}\label{eq:profile_terminal}
 \Var_{\mu_0}[f]\le e^\tau \cdot \Gamma^2 \cdot \ell_r(\tau^{-1})^2 \cdot 
 (G_1(\tau)+\tau^{-1}\int_0^\tau G_1(s)\d s).
\end{equation}
Eq.~\eqref{eq:profile_gradient}, Eq.~\eqref{eq:iterated_log_scaling}, and $\tau<1$ bound the right side of Eq.~\eqref{eq:profile_terminal} by
\begin{equation*}
 2 \cdot e^2\cdot257\cdot100 \cdot \Gamma^2 \cdot \ell_r(d)^2 \cdot 
 (d!)^2 \cdot \|T\|_{\mathrm{HS}}^2 \cdot R^{2d-2} \cdot b_{d-1}^2.
\end{equation*}
Since $b_d/b_{d-1}\ge\ell_r(d)/4$ and
$16\cdot2 \cdot e^2\cdot257\cdot100<2^{24}=R^2/\Gamma^2$, this closes the induction.
Conditioning on growing balls, centering, and whitening removes compact support at each degree, exactly as in Theorem~\ref{thm:polynomial}.
Moments through degree $2d$ and the corresponding Appell coefficients converge.
The affine contraction described in Section~\ref{tg:sec:families} then gives the assertion for covariance at most $I$.
\end{proof}

\subsection{A starting family from the restricted operator}\label{sec:restricted_operator}\label{sec:compensated_restart}

We construct a starting family whose energy is close to the restricted spectral scale and whose first centering losses are small. These two estimates limit the cost of starting the inverse-gradient hierarchy, which will be used repeatedly in the curvature comparison.

For the restricted operator and the initial estimates, $\mu$ is centered and regular with
$\Sigma:=\operatorname{Cov}(\mu)\preceq I$. We write
$P:=C_P(\mu)=\lambda^{-1}$ and define
\begin{equation}\label{eq:restricted_definition}
 R_{\mathrm{res}}:=\sup_{\E[\nabla f]=0}
       \frac{\Var[f]}{\E[|\nabla f|^2]},
\end{equation}
where the supremum is over nonconstant tests of finite energy.

\begin{lemma}\label{lem:restricted_operator}
We define $A:=LH^{1/2}$, $B:=H^{-1}$, and $E:=I-A^*A$ on centered
$L^2$. The operator $A$ extends boundedly, $AA^*=I$, and $E$ is an
orthogonal projection. Moreover,
\begin{equation}\label{eq:restricted_compressions}
 D^*D=E,\qquad
 R_{\mathrm{res}}=\|EBE\|_{\mathrm{op}}=\|B-L^*L\|_{\mathrm{op}},\qquad
 \|DH^{-1/2}\|_{\mathrm{op}}^2=R_{\mathrm{res}},
\end{equation}
and
\begin{equation}\label{eq:restricted_gap_one}
 R_{\mathrm{res}}\le P\le R_{\mathrm{res}}+1.
\end{equation}
There is a centered unit restricted extremizer $f$ such that
\begin{equation}\label{eq:restricted_extremizer}
 \E[\nabla f]=0,\qquad
 \E[|\nabla f|^2]=R_{\mathrm{res}}^{-1},\qquad
 H^{-1}f=R_{\mathrm{res}}f+q\cdot X,\qquad q:=\E[Xf].
\end{equation}
\end{lemma}

\begin{proof}
For $z\in\mathbb R^n$, $A^*z=H^{1/2}(z\cdot X)$, whose squared
norm is $|z|^2$ by the Dirichlet identity. Thus $AA^*=I$ and $E$ is a
projection. The uncentered gradient of $H^{-1/2}h$ has squared norm
$\|h\|_2^2$, and its removed mean is $Ah$. This proves $D^*D=E$.
Substituting $h:=H^{1/2}f$ in Eq.~\eqref{eq:restricted_definition} gives
$R_{\mathrm{res}}=\|EBE\|_{\mathrm{op}}$. The change of variables is onto $\operatorname{ran}E$.
Since
\begin{equation*}
 B^{1/2}EB^{1/2}=B-L^*L,
\end{equation*}
the equality of the norms of $TT^*$ and $T^*T$, with
$T:=B^{1/2}E$, proves the remaining norm identities in
Eq.~\eqref{eq:restricted_compressions}.

Relative to $\operatorname{ran}(A^*A)$ and
$\operatorname{ran}E$, the first diagonal block of $B$ is identified
with $ABA^*=\Sigma\preceq I$, and its second is at most $R_{\mathrm{res}}I$.
Positivity bounds the cross term by $\sqrt{R_{\mathrm{res}}} \cdot \|x\| \cdot \|y\|$.
Hence
\begin{equation*}
 \langle(x,y),B(x,y)\rangle
 \le(\|x\|+\sqrt{R_{\mathrm{res}}} \cdot \|y\|)^2
 \le(1+R_{\mathrm{res}}) \cdot (\|x\|^2+\|y\|^2).
\end{equation*}
Taking the norm proves Eq.~\eqref{eq:restricted_gap_one}.

The compact positive operator $B-L^*L$ has positive top eigenvalue $R_{\mathrm{res}}$.
Choose a centered unit eigenfunction $f$. Its eigenvector equation is
the last identity in Eq.~\eqref{eq:restricted_extremizer}. Since
$Bf$ and coordinate functions are in $\operatorname{Dom}(H)$, so is $f$.
Taking gradient means in that identity gives
$q=R_{\mathrm{res}} \cdot \E[\nabla f]+q$, hence $\E[\nabla f]=0$.
Applying $H$ and pairing with $f$ gives
$1=R_{\mathrm{res}} \cdot \E[|\nabla f|^2]$, as the coordinate term has zero pairing.
This proves all the assertions with the required domains.
\end{proof}

The extremizer of the restricted problem supplies a centered gradient whose energy is close to the spectral lower bound.

\begin{lemma}[A low-energy gradient]\label{lem:restricted_gradient_energy}
For $f$ in Eq.~\eqref{eq:restricted_extremizer}, the family
$U:=\sqrt{R_{\mathrm{res}}}\nabla f$ is centered, has unit norm, and has energy
\begin{equation}\label{eq:restricted_gradient_energy}
 \lambda\le e:=\|H^{1/2}U\|_2^2\le R_{\mathrm{res}}^{-1}.
\end{equation}
\end{lemma}
\begin{proof}
Apply Eq.~\eqref{tg:eq:bochner} to
$h:=H^{-1}f=R_{\mathrm{res}}f+q\cdot X$. Since $Hh=f$ and $D^2h=R_{\mathrm{res}}D^2f$,
\begin{equation*}
 1=R_{\mathrm{res}}^2 \cdot \E[\|D^2f\|_{\mathrm{HS}}^2]
   +\E[(R_{\mathrm{res}}\nabla f+q)^\top D^2W(R_{\mathrm{res}}\nabla f+q)].
\end{equation*}
The second term is nonnegative. The first, divided by $R_{\mathrm{res}}$, is exactly
$e$. The lower bound follows from the spectral gap for the centered
unit family. The calculation keeps the affine term in $h$ and the full
contribution of Fisher information in the Bochner identity for $f$.
\end{proof}

For a centered unit starting family of energy $\nu$, the construction in
Section~\ref{tg:sec:families} uses the same scalar normalizers. At a zero
raw successor we set all later families to zero, with $\beta_j:=\lambda$
and $\chi_j:=0$ at zero steps. With $\delta_0:=\nu-\lambda\ge0$, the
energy calculation gives
\begin{align}
 v_{j+1}&=v_j-p_j, & e_{j+1}&\le e_j-av_j-\chi_j,
 &\chi_j&\ge0,\label{eq:restricted_general_energy}\\
 v_N&=1-P_N, &a \cdot V_N+X_N&\le\delta_0+\lambda \cdot  P_N,
 \label{eq:restricted_general_budget}\\
 p_j&\le\min\{v_j,e_j\}\le\nu,&
 \lambda&\le\beta_j\le e_{j-1}/v_j,\ \forall j\ge1\text{ with }v_j>0.
 \label{eq:restricted_general_bounds}
\end{align}
Here $P_N:=\sum_{j<N}p_j$, $X_N:=\sum_{j<N}\chi_j$, and
$V_N:=\sum_{j<N}v_j$. The budget follows by subtracting $\lambda v_j$
from the energy recurrence and summing. All operator and form domains
are unchanged.

To keep this advantage after centering, we compare the skew part of the lost mass with the gain from normalization.

\begin{lemma}[Skew loss paid by normalization]\label{lem:skew_normalization_credit}
Let $U$ be a centered unit vector family in the form domain with
$\E[\nabla U]$ symmetric. We define
\begin{equation*}
 e:=\|H^{1/2}U\|_2^2,\qquad
 \beta:=\|H^{-1/2}U\|_2^{-2},\qquad
 F:=\sqrt\beta H^{-1/2}U.
\end{equation*}
The skew part of $p:=\|LH^{1/2}F\|^2$ satisfies
\begin{equation}\label{eq:restricted_skew_credit}
 p_{\mathrm{skew}}
 :=\beta \cdot \|\operatorname{Skew}(LU)\|_{\mathrm{HS}}^2
 \le1-\frac\beta e.
\end{equation}
\end{lemma}
\begin{proof}
For a skew matrix $C$ of unit Hilbert--Schmidt norm, we define $V:=CX$.
The Dirichlet identity gives $\|H^{1/2}V\|_2=1$ and
\begin{equation*}
 \langle H^{1/2}V,H^{1/2}U\rangle
   =\E[\langle C,\nabla U\rangle]=0.
\end{equation*}
In the full direct-sum Hilbert space, the projection of $H^{-1/2}U$
orthogonal to $H^{1/2}U$ has squared norm
$\beta^{-1}-e^{-1}$, since their inner product is one.
Cauchy--Schwarz therefore gives
$|\langle CX,U\rangle|^2\le\beta^{-1}-e^{-1}$.
Take the supremum over $C$ and multiply by $\beta$.
Every pairing is a form pairing; $HU$ need not belong to $L^2$.
\end{proof}

The remaining initial losses involve three slots. The following elementary frame estimate controls them by a symmetric component and the defect under an adjacent swap.

\begin{lemma}[A three-slot frame]\label{lem:startup_three_slot_frame}
On any real orthogonal representation of the symmetric group on three
slots, let $\mathsf P:=(I+s_{12})/2$ and $\mathsf S:=(I+s_{23})/2$.
Then
\begin{equation}\label{eq:restricted_three_slot}
 \|T\|^2\le8(\|\mathsf PT\|^2+\|(I-\mathsf S)T\|^2).
\end{equation}
\end{lemma}
\begin{proof}
On the trivial and sign representations,
$\mathsf P+(I-\mathsf S)=I$. On the standard two-dimensional
representation its eigenvalues are $1\pm\sqrt3/2$.
The smaller is greater than $1/8$. Decomposition into irreducible
orthogonal representations proves Eq.~\eqref{eq:restricted_three_slot}, including
all finite output direct sums.
\end{proof}

Combining the low-energy gradient with these two estimates gives a restart whose first three centering losses are small enough for the block argument.

\begin{proposition}[Compensated restart]\label{prop:compensated_restart}
For every centered regular law of covariance at most $I$ with $R_{\mathrm{res}}\ge16$,
there is a centered unit finite starting family $G_0$ whose energy $\nu$
and first three hierarchy losses satisfy
\begin{equation}\label{eq:restricted_restart_bounds}
 \nu\le\frac1{R_{\mathrm{res}}-K_2/(4R_{\mathrm{res}})},\qquad
 0\le\nu-\lambda\le2R_{\mathrm{res}}^{-2},\qquad
 p'_0+p'_1+p'_2\le160\lambda^2.
\end{equation}
\end{proposition}
\begin{proof}
We define $\sigma:=1/R_{\mathrm{res}}$ and use the gradient $U$ of
Lemma~\ref{lem:restricted_gradient_energy}. Define
$\beta:=\|H^{-1/2}U\|_2^{-2}$,
$F_0:=\sqrt\beta H^{-1/2}U$, and $\chi_{\mathrm{pre}}:=e-\beta$.
Equations~\eqref{eq:restricted_gap_one} and~\eqref{eq:restricted_gradient_energy} imply
\begin{equation}\label{eq:restricted_pre_defect}
 \lambda\le\beta\le e\le\sigma,\qquad
 \sigma-\lambda\le\sigma^2,\qquad
 \chi_{\mathrm{pre}}\le\sigma^2,\qquad
 \sigma/\lambda\le17/16.
\end{equation}
Run Eq.~\eqref{dt_eq_inverse} from $F_0$ and later restart
at $G_0:=F_1/\sqrt{v_1}$.

We first bound its energy. For $M:=\E[X\otimes\nabla f]$, quadratic
testing of $Hf=(f-H(q\cdot X))/R_{\mathrm{res}}$ gives
\begin{equation}\label{eq:restricted_symmetric_quadratic}
 2\operatorname{Sym} M=R_{\mathrm{res}}^{-1}Q_2f,\qquad
 \|\operatorname{Sym}(LU)\|_{\mathrm{HS}}^2\le\frac{K_2}{4R_{\mathrm{res}}}.
\end{equation}
Indeed the pairing of $H(q\cdot X)$ with the centered quadratic
Appell tensor is zero by integration by parts and $\E[X]=0$.
Lemma~\ref{lem:skew_normalization_credit} and Eq.~\eqref{eq:restricted_symmetric_quadratic}
therefore yield
\begin{equation}\label{eq:restricted_surviving_mass}
 1-p_0\ge\frac\beta e-\frac{K_2 \cdot \beta}{4R_{\mathrm{res}}}
 \ge\beta \cdot (R_{\mathrm{res}}-\frac{K_2}{4R_{\mathrm{res}}}).
\end{equation}
The denominator is positive. Energy decrease in the next normalization
gives
\begin{equation}\label{eq:restricted_restart_energy}
 \nu=\frac{e_1}{v_1}=\beta_1
 \le\frac\beta{1-p_0}\le\frac1{R_{\mathrm{res}}-K_2/(4R_{\mathrm{res}})}.
\end{equation}
Since $\lambda\ge\sigma/(1+\sigma)$ and $K_2/4\le2$,
\begin{equation*}
 \nu-\lambda
 \le\frac\sigma{1-2\sigma^2}-\frac\sigma{1+\sigma}
 =\frac{\sigma^2 \cdot (1+2\sigma)}{(1-2\sigma^2) \cdot (1+\sigma)}
 \le2\sigma^2.
\end{equation*}

It remains to bound the actual restarted losses. At the original loss
$p_1$, the tensor $T:=\sqrt{\beta_1}Q_1u^0$ has a coordinate slot,
a newest derivative slot, and the gradient slot inherited from $U$.
Polynomial testing and Eq.~\eqref{dt_eq_defect} give
\begin{align*}
 \mathsf PT&=\frac{\sqrt{\beta \cdot \beta_1}}2Q_2U,\\
 z_{\mathrm{pre}}&:=H^{-1/2}F_0-\beta^{-1/2}U,\\
 u^0&=\beta^{-1/2}P_+\nabla U+P_+\nabla z_{\mathrm{pre}},\qquad
 \|\nabla z_{\mathrm{pre}}\|_2^2=\chi_{\mathrm{pre}}/\beta.
\end{align*}
The first term in $u^0$ is symmetric in its two gradient slots.
Consequently Eq.~\eqref{eq:restricted_three_slot} gives
\begin{equation}\label{eq:restricted_inherited_loss}
 p_1\le2 \cdot K_2 \cdot \beta \cdot \beta_1
       +8(\beta_1/\beta) \cdot \chi_{\mathrm{pre}}.
\end{equation}
At $j=2,3$, apply the same frame to the coordinate and two newest
derivative slots; the oldest inherited slot is an output index.
The ordinary defect identity and polynomial testing give
\begin{equation}\label{eq:restricted_next_startup_losses}
 p_j\le2 \cdot K_2 \cdot \beta_j \cdot \beta_{j-1} \cdot v_{j-1}
       +8(\beta_j/\beta_{j-1}) \cdot \chi_{j-2}.
\end{equation}
Here $\beta_0$ means $\beta$ when needed.

All original energies and losses are at most $\beta\le\sigma$.
Through $j=3$ we therefore have
\begin{equation*}
 v_j\ge1-j\sigma\ge13/16,\quad
 \beta_j\le4\sigma/3,\quad
 \frac{\beta_j}{\beta_{j-1}}\le\frac43,\quad
 \beta_j \cdot \beta_{j-1} \cdot v_{j-1}\le4\sigma^2/3.
\end{equation*}
The ratio bound uses $e_j\le e_{j-1}$ and
$p_{j-1}\le e_{j-1}=\beta_{j-1} \cdot v_{j-1}$.
The initial excess is $\beta-\lambda\le\sigma^2$;
Eq.~\eqref{eq:restricted_general_budget} gives
$\chi_0\le2\sigma^2$ and $\chi_1\le3\sigma^2$.
Thus Eqs.~\eqref{eq:restricted_inherited_loss}--\eqref{eq:restricted_next_startup_losses}
yield
\begin{equation*}
 p_1\le32\sigma^2,\qquad
 p_2\le\frac{128}3 \cdot \sigma^2,\qquad
 p_3\le\frac{160}3 \cdot \sigma^2,
 \qquad p_1+p_2+p_3\le128\sigma^2.
\end{equation*}
All successors used here are nonzero by the displayed mass bounds.
A common scaling of a family leaves every scalar normalizer unchanged.
Hence the first three losses starting from $G_0$ are exactly
$(p_1,p_2,p_3)/v_1$. Since $v_1\ge1-\sigma$,
\begin{equation*}
 p'_0+p'_1+p'_2
 \le\frac{128\sigma^2}{1-\sigma}
 \le\frac{2048}{15} \cdot \sigma^2
 \le\frac{2312}{15} \cdot \lambda^2<160\lambda^2.
\end{equation*}
This proves Eq.~\eqref{eq:restricted_restart_bounds}.
\end{proof}

\subsection{Coefficient transfer with a small multiplicative loss}\label{sec:quadratic_coefficient}

We sharpen the transfer in Theorem~\ref{thm:profile_coefficients} so that repeated use does not accumulate a fixed multiplicative loss at every depth. By varying the localization parameters, we obtain losses whose product over successive depths is finite for profiles defined by iterated logarithms.

We use $K_d(\mu)$ from Section~\ref{sec:preliminaries_appell} and define its normalized coefficient
\begin{equation*}
 c_d(\mu):=\frac{\sqrt{K_d(\mu)}}{d!}.
\end{equation*}
All uniform assertions below may be made either in every dimension or in the fixed range $1\le m\le n$. The chosen range is kept fixed throughout each induction. A profile $\Phi:(0,\infty)\to(0,\infty)$ is admissible if it is finite, continuous, and nondecreasing, if $\Phi(x)\ge1$ for $x\ge1$, and if
\begin{equation}\label{eq:constant_profile_scaling}
 \Phi(x)\le\Phi(y)\le\frac yx \cdot \Phi(x),\ \forall 0<x\le y.
\end{equation}
When using such a profile for measures, we additionally require that every centered regular log-concave measure $\nu$ in the chosen dimension range with covariance at most $I$ satisfies
\begin{equation}\label{eq:constant_profile_input}
 D^2W\succeq x^{-1}I\quad\Longrightarrow\quad C_P(\nu)\le\Phi(x).
\end{equation}
The scalar requirements and the assertion about measures are distinct.

For example, let $\Phi_*(x)$ be the supremum of $C_P(\nu)$ over the regular measures in Eq.~\eqref{eq:constant_profile_input}. The matrix Poincar\'e inequality bounds $\Phi_*(x)$ by $x$. Contracting a measure by $\sqrt{x/y}$ proves $\Phi_*(y)\le(y/x) \cdot \Phi_*(x)$ for $x\le y$, while inclusion of the classes proves monotonicity. These inequalities give continuity. Scaled Gaussian measures show that $\Phi_*(x)=x$ for $x\le1$ and $\Phi_*(x)\ge1$ for $x\ge1$. Thus the intrinsic supremum is an admissible profile.

For $0<\epsilon<1$, we define
\begin{align}
 A_0(\epsilon)&:=9\epsilon \cdot  e^{13\epsilon}
       +(e^{2\epsilon}+3\sqrt{\epsilon} \cdot e^{13\epsilon/2})^2,
       \label{eq:constant_error_factor}\\
 \mathcal E_\Phi(x)&:=\inf_{0<\epsilon<1,\ s>0}
       e^\epsilon \cdot  A_0(\epsilon) \cdot (1+s) \cdot \Phi(x/(\epsilon s))
       ,\ \forall x\ge1.\label{eq:constant_quadratic_profile_transform}
\end{align}
This function is finite and nondecreasing. Moreover,
\begin{equation}\label{eq:constant_transform_lower}
 \mathcal E_\Phi(x)\ge\Phi(x)\ge1,\ \forall x\ge1.
\end{equation}
Indeed, $(1+s) \cdot \Phi(t/s)\ge\Phi(t)$: for $s\le1$ this follows from monotonicity, and for $s\ge1$ it follows from Eq.~\eqref{eq:constant_profile_scaling}. Apply this with $t:=x/\epsilon$ and use $e^\epsilon \cdot  A_0(\epsilon)\ge1$.

We define
\begin{equation}\label{eq:constant_product_weights}
 a_1:=1,\qquad a_j:=\sqrt{\mathcal E_\Phi(j^6)},\ \forall j\ge2,
 \qquad \mathcal B_0:=\mathcal B_1:=1,\qquad
 \mathcal B_d:=\prod_{j=2}^d a_j,\ \forall d\ge2.
\end{equation}
The sequence $(a_j)$ is nondecreasing, so
\begin{equation}\label{eq:constant_product_comparison}
 \mathcal B_k \cdot \mathcal B_{m-k}\le\mathcal B_m
 ,\ \forall 0\le k\le m.
\end{equation}

The product weights $\mathcal B_d$ bound the normalized polynomial coefficients under the curvature estimate for $\Phi$. Their dependence on the localization parameters makes it possible to control the loss in this transfer.

\begin{proposition}[Transfer from curvature to coefficients]
\label{prop:constant_quadratic_coefficient_transfer}
Suppose $\Phi$ is admissible and satisfies Eq.~\eqref{eq:constant_profile_input}. Then every isotropic log-concave measure $\mu$ in the chosen dimension range satisfies
\begin{equation}\label{eq:constant_coefficient_transfer}
 c_d(\mu)\le\mathcal B_d,\ \forall d\ge1.
\end{equation}
The same upper bound holds for centered log-concave measures with covariance at most $I$.
\end{proposition}
\begin{proof}
We induct simultaneously over all isotropic measures in the chosen dimension range. Degree one follows from isotropy, with $c_1(\mu)=1$. Assume the claim through degree $d-1$, where $d\ge2$, and first take a compactly supported isotropic initial measure $\mu_0$. We apply Lemma~\ref{lem:construction_localization_process} to this initial measure. Fix a nonzero symmetric $d$-tensor $T$, and define
\begin{equation*}
 f:=P_d^{\mu_0}[T],\qquad Q:=(d!)^2 \cdot \|T\|_{\mathrm{HS}}^2,
 \qquad M(t):=\max_{1\le j\le d}\frac{N_j(t)}{Q \cdot \mathcal B_{d-j}^2}.
\end{equation*}
Here $N_j,L_j$ are the covariance-weighted derivative quantities of Section~\ref{sec:polynomial_variance}. In particular, Eq.~\eqref{eq:appell_centering} gives
\begin{equation*}
 N_d(0)=Q,\qquad N_j(0)=0,\ \forall j<d,\qquad L_d=0.
\end{equation*}
Whitening each posterior and applying the inductive assertion in Eq.~\eqref{eq:polynomial_minkowski} gives
\begin{equation*}
 \sqrt{L_j(t)}\le\sum_{k=1}^{d-j}\mathcal B_k \cdot \sqrt{N_{j+k}(t)}.
\end{equation*}
This inequality is on the full output tensor direct sum, with one covariance factor in every derivative slot. Eq.~\eqref{eq:constant_product_comparison} therefore implies
\begin{equation*}
 L_j(t)\le d^2 \cdot Q \cdot \mathcal B_{d-j}^2 \cdot M(t).
\end{equation*}
Lemma~\ref{lem:hierarchy} now gives, after integration and taking the maximum,
\begin{equation*}
 M(t)\le1+13d^4\int_0^t M(u)\d u\le e^{13d^4 \cdot t}.
\end{equation*}
For the first inequality we used $4j^2+9j^2 \cdot d^2\le13d^4$; the second follows from Lemma~\ref{lem:gronwall}. Using the zero initial values for $j<d$ gives the sharper bounds
\begin{align}
 N_j(t)&\le9d^4 \cdot t \cdot  e^{13d^4 \cdot t} \cdot Q \cdot \mathcal B_{d-j}^2,\ \forall j<d,
       \label{eq:constant_zero_initial}\\
 N_d(t)&\le Q \cdot e^{4d^2 \cdot t}.\label{eq:constant_top_derivative}
\end{align}
For example, integrate the inequality for $N_j$ with the factor $e^{-4j^2t}$ and use $4j^2\le13d^4$ to obtain Eq.~\eqref{eq:constant_zero_initial}.

In the estimate for $L_1$, we separate the term $k=d-1$, which involves $N_d$, from the other terms. Eqs.~\eqref{eq:constant_product_comparison}--\eqref{eq:constant_top_derivative} give
\begin{equation*}
 \sqrt{L_1(t)}\le\sqrt Q \cdot \mathcal B_{d-1} \cdot 
       (e^{2d^2 \cdot t}+3d^3 \cdot \sqrt t \cdot  e^{13d^4 \cdot t/2}).
\end{equation*}
The second summand bounds at most $d-2$ terms, each bounded by
$3d^2 \cdot \sqrt t \cdot  e^{13d^4 \cdot t/2} \cdot \sqrt Q \cdot \mathcal B_{d-1}$.
We choose $\tau:=\epsilon/d^6$. For $0\le t\le\tau$,
\begin{equation}\label{eq:constant_gradient_energy}
 G_1(t):=N_1(t)+L_1(t)
 \le A_0(\epsilon) \cdot Q \cdot \mathcal B_{d-1}^2.
\end{equation}
Here $d^4 \cdot t\le\epsilon$, $d^2 \cdot t\le\epsilon$, and $d^3 \cdot \sqrt t\le\sqrt\epsilon$.

Write $w:=\nabla f$. The lower bound on the remaining variance from Theorem~\ref{thm:polynomial} gives
\begin{equation*}
 \Var_{\mu_0}[f]\le e^\tau \cdot \E_{\mathrm{loc}}[\Var_{\mu_\tau}[f]].
\end{equation*}
Apply Lemma~\ref{lem:profile_inflation} with the continuous curvature profile
$F(\delta):=\Phi(\delta^{-1})$ and with $\delta:=s\tau$. This is valid for every $s>0$, including $\delta>1$. The covariance and curvature matrices are $\Gamma_\tau$ and $\Lambda_\tau$. As in the proof of Eq.~\eqref{eq:profile_terminal},
\begin{equation*}
 \Lambda_\tau^{-1}\preceq\tau^{-2}\int_0^\tau\Gamma_t\d t,
 \qquad
 \E_{\mathrm{loc}}[\E_\tau[w^\top\Gamma_t w]]=G_1(t),\ \forall t\le\tau.
\end{equation*}
The second identity follows from the tower property for the posterior and includes its correlation with $\Gamma_t$. Consequently,
\begin{align*}
 \Var_{\mu_0}[f]
 &\le e^\tau \cdot \Phi(\delta^{-1}) \cdot 
       (G_1(\tau)+\delta \cdot \tau^{-2}\int_0^\tau G_1(t)\d t)\\
 &\le e^\epsilon \cdot  A_0(\epsilon) \cdot (1+s) \cdot 
       \Phi(d^6/(\epsilon s)) \cdot Q \cdot \mathcal B_{d-1}^2.
\end{align*}
Taking the supremum over $T$ gives this inequality for $c_d(\mu_0)^2$. It holds for every $\epsilon,s$, so the infimum in Eq.~\eqref{eq:constant_quadratic_profile_transform} closes the induction. No minimizer is required.

For each fixed degree, conditioning on growing balls, centering, and whitening removes compact support as in Theorem~\ref{thm:polynomial}. Moments through order $2d$ and the coefficients of the corresponding Appell polynomials converge. Finally, whitening and contraction of each tensor slot transfer the estimate to covariance at most $I$; a singular covariance can be treated on its supporting subspace. The approximation is carried out separately at each fixed degree.
\end{proof}

For a profile defined by iterated logarithms, the contraction in Eq.~\eqref{eq:iterated_log_contraction} absorbs the change of scale in the profile transform. This gives the small multiplicative loss needed in the depth induction.

\begin{proposition}[Specialization to iterated logarithms]
\label{prop:constant_log_specialization}
There is a universal integer $r_0$ such that the following holds for $r\ge r_0$ and $\Gamma\ge1$. If every centered regular measure with covariance at most $I$ satisfies
$C_P\le\Gamma^2 \cdot \ell_r(a^{-1})^2$ at curvature at least $aI$, then
\begin{equation}\label{eq:constant_unweighted_log}
 c_d(\mu)\le[(1+r^{-2}) \cdot \Gamma]^{d-1} \cdot 
        \prod_{j=2}^d\ell_r(j)
 \le[(1+r^{-2}) \cdot \Gamma \cdot \ell_r(d)]^{d-1}
\end{equation}
for every isotropic log-concave measure and every $d\ge1$.
\end{proposition}
\begin{proof}
The function $\Phi(x):=\Gamma^2 \cdot \ell_r(x)^2$ satisfies the profile conditions: Eq.~\eqref{eq:iterated_log_contraction} bounds its logarithmic elasticity by $2^{1-r}\le1$. We set $\alpha:=r^{-2}$ and $\epsilon:=c \cdot \alpha^2$, where $c>0$ is a sufficiently small universal constant. Since $A_0(\epsilon)=1+O(\sqrt\epsilon)$,
\begin{equation*}
 \sqrt{e^\epsilon \cdot (1+\epsilon) \cdot A_0(\epsilon)}\le e^{\alpha/4}.
\end{equation*}
For $j\ge1$, the first logarithm satisfies
$g(j^6/\epsilon^2)\le[6+2\log(1/\epsilon)] \cdot g(j)$.
Applying Eq.~\eqref{eq:iterated_log_contraction} to the remaining $r-1$ compositions gives
\begin{equation*}
 \log\frac{\ell_r(j^6/\epsilon^2)}{\ell_r(j)}
 \le2^{-(r-1)} \cdot \log[6+2\log(1/\epsilon)]\le\alpha/4
\end{equation*}
for sufficiently large universal $r_0$. In the infimum defining $\mathcal E_\Phi$, use $s:=\epsilon$. Then $a_j\le e^{\alpha/2} \cdot \Gamma \cdot \ell_r(j)\le(1+\alpha) \cdot \Gamma \cdot \ell_r(j)$ for $j\ge2$. Multiplication and Proposition~\ref{prop:constant_quadratic_coefficient_transfer} prove the result.
\end{proof}

Eq.~\eqref{eq:constant_unweighted_log} is an unweighted coefficient estimate. We combine it with the joint comparison in Section~\ref{sec:joint_comparison}; its degree-one normalization will control the initial terms of the error kernel.

\subsection{Controlling centering losses at all degrees}\label{sec:joint_frame}\label{sec:joint_comparison}\label{sec:restricted_joint_comparison}

We control centering losses across all dyadic degrees with a single frame estimate whose constant is independent of the number of degrees involved. Applied to the inverse-gradient hierarchy, it separates the initial losses, the terminal polynomial contribution, and the accumulated energy defects. This gives the loss estimate needed for the finite curvature comparison.

In an orthogonal representation of the symmetric group on $2d$ letters, we define $\mathsf P_k$ to average permutations of letters $1,\ldots,k$. For dyadic $k<d$, we define $\mathsf S_k$ to average letters $k+1,\ldots,4k$. The two averages commute. All averages below are orthogonal projections.

We define $V_s(N)$ to be the Euclidean permutation module on the $s$-subsets of $\{1,\ldots,N\}$. As in the proof of Eq.~\eqref{dt_eq_incidence}, we define the inclusion map $U_s:V_s(N)\to V_{s+1}(N)$ and its adjoint $D_{s+1}:=U_s^*$. The same counting identity gives
\begin{equation}\label{eq:joint_up_down}
 D_{s+1}U_s-U_{s-1}D_s=(N-2s)I.
\end{equation}
For $0\le j\le N/2$, the harmonic space $\ker D_j$ and its lifts to higher subset levels will be denoted by $\mathcal U_j(N)$. We take $\mathcal U_0(N)$ to be the constant representation.

\begin{lemma}[Subset decomposition]\label{lem:joint_subset_decomposition}
For $s\le N/2$,
\begin{equation*}
 V_s(N)=\bigoplus_{j=0}^s\mathcal U_j(N),\qquad
 \dim\mathcal U_j(N)=\binom Nj-\binom N{j-1}.
\end{equation*}
The summands are irreducible real orthogonal representations and are pairwise inequivalent. Averaging the subgroup $\mathfrak S_s\times\mathfrak S_{N-s}$ has rank one on each $\mathcal U_j(N)$, $0\le j\le s$, and is zero on every other irreducible representation.
\end{lemma}
\begin{proof}
Eq.~\eqref{eq:joint_up_down} shows that $U_s$ is injective for $s<N/2$. If $v\in\ker D_j$, repeated use of the same identity gives
\begin{equation*}
 DU^t v=t(N-2j-t+1)U^{t-1}v,
 \qquad
 \|U^tv\|^2=t! \cdot \frac{(N-2j)!}{(N-2j-t)!} \cdot \|v\|^2.
\end{equation*}
These identities give the stated decomposition and dimension formula. Orthogonality follows by moving lowering operators across the inner product until a harmonic vector is paired with the range of an upper operator.

The commutant on $V_s(N)$ has dimension $s+1$: an invariant matrix entry depends only on the intersection size of its two indexing subsets. The decomposition already has $s+1$ nonzero invariant orthogonal summands. A further splitting, a commutant that is a nontrivial real division algebra, or equivalent summands would increase that dimension. Thus the summands are irreducible and pairwise inequivalent.

Finally, an invariant vector $v$ for the subgroup fixing an $s$-subset $A$ gives an equivariant map $V_s(N)\to V$ by $\delta_{gA}\mapsto gv$. Conversely, evaluation on $\delta_A$ identifies such maps with subgroup-invariant vectors. The multiplicity-free decomposition proves the last assertion.
\end{proof}

The subset decomposition reduces the frame estimate to its harmonic components. The formula below describes the action of averaging on the invariant vector in each component.

\begin{lemma}[Averaging coefficient]\label{lem:joint_averaging_coefficient}
Let $z_j$ be a unit vector in $\mathcal U_j(N)$ invariant under $\mathfrak S_a\times\mathfrak S_{N-a}$, where $a\le N/2$ and $j\le a$. For $a\le b\le N$,
\begin{equation}\label{eq:joint_averaging_coefficient}
 \|\mathsf P_bz_j\|^2=
 \frac{\binom{b-j}{a-j} \cdot \binom{N-a-j}{b-a}}
 {\binom ba \cdot \binom{N-a}{b-a}}.
\end{equation}
Binomial coefficients outside their natural range are zero.
\end{lemma}
\begin{proof}
The inclusion operator from $a$-subsets to $b$-subsets is $U^{b-a}/(b-a)!$. The preceding up/down calculation gives its squared singular value on $\mathcal U_j(N)$ as
\begin{equation*}
 \lambda_j:=\binom{b-j}{a-j} \cdot \binom{N-a-j}{b-a}.
\end{equation*}
For $b>N/2$, the decomposition at level $b$ follows by complementing subsets; the displayed value is zero for a harmonic that does not occur there.

Let $\Pi_j$ be the harmonic projection at level $a$, and let $A:=\{1,\ldots,a\}$. The vector $\Pi_j\delta_A$ spans the invariant line in the subset realization of $\mathcal U_j(N)$, so we may take $z_j$ proportional to it. Transitivity gives
\begin{equation*}
 \|\Pi_j\delta_A\|^2=\frac{\dim\mathcal U_j(N)}{\binom Na}.
\end{equation*}
Averaging $\delta_A$ over the first $b$ letters gives the normalized sum of the $a$-subsets of that block. This is the adjoint inclusion operator applied to the block indicator, divided by $\binom ba$. Applying the singular-value formula and the same transitivity identity at level $b$ gives $\lambda_j/\lambda_0$, which is Eq.~\eqref{eq:joint_averaging_coefficient}.
\end{proof}

To compare consecutive dyadic scales, we also need the decomposition obtained by restricting a harmonic component to half of the letters.

\begin{lemma}[Half-block restriction]\label{lem:joint_half_restriction}
For $j\le m$, the vectors of $\mathcal U_j(4m)$ invariant under permutations of the last $2m$ letters form, as a representation on the first $2m$ letters,
\begin{equation}\label{eq:joint_half_restriction}
 \mathcal U_0(2m)\oplus\cdots\oplus\mathcal U_j(2m),
\end{equation}
with every summand occurring once.
\end{lemma}
\begin{proof}
The part of $V_r(4m)$ invariant under permutations of the last $2m$ letters is $\bigoplus_{t=0}^rV_t(2m)$: the index $t$ records the number of selected letters in the first half, while the selected letters in the second half form one orbit. Consequently a local harmonic $\mathcal U_\ell(2m)$ occurs $r-\ell+1$ times. Subtracting the invariant part of the lifted module $V_{r-1}(4m)$ leaves one copy of every $\ell\le r$. Use $r:=j$.
\end{proof}

Together, the averaging formula and the half-block decomposition determine the weights that enter the comparison between two consecutive scales.

On $\mathcal U_j(2m)$, we choose a unit vector $z_{m,j}$ invariant under permutations of the first $m$ and last $m$ letters. Its sign is immaterial.

On a representation of $\mathfrak S_{4m}$, we define $E_m:=\mathsf P_m\mathsf S_m$. Lemma~\ref{lem:joint_subset_decomposition} shows that on $\mathcal U_j(4m)$, $0\le j\le m$, this projection has rank one. Write $E_m=qq^*$ with $\|q\|=1$. The vector $q$ is invariant under the first $m$ and the last $3m$ letters. By Lemma~\ref{lem:joint_half_restriction}, its decomposition into local spherical vectors is
\begin{equation}\label{eq:joint_branch_weights}
 q=\sum_{\ell=0}^j\sqrt{w_\ell}z_{m,\ell},\qquad
 w_\ell\ge0,\quad\sum_{\ell=0}^jw_\ell=1,
\end{equation}
after choosing the signs of the local spherical vectors.

\begin{lemma}[Branching estimates]\label{lem:joint_branch_weights}
For $1\le j\le m$, the weights in Eq.~\eqref{eq:joint_branch_weights} satisfy
\begin{equation}\label{eq:joint_first_two_weights}
 w_0=\frac{(m)_{\underline j}}{(3m)_{\underline j}},\qquad
 \frac{w_1}{w_0}
 =\frac{j \cdot (2m-1) \cdot (4m-j+1)}{(2m-j) \cdot (2m-j+1)}.
\end{equation}
Here $(x)_{\underline j}$ is a falling factorial. For $j=1$, $w_0=1/3$ and $w_1=2/3$. For $j\ge2$, if $h:=\sum_{\ell\ge2}w_\ell$, then
\begin{equation}\label{eq:joint_branch_dominance}
 h\ge4w_0,\qquad h\ge w_1.
\end{equation}
\end{lemma}
\begin{proof}
The first identity is Eq.~\eqref{eq:joint_averaging_coefficient} with $a:=m$, $b:=2m$, and $N:=4m$. Only local harmonics zero and one survive averaging the first $2m-1$ letters. On the local standard spherical vector, this averaging has squared norm $1/(2m-1)$, by the same incidence formula or a direct calculation with a centered vector constant on each half. Hence
\begin{equation*}
 \|\mathsf P_{2m-1}q\|^2=w_0+\frac{w_1}{2m-1}.
\end{equation*}
Applying Eq.~\eqref{eq:joint_averaging_coefficient} again and simplifying proves the second identity.

For $j=2$, direct substitution gives
\begin{equation*}
 w_0=\frac{m-1}{3 \cdot (3m-1)},\qquad
 w_1=h=\frac{4m-1}{3 \cdot (3m-1)}.
\end{equation*}
For $j\ge3$, $w_0\le3^{-j}$ and $w_1/w_0<6j$. The latter follows because the two factors in Eq.~\eqref{eq:joint_first_two_weights}, apart from $j$, are less than two and three. Therefore
\begin{equation*}
 \frac h{w_0}=\frac1{w_0}-1-\frac{w_1}{w_0}
 \ge3^j-1-6j\ge8.
\end{equation*}
For $j\ge4$, $w_0+2w_1\le3^{-j} \cdot (1+12j)<1$, proving $h\ge w_1$. For $j=3$, substitution gives
\begin{equation*}
 h-w_1=\frac{2m \cdot (m-1) \cdot (14m-13)}{3 \cdot (2m-3) \cdot (3m-2) \cdot (3m-1)}>0.
\end{equation*}
This proves Eq.~\eqref{eq:joint_branch_dominance}.
\end{proof}

These weights allow the frame bound to pass from one dyadic scale to the next with a constant independent of the number of scales.

\begin{theorem}[Joint partial symmetrization]\label{thm:joint_uniform_frame}
For every dyadic $d\ge1$, every finite-dimensional real orthogonal representation of the symmetric group on $2d$ letters, and every vector $T$,
\begin{equation}\label{eq:joint_uniform_frame}
 \|T\|^2\le C_F \cdot 
 [ d\|\mathsf P_dT\|^2+
 \sum_{\substack{k<d\\k\ \mathrm{dyadic}}}
 k^2\|(I-\mathsf S_k)\mathsf P_kT\|^2],
 \qquad C_F:=10^4.
\end{equation}
The same assertion holds for tensor-slot permutations with arbitrary finite output direct sums.
\end{theorem}

\begin{proof}

We define
\begin{equation*}
 \mathsf F_k:=\mathsf P_k(I-\mathsf S_k),\qquad
 \mathsf H_d:=\sum_{\substack{k<d\\k\ \mathrm{dyadic}}}k^2\mathsf F_k,
 \qquad \mathsf M_d:=\mathsf H_d+d\mathsf P_d,\qquad c:=10^{-4}.
\end{equation*}
Since $\mathsf F_k$ is an orthogonal projection, proving $\mathsf M_d\succeq cI$ proves Theorem~\ref{thm:joint_uniform_frame}.

The operator $\mathsf H_m$ commutes with permutations of the last $m$ letters: its largest summand averages a block containing those letters, and earlier summands are supported in the first $m$ letters.

Assuming $\mathsf M_m-cI$ is positive definite, we define the shifted boundary trace
\begin{equation}\label{eq:joint_boundary_trace}
 a_{m,j}:=\inf_{\substack{T\ \mathrm{last}\text{-}m\ \mathrm{invariant}\\
                         \mathsf P_mT=z_{m,j}}}
          \langle T,(\mathsf H_m-cI)T\rangle.
\end{equation}
The infimum is finite and attained. On the kernel of $\mathsf P_m$, the quadratic form equals that of $\mathsf M_m-cI$, which is positive definite. In the sector invariant under the last $m$ letters, $\mathsf P_m$ has rank one. Thus Eq.~\eqref{eq:joint_boundary_trace} is its scalar Schur complement.

We prove simultaneously, for every dyadic $m$,
\begin{align}
 &\mathsf M_m-cI\text{ is positive definite in every representation},
                                                        \label{eq:joint_induction_positive}\\
 &a_{m,0}=-c,                                           \label{eq:joint_induction_trivial}\\
 &0<b_m:=-a_{m,1}\le4cm-2c,                             \label{eq:joint_induction_standard}\\
 &a_{m,j}\ge m/8,\ \forall 2\le j\le m.                   \label{eq:joint_induction_higher}
\end{align}
At $m=1$, $\mathsf H_1=0$ and $\mathsf M_1=I$. Therefore $a_{1,0}=a_{1,1}=-c$, and the assertion for higher harmonics is empty.

On a representation of $\mathfrak S_{4m}$,
\begin{equation}\label{eq:joint_rank_one_step}
 \mathsf M_{2m}=\mathsf H_m+m^2\mathsf P_m-m^2E_m+2m\mathsf P_{2m},
 \qquad E_m:=\mathsf P_m\mathsf S_m.
\end{equation}
By Lemma~\ref{lem:joint_subset_decomposition}, $E_m$ is zero except on $\mathcal U_j(4m)$ with $0\le j\le m$, where it has rank one. All sectors on which $E_m$ vanishes satisfy
\begin{equation*}
 \mathsf M_{2m}-cI\succeq\mathsf M_m-cI\succ0,
\end{equation*}
because $m^2\ge m$. Restriction to the first $2m$ letters may have multiplicities, which do not affect the induction hypothesis.

In an affected sector, use $E_m=qq^*$ and the weights from Eq.~\eqref{eq:joint_branch_weights}. The vector $q$ is invariant under the last $3m$ letters. The remaining operator in Eq.~\eqref{eq:joint_rank_one_step} commutes with permutations of the last $2m$ letters. Its noninvariant sector is unaffected by $qq^*$ and remains positive.

Each local component of $q$ is invariant under the next $m$ letters. The operator preserves that subspace. The other local subspaces contain no component of $q$ and remain positive. Eliminating the variables away from the boundary in the relevant subspaces using Eq.~\eqref{eq:joint_boundary_trace} gives
\begin{equation}\label{eq:joint_reduced_matrix}
 \operatorname{diag}(m^2+a_{m,\ell})_{\ell=0}^j
       -m^2bb^*+2m e_0e_0^*,\qquad b_\ell:=\sqrt{w_\ell}.
\end{equation}
Every eliminated block is positive by Eq.~\eqref{eq:joint_induction_positive}, since
$\mathsf H_m+m^2\mathsf P_m-cI\succeq\mathsf M_m-cI$. The boundary indexed by $\ell=0$ is the line invariant under the first $2m$ letters.

\paragraph{The scalar recurrence.}

For an affected global harmonic $j\le m$, we define
\begin{equation}\label{eq:joint_schur_denominator}
 D_{m,j}:=1-m^2\sum_{\ell=1}^j\frac{w_\ell}{m^2+a_{m,\ell}}.
\end{equation}
If $D_{m,j}>0$, the block of Eq.~\eqref{eq:joint_reduced_matrix} on $\ell\ge1$ is positive definite by the rank-one criterion. Its scalar Schur complement gives the exact recurrence
\begin{equation}\label{eq:joint_trace_recurrence}
 a_{2m,j}=m^2-c-\frac{m^2 \cdot w_0}{D_{m,j}}.
\end{equation}
The final term $2m\mathsf P_{2m}$ adds $2m$ to this boundary trace. For $m<j\le2m$, the projection $E_m$ is zero. The boundary vector invariant under the first $2m$ letters lies in the trivial local representation, so it is uncoupled and
\begin{equation}\label{eq:joint_absent_branch_trace}
 a_{2m,j}=m^2-c.
\end{equation}

For $j=1$, Eqs.~\eqref{eq:joint_first_two_weights}--\eqref{eq:joint_trace_recurrence} give
\begin{equation}\label{eq:joint_standard_recurrence}
 b_{2m}=c+\frac{2b_m}{1-3b_m/m^2}.
\end{equation}
The induction hypothesis gives $3b_m/m^2\le12c<1$ and
\begin{equation*}
 \frac{2b_m}{1-3b_m/m^2}
 \le2b_m+\frac{96c^2}{1-12c}\le2b_m+c.
\end{equation*}
Consequently $0<b_{2m}\le8cm-2c$. The denominator in Eq.~\eqref{eq:joint_schur_denominator} is positive because $m^2-3b_m>0$. The final boundary value satisfies $2m-b_{2m}>0$.

Now let $2\le j\le m$, so $m\ge2$. The induction hypotheses and Lemma~\ref{lem:joint_branch_weights} imply
\begin{align}
 D_{m,j}
 &=w_0+\sum_{\ell\ge1}\frac{w_\ell \cdot  a_{m,\ell}}{m^2+a_{m,\ell}}\notag\\
 &\ge w_0+\frac h{8m+1}-\frac{4c \cdot  w_1}{m-4c}
 \ge w_0+\frac h{9m}>0.                            \label{eq:joint_positive_denominator}
\end{align}
For the last inequality, use $h\ge w_1$,
$1/(8m+1)\ge2/(17m)$, $4c/(m-4c)\le8c/m$, and
$2/17-8c>1/9$. Eq.~\eqref{eq:joint_trace_recurrence} now gives
\begin{equation*}
 a_{2m,j}\ge\frac{m^2 \cdot h}{9m \cdot w_0+h}-c
 \ge\frac{4m^2}{9m+4}-c\ge\frac m4.
\end{equation*}
The last difference before subtracting $c$ is
$m \cdot (7m-4)/(4 \cdot (9m+4))\ge5m/44$ for $m\ge2$.
Eq.~\eqref{eq:joint_absent_branch_trace} gives the same bound when $j>m$, including $m=1$.

The trivial trace remains $-c$, and its full frame value is $2m-c>0$. The standard and higher boundary values are positive after adding $2m$. All eliminated blocks are positive, and Eq.~\eqref{eq:joint_positive_denominator} or the standard calculation gives positivity of the remaining rank-one block. The unaffected representation sectors were already positive. This closes Eqs.~\eqref{eq:joint_induction_positive}--\eqref{eq:joint_induction_higher} and proves Theorem~\ref{thm:joint_uniform_frame}.
\end{proof}

We apply the joint frame to the tensors representing centering losses. The resulting inequality separates the initial losses, the contribution of the terminal degree, and the energy defects accumulated along the hierarchy.

Use the normalized families and coefficients from Section~\ref{sec:centering_loss}. For a finite prefix, we define
\begin{equation}\label{eq:joint_budget}
 P_N:=\sum_{j<N}p_j,\qquad X_N:=\sum_{j<N}\chi_j,\qquad
 V_N:=\sum_{j<N}v_j,\qquad
 v_N=1-P_N,\quad a \cdot V_N+X_N\le\lambda \cdot  P_N.
\end{equation}
The last two relations follow from Eq.~\eqref{dt_eq_excess}.

\begin{lemma}[Joint local loss]\label{lem:joint_local_loss}
Let $d\ge2$ be dyadic. On a prefix with $\beta_j\le b\lambda$, where $b\ge1$, define
\begin{align}
 \theta_d&:=324 \cdot C_F \cdot 
       \sum_{\substack{k<d\\k\ \mathrm{dyadic}}}
       k^6 \cdot (b^4\lambda)^k \cdot c_k^2,                 \label{eq:joint_kernel}\\
 \Delta_d&:=3\lambda+2 \cdot C_F \cdot 
       \sum_{\substack{2\le k<d\\k\ \mathrm{dyadic}}}
       k^2 \cdot (b\lambda)^k \cdot c_k^2.                  \label{eq:joint_startup}
\end{align}
Then, whenever the prefix hypotheses are justified,
\begin{equation}\label{eq:joint_prefix_loss}
 P_N\le\Delta_d+N \cdot C_F \cdot d \cdot (b\lambda)^d \cdot c_d^2
              +\frac{\theta_d}{\lambda} \cdot X_N.
\end{equation}
\end{lemma}
\begin{proof}
For $j\ge2d-1$, we define $T_j:=\sqrt{\beta_j}Q_1u^{j-1}$, so $\|T_j\|^2=p_j$. Repeated use of Eq.~\eqref{eq:dyadic_testing} gives the exact compatibility identity
\begin{equation}\label{eq:joint_coherent_projections}
 \mathsf P_kT_j=\frac{\prod_{i=0}^{k-1}\sqrt{\beta_{j-i}}}{k!}
                      Q_ku^{j-k}.
\end{equation}
Thus all degrees are projections of one tensor. The old derivative slots beyond the first $2d$ slots form an output direct sum. In particular,
\begin{equation}\label{eq:joint_terminal_tensor}
 \|\mathsf P_dT_j\|^2\le(b\lambda)^d \cdot c_d^2.
\end{equation}
The derivative block averaged by $\mathsf S_k$ has $3k$ slots. Eqs.~\eqref{dt_eq_sort}, \eqref{eq:dyadic_swap}, and~\eqref{eq:joint_coherent_projections} give
\begin{equation}\label{eq:joint_frame_defect}
 \|(I-\mathsf S_k)\mathsf P_kT_j\|
 \le6k \cdot (b\lambda)^{k/2} \cdot c_k \cdot \lambda^{-1/2} \cdot 
       \sum_{a=1}^{3k-1}b^{(a-1)/2} \cdot 
                    \sqrt{\chi_{j-k-a}}.
\end{equation}
All indices are nonnegative, since $j\ge4k-1$ for the relevant $k<d$.

The $\ell^2$ convolution norm of the finite inner kernel is at most its $\ell^1$ norm, which is at most $3k \cdot b^{3k/2}$. Squaring Eq.~\eqref{eq:joint_frame_defect}, summing in $j$, and multiplying by the frame weight $k^2$ therefore gives
\begin{equation*}
 324 \cdot C_F \cdot  k^6 \cdot b^{4k} \cdot \lambda^{k-1} \cdot c_k^2 \cdot X_N.
\end{equation*}
Negative indices may be assigned zero. Summing the inequality in Theorem~\ref{thm:joint_uniform_frame} yields the kernel in Eq.~\eqref{eq:joint_kernel}.

At each generation we use the largest available dyadic degree, capped at $d$. Degree two becomes available at $j=3$; the first three losses are at most $\lambda$ each. A degree $k<d$ is used for exactly $2k$ generations before the next doubling, and its contribution apart from the error at each generation is $C_F \cdot k \cdot (b\lambda)^k \cdot c_k^2$. This gives Eq.~\eqref{eq:joint_startup}. The terminal contribution is bounded by the middle term in Eq.~\eqref{eq:joint_prefix_loss}. The same inequality holds before the terminal degree becomes available, since adding that term is harmless.
\end{proof}

For a general starting family, the restricted operator gives a sharper bound on the propagation of an existing defect. This bound gives the local inequality needed for the block construction.

\begin{lemma}[Refined joint local inequality]\label{lem:restricted_joint_local_loss}
For a unit starting family of energy $\nu$, suppose all normalizers
needed in a prefix satisfy $\beta_j\le B_*$, and define
$t_*:=\max\{1,B_* \cdot R_{\mathrm{res}}\}$. For a dyadic $d\ge2$, define
\begin{align}
 \theta&:=324 \cdot C_F \cdot \sum_{\substack{k<d\\k\ \mathrm{dyadic}}}
       k^6 \cdot (B_* \cdot t_*^3)^k \cdot c_k^2,\label{eq:restricted_refined_kernel}\\
 \Delta&:=S_3+2 \cdot C_F \cdot \sum_{\substack{2\le k<d\\k\ \mathrm{dyadic}}}
       k^2 \cdot B_*^k \cdot c_k^2,\qquad S_3:=p_0+p_1+p_2,\label{eq:restricted_refined_startup}\\
 \tau_d&:=C_F \cdot d \cdot B_*^d \cdot c_d^2.\label{eq:restricted_terminal_coefficient}
\end{align}
Theorem~\ref{thm:joint_uniform_frame} yields
\begin{equation}\label{eq:restricted_refined_local_loss}
 P_N\le\Delta+N \cdot \tau_d+\frac\theta\lambda \cdot  X_N.
\end{equation}
\end{lemma}
\begin{proof}
For $T_j:=\sqrt{\beta_j}Q_1u^{j-1}$, repeated polynomial testing
gives the coherent projections
\begin{equation*}
 \mathsf P_kT_j
 =\frac{\prod_{i=0}^{k-1}\sqrt{\beta_{j-i}}}{k!}Q_ku^{j-k},
 \qquad \|\mathsf P_dT_j\|^2\le B_*^d \cdot c_d^2.
\end{equation*}
A fresh adjacent swap has norm at most $2\sqrt{\chi/\lambda}$ by
Eq.~\eqref{dt_eq_defect}. Every subsequent map is
$\sqrt\beta DH^{-1/2}$ and has norm at most $\sqrt{t_*}$ by
Eq.~\eqref{eq:restricted_compressions}. These maps commute with
permutations of previously introduced slots. Eq.~\eqref{dt_eq_sort} for the derivative block of $3k$ slots consequently becomes
\begin{equation*}
 \|(I-\mathsf S_k)\mathsf P_kT_j\|
 \le6k \cdot B_*^{k/2} \cdot c_k \cdot \lambda^{-1/2} \cdot 
       \sum_{a=1}^{3k-1}t_*^{(a-1)/2} \cdot \sqrt{\chi_{j-k-a}}.
\end{equation*}
The finite inner kernel has $\ell^1$ norm at most $3k \cdot t_*^{3k/2}$.
Square, sum in $j$, and multiply by the frame weight $C_F \cdot k^2$.
This gives the kernel in Eq.~\eqref{eq:restricted_refined_kernel}. The fresh
error still has denominator $\lambda$; only its later propagation
has been sharpened.

At each generation use the largest available dyadic degree, capped at
$d$. The first three losses contribute $S_3$; degree $k<d$ is
then used for $2k$ generations, with coherent contribution at most
$C_F \cdot k \cdot B_*^k \cdot c_k^2$ each. This gives
Eq.~\eqref{eq:restricted_refined_startup}; the terminal contribution is
$N \cdot \tau_d$. Older output indices do not change the uniform frame.
All these estimates hold before the terminal degree is available by
adding the nonnegative terminal bound.
\end{proof}

\subsection{Coarse initialization and a finite comparison}\label{sec:coarse_profiles}\label{sec:delayed_block_comparison}

We establish the coarse bounds needed to start the refined argument: a curvature profile at every fixed depth and one coefficient estimate valid in every degree. The constants in these initial profiles may grow with depth, since only a fixed depth is needed to start the final induction. We then prove a finite energy comparison that converts bounds on centering losses into a curvature estimate.

\begin{lemma}\label{lem:coarse_profiles}
There is a universal constant $C_{\mathrm{seed}}\ge1$ such that
\begin{equation*}
 c_k\le[C_{\mathrm{seed}} \cdot \log(e+k)]^{k-1},\ \forall k\ge1
\end{equation*}
for every centered log-concave measure of covariance at most $I$.
For every fixed integer $r\ge1$, every centered regular measure of covariance at most $I$ and curvature at least $aI$ also satisfies
\begin{equation*}
 C_P\le\widehat\Gamma_r^2 \cdot \ell_r(a^{-1})^2,
 \qquad \widehat\Gamma_r:=2^{48+30(r-1)}.
\end{equation*}
\end{lemma}
\begin{proof}
Theorem~\ref{thm:polynomial} and affine contraction give
$c_k\le32^k \cdot k!\le128^k \cdot k^k/(k+1)^2$, since $k!\le k^k$ and $(k+1)^2\le4^k$.
Apply Theorem~\ref{thm:dyadic_curvature} with $\epsilon:=1$, $R:=2^{40}$, and $\ell(k):=k$.
Define $h:=g(a^{-1})\ge1$, and choose a dyadic $d$ with
$\max\{2,h\}\le d<2\max\{2,h\}\le6h$.
Then $\max\{1,a^{-1/(d+1)}\}\le e$, and
\begin{equation*}
 C_P(\nu)\le36 \cdot e \cdot 2^{85} \cdot g(a^{-1})^2
 \le(2^{48})^2 \cdot \ell_1(a^{-1})^2.
\end{equation*}
This establishes the first curvature profile.

Suppose a profile with constant $\Gamma\ge2^{48}$ is established at depth $r$.
Theorem~\ref{thm:profile_coefficients} gives the coefficient bound parameter $R:=2^{12} \cdot \Gamma\ge2^{40}$.
Apply Theorem~\ref{thm:dyadic_curvature} with $\epsilon:=1$ and the same dyadic choice of $d$.
Eq.~\eqref{eq:iterated_log_scaling} gives
$\ell_r(d)\le6 \cdot \ell_r(h)=6 \cdot \ell_{r+1}(a^{-1})$.
Therefore
\begin{equation*}
 C_P(\nu)\le36 \cdot e \cdot 2^{29} \cdot \Gamma^2 \cdot \ell_{r+1}(a^{-1})^2
 \le(2^{30} \cdot \Gamma)^2 \cdot \ell_{r+1}(a^{-1})^2.
\end{equation*}
Finite induction gives a profile at every fixed depth, with constant
$\widehat\Gamma_r:=2^{48+30(r-1)}$.

The first profile also supplies a fixed coefficient estimate. Apply Proposition~\ref{prop:constant_quadratic_coefficient_transfer} to
$\Phi(x):=(2^{48})^2 \cdot g(x)^2$, with fixed numerical choices of the two parameters in Eq.~\eqref{eq:constant_quadratic_profile_transform}. The elasticity estimate~\eqref{eq:iterated_log_contraction} makes this profile admissible. Since $g(C \cdot j^6)\le C' \cdot g(j)$ for fixed $C$, its consecutive product factors satisfy $a_j\le C_{\mathrm{seed}} \cdot g(j)$ for $j\ge2$. Hence
\begin{equation}\label{eq:joint_logarithmic_seed}
 c_k\le[C_{\mathrm{seed}} \cdot \log(e+k)]^{k-1},\ \forall k\ge1.
\end{equation}
This estimate is uniform over log-concave measures of covariance at most $I$. At degree one, isotropy gives $c_1=1$, and affine contraction gives $c_1\le1$ in the larger class.
\end{proof}

The local loss bounds must be matched to the energy budget. The restricted operator permits a one-step delay in that budget, which is useful when the starting energy is close to its restricted spectral scale.

For the delayed comparison and the initial block estimates, $\mu$ is centered and regular,
$\operatorname{Cov}(\mu)\preceq I$, and its potential has Hessian at least
$aI$. We write $R:=R_{\mathrm{res}}$, $\sigma:=R^{-1}$, and $B:=H^{-1}$,
and define
\begin{equation}\label{eq:block_operator}
 \mathcal T:=P_+\nabla B=DH^{-1/2}.
\end{equation}
The map acts componentwise on finite output families. By
Lemma~\ref{lem:restricted_operator},
\begin{equation}\label{eq:block_operator_norm}
 \mathcal T^*\mathcal T=B-L^*L,\qquad
 \|\mathcal T\|_{\mathrm{op}}^2=R,\qquad
 \frac\sigma{1+\sigma}\le\lambda\le\sigma.
\end{equation}
For a centered unit starting family of energy $\nu$, we use the
notation $v_j,e_j,p_j,\chi_j,P_N,V_N,X_N$ from
Eqs.~\eqref{eq:restricted_general_energy}--\eqref{eq:restricted_general_bounds}.
The restricted operator also gives a delayed version of the energy
budget. Indeed, $v_{N+1}=\|DF_N\|_2^2\le R \cdot e_N$, and hence
\begin{equation}\label{eq:delayed_block_budget}
 a \cdot V_N+X_N\le\nu-e_N
 \le\nu-\sigma+\sigma \cdot  P_{N+1}.
\end{equation}

The delayed energy budget bounds how long the hierarchy can stay within a prescribed cumulative centering loss. The resulting finite comparison gives the following bound on curvature.

\begin{lemma}[Delayed finite comparison]\label{lem:delayed_block_comparison}
We define $D_0:=\max\{\nu-\sigma,0\}$. Suppose $0<p\le1/8$, $\nu\le1/4$,
and $D_0\le\sigma$. We define $B_*:=\nu/(1-p)$. Assume that whenever the normalizers required
in a prefix are at most $B_*$, including a first exiting prefix, we have
\begin{equation}\label{eq:delayed_block_local}
 P_N\le\Delta+N\tau+
       \frac\theta\lambda \cdot  X_{\max\{N-2,0\}},
 \qquad p_0+p_1\le\Delta,
\end{equation}
where $\Delta,\tau,\theta$ are nonnegative and independent of $N$.
If
\begin{equation}\label{eq:delayed_block_conditions}
 \frac{\theta \cdot \sigma}{\lambda}\le\frac18,\qquad
 \Delta+\frac\theta\lambda \cdot  D_0\le\frac p{64},
\end{equation}
then
\begin{equation}\label{eq:delayed_block_comparison}
 a<12(\sigma+D_0/p) \cdot \tau.
\end{equation}
\end{lemma}
\begin{proof}
For $N\ge3$, Eq.~\eqref{eq:delayed_block_budget} at $N-2$ gives
$X_{N-2}\le D_0+\sigma \cdot  P_{N-1}\le D_0+\sigma \cdot  P_N$.
For $N\le2$, the startup bound in
Eq.~\eqref{eq:delayed_block_local} suffices. Thus every justified prefix
satisfies
\begin{equation}\label{eq:delayed_block_absorption}
 (1-\theta \cdot \sigma/\lambda) \cdot P_N
 \le\Delta+N\tau+(\theta/\lambda) \cdot D_0.
\end{equation}
At $N=1$, the energy budget gives
$a\le D_0+\sigma \cdot (p_0+p_1)\le D_0+\sigma p$.
The first two successors exist since each loss is at most $\nu\le1/4$.

Suppose the reverse of Eq.~\eqref{eq:delayed_block_comparison} holds,
and set
\begin{equation*}
 m:=\lceil(D_0+\sigma p)/a\rceil,\qquad M:=m+1.
\end{equation*}
Then $M\le3(D_0+\sigma p)/a$ and $M\tau\le p/4$.
Before a first $N\le M$ with $P_N>p$, every normalizer needed for
its losses uses a preceding mass at least $1-p$ and is at most $B_*$.
Its successor mass is at least $1-p-\nu>0$.
Equations~\eqref{eq:delayed_block_conditions} and
\eqref{eq:delayed_block_absorption} therefore give
\begin{equation*}
 P_N\le p/64+p/4+P_N/8,\qquad
 P_N\le cp,\qquad c:=17/56<1.
\end{equation*}
There is no exit, and the same bound holds at every earlier prefix.
At $m=M-1$, Eq.~\eqref{eq:delayed_block_budget} yields
\begin{equation*}
 a \cdot V_m\le D_0+c\sigma p,\qquad
 a \cdot V_m\ge(1-cp)am\ge(1-cp) \cdot (D_0+\sigma p).
\end{equation*}
The difference between the last lower bound and the upper bound is
\begin{equation*}
 \sigma p \cdot (1-c-cp-cD_0/\sigma)
 \ge\sigma p \cdot (1-2c-cp)>0.
\end{equation*}
This contradiction proves the assertion. Every horizon used is finite.
\end{proof}

\subsection{Averaged blocks and a static coefficient radius}\label{sec:initial_iterated_profile}\label{sec:mesoscopic_blocks}\label{sec:block_averaged_start}

We average a block of the inverse-gradient hierarchy to make its initial centering losses and excess energy cubic in the inverse block radius. This decay is the source of the squared-amplitude cost $C \cdot (r+1)^{1/3}$ in the final iteration. The same radius also controls every polynomial coefficient and the Poincar\'e constant up to a fixed additive error, so one parameter can be used throughout the iteration.

Throughout the block arguments, $\mu$ is centered and regular, its covariance is at most $I$, and $D^2W\succeq aI$ for $a>0$. We use $R:=R_{\mathrm{res}}$, $\sigma:=R^{-1}$, $B:=H^{-1}$, and $\mathcal T:=P_+\nabla B$.

We first bound the centering losses at degree three. Let
$\overline c_3\ge1$ be a universal upper bound for $c_3$, supplied by
Eq.~\eqref{eq:joint_logarithmic_seed}.

\begin{lemma}[A cubic local loss]\label{lem:block_cubic_loss}
Suppose the normalizers needed for a loss $p_j$, $j\ge3$, satisfy
$\beta_i\le3\sigma$ and $\beta_i/\lambda\le6$. Then
\begin{equation}\label{eq:block_cubic_loss}
 p_j\le A_0 \cdot \sigma^3+B_0 \cdot (\chi_{j-2}+\chi_{j-3}),\qquad
 A_0:=486 \cdot \overline c_3^2,\quad B_0:=50000.
\end{equation}
\end{lemma}
\begin{proof}
Define $T_j:=\sqrt{\beta_j}Q_1u^{j-1}$, so that
$\|T_j\|^2=p_j$. We consider its coordinate slot and its three newest
derivative slots, keeping all older slots as an output direct sum.
Let $\mathsf S$ average the three derivative slots. Block recovery,
Lemma~\ref{dt_lem_block} with $s:=1$, $\ell:=3$, and $q:=3$, gives
\begin{equation*}
 \|T_j\|\le3\|\mathsf P_3T_j\|
                  +4\|(I-\mathsf S)T_j\|.
\end{equation*}
Here $\mathsf P_3$ averages the coordinate slot and the two newest
derivative slots. Repeated polynomial testing gives
\begin{equation*}
 \mathsf P_3T_j=
 \frac{\sqrt{\beta_j \cdot \beta_{j-1} \cdot \beta_{j-2}}}{3!}
 Q_3u^{j-3},\qquad
 \|\mathsf P_3T_j\|^2\le27 \cdot \overline c_3^2 \cdot \sigma^3.
\end{equation*}
The fresh adjacent swap uses $\chi_{j-2}$. The older swap uses
$\chi_{j-3}$ and is propagated once by a map of norm at most
$\sqrt{3}$, by Eq.~\eqref{eq:block_operator_norm}.
The defect identity and Eq.~\eqref{dt_eq_sort} therefore give
\begin{equation*}
 \|(I-\mathsf S)T_j\|
 \le6\sqrt{\beta_j/\lambda} \cdot 
       (\sqrt{\chi_{j-2}}+\sqrt{3\chi_{j-3}}).
\end{equation*}
Its square is at most $1296 \cdot (\chi_{j-2}+\chi_{j-3})$.
Squaring the recovery bound with coefficients $18$ and $32$ proves
Eq.~\eqref{eq:block_cubic_loss}, since $32\cdot1296<50000$.
\end{proof}

The cubic loss estimate controls the mass remaining after a block of the hierarchy. Combined with the compensated restart, it gives upper and lower bounds for powers of the inverse-gradient map.

\begin{lemma}[Mesoscopic block norm]\label{lem:mesoscopic_block_norm}
There are universal constants $R_0,C_*>0$ such that, for $R\ge R_0$,
\begin{equation}\label{eq:mesoscopic_block_norm}
 R-C_*/R\le\|\mathcal T^m\|_{\mathrm{op}}^{2/m}\le R
 ,\ \forall 1\le m\le\lfloor R\rfloor+2.
\end{equation}
The map $\mathcal T^m$ includes all ordered output indices.
\end{lemma}
\begin{proof}
We define $K:=512+4 \cdot A_0$, $M:=\lfloor R\rfloor$, and $p_*:=K \cdot \sigma^2$.
Increase $R_0$ so that
\begin{equation*}
 R_0\ge\max\{16,4 \cdot B_0,\sqrt{8K}\}.
\end{equation*}
Use the compensated unit start from
Proposition~\ref{prop:compensated_restart}. Its energy and first three
losses satisfy
\begin{equation*}
 \nu\le\frac\sigma{1-2\sigma^2}\le\sigma+3\sigma^3,
 \qquad P_3\le160\sigma^2.
\end{equation*}
Before a first exit of $P_N$ above $p_*$, every required normalizer is
at most $\nu/(1-p_*)\le3\sigma$, and its ratio to $\lambda$ is at
most $6$. The successor at an exit is still positive because its mass
is at least $1-p_*-\nu$. Thus Lemma~\ref{lem:block_cubic_loss} applies
at the exit as well as at every earlier prefix.
For $4\le N\le M$, it gives
\begin{equation*}
 P_N\le160\sigma^2+A_0 \cdot N \cdot \sigma^3+2 \cdot B_0 \cdot X_{N-2}.
\end{equation*}
The delayed budget bounds $X_{N-2}$ by $3\sigma^3+\sigma \cdot  P_N$.
Since $N\le\sigma^{-1}$ and $\sigma\le(4 \cdot B_0)^{-1}$, we obtain
\begin{equation}\label{eq:mesoscopic_mass_bound}
 P_N\le(323+2 \cdot A_0) \cdot \sigma^2<3K \cdot \sigma^2/4.
\end{equation}
The first three prefixes satisfy the same conclusion directly.
Consequently no exit occurs through $M$, and all these bounds hold on
the actual hierarchy.

The hierarchy can be expressed in terms of powers of $\mathcal T$.
Let $f$ be the centered unit extremizer in
Eq.~\eqref{eq:restricted_extremizer}, and define
\begin{equation*}
 W_k:=\mathcal T^kf,\qquad a_k:=\|W_k\|_2^2,
 \qquad b_k:=\|LW_k\|^2.
\end{equation*}
Then $a_0=1$, $a_1=R$, and $W_1/\sqrt R=\sqrt R\nabla f$.
For each nonzero $W_k$, define
$\eta_k:=a_k/\langle W_k,BW_k\rangle$.
The compensated start used above is exactly
\begin{equation}\label{eq:mesoscopic_alignment}
 F_0=\frac{\sqrt{\eta_2}B^{1/2}W_2}{\sqrt{a_2}},
 \qquad F_j=s_jB^{1/2}W_{j+2},\qquad s_j>0.
\end{equation}
The first identity follows by taking the inverse normalization of
$W_1/\sqrt R$, centering its gradient, and normalizing once more, as
in Proposition~\ref{prop:compensated_restart}. The second follows
inductively from Eq.~\eqref{dt_eq_inverse}.
In particular $e_j/v_j=\eta_{j+2}$ and
\begin{equation*}
 \frac{b_{j+2}}{a_{j+2}}=
 \frac{p_j}{v_j \cdot \eta_{j+2}},\qquad \eta_{j+2}\ge\lambda.
\end{equation*}
Equation~\eqref{eq:mesoscopic_mass_bound} therefore implies
\begin{equation}\label{eq:mesoscopic_raw_loss}
 \sum_{k=2}^{M+1}\frac{b_k}{a_k}
 \le\frac{K \cdot \sigma^2}{\lambda \cdot (1-p_*)}\le2K \cdot \sigma.
\end{equation}

For the first raw ratio after $a_1$, the same skew credit used in
Eq.~\eqref{eq:restricted_surviving_mass} gives
\begin{equation}\label{eq:mesoscopic_first_ratio}
 \frac{a_2}{a_1}\ge R-\frac{K_2}{4R}\ge R-2\sigma.
\end{equation}
Indeed, for $U:=W_1/\sqrt R$ the ratio equals
$\|\mathcal TU\|_2^2=\eta_1^{-1}-\|LU\|^2$;
Lemma~\ref{lem:skew_normalization_credit} and
Eq.~\eqref{eq:restricted_symmetric_quadratic} bound it below by
$e_U^{-1}-K_2/(4R)$.

Bochner also gives $\mathcal T\mathcal T^*\preceq B$ on vector
families: the form norm of $\mathcal Th$ is at most $\|h\|_2$.
Hence, for $k\ge1$,
\begin{equation*}
 a_{k+1}+b_k=\langle W_k,BW_k\rangle
 \ge\|\mathcal T^*W_k\|_2^2
 \ge a_k^2/a_{k-1}.
\end{equation*}
Thus the positive ratios $r_k:=a_k/a_{k-1}$ obey
$r_{k+1}\ge r_k-b_k/a_k$.
All terms through $M+2$ are positive by the finite hierarchy above and
Eq.~\eqref{eq:mesoscopic_first_ratio}.
Combining Eqs.~\eqref{eq:mesoscopic_raw_loss} and
\eqref{eq:mesoscopic_first_ratio} gives
$r_k\ge R-(2+2K)/R$ for $2\le k\le M+2$.
Multiply the ratios and use $r_1=R$ to obtain the lower bound in
Eq.~\eqref{eq:mesoscopic_block_norm} with $C_*:=2+2K$.
The upper bound follows from $\|\mathcal T\|_{\mathrm{op}}^2=R$.
\end{proof}

Averaging along a block with the norm bounds in Lemma~\ref{lem:mesoscopic_block_norm} reduces the initial losses and excess energy. It also bounds the propagation of defects through subsequent blocks.

\begin{proposition}[Averaged block startup]\label{prop:block_averaged_start}
There are universal constants $R_1,C_S,C_7,C_0>0$ such that the
following holds for $R\ge R_1$. We define
\begin{equation}\label{eq:block_radius}
 m:=\lfloor R\rfloor,\qquad z:=\|\mathcal T^m\|_{\mathrm{op}}^{2/m}.
\end{equation}
There is a centered unit finite starting family with energy $\nu$ and
first seven losses satisfying
\begin{equation}\label{eq:block_startup_bounds}
 \nu\le\frac{z^{-1}}{1-C_S \cdot z^{-3}},\qquad
 S_7:=\sum_{j=0}^6p_j\le C_7 \cdot z^{-3},\qquad
 D_0:=\max\{\nu-R^{-1},0\}\le C_0 \cdot z^{-3}.
\end{equation}
Moreover,
\begin{equation}\label{eq:block_global_propagation}
 \|\mathcal T^\ell\|_{\mathrm{op}}
 \le C_{\mathrm{blk}} \cdot z^{\ell/2},\ \forall \ell\ge0,\qquad
 C_{\mathrm{blk}}:=(R/z)^{(m-1)/2}\le2.
\end{equation}
\end{proposition}
\begin{proof}
The operator $\mathcal T$ is compact, since
$\mathcal T=DH^{-1/2}$ with $D$ bounded and $H^{-1/2}$ compact.
The positive lower bound in Lemma~\ref{lem:mesoscopic_block_norm}
shows that $\mathcal T^m$ is nonzero. Choose a centered unit top right
singular vector $f$ and define
\begin{equation*}
 w_j:=z^{-j/2}\mathcal T^jf,\qquad
 b_j:=\|w_j\|_2^2,\qquad \ell_j:=\|Lw_j\|^2.
\end{equation*}
Then $b_0=b_m=1$ and
\begin{equation*}
 b_{j+1}\le(R/z) \cdot b_j,\qquad
 b_j\ge(z/R)^{m-j},\ \forall 0\le j\le m.
\end{equation*}
Lemma~\ref{lem:mesoscopic_block_norm} gives
$\log(R/z)\le2 \cdot C_*/R^2$ at a fixed large threshold. Consequently
\begin{equation}\label{eq:block_orbit_envelope}
 b_j\le2,\ \forall 0\le j\le m+9,\qquad
 b_j\ge1/2,\ \forall 0\le j\le m.
\end{equation}
No positive lower bound beyond $m$ is required.

For $j\ge1$, define
\begin{equation*}
 D_j:=\|\nabla[(B-zI)w_j]\|_2^2.
\end{equation*}
Bochner gives $\|H^{1/2}w_j\|_2^2\le b_{j-1}/z$, while
$\langle w_j,Bw_j\rangle=zb_{j+1}+\ell_j$.
Expanding the square proves
\begin{equation}\label{eq:block_fixed_defect}
 D_j\le z(b_{j+1}-2b_j+b_{j-1})+\ell_j.
\end{equation}
We define $N:=m+8$. Summing and using the orbit bounds gives
\begin{equation}\label{eq:block_defect_telescope}
 \sum_{j=1}^ND_j
 \le\sum_{j=1}^N\ell_j+
 z[(b_{N+1}-b_N)+(1-b_1)]
 \le\sum_{j=1}^N\ell_j+2 \cdot C_*.
\end{equation}
For the last inequality,
$z(b_{N+1}-b_N)\le(R-z)b_N\le2 \cdot C_*/R$, and
$z(1-b_1)\le z[1-(z/R)^{m-1}]\le C_*$.

The full linear moments satisfy the following estimates. Each $w_j$, $j\ge1$, is
a genuine gradient in its newest derivative slot. The identity
\begin{equation}\label{eq:block_fixed_decomposition}
 w_{j+1}=\sqrt z P_+\nabla w_j
             +z^{-1/2}P_+\nabla[(B-zI)w_j]
\end{equation}
bounds its fresh adjacent swap by $2\sqrt{D_j/z}$.
A previously created swap propagates by a map of norm at most
$\sqrt{R/z}$. Polynomial testing gives
\begin{equation*}
 \mathsf P_2Q_1w_2=\frac{Q_2w_1}{2\sqrt z}.
\end{equation*}
Lemma~\ref{lem:startup_three_slot_frame} and
Eq.~\eqref{eq:block_fixed_decomposition} therefore imply
\begin{equation}\label{eq:block_second_moment}
 \ell_2\le\frac{2 \cdot K_2 \cdot b_1+8D_1}{z}.
\end{equation}
For $j\ge3$, use the coordinate slot and the three newest derivative
slots. Block recovery again has constant $3$, and its coherent term is
\begin{equation*}
 \mathsf P_3Q_1w_j=\frac{Q_3w_{j-2}}{3! \cdot z}.
\end{equation*}
The error in symmetrizing the derivative slots is at most
\begin{equation*}
 6z^{-1/2}(\sqrt{D_{j-1}}+\sqrt{R/z} \cdot \sqrt{D_{j-2}}).
\end{equation*}
Squaring as in Lemma~\ref{lem:block_cubic_loss}, and using $R/z\le2$,
gives
\begin{equation}\label{eq:block_raw_cubic_moment}
 \ell_j\le18 \cdot \overline c_3^2 \cdot z^{-2} \cdot b_{j-2}
          +10^4 \cdot z^{-1} \cdot (D_{j-1}+D_{j-2}),\ \forall j\ge3.
\end{equation}
These inequalities control the full linear moments, not only their
symmetric parts.

Let $D_{\mathrm{tot}}:=\sum_{j=1}^ND_j$. Since $\ell_1\le b_1\le2$,
Eqs.~\eqref{eq:block_second_moment} and
\eqref{eq:block_raw_cubic_moment} show
\begin{equation*}
 \sum_{j=1}^N\ell_j
 \le2+4 \cdot K_2/z+36 \cdot \overline c_3^2 \cdot N/z^2
                +(8+2\cdot10^4) \cdot D_{\mathrm{tot}}/z.
\end{equation*}
Here $N\le4z$ at a fixed large threshold. Absorbing the last term in
Eq.~\eqref{eq:block_defect_telescope} yields universal constants
$C_D,C_L$ such that
\begin{equation}\label{eq:block_linear_moment_sum}
 D_{\mathrm{tot}}\le C_D,\qquad
 \sum_{j=2}^{m+8}\ell_j\le C_L/z.
\end{equation}

Define the finite family and its normalized gradient by
\begin{equation*}
 Z:=\bigoplus_{j=0}^{m-1}w_j,\qquad
 S_0:=\|Z\|_2^2=\sum_{j=0}^{m-1}b_j,\qquad
 U:=\frac{\mathcal TZ}{\sqrt{z \cdot S_0}}.
\end{equation*}
Since $b_m=b_0$, we have $\|\mathcal TZ\|_2^2=z \cdot S_0$.
Thus $U$ is centered, has unit norm, and is a genuine gradient in each
output column. Bochner gives $e_U\le z^{-1}$.
For every $j\ge1$, quadratic testing also gives
\begin{equation*}
 \|\operatorname{Sym}(Lw_j)\|_{\mathrm{HS}}^2
 \le\frac{K_2 \cdot b_{j-1}}{4z}.
\end{equation*}
Use this for $j=1,2$ and Eq.~\eqref{eq:block_linear_moment_sum} for the
remaining indices. Since $S_0\ge m/2$ and $m$ is comparable to $z$,
\begin{equation}\label{eq:block_symmetric_moment}
 \|\operatorname{Sym}(LU)\|_{\mathrm{HS}}^2\le C_S \cdot z^{-2}.
\end{equation}
The symmetrization exchanges the coordinate slot and the newest gradient
slot within each output column.

Start the hierarchy at $F_0:=B^{1/2}Z/\|B^{1/2}Z\|_2$.
Its first successor under inverse normalization, after rescaling to unit norm,
is the inverse normalization of $U$. Apply
Lemma~\ref{lem:skew_normalization_credit} to this unit gradient and
then apply centering and inverse normalization once more.
The unit family $G_0:=F_2/\|F_2\|_2$ has energy
\begin{equation*}
 \nu\le\frac1{e_U^{-1}-\|\operatorname{Sym}(LU)\|_{\mathrm{HS}}^2}
 \le\frac1{z-C_S \cdot z^{-2}}.
\end{equation*}
This proves the energy assertion in Eq.~\eqref{eq:block_startup_bounds}.
All denominators are positive at the fixed threshold.

To identify its actual losses, define
\begin{equation*}
 S_k:=\sum_{j=k}^{m+k-1}b_j,\qquad
 C_k:=\sum_{j=k}^{m+k-1}\ell_j.
\end{equation*}
Every $F_k$ is a positive scalar multiple of $B^{1/2}\mathcal T^kZ$, and
\begin{equation*}
 \|\mathcal T^kZ\|_2^2=z^k \cdot S_k,\qquad
 \langle\mathcal T^kZ,B\mathcal T^kZ\rangle
       =z^k \cdot (z \cdot S_{k+1}+C_k).
\end{equation*}
Hence its loss fraction is exactly
\begin{equation}\label{eq:block_loss_fraction}
 \frac{p_k}{v_k}=\frac{C_k}{z \cdot S_{k+1}+C_k}.
\end{equation}
For $2\le k\le8$, Eq.~\eqref{eq:block_linear_moment_sum} gives
$C_k\le C_L/z$. The interval in $S_{k+1}$ includes $k+1,\ldots,m$,
so Eq.~\eqref{eq:block_orbit_envelope} gives
$S_{k+1}\ge m/4\ge z/8$ when $m\ge16$.
Thus $p_k/v_k\le8 \cdot C_L \cdot z^{-3}$ for these seven indices.
They are exactly the first seven losses of the unit restart $G_0$;
its masses are at most one. We may therefore take $C_7:=56 \cdot C_L$.
Every family used, including $F_9$, is nonzero because its window
contains positive interior terms from Eq.~\eqref{eq:block_orbit_envelope}.
All required domains follow from the same finite Bochner hierarchy.

Lemma~\ref{lem:mesoscopic_block_norm} and the energy bound imply
$D_0\le C_0 \cdot z^{-3}$: the difference $z^{-1}-R^{-1}$ is $O(z^{-3})$,
and the remaining energy error is $O(z^{-4})$.
Finally, for $\ell=qm+r$ with $0\le r<m$,
\begin{equation*}
 \|\mathcal T^\ell\|_{\mathrm{op}}
 \le\|\mathcal T^r\|_{\mathrm{op}} \cdot \|\mathcal T^m\|_{\mathrm{op}}^q
 \le z^{\ell/2} \cdot (R/z)^{r/2}.
\end{equation*}
This proves Eq.~\eqref{eq:block_global_propagation} after another fixed
increase of $R_1$.
\end{proof}

The frame estimate and the propagation bound for a whole block control the local losses of the averaged starting family. This is the form needed for the delayed comparison.

\begin{lemma}[Joint loss for an averaged block]\label{lem:block_joint_local_loss}
Use the unit starting family in Proposition~\ref{prop:block_averaged_start}.
On a prefix whose required normalizers are at most $B_*$, we define
$t:=\max\{1,B_* \cdot z\}$. For a dyadic $d\ge4$, define
\begin{align}
 \theta&:=324 \cdot C_F \cdot C_{\mathrm{blk}}^2 \cdot 
 \sum_{\substack{k<d\\k\ \mathrm{dyadic}}}
       k^6 \cdot (B_* \cdot t^3)^k \cdot c_k^2,\label{eq:block_joint_kernel}\\
 \Delta&:=S_7+2 \cdot C_F \cdot 
 \sum_{\substack{4\le k<d\\k\ \mathrm{dyadic}}}
       k^2 \cdot B_*^k \cdot c_k^2,\label{eq:block_joint_startup}\\
 \tau_d&:=C_F \cdot d \cdot B_*^d \cdot c_d^2.\label{eq:block_joint_terminal}
\end{align}
Then every justified prefix satisfies
\begin{equation}\label{eq:block_joint_local_loss}
 P_N\le\Delta+N \cdot \tau_d+
          \frac\theta\lambda \cdot  X_{\max\{N-2,0\}}.
\end{equation}
In particular Lemma~\ref{lem:delayed_block_comparison} applies with
$B_*:=\nu/(1-p)$ whenever its stated conditions hold.
\end{lemma}
\begin{proof}
The coherent projection identity
Eq.~\eqref{eq:joint_coherent_projections} is unchanged and gives
$\|\mathsf P_kT_j\|^2\le B_*^k \cdot c_k^2$.
A fresh adjacent swap has norm at most $2\sqrt{\chi/\lambda}$.
After $\ell$ further steps, Eq.~\eqref{eq:block_global_propagation}
bounds it by $2 \cdot C_{\mathrm{blk}} \cdot t^{\ell/2} \cdot \sqrt{\chi/\lambda}$.
The factor $C_{\mathrm{blk}}$ occurs once for the entire propagation.
These maps commute with permutations of all previously introduced slots.
Consequently the sorting argument in
Lemma~\ref{lem:restricted_joint_local_loss} gives
\begin{equation*}
 \|(I-\mathsf S_k)\mathsf P_kT_j\|
 \le6k \cdot C_{\mathrm{blk}} \cdot B_*^{k/2} \cdot c_k \cdot \lambda^{-1/2} \cdot 
       \sum_{h=1}^{3k-1}t^{(h-1)/2} \cdot \sqrt{\chi_{j-k-h}}.
\end{equation*}
This display is used when the degree-$k$ frame term occurs, so its
required derivative slots are present and the defect indices are
nonnegative. The inner convolution kernel has $\ell^1$ norm at most $3k \cdot t^{3k/2}$.
Square, sum in $j$, and multiply by the frame weight $C_F \cdot k^2$.
This gives Eq.~\eqref{eq:block_joint_kernel}.
In a prefix $P_N$, every defect index is at most $N-3$, because
$j\le N-1$ and $k,h\ge1$. Its contribution is therefore bounded by
$X_{\max\{N-2,0\}}$.

We sum the first seven actual losses separately. Degree four is available at
$j=7$; thereafter use the largest available dyadic degree, capped at $d$.
Each $k<d$ is used for $2k$ generations, with coherent contribution
at most $C_F \cdot k \cdot B_*^k \cdot c_k^2$ per generation. Theorem~\ref{thm:joint_uniform_frame}
now gives Eqs.~\eqref{eq:block_joint_startup}--\eqref{eq:block_joint_local_loss}.
For shorter prefixes, the initial losses and the nonnegative terminal term suffice.
\end{proof}

The block radius also controls polynomial coefficients: polynomial testing relates successive degrees to the same inverse-gradient map. We record this as a single parameter that bounds every degree and controls the Poincar\'e constant.

\begin{lemma}[A static block radius for every degree]\label{lem:block_radius_static}
There is a universal constant $H_0$ and a parameter $Z(\mu)$ for every
centered regular law of covariance at most $I$ such that
\begin{equation}\label{eq:block_static_radius}
 c_d(\mu)\le Z(\mu)^{(d-1)/2},\ \forall d\ge1,\qquad
 C_P(\mu)\le Z(\mu)+2,\qquad
 Z(\mu)\le\max\{H_0,C_P(\mu)\}.
\end{equation}
Whenever $Z(\mu)>H_0$, it equals the block radius $z$ in
Eq.~\eqref{eq:block_radius}, and Proposition~\ref{prop:block_averaged_start}
applies with universal constants.
\end{lemma}
\begin{proof}
Choose a universal $R_*$ above all preceding thresholds and large enough
that $R\ge R_*$ implies $z\ge R/2$ and
\begin{equation*}
 C_{\mathrm{seed}}^2 \cdot \log(e+R)^2\le R/2\le z.
\end{equation*}
Such a choice is possible because $\log(e+R)^2=o(R)$.
Set $H_0:=R_*+1$. For $R\ge R_*$, use $m:=\lfloor R\rfloor$ and $z$
as above and define $Z(\mu):=\max\{H_0,z\}$.
For $R<R_*$, define $Z(\mu):=H_0$.

Integration by parts against a degree-$(k+1)$ Appell tensor gives
\begin{equation}\label{eq:block_polynomial_testing}
 \frac{Q_{k+1}}{k+1}=\mathsf P_{k+1}Q_k\mathcal T,\ \forall k\ge1.
\end{equation}
Indeed its derivative is $(k+1)$ times the corresponding degree-$k$
Appell tensor, whose mean is zero; subtracting the mean gradient in
$\mathcal T$ changes no pairing. Iterating this identity $m$ times,
using contraction of every symmetrization and the factorial convention
for $c_d$, gives
\begin{equation}\label{eq:block_coefficient_recursion}
 c_d\le z^{m/2} \cdot c_{d-m},\ \forall d\ge m+1.
\end{equation}
For $1\le d\le m$, the fixed seed in
Eq.~\eqref{eq:joint_logarithmic_seed} gives
$c_d\le z^{(d-1)/2}$; at $d=1$ this is covariance normalization.
Repeatedly applying Eq.~\eqref{eq:block_coefficient_recursion} proves
the same bound at every degree, including the degree-one remainder.

When $R<R_*$, Lemma~\ref{lem:restricted_operator} gives
$C_P\le R+1\le H_0$. Applying ordinary Poincar\'e repeatedly to the
Appell derivatives yields $c_d\le C_P^{(d-1)/2}$, hence the first
claim in Eq.~\eqref{eq:block_static_radius} in this case as well.
For $R\ge R_*$, enlarge $R_*$ so that $C_*/R\le1$.
Then $C_P\le R+1\le z+2\le Z+2$.
The final inequality follows from $z\le R\le C_P$ in the large-$R$
case and from the definition in the remaining case.
The coefficient inequalities for $Z$ at each finite degree are sufficient for the
transfer below, even when $m$, $z$, and $Z$ vary discontinuously with the law.
\end{proof}

\subsection{From the block radius to polynomial coefficients}\label{sec:block_profile_transfer}

We transfer a curvature profile for the common radius $Z$ from Lemma~\ref{lem:block_radius_static} to a degree-dependent coefficient bound, with multiplicative losses whose product over the iteration depths is finite. Working directly with this radius keeps its fixed additive comparison with the Poincar\'e constant out of the repeated transfer. The proof uses an accumulated covariance metric to control all derivative orders along localization.

The posterior estimate below transfers a continuous bound on all polynomial coefficients through an affine normalization of a Gaussian posterior.

\begin{lemma}[Posterior Appell estimate]\label{lem:block_posterior_appell}
Fix a dimension range, either all dimensions or $1\le m\le n$.
Let $\Psi:(0,\infty)\to[1,\infty)$ be finite, continuous, and nondecreasing.
Suppose that, for every $a>0$, every centered regular law $\mu$ in this range with covariance at most $I$ and curvature at least $aI$ satisfies
$c_k(\mu)\le\Psi(a^{-1})^{(k-1)/2}$ for every integer $k\ge1$.
Let a law in the same dimension range have covariance $A\succ0$ and density proportional to
$\exp(-x^\top\Lambda x/2-V(x))$, where $\Lambda\succ0$ and $V$ is extended-valued and convex.
If $\delta>0$, $M\succeq A$, and $M\succeq\delta\Lambda^{-1}$, then every degree-$d$ polynomial $f$, $d\ge1$, satisfies
\begin{equation}\label{eq:block_posterior_appell}
 \sqrt{\Var[f]}
 \le\sum_{k=1}^d J^{k-1} \cdot 
       \|\E[D^kf]\|_{M^{\otimes k}},
 \qquad J:=\sqrt{\Psi(\delta^{-1})}.
\end{equation}
The expectations are under this law, and $\|v\|_K^2:=\langle v,Kv\rangle$.
The estimate also holds on finite output tensor direct sums.
\end{lemma}
\begin{proof}
Let $X$ have the law in the statement. The law of
$Y:=M^{-1/2}(X-\E[X])$ has covariance at most $I$ and curvature at least
$\delta I$. Indeed, inversion of the second matrix inequality gives
$\Lambda\succeq\delta M^{-1}$, so
$M^{1/2}\Lambda M^{1/2}\succeq\delta I$.
Applying the full Appell expansion in Eq.~\eqref{eq:polynomial_minkowski}
to the polynomial $f$ therefore gives Eq.~\eqref{eq:block_posterior_appell}.
We use $\|v\|_K^2:=\langle v,Kv\rangle$. The factorial $1/k!$ in the
Appell expansion is already included in $c_k$.

To justify this application to a nonsmooth posterior, apply the
convolution and rescaling argument of Lemma~\ref{lem:profile_inflation}
to the law of $Y$. Convolution with $N(0,\epsilon I)$ followed by division
by $\sqrt{1+\epsilon}$ gives a regular measure with covariance at most
$I$ and curvature at least
$\delta_\epsilon I$, where
$\delta_\epsilon:=\delta(1+\epsilon)/(1+\delta\epsilon)\to\delta$.
At each fixed degree, its polynomial moments and Appell coefficients
converge. Continuity of $\Psi$ passes Eq.~\eqref{eq:block_posterior_appell}
to the original posterior. This passes a finite-degree polynomial
inequality; no continuity of an auxiliary radius or block length is required.
The same inequality holds on a finite output tensor direct sum by
Minkowski's inequality.
\end{proof}

To use the posterior estimate along localization, we track derivative moments in a metric containing both the current covariance and its accumulated past values.

Apply Lemma~\ref{lem:construction_localization_process} to a compactly supported isotropic log-concave measure $\mu_0$. We use its notation $\Gamma_t$, $U_{i,t}$, and $S_{i,t}$. Fix $\tau,s>0$, and define
\begin{equation}\label{eq:block_accumulated_metric}
 \rho:=s/\tau,\qquad
 M_t:=\Gamma_t+\rho\int_0^t\Gamma_u\d u.
\end{equation}
For a degree-$d$ polynomial $f$ and $1\le j\le d$, we define
\begin{align*}
 h_j(t)&:=\E_t[D^jf],\\
 N_j(t)&:=\E_{\mathrm{loc}}[
       \langle h_j(t),M_t^{\otimes j}h_j(t)\rangle],\\
 L_j(t)&:=\E_{\mathrm{loc}}[\E_t[
       \langle D^jf-h_j(t),M_t^{\otimes j}(D^jf-h_j(t))\rangle]].
\end{align*}
These quantities use the accumulated metric $M_t$, rather than only
$\Gamma_t$. 
\begin{lemma}[Derivative hierarchy in the accumulated metric]\label{lem:block_metric_hierarchy}
Under the preceding assumptions, for $1\le j\le d$ and almost every $t\ge0$,
\begin{equation}\label{eq:block_metric_hierarchy}
 N_j'(t)\le(4j^2+j\rho) \cdot N_j(t)+9j^2 \cdot L_j(t).
\end{equation}
\end{lemma}
\begin{proof}
By definition, $M_0=I$ and $M_t\succeq\Gamma_t$.
The drift of $M_t$ is $(\rho-1)\Gamma_t\preceq\rho M_t$, and its noise
coefficients are the matrices $U_{i,t}$ from
Eq.~\eqref{eq:inverse_covariance_noise}. We define
\begin{equation*}
 C_t:=M_t^{-1/2}\Gamma_t^{1/2},\qquad
 \widetilde S_{i,t}:=M_t^{-1/2}U_{i,t}M_t^{-1/2}
                  =C_t S_{i,t}C_t^\top.
\end{equation*}
Since $C_tC_t^\top\preceq I$ and $C_t^\top C_t\preceq I$,
Eq.~\eqref{eq:inverse_covariance_noise} implies
\begin{equation}\label{eq:block_metric_noise}
 \sum_i\widetilde S_{i,t}^2
 =C_t(\sum_iS_{i,t}C_t^\top C_tS_{i,t})C_t^\top
 \preceq C_t(\sum_iS_{i,t}^2)C_t^\top
 \preceq8I.
\end{equation}

Write $K:=M_t^{\otimes j}$. The tensor product rule
and Eq.~\eqref{eq:block_metric_noise} bound its drift by
$[j\rho+8\binom j2]K$. After congruence by $K^{-1/2}$, its noise
coefficients are $V_i:=\sum_{h=1}^j\widetilde S_{i,t}^{(h)}$, where
the superscript selects a tensor slot. The operator Cauchy--Schwarz
inequality gives $\sum_iV_i^2\preceq8j^2I$.
Write
$\d h_j=\sum_i c_i\d\beta_{i,t}$, where
$c_i:=\E_t[(D^jf-h_j)\xi_{i,t}]$ and
$\xi_t:=\Gamma_t^{-1/2}(X-m_t)$.
The weighted Bessel inequality gives
\begin{equation*}
 Q_j:=\sum_i\langle c_i,Kc_i\rangle
 \le\E_t[\langle D^jf-h_j,K(D^jf-h_j)\rangle].
\end{equation*}
The cross variation between $K$ and $h_j$ is bounded by
\begin{equation*}
 2\sum_i\langle K^{1/2}h_j,V_iK^{1/2}c_i\rangle
 \le\langle h_j,Kh_j\rangle+8j^2 \cdot Q_j.
\end{equation*}
Including the quadratic variation of $h_j$, the full drift of
$\langle h_j,Kh_j\rangle$ is therefore at most
\begin{equation*}
 [j\rho+8\binom j2+1]\langle h_j,Kh_j\rangle
       +(1+8j^2) \cdot Q_j.
\end{equation*}
Since $8\binom j2+1\le4j^2$ and $1+8j^2\le9j^2$, taking expectations
proves Eq.~\eqref{eq:block_metric_hierarchy}. These calculations may
first be stopped on compact sets of localization parameters. The fixed
spatial support bounds the derivative tests and $\Gamma_t$, and the
integral in Eq.~\eqref{eq:block_accumulated_metric} is bounded on every
finite deterministic time interval. Bessel's inequality and
Eq.~\eqref{eq:block_metric_noise} then permit removal of the stops,
as in Lemma~\ref{lem:hierarchy}.
\end{proof}

The posterior estimate and derivative hierarchy bound polynomial coefficients in terms of the curvature profile for $Z$. The factor $1+r^{-2}$ in this bound has a finite product over the depths of the iteration.

\begin{proposition}[Polynomial coefficients from the block radius]
\label{prop:block_profile_coefficients}
There is a universal integer $r_0$ with the following property. Let
$r\ge r_0$ and $\Gamma\ge\sqrt{H_0}$. Suppose that every centered regular
measure with covariance at most $I$ in the chosen dimension range
satisfies
\begin{equation}\label{eq:block_profile_input}
 D^2W\succeq aI
 \quad\Longrightarrow\quad
 Z(\mu)\le\Gamma^2 \cdot \ell_r(a^{-1})^2
 ,\ \forall a>0.
\end{equation}
Then every centered log-concave measure with covariance at most $I$ in
that dimension range satisfies
\begin{equation}\label{eq:block_profile_coefficients}
 c_d(\mu)\le[(1+r^{-2}) \cdot \Gamma \cdot \ell_r(d)]^{d-1}
 ,\ \forall d\ge1.
\end{equation}
Here $H_0$ and $Z$ are the constant and parameter from
Lemma~\ref{lem:block_radius_static}. The dimension range is either all
dimensions or $1\le m\le n$, fixed throughout the argument.
\end{proposition}

\begin{proof}
We first record the polynomial consequence of a general profile. Suppose
$\Psi:(0,\infty)\to[H_0,\infty)$ is finite, continuous, and nondecreasing,
and bounds $Z(\mu)$ by $\Psi(a^{-1})$ whenever a centered regular measure
has covariance at most $I$ and curvature at least $aI$.
Lemma~\ref{lem:block_radius_static} then gives
$c_k(\mu)\le\Psi(a^{-1})^{(k-1)/2}$ for every $k\ge1$.

Lemma~\ref{lem:block_posterior_appell} therefore gives Eq.~\eqref{eq:block_posterior_appell} for the inflated posterior metrics used below.

Apply Lemma~\ref{lem:construction_localization_process} to a compactly supported isotropic log-concave measure $\mu_0$ and fix $\tau,s>0$. For a degree-$d$ polynomial $f$, use the accumulated metric $M_t$ and derivative moments $N_j,L_j$ defined before Lemma~\ref{lem:block_metric_hierarchy}. That lemma supplies Eq.~\eqref{eq:block_metric_hierarchy}.

Suppose, for the moment, that every smaller degree and every law in the
dimension range satisfies $c_k\le A^{k-1}$, where $A\ge1$ is fixed.
Whitening a posterior and applying its Appell expansion to $D^jf$ gives
\begin{equation}\label{eq:block_derivative_fluctuations}
 \sqrt{L_j(t)}\le\sum_{k=1}^{d-j}A^{k-1} \cdot \sqrt{N_{j+k}(t)}
 ,\ \forall 1\le j<d.
\end{equation}
Here the old output slots carry the weight $M_t$. Whitening initially
puts $\Gamma_t$ in each new derivative slot; the inequality
$\Gamma_t\preceq M_t$ enlarges those weights to $M_t$.
Minkowski's inequality on the localization probability space and the
complete output tensor direct sum proves the displayed bound. All
weights in this estimate are taken at the same time.

Take $f:=P_d^{\mu_0}[T_d]$ for a nonzero symmetric tensor $T_d$, and define
\begin{equation*}
 Q:=(d!)^2 \cdot \|T_d\|_{\mathrm{HS}}^2,\qquad
 w_l:=A^{2\max(l-1,0)},\qquad
 M_*(t):=\max_{1\le j\le d}\frac{N_j(t)}{Q \cdot w_{d-j}}.
\end{equation*}
For $1\le k\le l$, one has
$A^{k-1} \cdot \sqrt{w_{l-k}}\le\sqrt{w_l}$. Thus
Eq.~\eqref{eq:block_derivative_fluctuations} implies
$L_j(t)\le d^2 \cdot Q \cdot w_{d-j} \cdot M_*(t)$.
Eq.~\eqref{eq:appell_centering} and $M_0=I$ give
$N_j(0)=0$ for $j<d$ and $N_d(0)=Q$; also $L_d=0$.
Integrating Eq.~\eqref{eq:block_metric_hierarchy}, taking the maximum,
and applying Lemma~\ref{lem:gronwall} therefore give
\begin{align}
 M_*(t)&\le e^{(13d^4+d\rho) \cdot t},
       \label{eq:block_weighted_maximum}\\
 N_j(t)&\le9d^4 \cdot t \cdot  e^{(13d^4+d\rho) \cdot t} \cdot Q \cdot w_{d-j}
       ,\ \forall j<d,\label{eq:block_zero_initial}\\
 N_d(t)&\le Q \cdot e^{(4d^2+d\rho) \cdot t}.
       \label{eq:block_top_derivative}
\end{align}
For Eq.~\eqref{eq:block_weighted_maximum}, use
$4j^2+j\rho+9j^2 \cdot d^2\le13d^4+d\rho$.
For Eq.~\eqref{eq:block_zero_initial}, use the zero initial value and
integrate with the factor $e^{-(4j^2+j\rho)t}$. The last inequality uses
$L_d=0$.

Choose $0<\epsilon<1$, and set
\begin{equation}\label{eq:block_transfer_times}
 \tau:=\epsilon/d^6,\qquad s:=\epsilon/d,\qquad
 \delta:=s\tau=\epsilon^2/d^7.
\end{equation}
The estimate for the inverse precision matrix in the proof of
Theorem~\ref{thm:polynomial} gives
\begin{equation*}
 \delta\Lambda_\tau^{-1}
 \preceq\delta \cdot \tau^{-2}\int_0^\tau\Gamma_t\d t
 =\rho\int_0^\tau\Gamma_t\d t
 \preceq M_\tau.
\end{equation*}
For completeness, the underlying estimate follows for every vector $v$
from $\Lambda_\tau=\int_0^\tau\Gamma_t^{-1}\d t$ and
\begin{equation*}
 \int_0^\tau v^\top\Gamma_t v\d t
       -\tau^2v^\top\Lambda_\tau^{-1}v
 =\int_0^\tau
       |\Gamma_t^{1/2}v-\tau\Gamma_t^{-1/2}\Lambda_\tau^{-1}v|^2\d t
 \ge0.
\end{equation*}
We may therefore use Eq.~\eqref{eq:block_posterior_appell} at time $\tau$
with
$J:=\sqrt{\Psi(d^7/\epsilon^2)}$.
Since $(13d^4+d\rho) \cdot \tau\le14\epsilon$ and
$(4d^2+d\rho) \cdot \tau\le5\epsilon$, Minkowski's inequality and
Eqs.~\eqref{eq:block_zero_initial}--\eqref{eq:block_top_derivative} give,
provided $J\le A$,
\begin{equation}\label{eq:block_terminal_appell}
 \sqrt{\E_{\mathrm{loc}}[\Var_\tau[f]]}
 \le\sqrt Q \cdot 
       [e^{5\epsilon/2} \cdot J^{d-1}
          +3d^3 \cdot \sqrt\tau \cdot  e^{7\epsilon} \cdot A^{d-2}].
\end{equation}
The first term is the top Appell term $k=d$. Every lower term has the
factor $A^{(k-1)+(d-k-1)}=A^{d-2}$; this includes $k=d-1$, since
$w_1=1$. There are at most $d-1$ lower terms.

For the process in Lemma~\ref{lem:construction_localization_process}, the posterior variance satisfies
\begin{equation*}
 \frac{\d}{\d t}\E_{\mathrm{loc}}[\Var_t[f]]
 =-\E_{\mathrm{loc}}[|\operatorname{Cov}_t(f,\xi_t)|^2]
 \ge-\E_{\mathrm{loc}}[\Var_t[f]]
\end{equation*}
by Bessel's inequality. Hence
$\Var_{\mu_0}[f]\le e^\tau \cdot \E_{\mathrm{loc}}[\Var_\tau[f]]$.
Taking square roots in this estimate, using
Eq.~\eqref{eq:block_terminal_appell}, and taking the supremum over $T_d$
prove the finite-degree bound
\begin{equation}\label{eq:block_finite_transfer}
 c_d(\mu_0)\le e^{3\epsilon} \cdot J^{d-1}
                 +3\sqrt\epsilon \cdot  e^{8\epsilon} \cdot A^{d-2}.
\end{equation}

We specialize to $\Psi(x):=\Gamma^2 \cdot \ell_r(x)^2$ and prove
Eq.~\eqref{eq:block_profile_coefficients} by induction on $d$,
simultaneously over all isotropic laws in the dimension range. Degree
one is isotropy. For $d\ge2$, we define
\begin{equation*}
 \alpha:=r^{-2},\qquad
 \epsilon:=2^{-20} \cdot \alpha^2,\qquad
 A:=(1+\alpha) \cdot \Gamma \cdot \ell_r(d).
\end{equation*}
Monotonicity of $\ell_r$ makes the smaller-degree induction hypotheses
sufficient for Eq.~\eqref{eq:block_finite_transfer} with this $A$.
The first logarithm satisfies
\begin{equation*}
 g(d^7/\epsilon^2)\le[7+2\log(1/\epsilon)] \cdot g(d).
\end{equation*}
Eq.~\eqref{eq:iterated_log_contraction}, applied to the remaining
$r-1$ compositions, gives
\begin{equation*}
 \log\frac{\ell_r(d^7/\epsilon^2)}{\ell_r(d)}
 \le2^{-(r-1)} \cdot \log[7+2\log(1/\epsilon)]
 \le\alpha/8
\end{equation*}
for every $r\ge r_0$, after fixing a sufficiently large universal $r_0$.
This choice is uniform in $d$ and $\Gamma$.
Since $\log(1+\alpha/4)\ge\alpha/8$ for $0<\alpha\le1$,
\begin{equation*}
 J\le(1+\alpha/4) \cdot \Gamma \cdot \ell_r(d),\qquad
 \frac JA\le e^{-3\alpha/8}.
\end{equation*}
The second inequality uses
$\log(1+\alpha)-\log(1+\alpha/4)\ge3\alpha/8$.
Dividing the two terms of Eq.~\eqref{eq:block_finite_transfer} by
$A^{d-1}$, respectively, bounds them by
\begin{equation*}
 e^{3\epsilon-3\alpha/8}\le e^{-\alpha/4}\le1-\alpha/8,
 \qquad
 \frac{3\sqrt\epsilon \cdot  e^{8\epsilon}}{A}\le\alpha/128.
\end{equation*}
Their sum is less than one, so the induction closes for compactly
supported isotropic measures.

At each fixed degree, conditioning on growing balls, centering, and
whitening removes compact support as in Theorem~\ref{thm:polynomial}.
Moments through order $2d$ and the coefficients of the corresponding
Appell polynomials converge. Whitening and contraction in every tensor
slot then give the conclusion for covariance at most $I$; singular covariance
is treated on its supporting subspace. This completes the simultaneous induction at each fixed degree.
\end{proof}

\subsection{Proof of iterative logarithmic curvature refinement}\label{sec:sharp_iterated_profiles}

We combine the coefficient transfer and the comparison for averaged blocks to add one logarithm at each depth while proving the amplitude recurrence stated at the beginning of the section. The induction is carried out for the common radius $Z$; Lemma~\ref{lem:block_radius_static} converts it to a Poincar\'e bound only at the end. This yields one sequence of amplitudes valid at every depth and for every curvature lower bound.

\begin{proof}[Proof of Proposition~\ref{prop:sharp_iterated_curvature}]
We first prove the profile for the parameter $Z$ of
Lemma~\ref{lem:block_radius_static}. We use the common floor $H_0$ from
that lemma, so that
\begin{equation}\label{eq:block_profile_comparison}
 C_P(\mu)\le Z(\mu)+2,
 \qquad Z(\mu)\le\max\{H_0,C_P(\mu)\}.
\end{equation}
Let $r_0$ be the universal depth in
Proposition~\ref{prop:block_profile_coefficients}. Suppose that, uniformly
over centered regular measures of covariance at most $I$ and curvature
at least $aI$, we have
\begin{equation}\label{eq:block_induction_hypothesis}
 Z(\mu)\le\widetilde\Gamma_r^2 \cdot \ell_r(a^{-1})^2,
 \qquad \widetilde\Gamma_r\ge G \cdot (r+1)^{1/6}
 ,\ \forall r\ge r_0.
\end{equation}
The universal constant $G$ will be chosen below. In particular, we require
$G^2>2 \cdot H_0$, so the hypothesis of the coefficient transfer is satisfied.
We define
\begin{align}
 s&:=r+1,& Q&:=(1+r^{-2}) \cdot \widetilde\Gamma_r,&
 \delta&:=1/(32s),&p&:=1/(4096s),\notag\\
 b&:=(1-2p)^{-1},&q&:=e^{-\delta},&
 M_0&:=2^{44} \cdot (1+C_F),&
 k_0&:=\lceil64\delta^{-1} \cdot \log(M_0/\delta)\rceil.
 \label{eq:block_refinement_parameters}
\end{align}
The scalar $\delta$ controls the tail of the coefficient sum. By
Proposition~\ref{prop:block_profile_coefficients} and
Eq.~\eqref{eq:joint_logarithmic_seed}, the same coefficients satisfy
\begin{equation}\label{eq:block_two_coefficients}
 c_k\le[Q \cdot \ell_r(k)]^{k-1},
 \qquad c_k\le[C_{\mathrm{seed}} \cdot \log(e+k)]^{k-1}
 ,\ \forall k\ge1.
\end{equation}

Fix a dyadic integer $d\ge4$. We seek a contradiction in the regime
\begin{equation}\label{eq:block_large_radius_regime}
 Z(\mu)>b^4 \cdot e^\delta \cdot  Q^2 \cdot \ell_r(d)^2.
\end{equation}
The right side exceeds $H_0$, so $Z(\mu)=z$ is the actual block radius
from Lemma~\ref{lem:block_radius_static}. We write
$R:=R_{\mathrm{res}}(\mu)$, $u:=z^{-1}$, and $\sigma_R:=R^{-1}$.
Then
\begin{equation}\label{eq:block_inverse_radius_bounds}
 b^4 \cdot u \cdot Q^2 \cdot \ell_r(d)^2<q,
 \qquad u\le G^{-2} \cdot s^{-1/3},
 \qquad u^3\le G^{-6} \cdot s^{-1}.
\end{equation}
Proposition~\ref{prop:block_averaged_start} supplies a unit finite starting
family of energy $\nu$ such that
\begin{equation}\label{eq:block_refinement_startup}
 \nu\le\frac{u}{1-C_S \cdot u^3},
 \qquad S_7:=\sum_{j=0}^6p_j\le C_7 \cdot u^3,
 \qquad D_0:=\max\{\nu-\sigma_R,0\}\le C_0 \cdot u^3.
\end{equation}
Also $\lambda=C_P(\mu)^{-1}\ge1/(z+2)$. We enlarge $G$ so that
$C_S \cdot u^3\le p$ and all the fixed large-radius hypotheses hold. Before a
first exit of accumulated mass $p$, every needed normalizer is bounded by
\begin{equation}\label{eq:block_refinement_normalizers}
 B_*\le\frac{u}{(1-p)^2}\le bu,
 \qquad \max\{1,B_* \cdot z\}\le b.
\end{equation}
The factor governing propagation through a block is at most two. We may therefore use
the following upper bounds in Lemma~\ref{lem:block_joint_local_loss}:
\begin{align}
 \theta_d&:=1296 \cdot C_F \cdot \sum_{\substack{k<d\\k\ \mathrm{dyadic}}}
                 k^6 \cdot (b^4 \cdot u)^k \cdot c_k^2,\notag\\
 \Delta_d&:=C_7 \cdot u^3+2 \cdot C_F \cdot 
           \sum_{\substack{4\le k<d\\k\ \mathrm{dyadic}}}
                 k^2 \cdot (bu)^k \cdot c_k^2,\notag\\
 \tau_d&:=C_F \cdot d \cdot (bu)^d \cdot c_d^2.
 \label{eq:block_refinement_kernel}
\end{align}
The coherent sum begins at degree four because the first seven losses
are included in $S_7$.

The two sums admit bounds uniform in $r$ and $d$. We define
\begin{equation*}
 C_M:=2048 \cdot [\log(32 \cdot M_0)+1]+1.
\end{equation*}
For $s\ge1$, Eq.~\eqref{eq:block_refinement_parameters} gives
$k_0\le C_M \cdot s^2$, and hence
\begin{equation*}
 \frac{\log(e+k_0)}{s^{1/6}}
 \le\log(e+C_M)+12/e.
\end{equation*}
Thus one fixed choice of $G$ ensures
$\widetilde\Gamma_r\ge2 \cdot C_{\mathrm{seed}} \cdot \log(e+k_0)$ for all $r\ge r_0$.
For $k\le k_0$, Eqs.~\eqref{eq:block_two_coefficients} and
\eqref{eq:block_inverse_radius_bounds} imply
\begin{align*}
 b^4 \cdot u \cdot [C_{\mathrm{seed}} \cdot \log(e+k)]^2&\le1/4,\\
 (b^4 \cdot u)^k \cdot c_k^2&\le b^4 \cdot u \cdot 4^{-(k-1)},\\
 (bu)^k \cdot c_k^2&\le b^3 \cdot u^3 \cdot C_{\mathrm{seed}}^4 \cdot 
           \log(e+k)^4 \cdot  4^{-(k-3)},\ \forall k\ge4.
\end{align*}
The integer sums with weights $k^6$ and $k^2$, respectively, converge.
They bound the low-degree error by a universal multiple of $u$ and the
low-degree coherent contribution by a universal multiple of $u^3$.

For $k_0<k<d$, the first bound in
Eq.~\eqref{eq:block_two_coefficients} and monotonicity of $\ell_r$ give
\begin{equation*}
 (b^4 \cdot u)^k \cdot c_k^2\le b^4 \cdot u \cdot  q^{k-1},
 \qquad (bu)^k \cdot c_k^2\le q^k.
\end{equation*}
The elementary bounds for the exponential tails are
\begin{align*}
 \sum_{k>k_0}k^6 \cdot e^{-\delta k}
   &\le2^{30} \cdot \delta^{-7} \cdot e^{-\delta k_0/2},\\
 \sum_{k>k_0}k^2 \cdot e^{-\delta k}
   &\le2^{10} \cdot \delta^{-3} \cdot e^{-\delta k_0/2},
 \qquad e^{-\delta k_0/2}\le(\delta/M_0)^{32}.
\end{align*}
The choice of $M_0$ makes the high-degree contribution to $\theta_d$ at
most $u$ and the high-degree contribution to $\Delta_d$ at most $p/512$.
Replacing a dyadic sum by the corresponding integer sum only increases
it. Consequently there are universal constants $C_\theta,C_\Delta$ with
\begin{equation}\label{eq:block_split_bounds}
 \theta_d\le C_\theta \cdot  u,
 \qquad \Delta_d\le C_\Delta \cdot  u^3+p/512.
\end{equation}
For example, since $b\le2$, we may take
\begin{align*}
 C_\theta&:=1+20736 \cdot C_F \cdot \sum_{k\ge1}k^6 \cdot  4^{-(k-1)},\\
 C_\Delta&:=C_7+16 \cdot C_F \cdot C_{\mathrm{seed}}^4 \cdot 
        \sum_{k\ge4}k^2 \cdot \log(e+k)^4 \cdot  4^{-(k-3)}.
\end{align*}
In particular the error bound still includes the leading factor $u$.

We choose $G$ once so that all preceding requirements and
\begin{align}
 \frac{2 \cdot C_\theta}{G^2}&\le1/8,&
 \frac{C_\Delta+2 \cdot C_\theta \cdot  C_0}{G^6}&\le2^{-20},\notag\\
 \frac{C_S}{G^6}&\le1/4096,&
 G^{-2}+4096 \cdot C_0 \cdot G^{-6}&\le1,&G^4&\ge2 \cdot C_0
 \label{eq:block_fixed_g_conditions}
\end{align}
hold. We also require $G^2\ge16$. All these are compatible lower bounds
on one universal constant. Since $z\ge2$,
\begin{equation*}
 \frac{\theta_d \cdot D_0}{\lambda}\le2 \cdot C_\theta \cdot  C_0 \cdot u^3,
 \qquad \frac{\sigma_R}{\lambda}\le2.
\end{equation*}
Eqs.~\eqref{eq:block_inverse_radius_bounds},
\eqref{eq:block_split_bounds}, and \eqref{eq:block_fixed_g_conditions}
therefore yield
\begin{equation}\label{eq:block_refinement_absorption}
 \frac{\theta_d \cdot \sigma_R}{\lambda}\le1/8,
 \qquad \Delta_d+\frac{\theta_d \cdot D_0}{\lambda}\le p/64,
 \qquad \sigma_R+D_0/p\le1.
\end{equation}
For the middle inequality, the low-degree part is at most $p/256$ and
the high-degree part is at most $p/512$. For the last inequality we use
\begin{equation*}
 \sigma_R+D_0/p\le u+4096s \cdot C_0 \cdot u^3
 \le G^{-2}+4096 \cdot C_0 \cdot G^{-6}.
\end{equation*}
The fixed large-radius threshold also ensures $\nu\le1/4$ and
$D_0\le\sigma_R$: indeed $\sigma_R\ge u/2$ and $C_0 \cdot u^2\le1/2$.
Thus Lemma~\ref{lem:delayed_block_comparison} applies and gives
\begin{equation}\label{eq:block_terminal_radius_comparison}
 a<12(\sigma_R+D_0/p) \cdot \tau_d
    \le12 \cdot C_F \cdot d \cdot (bu)^d \cdot c_d^2.
\end{equation}
This comparison uses the original one-step restricted radius in its
energy budget. Its terminal power is $z^d$; an additional factor $u$ is
not required below.

For $0<a\le1$, we define $h:=\log(e+a^{-1})$ and choose a dyadic $d$ with
\begin{equation}\label{eq:block_terminal_degree}
 C_T \cdot s^2 \cdot h\le d<2 \cdot C_T \cdot s^2 \cdot h.
\end{equation}
Here $C_T$ is a sufficiently large universal constant. Taking
$C_F\ge1$ and $C_T\ge2^{24} \cdot (1+C_F)$, enlarged if necessary, ensures
\begin{equation}\label{eq:block_terminal_exponent}
 \log(12 \cdot C_F \cdot d/a)\le\delta d.
\end{equation}
Indeed the left side is at most
$\log(24 \cdot C_F \cdot C_T)+2\log s+2h$, while the right side is at least
$(C_T/32) \cdot sh$. The inequality follows for every $s,h\ge1$ from the fixed
condition $C_T/32\ge\log(24 \cdot C_F \cdot C_T)+4$. This choice also gives $d\ge4$.

We define $L_*:=Q^2 \cdot \ell_r(d)^2\ge1$. Under
Eq.~\eqref{eq:block_large_radius_regime}, we have $b^4 \cdot u \cdot L_*<e^{-\delta}$.
Together with $c_d^2\le L_*^{d-1}$, Eq.~\eqref{eq:block_terminal_radius_comparison}
would imply
\begin{equation*}
 a<12 \cdot C_F \cdot d \cdot L_*^{-1} \cdot (bu \cdot L_*)^d
   \le12 \cdot C_F \cdot d \cdot  e^{-\delta d}\le a,
\end{equation*}
a contradiction. Hence
\begin{equation}\label{eq:block_one_step_profile}
 Z(\mu)\le b^4 \cdot e^\delta \cdot  Q^2 \cdot \ell_r(d)^2.
\end{equation}
If $a\ge1$, the curvature bound gives $C_P(\mu)\le1$, and therefore
$Z(\mu)=H_0$. This case is already covered by the large universal
profile. The floor and bounded-radius cases were excluded from
Eq.~\eqref{eq:block_large_radius_regime} before the hierarchy was used.

It remains to control the amplitude. The estimates for iterated logarithms and
Eq.~\eqref{eq:block_terminal_degree} give
\begin{equation}\label{eq:block_profile_distortion}
 \ell_r(d)\le\zeta_r \cdot \ell_{r+1}(a^{-1}),
 \qquad
 \zeta_r:=\exp(2^{-r} \cdot \log(2 \cdot C_T \cdot (r+1)^2)).
\end{equation}
The product of the $\zeta_r$ is finite. Also
\begin{equation}\label{eq:block_amplitude_increment}
 \log(b^4 \cdot e^\delta)\le16p+\delta
    =\frac9{256s}
    \le\frac13 \cdot \log(1+1/s).
\end{equation}
For the last inequality, $\log(1+1/s)\ge1/(2s)$ suffices. We define
\begin{equation}\label{eq:block_amplitude_recurrence}
 B_{r_0}:=1,
 \qquad B_{r+1}:=B_r \cdot (1+r^{-2})^2 \cdot \zeta_r^2,
 \qquad \widetilde\Gamma_r^2:=G^2 \cdot (r+1)^{1/3} \cdot B_r
 ,\ \forall r\ge r_0.
\end{equation}
The $B_r$ are bounded above and below by positive universal constants.
Lemma~\ref{lem:coarse_profiles} gives a coarse ordinary profile with
constant $\widehat\Gamma_{r_0}$ at the fixed depth $r_0$.
By Eq.~\eqref{eq:block_profile_comparison}, enlarging $G$ initializes
Eq.~\eqref{eq:block_induction_hypothesis} at that depth. This is one
fixed enlargement, independent of every future depth. Eqs.~\eqref{eq:block_one_step_profile}--\eqref{eq:block_amplitude_recurrence} then give
\begin{equation*}
 Z(\mu)\le b^4 \cdot e^\delta \cdot (1+r^{-2})^2 \cdot \zeta_r^2 \cdot 
       \widetilde\Gamma_r^2 \cdot \ell_{r+1}(a^{-1})^2
   \le\widetilde\Gamma_{r+1}^2 \cdot \ell_{r+1}(a^{-1})^2.
\end{equation*}
The induction closes.

For the finitely many depths $1\le r<r_0$, we define
\begin{equation*}
 \widetilde\Gamma_r:=\max\{\sqrt{H_0},\widehat\Gamma_r\},
\end{equation*}
using the coarse profiles from Lemma~\ref{lem:coarse_profiles}.
Eq.~\eqref{eq:block_profile_comparison} then gives the same $Z$ profile
at these depths. We finally define
\begin{equation*}
 \Gamma_r^2:=\widetilde\Gamma_r^2+2,\ \forall r\ge1.
\end{equation*}
Since $\ell_r\ge1$, Eq.~\eqref{eq:block_profile_comparison} proves the
ordinary curvature profile for every $r\ge1$. Increasing a universal
constant to cover the finitely many smaller depths gives
$\Gamma_r^2\le C_0 \cdot (r+1)^{1/3}$, as required in
Eq.~\eqref{eq:sharp_iterated_curvature}. The additive two is absorbed once,
after the $Z$-profile induction.
\end{proof}

\Needspace{22\baselineskip}
\section{Dimension Bounds from Gaussian Localization}\label{sec:initial_dimension_bound}

This section gives our first dimensional KLS bound of order $(1+\log^*(n+2))^{1/6}$ for arbitrary isotropic log-concave measures, without a positive curvature assumption on the initial measure. Gaussian localization creates a posterior at time of order $1/\log(en)$ to which the curvature estimate of Proposition~\ref{prop:sharp_iterated_curvature} applies. Preserving a fixed fraction of a Poincar\'e witness's variance lets us transfer that estimate back to the initial measure.

\begin{theorem}[Initial dimension bounds]\label{thm:initial_dimension_bound}
There are universal constants $C,C'>0$ such that every isotropic log-concave probability measure $\mu$ on $\mathbb R^n$ satisfies, for every integer $r\ge1$,
\begin{equation*}
 C_P(\mu)\le C \cdot (r+1)^{1/3} \cdot \ell_r(L_n)^2,\qquad
 \psi_n\le C' \cdot (r+1)^{1/6} \cdot \ell_r(L_n).
\end{equation*}
Here $L_n:=\log(en)$, and $\ell_r$ is the iterated logarithm defined in Section~\ref{sec:iterated_curvature}. The constants are independent of both $r$ and $n$.
In particular, after enlarging them if necessary,
\begin{equation*}
 C_P(\mu)\le C \cdot (1+\log^*(n+2))^{1/3},\qquad
 \psi_n\le C' \cdot (1+\log^*(n+2))^{1/6}.
\end{equation*}
\end{theorem}

Throughout this section, $\mu_t$ denotes the conditional distribution produced by Gaussian localization, $m_t$ its mean, and $A_t$ its covariance matrix. The expectation $\E_{\mathrm{loc}}$ is over the localization process. We use the normalized bounded function and constant $B_0$ from Eq.~\eqref{tg:eq:witness}.

\subsection{Gaussian posterior identities}\label{sec:gaussian_localization}\label{tg:sec:localization}

We derive the posterior identities for the Gaussian observation construction from~\cite[Section~6.2 and Lemma~54]{klartag_lehec_2025_survey}. These identities will control the covariance fluctuations and the amount of the witness's variance retained after localization.

\begin{lemma}[Gaussian posterior identities]\label{lem:posterior_identities}
Let $\mu$ be a regular isotropic measure satisfying Eq.~\eqref{tg:eq:regular}.
Let $X\sim\mu$ be independent of a standard Brownian motion $B_t$, define $Y_t:=tX+B_t$, and let $\mathcal F_t$ be the usual augmentation of the observation filtration generated by $(Y_s)_{0\le s\le t}$.
The conditional distribution $\mu_t$ of $X$ given $\mathcal F_t$ is
\begin{equation*}
 \d\mu_t(x)=Z_t^{-1}\exp(Y_t\cdot x-t|x|^2/2)\d\mu(x),\qquad
 Z_t:=\int\exp(Y_t\cdot x-t|x|^2/2)\d\mu(x).
\end{equation*}
Write $m_t:=\E_{\mu_t}[X]$ and $A_t:=\operatorname{Cov}_{\mu_t}(X)$.
Then $W_t:=Y_t-\int_0^t m_s\d s$ is a standard Brownian motion in $(\mathcal F_t)$.
For every Borel $q$ of polynomial growth, $\E_{\mu_t}[q]$ is a square-integrable martingale and
\begin{equation}\label{eq:posterior_martingale}
 \d\E_{\mu_t}[q]=\operatorname{Cov}_{\mu_t}(q,X)\cdot\d W_t.
\end{equation}
In particular, $\d m_t=A_t\d W_t$ and
\begin{equation}\label{tg:eq:covsde}
 \d A_t=\sum_j U_{j,t}\d W_{j,t}-A_t^2\d t,\qquad
 U_{j,t}:=\E_{\mu_t}[(X_j-m_{j,t})(X-m_t)(X-m_t)^\top].
\end{equation}
Every posterior has a positive smooth log-concave density with potential Hessian $D^2W+tI$, bounded between $(a+t)I$ and $(b+t)I$.
\end{lemma}
\begin{proof}
For fixed $x$, the distribution of $(sx+B_s)_{0\le s\le t}$ has density $\exp(x\cdot Y_t-t|x|^2/2)$ relative to Wiener measure on paths up to time $t$.
Indeed, multiplying the Gaussian likelihood ratios of the independent increments along any finite partition gives this endpoint expression. Cylinder sets generate the path sigma-field, so the same identity holds on that sigma-field.
Integrating against $\mu(\d x)$ and applying Bayes' formula proves the posterior formula and $\E_{\mu_t}[q]=\E[q(X)\mid\mathcal F_t]$.
Strong convexity gives Gaussian tails for $\mu$, hence all polynomial moments are finite; the conditional form of Jensen's inequality therefore makes these martingales square integrable.

For $0\le s<t$, Brownian increments after $s$ are independent of $X$ and of the Brownian path up to $s$.
Moreover, the tower property gives $\E[m_r\mid\mathcal F_s]=\E[X\mid\mathcal F_s]=m_s$ for $r\ge s$.
Consequently $\E[W_t-W_s\mid\mathcal F_s]=0$.
All terms are integrable because $\E[|m_r|^2]\le\E[|X|^2]$ by the conditional form of Jensen's inequality.
The process $W_t$ is continuous and its quadratic covariation is $tI$, since subtracting the finite-variation drift does not change quadratic covariation.
L\'evy's characterization (Lemma~\ref{lem:martingale_tools}) now shows that $W_t$ is Brownian in the observation filtration.

For the stochastic differential, define
\begin{equation*} I_q(t,y):=\int q(x)\exp(y\cdot x-t|x|^2/2)\d\mu(x). \end{equation*}
Gaussian tails permit differentiation under this integral locally uniformly in $(t,y)$, including a neighborhood of $t=0$.
In particular $\partial_t I_q+\Delta_y I_q/2=0$ and $\nabla_y I_q=I_{qx}$, where the vector integral $I_{qx}$ is defined componentwise.
Since $Z_t=I_1(t,Y_t)$, It\^o's formula gives $\d I_q(t,Y_t)=\nabla_y I_q(t,Y_t)\cdot\d Y_t$ and $\d Z_t=Z_t m_t\cdot\d Y_t$.
The quotient rule, including its quadratic-covariation terms, yields
\begin{equation*}
 \d(I_q(t,Y_t)/Z_t)
 =\operatorname{Cov}_{\mu_t}(q,X)\cdot(\d Y_t-m_t\d t)
 =\operatorname{Cov}_{\mu_t}(q,X)\cdot\d W_t.
\end{equation*}
The first equality uses $I_{qx}/Z_t-(I_q/Z_t)m_t=\operatorname{Cov}_{\mu_t}(X,q)$; the second uses the definition of $W_t$.
These calculations may first be stopped when $|Y_t|$ reaches $R$.
For every finite $T$, conditional Cauchy--Schwarz followed by Cauchy--Schwarz and the conditional form of Jensen's inequality gives
\begin{equation*} \E_{\mathrm{loc}}[\int_0^T|\operatorname{Cov}_{\mu_t}(q,X)|^2\d t]\le T \cdot (\E_\mu[q^4] \cdot \E_\mu[|X|^4])^{1/2}<\infty. \end{equation*}
Here conditional Cauchy--Schwarz bounds the squared covariance by $\E_{\mu_t}[q^2] \cdot \E_{\mu_t}[|X|^2]$, and the other two inequalities bound the expectation of this product by the displayed fourth moments.
Continuity of $Y_t$ on finite intervals and the It\^o isometry now let $R\to\infty$, proving Eq.~\eqref{eq:posterior_martingale} for the unbounded tests as well.

Apply Eq.~\eqref{eq:posterior_martingale} to $x_i$ and $x_ix_k$.
The first choice gives $\d m_t=A_t\d W_t$.
Writing $Q_t:=\E_{\mu_t}[XX^\top]$ and $A_t=Q_t-m_tm_t^\top$, the product rule contributes the drift $-\d m_t\d m_t^\top=-A_t^2\d t$.
For each Brownian coordinate $j$, its remaining coefficient is
\begin{equation*}
 \E_{\mu_t}[XX^\top(X_j-m_{j,t})]-(A_te_j)m_t^\top-m_t(A_te_j)^\top
 =U_{j,t},
\end{equation*}
where $e_j$ is the $j$th coordinate vector; expanding $X=(X-m_t)+m_t$ proves the equality.
This establishes Eq.~\eqref{tg:eq:covsde}.
Finally, the posterior potential is $W(x)-Y_t\cdot x+t|x|^2/2+\log Z_t$, which gives the asserted Hessian bounds.
\end{proof}

\subsection{Covariance control}

We show that the posterior covariance stays in a fixed neighborhood of the identity with high probability up to a sufficiently small constant multiple of $1/\log(en)$. Its lower bound preserves the curvature gained from localization under whitening, while its upper bound controls the transformed witness's energy.

\begin{lemma}[Covariance control]\label{tg:lem:stop}
Let $\tau$ be the first exit of $A_t$ from $(I/2,2I)$. For
$0<t\le1/16$,
\begin{equation}\label{tg:eq:exit}
 \Pr[\tau\le t]\le2n \cdot \exp(-1/(2048t)).
\end{equation}
\end{lemma}
\begin{proof}
For any centered log-concave $Y$ with positive covariance $A$,
applying Lemma~\ref{lem:quadratic_forms} after whitening gives,
for every symmetric $D$,
\begin{equation*}
 \Var[Y^\top D Y]\le8\|A^{1/2}DA^{1/2}\|_\mathrm{HS}^2
 \le8\|A\|_\mathrm{op}^2 \cdot \|D\|_\mathrm{HS}^2.
\end{equation*}
We define $S_z:=\E[(z\cdot Y)YY^\top]$. Trace duality and Cauchy--Schwarz give
\begin{equation*}
 |\langle D,S_z\rangle_\mathrm{HS}|
 =|\E[(z\cdot Y)(Y^\top D Y-\tr[DA])]|\le
 \sqrt{8z^\top Az} \cdot \|A\|_\mathrm{op} \cdot \|D\|_\mathrm{HS}.
\end{equation*}
Taking the supremum over symmetric $D$ of unit norm shows
$\|S_z\|_\mathrm{HS}^2\le8\|A\|_\mathrm{op}^2 \cdot z^\top Az$.
Symmetry of the third moment tensor gives
$\sum_j|U_jz|^2=\|S_z\|_\mathrm{HS}^2$.
Apply this to each centered posterior to conclude
$\sum_jU_{j,s}^2\preceq64I$ before $\tau$.

We define $M_u:=\int_0^{u\wedge\tau}\sum_jU_{j,s} \d W_{j,s}$.
Fix $\theta>0$. For symmetric matrices $M,T$,
\begin{equation*}
 D^2\tr[e^{\theta M}][T,T]\le
 \theta^2\tr[e^{\theta M}T^2].
\end{equation*}
To verify it, diagonalize $M$. The coefficient of $T_{ij}^2$ on the
left is $\theta^2$ times the logarithmic mean of
$e^{\theta m_i}$ and $e^{\theta m_j}$; this is at most their arithmetic
mean, which gives the right side after summation.
It\^o's formula now makes
$Z_u:=e^{-32\theta^2u}\tr[e^{\theta M_u}]$ a nonnegative local
supermartingale, hence a supermartingale by Lemma~\ref{lem:martingale_tools}. At the first upper eigenvalue
crossing of level $r$ before time $t$, $Z_u\ge e^{\theta r-32\theta^2t}$.
Optional stopping, $Z_0=n$, optimization in $\theta>0$, and the same
argument for $-M$ yield
\begin{equation*}
 \Pr[\sup_{u\le t}\|M_u\|_\mathrm{op}\ge r]
 \le2n \cdot \exp(-r^2/(128t)).
\end{equation*}
On $\tau\le t$, Eq.~\eqref{tg:eq:covsde} gives
$M_\tau=A_\tau-I+\int_0^\tau A_s^2\d s$.
An upper exit forces an eigenvalue at least one. A lower exit, tested
against its unit eigenvector, forces an eigenvalue at most
$-1/2+4t\le-1/4$. Use $r:=1/4$ to conclude.
\end{proof}

\subsection{Proof of the dimension bounds}\label{tg:sec:witness}

We transfer the curvature estimate to the initial measure by selecting a posterior that retains a fixed fraction of the witness's variance and has controlled covariance. Choosing the refinement depth gives the stated log-star bounds for regular measures. Approximation then removes the regularity and positive curvature assumptions.

\begin{proof}[Proof of Theorem~\ref{thm:initial_dimension_bound}]

\emph{The bound for regular measures.}
Let $\mu$ be regular and isotropic, with potential $W$. We define $k:=C_P(\mu)$.
Choose $f$ from Lemma~\ref{lem:bounded_witness},
so that Eq.~\eqref{tg:eq:witness} holds.
Let $\tau$ be the exit time from Lemma~\ref{tg:lem:stop}.
Lemma~\ref{lem:posterior_identities} shows that the posterior expectations of $f$ and $f^2$ are bounded martingales.
For $V_s:=\Var_{\mu_{s\wedge\tau}}[f]$, Eq.~\eqref{eq:posterior_martingale}, It\^o's formula, and optional stopping (Lemmas~\ref{lem:ito_calculus}--\ref{lem:martingale_tools}) give
\begin{equation*}
 F(t):=\E_{\mathrm{loc}}[V_t]
 =1-\int_0^t\E_{\mathrm{loc}}
 [\mathbf1_{\{s<\tau\}}|\operatorname{Cov}_{\mu_s}(X,f)|^2]\d s.
\end{equation*}
Before exit, directional Cauchy--Schwarz gives
$|\operatorname{Cov}_{\mu_s}(X,f)|^2\le\|A_s\|_\mathrm{op} \cdot  V_s\le2V_s$.
Therefore $F'(t)\ge-2F(t)$ almost everywhere, and
\begin{equation}\label{tg:eq:survival}
 F(t)\ge e^{-2t},\qquad0\le V_t\le B_0^2.
\end{equation}
Choose
\begin{equation*}
 c_0:=\frac1{2048 \cdot \log(16 \cdot B_0^2)},
 \qquad t:=\frac{c_0}{\log(en)}.
\end{equation*}
Since $B_0\ge1$ and $\log(en)\ge1$, we have $0<t\le c_0<1/16$.
Writing $q:=1/(2048 \cdot c_0)=\log(16 \cdot B_0^2)>1$, Lemma~\ref{tg:lem:stop}, in the form of Eq.~\eqref{tg:eq:exit}, gives
$\Pr[\tau\le t]\le2n \cdot  e^{-1/(2048t)}=2 \cdot e^{-q} \cdot n^{1-q}\le1/(8 \cdot B_0^2)$.
The equality uses $t=c_0/\log(en)$, and the last inequality uses $q>1$ and $e^{-q}=1/(16 \cdot B_0^2)$.
Also $e^{-2t}\ge1-2t\ge7/8$, first by convexity of the exponential and then by $t\le1/16$.
Using Eq.~\eqref{tg:eq:survival} and subtracting the contribution of the exit event, we obtain
\begin{equation*}
 \E_{\mathrm{loc}}[\mathbf1_{\{\tau>t\}}\Var_{\mu_t}[f]]\ge3/4.
\end{equation*}
There is a realization with $\tau>t$ and
$\Var_{\mu_t}[f]>1/2$. Fix this realization and its mean $m$ and covariance
$A$, and let $\eta$ be the distribution of $A^{-1/2}(X-m)$ under $\mu_t$.
The whitened potential has Hessian
\begin{equation*}
 A^{1/2}(D^2W+tI)A^{1/2}\succeq(t/2)I,
\end{equation*}
and it has a finite upper Hessian bound.
The measure $\eta$ is regular and isotropic, so Proposition~\ref{prop:sharp_iterated_curvature} applies with curvature at least $t/2$.
On the other hand, $z\mapsto f(A^{1/2}z+m)$ has variance greater than $1/2$ and Dirichlet energy at most
$\|A\|_{\mathrm{op}} \cdot \operatorname{Lip}(f)^2\le2 \cdot B_0^2/k$.
Hence $k\le4 \cdot B_0^2 \cdot C_P(\eta)$.
For every integer $r\ge1$, Eq.~\eqref{eq:sharp_iterated_curvature} gives
\begin{equation*}
 k\le4 \cdot B_0^2 \cdot \Gamma_r^2 \cdot \ell_r(2/t)^2
 =4 \cdot B_0^2 \cdot \Gamma_r^2 \cdot \ell_r((2/c_0) \cdot L_n)^2.
\end{equation*}
Since $2/c_0\ge1$, Eq.~\eqref{eq:iterated_log_scaling} and the bound on $\Gamma_r$ in Eq.~\eqref{eq:sharp_iterated_curvature} give
\begin{equation}\label{eq:uniform_dimension}
 C_P(\mu)\le4 \cdot B_0^2 \cdot (2/c_0)^2 \cdot \Gamma_r^2 \cdot \ell_r(L_n)^2
 \le C \cdot (r+1)^{1/3} \cdot \ell_r(L_n)^2.
\end{equation}
The constant $C$ is independent of $r$ and $n$.
This transfer applies uniformly to every regular isotropic initial measure.

We choose the iteration depth so that the factor involving the iterated logarithm is bounded by a universal constant.
For $y\ge4$,
\begin{equation}\label{eq:shifted_log_step}
 g(y+1)=\log y+\log(1+(e+1)/y)\le\log y+1.
\end{equation}
The inequality follows from $\log(1+t)\le t$ and $(e+1)/y<1$.
Let $m$ be the least number of ordinary logarithms needed to bring $n+2$ to at most four.
Since $L_n\le1+\log(n+2)$, induction using Eq.~\eqref{eq:shifted_log_step} gives
$\ell_{m-1}(L_n)\le5$ when $m\ge1$.
Indeed, every ordinary iterate before the last is greater than four, so the estimate applies at every induction step.
When $m=0$, one has $L_n<5$ directly.
Also $g$ maps $[0,5]$ into $[0,5]$.
Thus
\begin{equation}\label{eq:chosen_depth}
 r:=\max\{1,m-1\}\le\log^*(n+2),\qquad \ell_r(L_n)\le5.
\end{equation}
Here $m\le\log^*(n+2)$ and $\log^*(n+2)\ge1$ because $n\ge1$.
Since $r+1\le1+\log^*(n+2)$, substitution into Eq.~\eqref{eq:uniform_dimension} gives
$C_P(\mu)\le C \cdot (1+\log^*(n+2))^{1/3}$ after enlarging the universal constant.

Equation~\eqref{eq:uniform_dimension} holds for every regular isotropic measure with one constant independent of $r$ and $n$.
Fix $n$ and a finite integer $r\ge1$.
For an arbitrary isotropic log-concave measure $\mu$, take the regular isotropic sequence from Lemma~\ref{lem:regular_isotropic_approximation}.
An isotropic log-concave measure is full-dimensional and hence has a density.
Lemma~\ref{lem:poincare_limit}, with $K:=C \cdot (r+1)^{1/3} \cdot \ell_r(L_n)^2$, passes the same scalar Poincar\'e inequality to $\mu$ and proves that every locally Lipschitz test of finite energy belongs to $L^2(\mu)$.

Lemma~\ref{lem:cheeger_poincare_comparison}, in the form of Eq.~\eqref{tg:eq:cheeger}, now gives
$\psi_\mu\le\sqrt{\pi \cdot  C} \cdot (r+1)^{1/6} \cdot \ell_r(L_n)$.
Taking the supremum over isotropic log-concave measures proves the asserted family of bounds for $\psi_n$.
Finally choose the finite depth from Eq.~\eqref{eq:chosen_depth}.
Then $\ell_r(L_n)\le5$ and $r\le\log^*(n+2)$, giving
$C_P(\mu)\le C \cdot (1+\log^*(n+2))^{1/3}$ and $\psi_n\le C' \cdot (1+\log^*(n+2))^{1/6}$ after enlarging the universal constants.
The chosen depth is finite for each dimension.
\end{proof}

%% file: 30_upper_cont.tex
\section{Repeated Height Reduction}
\label{sec:outer_refinement}
\label{sec:height_reduction}
\label{sec:repeated_height_reduction}

We now use the polynomial bounds from the inner refinement to repeat
the whole procedure, applying another log-star at each outer round.
Higher-order blocks reduce the cost of increasing the inner depth,
while the coefficient transfer carries the improvement from one round
to the next. The resulting curvature and coefficient bounds incur one
universal multiplicative cost per round, with starting depths independent
of the number of rounds.

We use $\ell_r:=g^{\circ r}$, where $g(x):=\log(e+x)$, and write
$\mathfrak t(x):=1+\log^*(x+2)$ and
$\mathfrak t_j:=\mathfrak t^{\circ j}$, with $\mathfrak t_0(x):=x$.
Here $\log^*x$ is the number of successive natural logarithms needed
to reach a value at most one. We call $\mathfrak t(x)$ the logarithmic
height of $x$; apart from the fixed shifts, it counts the logarithms
needed to reduce $x$ to a bounded value. The normalized polynomial
coefficients are $c_d:=\sqrt{K_d}/d!$. The integer inner depth $r$ counts
applications of $g$ in iterative logarithmic curvature refinement, while $j$ counts the
outer rounds of height reduction.
The starting depths have the form
$r_Q:=\max\{r_*,\lceil C_b \cdot \mathfrak t(Q)\rceil\}$, with universal
constants chosen in Lemma~\ref{lem:repeated_height_distortion}.

\begin{proposition}[Repeated height reduction]
\label{prop:repeated_height_outer_profiles}
There are universal constants $A_1\ge1$ and $C_*\ge1$ with the following
property. We define $A_j:=A_1 \cdot C_*^{j-1}$ for $j\ge1$.
For every $a>0$, every centered regular measure of covariance at most
$I$ and curvature at least $aI$ satisfies
\begin{equation}\label{eq:repeated_height_outer_profile}
 C_P\le A_j \cdot \mathfrak t_j(Q)^{1/3} \cdot (r+1)^{1/Q} \cdot
                \ell_r(a^{-1})^2
 ,\ \forall Q\ge3\text{ odd},\ r\ge r_Q.
\end{equation}
The depths $r_Q$ are exactly those in
Eq.~\eqref{eq:second_logstar_profile_depth}; they do not depend on
$j$. There is also a universal $C_e$ such that, for every $j\ge1$,
all centered log-concave measures of covariance at most $I$ satisfy
\begin{equation}\label{eq:repeated_height_outer_coefficients}
 c_d\le[C_e \cdot \sqrt{A_j} \cdot \mathfrak t_{j+1}(d)^{1/6}]^{d-1}
 ,\ \forall d\ge1.
\end{equation}
\end{proposition}

The coefficient estimate in this proposition makes another repetition
possible: it supplies polynomial bounds at every degree, to which we
can apply the curvature comparison again. We first derive the
coefficient bound from Proposition~\ref{prop:sharp_iterated_curvature},
then establish the block estimates that make this argument uniform in
the number of repetitions. The proof of
Proposition~\ref{prop:repeated_height_outer_profiles} concludes the section.

\subsection{Turning curvature profiles into coefficient bounds}
\label{sec:quantitative_order_profiles}

To begin the outer iteration, we turn the full inner curvature profile
into coefficient bounds at every degree. Choosing the inner depth
separately for each degree gives the first log-star coefficient seed.
We also establish the transfer for a general static coefficient bound,
so that it remains available when the block length and order change.

\begin{lemma}\label{lem:iterated_coefficient_seed}
There is a universal constant $C_6$ such that every centered log-concave
measure of covariance at most $I$ satisfies
\begin{equation}\label{eq:iterated_coefficient_seed}
 c_d\le[C_6 \cdot (1+\log^*(d+2))^{1/6}]^{d-1}
 ,\ \forall d\ge1.
\end{equation}
\end{lemma}
\begin{proof}
The case $d=1$ follows from the covariance bound. For $d\ge2$, let $m$
be the least number of ordinary logarithms needed to bring $d+2$ to at
most four. For $y\ge4$,
\begin{equation*}
 g(y+1)=\log y+\log(1+(e+1)/y)\le\log y+1.
\end{equation*}
Since $g(d)\le1+\log(d+2)$, induction gives $\ell_m(d)\le5$ when
$m\ge1$. If $m=0$, then $d\le2$. The map $g$ preserves $[0,5]$,
so $\ell_r(d)\le5$ for $r:=\max\{r_*,m,1\}$, where $r_*$ is the
universal threshold in Proposition~\ref{prop:constant_log_specialization}.
Proposition~\ref{prop:sharp_iterated_curvature} and
Eq.~\eqref{eq:constant_unweighted_log} now give
\begin{equation*}
 c_d\le[5(1+r^{-2}) \cdot \Gamma_r]^{d-1}
     \le[C \cdot (r+1)^{1/6}]^{d-1}.
\end{equation*}
Here $r+1\le C \cdot (1+\log^*(d+2))$. This proves the assertion for isotropic
measures. Whitening and contraction of each coefficient tensor give the
case of covariance at most $I$, as in Lemma~\ref{lem:coarse_profiles}; a
singular covariance is treated on its supporting subspace.
\end{proof}

The following direct coefficient transfer has constants independent of
the chosen block length and order.

\begin{proposition}[Transfer from a static coefficient bound]
\label{prop:static_profile_coefficients}
There is a universal integer $r_*$ with the following property. Let
$H\ge1$, and assign to each centered regular measure $\mu$ of covariance
at most $I$ a number $Z(\mu)\ge H$ such that
\begin{equation*}
 c_d(\mu)\le Z(\mu)^{(d-1)/2},\ \forall d\ge1.
\end{equation*}
Suppose that $r\ge r_*$, $\Gamma\ge\sqrt H$, and, uniformly over this
class, curvature at least $aI$ implies
$Z(\mu)\le\Gamma^2 \cdot \ell_r(a^{-1})^2$. Then every centered log-concave
measure of covariance at most $I$ satisfies
\begin{equation}\label{eq:static_profile_coefficients}
 c_d\le[(1+r^{-2}) \cdot \Gamma \cdot \ell_r(d)]^{d-1}
 ,\ \forall d\ge1.
\end{equation}
The constant $r_*$ is independent of $H$ and of the assignment $Z$.
\end{proposition}
\begin{proof}
We apply the proof of Proposition~\ref{prop:block_profile_coefficients}
with $\Psi(x):=\Gamma^2 \cdot \ell_r(x)^2$. Its use of the particular block
radius is confined to the implication
\begin{equation*}
 \text{curvature at least }aI
 \quad\Longrightarrow\quad
 c_d\le\Psi(a^{-1})^{(d-1)/2}
 ,\ \forall d\ge1.
\end{equation*}
This implication holds by the present hypotheses, and $\Psi(x)\ge H$
for every $x>0$. Thus the posterior Appell estimate in
Lemma~\ref{lem:block_posterior_appell}, the derivative hierarchy in
Lemma~\ref{lem:block_metric_hierarchy}, and the degree induction in the
proof of Proposition~\ref{prop:block_profile_coefficients} apply unchanged. In particular the
choices of the localization time and of the small parameter depend only
on $r$ and the polynomial degree, not on $H$ or $Z$.

For completeness, the application to a nonsmooth posterior passes only
coefficient inequalities at finitely many degrees. After whitening with an inflated
covariance matrix, Gaussian convolution and rescaling give regular laws
with covariance at most $I$ and curvature tending to the required value.
Their moments through the required finite order converge. Continuity of
$\Psi$ passes the displayed coefficient bound to the posterior, exactly
as in Lemma~\ref{lem:block_posterior_appell}.
No continuity of $Z$, or of an integer block length defining it, is used.
The same conditioning and affine contraction at the end of that proof
establish the stated class of measures.
\end{proof}

\subsection{From symmetry errors to a block comparison}
\label{sec:joint_block_estimates}
\label{sec:quantitative_block_geometry}

Higher-order blocks require estimates for the mass lost when derivative
tensors are centered, including the part missed by symmetric polynomial
tests. We control that part through adjacent swaps, then use the resulting
estimates to initialize a block and compare its cumulative losses with
its energy budget. This prepares the curvature comparison for higher
block orders.

For the block estimates below the law is centered and regular, its covariance
is at most $I$, and $D^2W\succeq aI$. We write $B:=H^{-1}$,
$R:=R_{\mathrm{res}}$, and $\mathcal T:=P_+\nabla B$. We define
\begin{equation}\label{eq:quarter_order_seed}
 G_*(k):=\max\{1,C_6 \cdot (1+\log^*(k+2))^{1/6}\}.
\end{equation}
Lemma~\ref{lem:iterated_coefficient_seed} gives
$c_k\le G_*(k)^{k-1}$. For an odd integer $p\ge3$, let $d_p$ be the
least dyadic integer at least $p+2$, and define
\begin{equation}\label{eq:quarter_order_coherent_constant}
 p+3\le d_p\le2p+2,\qquad
 \mathsf A_p:=C_F \cdot d_p \cdot G_*(d_p)^{2d_p-2}.
\end{equation}
All tensor families below include every ordered output index.

We first isolate the permutation representation generated by tensors
that are symmetric in all but one slot.

\begin{lemma}[Point-stabilizer reduction]\label{lem:point_stabilizer_reduction}
Let $N\ge2$, let $\rho$ be a finite-dimensional real orthogonal representation of the permutations of $N$ letters, and let $K$ be a point stabilizer. We define $M:=\operatorname{Fix}(K)$ and let $W$ be the span of the translates of $M$.
Then $W$ is equivariantly isometric to an invariant subspace of $\mathbb R^N\otimes M$, on which the group acts by permuting the natural coordinates.
Every point-stabilizer average vanishes on $W^\perp$.
\end{lemma}
\begin{proof}
Choose $g_i$ sending the distinguished letter to $i$. The map
\begin{equation*}
 \Phi:\mathbb R^N\otimes M\longrightarrow W,\qquad
 \Phi(e_i\otimes v):=\rho(g_i)v
\end{equation*}
is well-defined and intertwines the natural permutation representation
with $\rho$. Its kernel is invariant. The polar factor of its restriction
to the orthogonal complement of the kernel is an equivariant isometry
onto $W$, since $\Phi^*\Phi$ commutes with every permutation. Thus an
inequality in permutation operators that holds on Hilbert-valued natural
coordinates also holds on $W$. On $W^\perp$, every point-stabilizer
average vanishes, because the corresponding fixed spaces are translates
of $M$.
\end{proof}

The reduction to a point-stabilizer fixed space makes the permutation
inequality accessible through the natural representation. The next
estimate converts the gap for all transpositions into errors for adjacent
transpositions.

\begin{lemma}[Permutation gaps]\label{lem:permutation_path_gaps}
Let $N\ge2$, and let $\rho$ be a finite-dimensional real orthogonal representation of the permutations of $1,\ldots,N$. Let $\mathsf S$ be the full average. We define $\tau_i:=(i,i+1)$ for $1\le i<N$ and $\Lambda_N:=\sum_{i<j}(I-\rho((ij)))$.
Then
\begin{equation}\label{eq:complete_transposition_gap}
 \Lambda_N\succeq N(I-\mathsf S).
\end{equation}
For every vector $Y$, one also has
\begin{equation}\label{eq:adjacent_path_poincare}
 \|(I-\mathsf S)Y\|^2\le\frac{N^2}4
       \sum_{t=1}^{N-1}\|(I-\rho(\tau_t))Y\|^2.
\end{equation}
\end{lemma}
\begin{proof}
For a point stabilizer, use the spaces $M,W$ from Lemma~\ref{lem:point_stabilizer_reduction}. We first prove Eq.~\eqref{eq:complete_transposition_gap} by induction on $N$.
For $N=2$ this follows directly from the two eigenspaces of the
transposition. On natural coordinates,
$\Lambda_N=N(I-\mathsf S)$, so the same identity holds on $W$.
On $W^\perp$, induction applied to each of the $N$ point stabilizers
gives a lower bound $(N-1)I$. Each transposition belongs to $N-2$ of
these stabilizers. Their sum therefore gives
$(N-2)\Lambda_N\succeq N(N-1)I$ on $W^\perp$, which proves the
claim for $N\ge3$.

Since a transposition is an orthogonal involution,
Eq.~\eqref{eq:complete_transposition_gap} implies
\begin{equation*}
 \|(I-\mathsf S)Y\|^2
 \le\frac1{2N}\sum_{i<j}\|(I-\rho((ij)))Y\|^2.
\end{equation*}
Number the letters $1,\ldots,N$. The word
\begin{equation*}
 (ij)=\tau_i\cdots\tau_{j-2}\tau_{j-1}
                   \tau_{j-2}\cdots\tau_i
\end{equation*}
has length $2(j-i)-1$ and uses each adjacent transposition at most twice.
Telescoping and Cauchy--Schwarz give
\begin{equation*}
 \|(I-\rho((ij)))Y\|^2
 \le4(j-i)\sum_{t=i}^{j-1}\|(I-\rho(\tau_t))Y\|^2.
\end{equation*}
For each $t$, the sum of $j-i$ over $i\le t<j$ is
$Nt(N-t)/2\le N^3/8$. Consequently Eq.~\eqref{eq:adjacent_path_poincare} follows.
\end{proof}

These permutation estimates recover the symmetric part from errors
between adjacent derivative slots. We combine recovery and symmetrization
in one inequality that also applies when the tensor has several output
slots. Its constants will remain uniform as the block construction is
repeated.

\begin{lemma}[Recovery from adjacent errors]\label{lem:combined_path_recovery}
Let $k\ge2$, and let $\rho$ be a finite-dimensional real orthogonal
representation of the permutations of $0,1,\ldots,k$. Let $\mathsf P$
average permutations of $0,\ldots,k-1$. We define $\tau_i:=(i,i+1)$ for
$1\le i<k$. Then
\begin{equation}\label{eq:combined_path_recovery}
 \|Y\|^2\le4k^2\|\mathsf PY\|^2+
 k^3\sum_{i=1}^{k-1}\|(I-\rho(\tau_i))Y\|^2.
\end{equation}
\end{lemma}
\begin{proof}
Apply Lemma~\ref{lem:point_stabilizer_reduction} with the subgroup permuting $1,\ldots,k$, and use its decomposition $W\oplus W^\perp$.
On $W^\perp$, its average vanishes, so
Eq.~\eqref{eq:adjacent_path_poincare} proves
Eq.~\eqref{eq:combined_path_recovery}, even with coefficient $k^2/4$
in place of $k^3$. It remains to check natural coordinates
$x:=(x_0,\ldots,x_k)$, which may take values in an auxiliary Hilbert
space. We define
\begin{equation*}
 a_0:=\frac1k\sum_{i=0}^{k-1}x_i,\qquad
 d_i:=x_i-x_{i+1},\qquad E_0:=\sum_{i=1}^{k-1}\|d_i\|^2.
\end{equation*}
Then $\|\mathsf Px\|^2=k\|a_0\|^2+\|x_k\|^2$ and
$x_0=ka_0-(k-1)x_k-\sum_{i=1}^{k-1}i d_i$. Since
$\sum_{i=1}^{k-1}i^2<k^3/3$, we have
\begin{equation*}
 \|x_0\|^2\le3k^2\|a_0\|^2+3k^2\|x_k\|^2+k^3E_0,
 \qquad
 \sum_{i=1}^k\|x_i\|^2\le2k\|x_k\|^2+2k^2E_0.
\end{equation*}
For $k\ge2$, their sum is at most
$4k^2\|\mathsf Px\|^2+2k^3E_0$.
Finally, $\sum_i\|(I-\rho(\tau_i))x\|^2=2E_0$.
This proves the inequality on natural coordinates. The equivariant
isometry transfers it to $W$, and the orthogonal decomposition
$W\oplus W^\perp$ completes the proof.
\end{proof}

For the higher-order recovery and initialization estimates below, $\mu$ is centered and regular,
$\operatorname{Cov}(\mu)\preceq I$, and $D^2W\succeq aI$.
We use $B:=H^{-1}$, $R:=R_{\mathrm{res}}$, and
$\mathcal T:=P_+\nabla B$, with every ordered output index included in
each power. We define
\begin{equation*}
 G_{\mathrm{seed}}(k):=\max\{1,C_{\mathrm{seed}} \cdot \log(e+k)\},\qquad
 \alpha_k:=4k^2 \cdot G_{\mathrm{seed}}(k)^{2k-2}.
\end{equation*}
Thus Eq.~\eqref{eq:joint_logarithmic_seed} gives $c_k\le G_{\mathrm{seed}}(k)^{k-1}$.

Applying the recovery estimate to a derivative tensor gives the frame at
an arbitrary degree. Its polynomial term is controlled by the coefficient
seed, while its remaining terms record the symmetry defects.

\begin{lemma}[A raw frame at arbitrary degree]\label{lem:quantitative_raw_frame}
Let $z>0$, let $w_j:=z^{-j/2}\mathcal T^jf$, and define
$b_j:=\|w_j\|_2^2$, $\ell_j:=\|Lw_j\|^2$, and
$D_j:=\|\nabla[(B-zI)w_j]\|_2^2$ for $j\ge1$.
For $j\ge k\ge2$,
\begin{equation}\label{eq:quantitative_raw_frame}
 \ell_j\le\alpha_k \cdot z^{-(k-1)} \cdot b_{j-k+1}
 +\frac{4k^3}{z}\sum_{h=1}^{k-1}(R/z)^{h-1} \cdot D_{j-h}.
\end{equation}
\end{lemma}
\begin{proof}
Apply Lemma~\ref{lem:combined_path_recovery} to $Q_1w_j$, with
the coordinate slot labeled $0$ and the newest $k$ derivative slots
labeled $1,\ldots,k$. Older indices form an orthogonal output direct
sum. Repeated use of Eq.~\eqref{eq:block_polynomial_testing} gives
\begin{equation*}
 \mathsf P Q_1w_j
 =\frac{Q_kw_{j-k+1}}{k! \cdot z^{(k-1)/2}}.
\end{equation*}
Its squared norm is at most
$c_k^2 \cdot z^{-(k-1)} \cdot b_{j-k+1}$.
By Eq.~\eqref{eq:block_fixed_decomposition}, a fresh adjacent swap
has norm at most $2\sqrt{D_i/z}$. After $h-1$ further steps its
norm is multiplied by at most $(R/z)^{(h-1)/2}$.
The map $Q_1$ is a contraction because the covariance is at most $I$.
The $k-1$ adjacent swaps therefore have total squared error at most
\begin{equation*}
 \frac4z\sum_{h=1}^{k-1}(R/z)^{h-1} \cdot D_{j-h}.
\end{equation*}
Substitution in Eq.~\eqref{eq:combined_path_recovery} proves the claim.
The domain justification is the same finite Bochner calculation as in
Eq.~\eqref{eq:block_fixed_decomposition}; no inverse operator is
commuted with a spatial derivative.
\end{proof}

We now apply the frame to an entire block. Averaging its starting family
produces the estimates for energy and losses over a prefix that are needed
for the comparison.

\begin{lemma}[Full-block initialization]\label{lem:full_block_initialization}
Fix positive universal constants $b_-,b_+$. There is a constant $C$,
depending only on these constants and the fixed coefficient seed,
with the following property. Let $p\ge3$ be an integer, let $\kappa\ge1$,
and let
\begin{equation*}
 z:=\|\mathcal T^m\|_{\mathrm{op}}^{2/m},\qquad
 \frac{z^{p-2}}{4\kappa}\le m\le\frac{4z^{p-2}}\kappa,
 \qquad m\ge40p,
 \qquad \alpha_p\le(1+\alpha_3)\kappa.
\end{equation*}
Assume $R-z\le2$, $z\ge C \cdot p^4$, and that a centered unit top right
singular vector $f$ gives $b_-\le b_j\le b_+$ for $0\le j\le m$.
There is a centered unit finite starting family with energy $\nu$ such
that, for $u:=z^{-1}$,
\begin{equation}\label{eq:full_block_startup}
 \nu\le\frac{u}{1-C \cdot \kappa \cdot p^4 \cdot u^p},\qquad
 P_J\le C \cdot J \cdot \kappa \cdot p^4 \cdot u^p,\ \forall 1\le J\le16p,
\end{equation}
provided $C \cdot \kappa \cdot p^4 \cdot u^p\le1/8$. The same family satisfies
\begin{equation}\label{eq:full_block_budget}
 a \cdot V_N+X_N\le\nu-u+u \cdot P_{N+1},\ \forall N\ge0.
\end{equation}
\end{lemma}
\begin{proof}
Use the notation of Lemma~\ref{lem:quantitative_raw_frame}.
Maximality of the $m$-step norm gives
\begin{equation}\label{eq:full_block_shift}
 b_0=b_m=1,\qquad b_{j+m}\le b_j,\ \forall j\ge0.
\end{equation}
Equation~\eqref{eq:block_fixed_defect} and the case $j=1$ of
Eq.~\eqref{eq:full_block_shift} therefore give the exact cancellation
\begin{equation*}
 D_{\mathrm{int}}:=\sum_{j=1}^mD_j
 \le\sum_{j=1}^m\ell_j.
\end{equation*}
For $k\le p$, the assumption $z\ge C \cdot p^4$ ensures
$(R/z)^{k-1}\le e$, so the error coefficient in
Eq.~\eqref{eq:quantitative_raw_frame} is at most $4ek^3$.
For $2\le j<p$ apply that equation with $k:=j$; for $j\ge p$ use $k:=p$.
The early coherent sum is $O(z^{-1})$ uniformly in $p$.
Indeed its $j=2$ term uses the fixed number $G_{\mathrm{seed}}(2)$, while for $j\ge3$
and $z\ge4G_{\mathrm{seed}}(p)^4$,
\begin{equation*}
 z \cdot \alpha_j \cdot z^{-(j-1)}
 \le4j^2\frac{G_{\mathrm{seed}}(p)^4}{z}
              (G_{\mathrm{seed}}(p)^2/z)^{j-3}
 \le j^2 \cdot 4^{-(j-3)}.
\end{equation*}
The sum of the last expression is finite. The late coherent sum is at
most $b_+ \cdot \alpha_p \cdot m/z^{p-1}\le4b_+ \cdot (1+\alpha_3)/z$.
The early errors cost at most $4e\sum_{k<p}k^3D_{\mathrm{int}}/z$;
the late errors cost at most $4e(p-1) \cdot p^3 \cdot D_{\mathrm{int}}/z$.
Since $\ell_1\le b_+$, absorption yields
\begin{equation}\label{eq:full_block_interior}
 D_{\mathrm{int}}\le C,\qquad
 \sum_{j=1}^m\ell_j\le C,\qquad
 \sum_{j=2}^m\ell_j\le C \cdot p^4/z.
\end{equation}

We also need a finite future tail. Let $L\le16p+2$ and write
$\delta_j:=b_j-b_{j-1}$. Equation~\eqref{eq:block_fixed_defect}
implies $\delta_{j+1}\ge\delta_j-\ell_j/z$.
Some $\delta_j$, $1\le j\le m$, is nonnegative, because their sum is
$b_m-b_0=0$. Hence $\delta_{m+1}\ge-C/z$ by
Eq.~\eqref{eq:full_block_interior}.
For $0\le h\le L+1$, the norm bound for one step gives
\begin{equation*}
 b_{m+h}\le b_+ \cdot (R/z)^h\le e \cdot b_+,
 \qquad
 \delta_{m+L+1}\le(R/z-1)b_{m+L}\le2e \cdot b_+/z.
\end{equation*}
Summing Eq.~\eqref{eq:block_fixed_defect} over this tail thus gives
$D_{\mathrm{tail}}\le\sum_{j=m+1}^{m+L}\ell_j+C$.
Apply Eq.~\eqref{eq:quantitative_raw_frame} at degree $p$ throughout
the tail. Its coherent sum is at most
$e \cdot b_+ \cdot \alpha_p \cdot L/z^{p-1}\le C/z$, since $L\le m$.
Its error coefficient is at most $4e(p-1) \cdot p^3/z$ times
$D_{\mathrm{int}}+D_{\mathrm{tail}}$. Another absorption proves
\begin{equation}\label{eq:full_block_tail}
 \sum_{j=1}^{m+L}D_j\le C,\qquad
 \sum_{j=2}^{m+L}\ell_j\le C \cdot p^4/z.
\end{equation}

We define
\begin{equation*}
 \mathcal W:=\bigoplus_{j=0}^{m-1}w_j,\qquad
 S_k:=\sum_{j=k}^{m+k-1}b_j,\qquad
 C_k:=\sum_{j=k}^{m+k-1}\ell_j,\qquad
 U:=\frac{\mathcal T\mathcal W}{\sqrt{z \cdot S_0}}.
\end{equation*}
Equation~\eqref{eq:full_block_shift} gives $S_0=S_1$.
Consequently $U$ is a centered unit genuine gradient, and Bochner gives
$e_U\le z^{-1}$. Quadratic testing for $w_1$ and
Eq.~\eqref{eq:full_block_interior} for all subsequent terms yield
\begin{equation*}
 \|\operatorname{Sym}(LU)\|_{\mathrm{HS}}^2
 \le\frac{C \cdot p^4}{zm}\le C \cdot \kappa \cdot p^4 \cdot z^{-(p-1)}.
\end{equation*}
As in Proposition~\ref{prop:block_averaged_start}, start at
$F_0:=B^{1/2}\mathcal W/\|B^{1/2}\mathcal W\|_2$ and normalize its
second successor to obtain $G_0:=F_2/\|F_2\|_2$.
Lemma~\ref{lem:skew_normalization_credit} gives
\begin{equation*}
 \nu\le\frac1{e_U^{-1}-\|\operatorname{Sym}(LU)\|_{\mathrm{HS}}^2}
 \le\frac{u}{1-C \cdot \kappa \cdot p^4 \cdot u^p}.
\end{equation*}
All denominators are positive under the stated smallness condition.

The exact loss fractions are still those of
Eq.~\eqref{eq:block_loss_fraction}. The first $J$ losses of $G_0$
correspond to $2\le k\le J+1$ in the original hierarchy.
With $L:=J+1$, Eq.~\eqref{eq:full_block_tail} gives $C_k\le C \cdot p^4/z$.
The window in $S_{k+1}$ contains at least $m-J-1\ge m/2$ interior
indices, so $S_{k+1}\ge b_- \cdot m/2$. Each of these loss fractions is at
most $C \cdot \kappa \cdot p^4 \cdot u^p$. The restarted masses are at most one, which
proves the bound on $P_J$. The same positive interior masses guarantee
every successor used here is nonzero.

Finally, Eq.~\eqref{eq:full_block_shift} implies $S_{k+1}\le S_k$
for every $k$. Writing the unnormalized hierarchy as
$F_j=t_jB^{1/2}\mathcal T^j\mathcal W$ gives
\begin{equation*}
 e_j=t_j^2 \cdot z^j \cdot S_j,\qquad
 v_{j+1}=t_j^2 \cdot z^{j+1} \cdot S_{j+1}\le z \cdot e_j.
\end{equation*}
This inequality survives the two-step shift and the common normalization.
Subtract it from the energy budget
$a \cdot V_N+X_N\le\nu-e_N$ to obtain
Eq.~\eqref{eq:full_block_budget}.
If a later family vanishes, the absorbing convention gives the same
inequality. All calculations use finite families and the form domains
already justified for the original block construction.
\end{proof}

The later comparison may inspect any permitted prefix of this family.
We therefore estimate its joint loss while keeping the prefix length
explicit.

\begin{lemma}[Joint loss with an arbitrary retained prefix]
\label{lem:order_uniform_joint_loss}
Suppose a centered unit starting family has
$\|\mathcal T^\ell\|_{\mathrm{op}}\le E_* \cdot z^{\ell/2}$ for every
$\ell\ge0$. Let $d_0\le d$ be dyadic integers. We define $J:=2d_0-1$.
On a prefix whose required normalizers are at most $B_*$, we define
$t:=\max\{1,B_* \cdot z\}$ and
\begin{align}
 \theta&:=324 \cdot C_F \cdot E_*^2 \cdot \sum_{\substack{k<d\\k\ \mathrm{dyadic}}}
                  k^6 \cdot (B_* \cdot t^3)^k \cdot c_k^2,\label{eq:order_uniform_kernel}\\
 \Delta&:=P_J+2 \cdot C_F \cdot \sum_{\substack{d_0\le k<d\\k\ \mathrm{dyadic}}}
                  k^2 \cdot B_*^k \cdot c_k^2,
 \qquad \tau:=C_F \cdot d \cdot B_*^d \cdot c_d^2.\label{eq:order_uniform_coefficients}
\end{align}
Then
\begin{equation}\label{eq:order_uniform_joint_loss}
 P_N\le\Delta+N \cdot \tau+\frac\theta\lambda
                       X_{\max\{N-2,0\}}.
\end{equation}
\end{lemma}
\begin{proof}
The proof of Lemma~\ref{lem:block_joint_local_loss} uses its initial
family only through the retained losses and the global bound for powers
of $\mathcal T$. With the present bound, a propagated fresh error has
norm at most $2 \cdot E_* \cdot t^{h/2} \cdot \sqrt{\chi/\lambda}$. The same sorting and
convolution calculation gives Eq.~\eqref{eq:order_uniform_kernel},
with $E_*$ charged once for each entire propagation.
We keep the first $2d_0-1$ losses explicitly. Beginning at that index,
degree $d_0$ is available; then use the largest available dyadic degree,
capped at $d$. Each nonterminal degree $k$ is used for at most $2k$
generations, giving the sum in Eq.~\eqref{eq:order_uniform_coefficients}.
Every defect index in a prefix is at most $N-3$, exactly as in the
original proof. This proves Eq.~\eqref{eq:order_uniform_joint_loss};
for a shorter prefix its retained losses already suffice.
\end{proof}

The initialization and joint loss bounds now give one energy comparison.
Matching the loss estimate to the family's budget gives the contradiction
that will close each curvature improvement.

\begin{lemma}[Comparison at a matched scale]\label{lem:matched_scale_comparison}
Lemma~\ref{lem:delayed_block_comparison} remains valid with $\sigma$
replaced by any $u>0$ for which the actual family satisfies
\begin{equation*}
 a \cdot V_N+X_N\le\nu-u+u \cdot P_{N+1},\ \forall N\ge0.
\end{equation*}
In particular, with $D_0:=\max\{\nu-u,0\}$, its hypotheses give
\begin{equation}\label{eq:matched_scale_comparison}
 a<12(u+D_0/p) \cdot \tau.
\end{equation}
\end{lemma}
\begin{proof}
In the proof of Lemma~\ref{lem:delayed_block_comparison}, the scale
$\sigma$ occurs only in its displayed energy budget and in the scalar
conditions $D_0\le\sigma$ and $\theta\sigma/\lambda\le1/8$.
The same budget here gives
$X_{N-2}\le D_0+u \cdot P_N$ and $a\le D_0+u(p_0+p_1)$.
The identical finite horizon
$m:=\lceil(D_0+up)/a\rceil$, $M:=m+1$, and first-exit argument therefore
apply with $u$ throughout. Their final contradiction is
$up(1-2c-cp)>0$, where $c:=17/56$ as before.
No operator inequality comparing $R$ with $u^{-1}$ is required.
\end{proof}

\subsection{Propagating losses along growing blocks}
\label{sec:quarter_order_blocks}
\label{sec:second_logstar_propagation}

We need loss estimates that remain valid as the block length grows and
inverse-gradient iterates propagate earlier errors. Joint bounds for raw
moments and averaging over starting positions provide starting families
with controlled losses. We then allow polynomial growth in the propagation
bound, as required when successive blocks are combined.

\begin{lemma}[A joint frame for raw moments]
\label{lem:quarter_order_raw_joint}
Suppose $z>0$ and
$\|\mathcal T^h\|_{\mathrm{op}}\le E_* \cdot z^{h/2}$ for every integer $h\ge0$.
Let $f$ be centered, and define
\begin{equation*}
 w_j:=z^{-j/2}\mathcal T^jf,\qquad
 b_j:=\|w_j\|_2^2,\qquad \ell_j:=\|Lw_j\|^2,\qquad
 D_j:=\|\nabla[(B-zI)w_j]\|_2^2,\ \forall j\ge1.
\end{equation*}
For every dyadic $d\ge2$ and $j\ge2d-1$,
\begin{align}
 \ell_j&\le C_F \cdot d \cdot c_d^2 \cdot z^{-(d-1)} \cdot b_{j-d+1}\notag\\
 &\quad+108 \cdot C_F \cdot E_*^2
 \sum_{\substack{k<d\\k\ \mathrm{dyadic}}}
 k^5 \cdot c_k^2 \cdot z^{-k} \cdot \sum_{h=1}^{3k-1}D_{j-k+1-h}.
 \label{eq:quarter_order_raw_joint}
\end{align}
In a sum of such inequalities, even with varying terminal degrees at most
$d$, the total coefficient of any one defect is at most
\begin{equation}\label{eq:quarter_order_raw_kernel}
 \mathcal K(z,d):=324 \cdot C_F \cdot E_*^2
 \sum_{\substack{k<d\\k\ \mathrm{dyadic}}}k^6 \cdot c_k^2 \cdot z^{-k}.
\end{equation}
If $z\ge64 \cdot G_*(d)^4$, then $\mathcal K(z,d)\le C/z$, where $C$ depends
only on $E_*$ and the universal constants.
\end{lemma}
\begin{proof}
We define $Y_j:=Q_1w_j$. Repeated polynomial testing, as in
Eq.~\eqref{eq:block_polynomial_testing}, gives
\begin{equation}\label{eq:quarter_order_raw_projection}
 \mathsf P_kY_j=
 \frac{Q_kw_{j-k+1}}{k! \cdot z^{(k-1)/2}}.
\end{equation}
Thus the coordinate slot and the $k-1$ newest derivative slots form the
first symmetric block. The remaining family $w_{j-k+1}$ has
$j-k+1$ derivative slots.

Equation~\eqref{eq:block_fixed_decomposition} bounds a fresh adjacent swap
in $w_i$ by $2\sqrt{D_{i-1}/z}$. The main gradient term has zero fresh
swap because $w_{i-1}$ is itself a genuine gradient in its newest slot.
A swap in positions $h,h+1$, counted from the newest slot of
$w_{j-k+1}$, was created using $D_{j-k+1-h}$ and subsequently propagated
through $h-1$ normalized copies of $\mathcal T$. The assumed global
bound charges a factor at most $E_*$ for this entire propagation.
These maps commute with permutations of every previously introduced
slot. Since $\|Q_k\|\le k! \cdot c_k$, Eq.~\eqref{dt_eq_sort} on the next
$3k$ derivative slots implies
\begin{equation}\label{eq:quarter_order_raw_swap}
 \|(I-\mathsf S_k)\mathsf P_kY_j\|
 \le6k \cdot E_* \cdot c_k \cdot z^{-k/2}
       \sum_{h=1}^{3k-1}\sqrt{D_{j-k+1-h}}.
\end{equation}
For each degree in this display, $k\le d/2$ and $j\ge2d-1\ge4k-1$.
The smallest defect index is therefore $j-4k+2\ge1$, and all the
required derivative slots are present.

The coherent projection in Eq.~\eqref{eq:quarter_order_raw_projection} has
squared norm at most $c_d^2 \cdot z^{-(d-1)} \cdot b_{j-d+1}$. Square
Eq.~\eqref{eq:quarter_order_raw_swap}, using at most $3k$ summands, and
apply Theorem~\ref{thm:joint_uniform_frame}. This proves
Eq.~\eqref{eq:quarter_order_raw_joint}. At a fixed degree $k$, any one
$D_i$ occurs at most $3k-1$ times when $j$ is summed. Omitting some
indices, or some degrees at those indices, only decreases that count.
This proves Eq.~\eqref{eq:quarter_order_raw_kernel}.
Finally, with $H:=G_*(d)\ge1$, the seed and $H^2/z\le1/64$ give
\begin{equation*}
 \mathcal K(z,d)
 \le\frac{324 \cdot C_F \cdot E_*^2}{z}\sum_{k\ge1}k^6 \cdot (H^2/z)^{k-1}
 \le\frac Cz.
\end{equation*}
The same finite Bochner domain argument used for
Lemma~\ref{lem:quantitative_raw_frame} justifies every displayed pairing.
\end{proof}

The raw frame identifies the loss at each starting position. Weighting
these positions before initialization gives a family whose prefix losses
are controlled by the corresponding weighted sums.

\begin{lemma}[Initialization with weighted raw losses]
\label{lem:quarter_order_initialization}
Fix positive universal constants $b_-,b_+,E_*$. There is a constant $C$
depending only on these constants and the coefficient seed with the
following property. Suppose $p\ge3$ is odd, $\kappa\ge1$, $m$ is a
positive integer, and
\begin{gather*}
 z:=\|\mathcal T^m\|_{\mathrm{op}}^{2/m},\qquad
 \frac{z^{p-2}}{4\kappa}\le m\le\frac{4z^{p-2}}\kappa,
 \qquad m\ge100p,\qquad
 \mathsf A_p\le(1+\mathsf A_3)\kappa,\\
 0\le R-z\le2,\qquad z\ge C,\qquad
 z\ge64 \cdot G_*(d_p)^4,\qquad 16p z^{p-d_p}\le1,\\
 \|\mathcal T^h\|_{\mathrm{op}}\le E_* \cdot z^{h/2},\ \forall h\ge0.
\end{gather*}
Suppose also that a centered unit top right singular vector of
$\mathcal T^m$ gives $b_-\le b_j\le b_+$ for $0\le j\le m$.
There is one centered unit finite starting family whose energy $\nu$
and actual losses satisfy, with $u:=z^{-1}$,
\begin{equation}\label{eq:quarter_order_initialization}
 \nu\le\frac{u}{1-C \cdot \kappa \cdot u^p},\qquad
 P_J\le C \cdot \kappa \cdot u^p,\ \forall 1\le J\le16p,
\end{equation}
provided $C \cdot \kappa \cdot u^p\le1/8$. The same family satisfies
\begin{equation}\label{eq:quarter_order_matched_budget}
 a \cdot V_N+X_N\le\nu-u+u \cdot P_{N+1},\ \forall N\ge0.
\end{equation}
\end{lemma}
\begin{proof}
Use the notation of Lemma~\ref{lem:quarter_order_raw_joint}. Compactness
supplies the stated singular vector, as in
Proposition~\ref{prop:block_averaged_start}. Maximality gives
\begin{equation}\label{eq:quarter_order_orbit_shift}
 b_0=b_m=1,\qquad b_{j+m}\le b_j,\ \forall j\ge0.
\end{equation}
By Eq.~\eqref{eq:block_fixed_defect},
\begin{equation}\label{eq:quarter_order_interior_telescope}
 D_{\mathrm{int}}:=\sum_{j=1}^mD_j
 \le\sum_{j=1}^m\ell_j-z(b_1-b_{m+1})
 \le\sum_{j=1}^m\ell_j.
\end{equation}
Indeed the boundary term is $z(b_{m+1}-b_m+1-b_1)\le0$.

For $j\ge3$, use Eq.~\eqref{eq:quarter_order_raw_joint} at the largest
dyadic degree at most $\min\{d_p,(j+1)/2\}$. This degree is available
at the actual index $j$. Every nonterminal degree $k$ is used for at
most $2k$ generations. Their combined coherent contribution is at most
\begin{equation}\label{eq:quarter_order_early_coherent}
 2 \cdot C_F \cdot b_+\sum_{\substack{2\le k<d_p\\k\ \mathrm{dyadic}}}
       k^2 \cdot c_k^2 \cdot z^{-(k-1)}\le C/z.
\end{equation}
Here the degree-two term uses the universal bound on $c_2$. For $k\ge4$,
we define $H:=G_*(d_p)$ and note that
\begin{equation*}
 z \cdot k^2 \cdot c_k^2 \cdot z^{-(k-1)}
 \le k^2\frac{H^6}{z^2}(H^2/z)^{k-4}.
\end{equation*}
Its sum is uniformly bounded by the hypothesis $z\ge64H^4$.
The late coherent sum, at degree $d_p\ge p+3$, is at most
\begin{equation*}
 \mathsf A_p \cdot b_+ \cdot m \cdot z^{-(d_p-1)}\le C/z.
\end{equation*}
The total defect contribution of every degree and every generation is
at most $C \cdot D_{\mathrm{int}}/z$ by
Eq.~\eqref{eq:quarter_order_raw_kernel}.

The remaining first two indices satisfy $\ell_1\le b_+$ and
\begin{equation}\label{eq:quarter_order_second_moment}
 \ell_2\le(2 \cdot K_2 \cdot b_1+8 \cdot D_1)/z,
\end{equation}
as in Eq.~\eqref{eq:block_second_moment}. Thus
Eq.~\eqref{eq:quarter_order_interior_telescope} and absorption at a fixed
large $z$ give
\begin{equation}\label{eq:quarter_order_interior}
 D_{\mathrm{int}}\le C,\qquad
 \sum_{j=1}^m\ell_j\le C,\qquad
 \sum_{j=2}^m\ell_j\le C/z.
\end{equation}
In particular there is no factor depending on $p$ in the last estimate.

We need the same estimate on a finite future tail. Let $1\le L\le16p+2$
be an integer and write $\delta_j:=b_j-b_{j-1}$. The defect identity gives
$\delta_{j+1}\ge\delta_j-\ell_j/z$. Since
$\sum_{j=1}^m\delta_j=0$, some interior $\delta_j$ is nonnegative,
and Eq.~\eqref{eq:quarter_order_interior} gives
$\delta_{m+1}\ge-C/z$. The global propagation bound and $\|f\|_2=1$
give $b_j\le E_*^2$ for every $j\ge0$. Since
$b_{j+1}\le(R/z)b_j$ and $R-z\le2$, we also have
\begin{equation*}
 \delta_{m+L+1}\le(R/z-1)b_{m+L}\le2 \cdot E_*^2/z.
\end{equation*}
The tail telescope therefore gives
\begin{equation*}
 D_{\mathrm{tail}}:=\sum_{j=m+1}^{m+L}D_j
 \le\sum_{j=m+1}^{m+L}\ell_j+C.
\end{equation*}
Degree $d_p$ is available throughout this tail because
$2d_p-1\le4p+3<m$.
Its coherent contribution is at most
$E_*^2 \cdot \mathsf A_p \cdot L/z^{d_p-1}\le C/z$, since $L\le m$.
The defect contribution is bounded by
$C(D_{\mathrm{int}}+D_{\mathrm{tail}})/z$; the exact indices in
Lemma~\ref{lem:quarter_order_raw_joint} include only this interior and tail.
A second absorption proves
\begin{equation}\label{eq:quarter_order_tail}
 \sum_{j=1}^{m+L}D_j\le C,\qquad
 \sum_{j=2}^{m+L}\ell_j\le C/z.
\end{equation}

We define the averaged family
\begin{equation*}
 \mathcal W:=\bigoplus_{j=0}^{m-1}w_j,\qquad
 S_k:=\sum_{j=k}^{m+k-1}b_j,\qquad
 C_k:=\sum_{j=k}^{m+k-1}\ell_j,\qquad
 U:=\frac{\mathcal T\mathcal W}{\sqrt{z \cdot S_0}}.
\end{equation*}
Equation~\eqref{eq:quarter_order_orbit_shift} gives $S_0=S_1$.
Thus $U$ is a centered unit genuine gradient and $e_U\le z^{-1}$.
Quadratic testing at its first summand gives
$\|\operatorname{Sym}(Lw_1)\|^2\le K_2/(4z)$.
Use Eq.~\eqref{eq:quarter_order_interior} at all subsequent summands and
$S_0\ge b_- \cdot m$ to obtain
\begin{equation*}
 \|\operatorname{Sym}(LU)\|_{\mathrm{HS}}^2
 \le C/(zm)\le C \cdot \kappa \cdot z^{-(p-1)}.
\end{equation*}
Start at $F_0:=B^{1/2}\mathcal W/\|B^{1/2}\mathcal W\|_2$ and take
the compensated unit restart $G_0:=F_2/\|F_2\|_2$.\linebreak[4]
Lemma~\ref{lem:skew_normalization_credit} gives
\begin{equation*}
 \nu\le\frac1{e_U^{-1}-\|\operatorname{Sym}(LU)\|_{\mathrm{HS}}^2}
 \le\frac{u}{1-C \cdot \kappa \cdot u^p}.
\end{equation*}
The denominators are positive under the stated smallness condition.

Fix $1\le J\le16p$ and define the window multiplicities
\begin{equation*}
 n_J(i):=\#\{k\in\{2,\ldots,J+1\}:k\le i\le m+k-1\},\qquad
 \mathcal E_J:=\sum_{i\ge1} n_J(i)\ell_i,\qquad
 \mathcal D_J:=\sum_{i\ge1} n_J(i)D_i.
\end{equation*}
Here $n_J$ is supported on $2,\ldots,m+J$, is $1$-Lipschitz on the
integer line, and satisfies $n_J(i)\le i-1$ for $i\ge1$ and
$\sum_i n_J(i)=Jm$. Summing Eq.~\eqref{eq:block_fixed_defect} over
the same $J$ windows gives
\begin{align}
 \mathcal D_J
 &\le\mathcal E_J+
 z\sum_{k=2}^{J+1}(\delta_{m+k}-\delta_k)\notag\\
 &=\mathcal E_J+
 z(b_{m+J+1}-b_{m+1}-b_{J+1}+b_1)\notag\\
 &\le\mathcal E_J+z(b_1-b_{m+1})
 \le\mathcal E_J+C.\label{eq:quarter_order_weighted_defect}
\end{align}
The penultimate inequality uses Eq.~\eqref{eq:quarter_order_orbit_shift};
the last follows from Eqs.~\eqref{eq:quarter_order_interior_telescope}
and~\eqref{eq:quarter_order_interior}.

To apply the raw moment estimate with these weights, set
$D_i:=0$ for $i\le0$. For $j\ge3$, let $h_j$ be its coherent term at
the largest available degree, capped at $d_p$, and set
$h_2:=2 \cdot K_2 \cdot b_1/z$. Eqs.~\eqref{eq:quarter_order_raw_joint}
and~\eqref{eq:quarter_order_second_moment} give nonnegative coefficients
$\alpha_r$ such that
\begin{equation}\label{eq:quarter_order_raw_convolution}
 \ell_j\le h_j+\sum_{r\ge1}\alpha_r \cdot D_{j-r}
 ,\ \forall 2\le j\le m+J.
\end{equation}
Explicitly, add $8/z$ at delay one, and, for each dyadic $k<d_p$,
add $108 \cdot C_F \cdot E_*^2 k^5 \cdot c_k^2 \cdot z^{-k}$ at each delay
$r=k,\ldots,4k-2$. Including unavailable degrees only increases the
right side. Therefore
\begin{align}
 \eta:=\sum_r\alpha_r
 &\le\frac8z+324 \cdot C_F \cdot E_*^2
    \sum_{\substack{k<d_p\\k\ \mathrm{dyadic}}}k^6 \cdot c_k^2 \cdot z^{-k}
 \le\frac Cz,\notag\\
 \eta_1:=\sum_r r\alpha_r
 &\le\frac8z+1296 \cdot C_F \cdot E_*^2
    \sum_{\substack{k<d_p\\k\ \mathrm{dyadic}}}k^7 \cdot c_k^2 \cdot z^{-k}
 \le\frac Cz.\label{eq:quarter_order_delay_moments}
\end{align}
Both bounds follow from $c_k^2\le H^{2k-2}$ and $H^2/z\le1/64$.
Since $n_J(i+r)\le n_J(i)+r$, multiplying
Eq.~\eqref{eq:quarter_order_raw_convolution} by $n_J(j)$ and summing
gives, by Eq.~\eqref{eq:quarter_order_tail},
\begin{equation}\label{eq:quarter_order_weighted_convolution}
 \mathcal E_J\le\sum_j n_J(j)h_j+
 \eta \cdot \mathcal D_J+\eta_1 \cdot \sum_{i=1}^{m+J}D_i
 \le\sum_j n_J(j)h_j+\frac Cz(\mathcal D_J+1).
\end{equation}
Every defect originally used here has index at least one and at most
$m+J-1$; thus the finite tail already proved is sufficient.

The $j=2$ coherent contribution is at most $C/z$, since $n_J(2)=1$.
A nonterminal degree $k<d_p$ is used only at
$j=2k-1,\ldots,4k-2$. There are $2k$ such indices, with
$n_J(j)<4k$. Using $b_i\le E_*^2$, their total contribution is at most
\begin{equation*}
 8 \cdot C_F \cdot E_*^2 \cdot \sum_{\substack{2\le k<d_p\\k\ \mathrm{dyadic}}}
 k^3 \cdot c_k^2 \cdot z^{-(k-1)}\le C/z.
\end{equation*}
For $k=2$ use $c_2^2\le2$; for $k\ge4$ use
$z \cdot k^3 \cdot c_k^2 \cdot z^{-(k-1)}\le k^3 \cdot H^6 \cdot z^{-2} \cdot (H^2/z)^{k-4}$.
The terminal coherent contribution is at most
\begin{equation*}
 E_*^2 \cdot \mathsf A_p \cdot Jm \cdot z^{-(d_p-1)}
 \le C \cdot J \cdot z^{p-d_p-1}\le C/z,
\end{equation*}
by the hypotheses on $m$, $\mathsf A_p$, and $16p z^{p-d_p}$.
Eqs.~\eqref{eq:quarter_order_weighted_defect}
and~\eqref{eq:quarter_order_weighted_convolution} now give, after
absorption at a fixed large $z$,
\begin{equation}\label{eq:quarter_order_weighted_loss}
 \mathcal E_J=\sum_{k=2}^{J+1}C_k\le C/z.
\end{equation}

The first $J$ losses of $G_0$ are the original losses at indices
$2,\ldots,J+1$, up to their common unit normalization. Their fractions
are $C_k/(z \cdot S_{k+1}+C_k)$, by Eq.~\eqref{eq:block_loss_fraction}.
The denominator window contains at least $m-J-1\ge m/2$ interior
indices, so $S_{k+1}\ge b_- \cdot m/2$. The unit restarted masses are at most
one, and Eq.~\eqref{eq:quarter_order_weighted_loss} yields
\begin{equation*}
 P_J\le\frac{2\mathcal E_J}{b_- \cdot zm}
 \le\frac C{z^2 \cdot m}\le C \cdot \kappa \cdot u^p.
\end{equation*}
This proves Eq.~\eqref{eq:quarter_order_initialization}. The same
positive interior entries justify every normalization through these
indices.

Finally, Eq.~\eqref{eq:quarter_order_orbit_shift} implies $S_{k+1}\le S_k$
for all $k$. Writing the unnormalized hierarchy as
$F_j=t_jB^{1/2}\mathcal T^j\mathcal W$ gives
\begin{equation*}
 e_j=t_j^2 \cdot z^j \cdot S_j,\qquad
 v_{j+1}=t_j^2 \cdot z^{j+1} \cdot S_{j+1}\le z \cdot e_j.
\end{equation*}
The inequality survives the two-step shift and common normalization.
Subtract it from $a \cdot V_N+X_N\le\nu-e_N$ to obtain
Eq.~\eqref{eq:quarter_order_matched_budget}. The absorbing zero convention
handles any later vanishing family, as in
Lemma~\ref{lem:full_block_initialization}. All families and all tails
used in this proof are finite.
\end{proof}

The same construction must work when the operator powers have a
polynomial growth factor. Each error propagates through a number of
steps proportional to its degree. We can therefore absorb this growth
into the degree weights in the error kernel, giving the loss bound
needed to extend the blocks.

\begin{lemma}[Joint loss under polynomial propagation]
\label{lem:second_logstar_polynomial_loss}
Suppose that, for a fixed $E\ge1$ and some $0\le\alpha\le1$,
\begin{equation}\label{eq:second_logstar_polynomial_propagation}
 \|\mathcal T^h\|_{\mathrm{op}}
 \le E \cdot (h+1)^\alpha \cdot z^{h/2},\ \forall h\ge0.
\end{equation}
Let a centered unit starting family, dyadic integers $d_0\le d$, and
$J:=2d_0-1$ be given. On a prefix whose required normalizers are at most
$B_*$, we define $t:=\max\{1,B_* \cdot z\}$ and
\begin{align}
 \theta&:=5184 \cdot C_F \cdot E^2 \cdot \sum_{\substack{k<d\\k\ \mathrm{dyadic}}}
          k^8 \cdot (B_* \cdot t^3)^k \cdot c_k^2,\notag\\
 \Delta&:=P_J+2 \cdot C_F \cdot \sum_{\substack{d_0\le k<d\\k\ \mathrm{dyadic}}}
          k^2 \cdot B_*^k \cdot c_k^2,
 \qquad \tau:=C_F \cdot d \cdot B_*^d \cdot c_d^2.
 \label{eq:second_logstar_loss_coefficients}
\end{align}
Then every justified prefix satisfies
\begin{equation}\label{eq:second_logstar_polynomial_loss}
 P_N\le\Delta+N \cdot \tau+\frac\theta\lambda X_{\max\{N-2,0\}}.
\end{equation}
The constants are uniform in $\alpha$ and the terminal degree.
\end{lemma}
\begin{proof}
A fresh adjacent swap has norm at most $2\sqrt{\chi/\lambda}$.
The next $h$ maps are scalar multiples of $\mathcal T$, with product
$\prod\sqrt{\beta_j}\mathcal T^h$. Their norm is at most
$E \cdot (h+1)^\alpha \cdot t^{h/2}$. These maps commute with permutations of
previously introduced output indices, and use the scalars of the
original complete family. Thus the sorting argument in
Lemma~\ref{lem:order_uniform_joint_loss} gives
\begin{equation*}
 \|(I-\mathsf S_k)\mathsf P_kT_j\|
 \le6kE \cdot B_*^{k/2} \cdot c_k \cdot \lambda^{-1/2}
 \sum_{h=1}^{3k-1}h^\alpha \cdot t^{(h-1)/2} \cdot \sqrt{\chi_{j-k-h}}.
\end{equation*}
The inner convolution has $\ell^1$ norm at most
$3k \cdot (3k+1)^\alpha \cdot t^{3k/2}$. Squaring, summing in $j$, and applying
the joint-frame weight $C_F \cdot k^2$ gives an error coefficient at most
\begin{equation*}
 324 \cdot C_F \cdot E^2 \cdot \sum_{\substack{k<d\\k\ \mathrm{dyadic}}}
 k^6 \cdot (3k+1)^{2\alpha} \cdot (B_* \cdot t^3)^k \cdot c_k^2\le\theta.
\end{equation*}
Here $(3k+1)^{2\alpha}\le16k^2$ for $k\ge1$.
The coherent projections, the retained first $J$ losses, and the at
most $2k$ generations at each nonterminal degree are unchanged.
Every defect in $P_N$ has index $j-k-h\le N-3$, giving the displayed
delay. Shorter prefixes are covered by the retained actual losses.
This proves Eq.~\eqref{eq:second_logstar_polynomial_loss}.
\end{proof}

\subsection{Extending blocks under a uniform coefficient cap}
\label{sec:repeated_height_scalars}

To repeat the block construction with successively better coefficients,
its constants must remain uniform in the outer round. We keep each new
coefficient bound below the original seed and prove uniform height and
moment estimates for these capped bounds. They allow us to extend a
finite block while keeping its radius close to the coefficient scale.

For $x\ge1$, we define
\begin{equation}\label{eq:second_logstar_seed_notation}
 \mathfrak t(x):=1+\log^*(x+2),\qquad
 \mathfrak t_j:=\mathfrak t^{\circ j},\ \forall j\ge0,\qquad
 t_Q:=\mathfrak t(Q).
\end{equation}
Here $\mathfrak t_0(x):=x$, and $\log^*x$ is the least number of successive natural logarithms needed to bring $x$ to at most one. The following bounds are uniform in the number of compositions.

\begin{lemma}[Uniform height and coefficient-seed estimates]
\label{lem:second_logstar_seed_estimates}
For $x,y\ge1$ and $j\ge1$,
\begin{equation}\label{eq:second_logstar_height_product}
 \mathfrak t_j(xy)\le\mathfrak t_j(x)+\mathfrak t_j(y).
\end{equation}
For each fixed $C,p\ge1$, there is an integer $b_{C,p}$ such that
\begin{equation}\label{eq:repeated_height_power_argument}
 \mathfrak t_j(C \cdot x^p)\le\mathfrak t_j(x)+b_{C,p},\qquad
 \mathfrak t_j(x^x)\le\mathfrak t_j(x)+2.
\end{equation}
For arbitrary $A\ge1$, we define
\begin{equation}\label{eq:repeated_height_capped_seed}
 \mathcal G(x):=\min\{G_*(x),A \cdot \mathfrak t_j(x)^{1/6}\},\qquad
 \mathfrak g(x):=\mathcal G(x)^2.
\end{equation}
This seed is nondecreasing, satisfies $1\le\mathcal G\le G_*$, and
obeys the following estimates for every integer $Q\ge3$, with constants
independent of $j,A$:
\begin{align}
 \mathcal G(C \cdot x^p)&\le C_{C,p} \cdot \mathcal G(x),\qquad
 \mathcal G(k)^6\le C_0 \cdot k,\label{eq:repeated_height_seed_scaling}\\
 \sup_{s\ge1}\frac{\mathfrak g(C \cdot (Qs)^2)}{s^{1/Q}}
 &\le C'_C \cdot \mathfrak g(Q),\label{eq:second_logstar_seed_cutoff}\\
 D_Q:=2 \cdot C_F \cdot 2^Q \cdot \sum_{k\ge Q}k^2 \cdot \mathfrak g(k)^{Q-1} \cdot 4^{-(k-Q)}
 &\le[C_1 \cdot \mathfrak g(Q)]^Q.\label{eq:second_logstar_coherent_root}
\end{align}
The same estimates hold for $\mathcal G=G_*$. If $m\le R^Q$ and
$R\ge1$, then also
\begin{equation}\label{eq:repeated_height_static_seed}
 \mathcal G(m)^2\le C_2 \cdot \mathcal G(8Q)^2+C_2 \cdot G_*(R)^2.
\end{equation}
\end{lemma}
\begin{proof}
The thresholds of the ordinary log-star function are
$a_0:=1$, $a_{i+1}:=e^{a_i}$. Consecutive thresholds are separated by
more than one. Consequently, for every integer $b\ge0$,
\begin{equation*}
 \mathfrak t(x+b)\le\mathfrak t(x)+b.
\end{equation*}
The same inequality holds after any finite number of compositions.
Also $\mathfrak t(x)\ge3$ for $x\ge1$. Since $\mathfrak t(3)=3$,
the integer-shift inequality implies $\mathfrak t(n)\le n$ for all
integers $n\ge3$. Therefore
\begin{equation*}
 3\le\mathfrak t_j(x)\le\mathfrak t(x)\le x+3\le4x.
\end{equation*}
Here $\log^*v\le v$ for $v\ge1$, by induction on its logarithm height.

Both $xy+2$ and $x+y+2$ are at most $(x+2)(y+2)$, which is at most
$\exp(\max\{x+2,y+2\})$. Thus
$\mathfrak t(xy)\le\mathfrak t(x)+\mathfrak t(y)$ and
$\mathfrak t(x+y)\le\mathfrak t(x)+\mathfrak t(y)$.
Induction proves Eq.~\eqref{eq:second_logstar_height_product}.
For fixed $C,p$, let $K:=p+\log(C+2)$. The inequality
$\log(C \cdot x^p+2)\le K(x+2)$ gives
$\mathfrak t(C \cdot x^p)\le\mathfrak t(x)+2+\log^*K$.
The integer-shift bound propagates this additive constant through every
later composition. Similarly,
$x^x+2\le(x+2)^{x+2}\le\exp(\exp(x+2))$ gives the second assertion
of Eq.~\eqref{eq:repeated_height_power_argument}.

Both branches in Eq.~\eqref{eq:repeated_height_capped_seed} satisfy
the fixed-power argument bound, and taking their minimum preserves it.
The cap and the original seed give $\mathcal G(k)^6\le C_0 \cdot k$.
For $k\ge Q$, writing $v:=k/Q$ gives
$\mathfrak t_j(k)\le\mathfrak t_j(Q)+4v\le5v\mathfrak t_j(Q)$.
Applying this to both branches yields
\begin{equation*}
 \mathfrak g(k)\le C \cdot \mathfrak g(Q) \cdot (k/Q)^{1/3}
 \le C \cdot \mathfrak g(Q) \cdot (k/Q).
\end{equation*}
The integer moment bound
\begin{equation*}
 \sum_{k\ge0}k^\beta \cdot 4^{-k}\le2(\beta/(e \cdot \log2))^\beta
\end{equation*}
with $\beta:=Q+1$ proves Eq.~\eqref{eq:second_logstar_coherent_root};
the remaining factor $Q^{2/Q}$ is uniformly bounded.

For the cutoff, we define $y:=\log(s)/Q$ and $X:=C \cdot Q^2 \cdot e^{2Qy}$. Then
$\log(X+2)\le C_3 \cdot Q \cdot (1+y)$, so
$\mathfrak t(X)\le1+\mathfrak t(C_3 \cdot Q \cdot (1+y))$.
The integer-shift, product, and fixed-power bounds imply
\begin{equation*}
 \mathfrak t_j(X)\le1+\mathfrak t_j(C_3 \cdot Q \cdot (1+y))
 \le1+\mathfrak t_j(Q)+b_{C_3,1}+4(1+y)
 \le C_4 \cdot \mathfrak t_j(Q) \cdot (1+y).
\end{equation*}
After division by $s^{1/Q}=e^y$, the factor
$(1+y)^{1/3} \cdot e^{-y}$ is bounded. Both branches obey this pointwise
estimate, proving Eq.~\eqref{eq:second_logstar_seed_cutoff} for their
minimum as well.

Finally, we define $M:=\max\{R,Q\}$. Then $m\le M^M$, so the bound for raising the argument to its own power gives $\mathcal G(m)\le C \cdot \mathcal G(M)$. If $R\le Q$, use
$\mathcal G(M)\le\mathcal G(8Q)$; otherwise use
$\mathcal G(M)\le G_*(R)$. This proves
Eq.~\eqref{eq:repeated_height_static_seed}.
\end{proof}

The scalar estimates control the growth of the cap. The next moment
bounds show that its exponentially damped tails have universal constants,
which will control the extension kernels.

\begin{lemma}[A uniform moment bound for capped seeds]
\label{lem:nearunit_seed_moments}
Fix nonnegative integers $M,r$ and constants $c,C_*>0$. We define
$\rho_k:=\min\{1,C_* \cdot k^{1/3}/F^2\}$. For $0<\epsilon\le1/16$,
\begin{equation}\label{eq:nearunit_seed_moments}
 \sum_{k\ge r+1}k^M \cdot \rho_k^{k-r} \cdot e^{-c \cdot \epsilon \cdot k}\le C
 \qquad\text{if }\epsilon \cdot F^2\ge C_0.
\end{equation}
Here $C_0,C$ depend only on the fixed parameters. The left side tends
uniformly to zero as the lower bound $C_0$ tends to infinity.
\end{lemma}
\begin{proof}
For $k\le F^3$, increasing $C_0$ gives
$\rho_k\le C_*/F\le F^{-1/2}$. The corresponding sum is bounded by
$\sum_{k\ge r+1}k^M \cdot F^{-(k-r)/2}$, which tends to zero as
$F\ge\sqrt{16 \cdot C_0}$ tends to infinity. For $k>F^3$, discard
$\rho_k\le1$ and take out $e^{-c \cdot \epsilon \cdot F^3/2}$. The remaining
exponential moment is at most $C \cdot \epsilon^{-M-1}$. Since
$\epsilon^{-1}\le F^2/C_0$ and $\epsilon \cdot F^3\ge C_0 \cdot F$, this part
is at most
\begin{equation*}
 C \cdot C_0^{-M-1} \cdot F^{2M+2} \cdot e^{-c \cdot C_0 \cdot F/2},
\end{equation*}
which tends uniformly to zero as $C_0$ increases. This proves the
assertion. We will use only finitely many fixed choices of $M,r,c$.
\end{proof}

We can now extend a block using its existing family and loss budget.
The moment bounds for capped seeds control the added terms, while the normalization
conditions keep the extension in the permitted regime.

\begin{lemma}[Block extension near the seed radius]
\label{lem:nearunit_block_extension}
Fix a universal $E\ge1$. Let $0<\epsilon\le1/16$, let $q\ge3$ be
odd, and suppose $\epsilon \cdot F^2\ge C_0$, where $C_0$ is a sufficiently
large universal constant. Let $\mathcal G$ be a nondecreasing coefficient
majorant satisfying $1\le\mathcal G\le G_*$ and
$c_k\le\mathcal G(k)^{k-1}$. We define
\begin{equation*}
 d:=\min\{2^j:2^j\ge8q\},\qquad
 \kappa_q:=F^{2q-2},\qquad \kappa_{q-2}:=F^{2q-6}.
\end{equation*}
Assume $\mathcal G(d)\le F$, $0\le R-z\le2$, and
\begin{equation*}
 \|\mathcal T^h\|_{\mathrm{op}}
 \le E \cdot (h+1)^\alpha \cdot z^{h/2},\ \forall h\ge0,\qquad0\le\alpha\le1.
\end{equation*}
Suppose one centered unit starting family satisfies, with $u:=z^{-1}$,
\begin{equation*}
 \nu\le\frac{u}{1-a_q \cdot u^q},\qquad
 P_J\le s_q \cdot u^q,\ \forall J\le32q,\qquad
 a \cdot V_N+X_N\le\nu-u+u \cdot P_{N+1},\ \forall N\ge0.
\end{equation*}
There is a universal $C_\theta$, depending only on $E$ and the
previously fixed constants, for which the following holds. We define
\begin{equation*}
 \mathsf K_q:=4(s_q+C_F \cdot \kappa_{q-2}+4 \cdot C_\theta \cdot a_q+1),\qquad
 C_q:=a_q+4 \cdot \mathsf K_q.
\end{equation*}
If
\begin{equation}\label{eq:nearunit_extension_thresholds}
 z\ge\max\{(1+\epsilon) \cdot F^2,8,8 \cdot C_\theta\},\qquad
 a_q \cdot u^q\le\epsilon/1024,\qquad
 \mathsf K_q \cdot u^q\le\epsilon/1024,
\end{equation}
then
\begin{equation}\label{eq:nearunit_block_extension}
 \|\mathcal T^m\|_{\mathrm{op}}^{2/m}
 \ge z-C_q \cdot z^{1-q}
 ,\ \forall 1\le m\le\lfloor z^q/\kappa_q\rfloor.
\end{equation}
In particular $C_q\le C \cdot (a_q+s_q+\kappa_{q-2}+1)$.
\end{lemma}
\begin{proof}
Use Lemma~\ref{lem:second_logstar_polynomial_loss} with $d_0:=d$.
The available range covers $2d-1<32q$. We define $p_*:=\mathsf K_q \cdot u^q$.
Before a first exit of this loss, the actual normalizers obey
\begin{equation}\label{eq:nearunit_extension_normalizers}
 B_*\le\frac\nu{1-p_*}
 \le\frac{u}{(1-a_q \cdot u^q)(1-p_*)}
 \le(1+\epsilon/128)u.
\end{equation}
Thus $t_*:=\max\{1,B_* \cdot z\}\le1+\epsilon/128$ and
$B_* \cdot t_*^3\le(1+\epsilon/16)u$. Also
$C_P\le R+1\le z+3$ gives $\lambda\ge u/2$.
We define $\rho_k:=\min\{1,C_* \cdot k^{1/3}/F^2\}$ for all $k\ge1$.
The original cap and the seed bound at $d$ give
\begin{equation*}
 \mathcal G(k)^2/F^2\le\rho_k,\qquad
 B_* \cdot t_*^3 \cdot \mathcal G(k)^2\le\rho_k \cdot e^{-\epsilon/2}
 ,\ \forall k\le d.
\end{equation*}
Lemma~\ref{lem:nearunit_seed_moments}, with the fixed term $k=1$
treated separately, therefore bounds the error kernel by
\begin{equation*}
 \theta\le C \cdot u\sum_{k\ge1}k^8 \cdot \rho_k^{k-1}
                    e^{-\epsilon(k-1)/2}\le C_\theta \cdot u.
\end{equation*}
The constant $C_\theta$ is independent of the order, the seed amplitude,
and $\epsilon$.

We use only the horizon required for the next block. Write
$D:=1+\epsilon$, $\beta:=1+\epsilon/128$, and
$\rho:=\mathcal G(d)^2/F^2$. Since $d\ge8q$,
$\beta^d \cdot D^{2q-d}\le e^{-\epsilon d/2}$. The terminal coefficient
thus satisfies
\begin{align*}
 \frac{\tau}{\kappa_q \cdot \kappa_{q-2} \cdot u^{2q}}
 &\le C_F \cdot d \cdot F^6 \cdot \rho^{d-1} \cdot \beta^d \cdot D^{2q-d}\\
 &\le C_F \cdot d \cdot \mathcal G(d)^6 \cdot \rho^{d-4} \cdot e^{-\epsilon d/2}
 \le C \cdot d^2 \cdot \rho^{d-4} \cdot e^{-\epsilon d/2}\le C_F.
\end{align*}
For the last inequality use $d\ge24$, the original cap
$\mathcal G(d)^6\le C_* \cdot d$, and the smallness assertion of
Lemma~\ref{lem:nearunit_seed_moments}, after increasing $C_0$.
The equality $F^6 \cdot \rho^{d-1}=\mathcal G(d)^6 \cdot \rho^{d-4}$ avoids
any growth of the terminal degree with $\log F$.

The initial excess is at most $2a_q \cdot u^{q+1}$. For
$N\le\lfloor z^q/\kappa_q\rfloor$, the delayed loss estimate
and matched budget give
\begin{equation*}
 (1-2 \cdot C_\theta \cdot u)P_N
 \le s_q \cdot u^q+C_F \cdot N \cdot \kappa_q \cdot \kappa_{q-2} \cdot u^{2q}
                        +4 \cdot C_\theta \cdot a_q \cdot u^{q+1}
 \le\mathsf K_q \cdot u^q/4.
\end{equation*}
The threshold $z\ge8 \cdot C_\theta$ excludes a first exit and gives
$P_N<p_*/2$ throughout this finite horizon. Shorter prefixes are
covered by their retained actual losses. The successor masses are
positive because $p_*\le\epsilon/1024$ and $\nu<1/4$.

For completeness, write the starting family as a normalized copy of
$B^{1/2}Y$, and set
\begin{equation*}
 Y_j:=\mathcal T^jY,\qquad A_j:=\|Y_j\|_2^2,\qquad B_j:=\|LY_j\|^2,
 \qquad\eta_j:=\frac{A_j}{\langle Y_j,BY_j\rangle}.
\end{equation*}
The hierarchy uses positive scalar multiples of $B^{1/2}Y_j$, so
$\eta_0=\nu$, $\eta_j\ge\lambda$, and
$B_j/A_j=p_j/(v_j \cdot \eta_j)$. The Bochner calculation above gives
\begin{equation*}
 \frac{A_1}{A_0}=\nu^{-1}-\frac{B_0}{A_0},\qquad
 \frac{A_{j+1}}{A_j}\ge\frac{A_j}{A_{j-1}}-\frac{B_j}{A_j}.
\end{equation*}
All denominators are positive on the horizon. Moreover,
\begin{equation*}
 \sum_{j<m}\frac{B_j}{A_j}
 \le\frac{P_m}{\lambda(1-p_*)}\le4 \cdot \mathsf K_q \cdot u^{q-1},\qquad
 \nu^{-1}\ge z-a_q \cdot u^{q-1}.
\end{equation*}
Every ratio is at least $z-C_q \cdot u^{q-1}$. Multiplication proves
Eq.~\eqref{eq:nearunit_block_extension}; finite orthogonal amplification
does not change the operator norm. The definition of $\mathsf K_q$
gives the asserted bound on $C_q$.
\end{proof}

\subsection{A common radius for all coefficient bounds}
\label{sec:orbits_radii}
\label{sec:restricted_radius_geometry}
\label{sec:coefficient_radius_static}

We introduce a single coefficient radius that is independent of the
block length and order. It controls every polynomial degree and bounds
the Poincar\'e constant, and we compare it with the radius of each finite
block. This gives later refinements one scalar quantity to improve,
while postponing the final conversion to the Poincar\'e constant.

We use $B:=H^{-1}$, $R:=R_{\mathrm{res}}=\|\mathcal T\|_{\mathrm{op}}^2$,
$\mathcal T:=P_+\nabla B$, and the coefficient seed $G_*$ in
Eq.~\eqref{eq:quarter_order_seed}. For the orbit and block estimates below, all laws are
centered and regular, have covariance at most $I$, and satisfy
$D^2W\succeq aI$. All families include every ordered output index. Let
$\mathcal G:[1,\infty)\to[1,\infty)$ be a nondecreasing coefficient
majorant satisfying
\begin{equation}\label{eq:second_logstar_majorant}
 1\le\mathcal G(k)\le G_*(k),\qquad
 c_k\le\mathcal G(k)^{k-1},\ \forall k\ge1.
\end{equation}
The bound $\mathcal G\le G_*$ implies, for every fixed $r>0$,
$\mathcal G(k)^r\le C_r \cdot k$ for $k\ge2$, with constants independent of
$\mathcal G$.

Polynomial testing bounds successive coefficients by the same restricted
operator norm. This gives the following estimate for centered regular
measures of covariance at most $I$.

\begin{lemma}[Static coefficients from the restricted radius]
\label{lem:restricted_radius_static}
Every such law satisfies
\begin{equation}\label{eq:restricted_radius_static_coefficients}
 c_d\le R^{(d-1)/2},\ \forall d\ge1.
\end{equation}
In particular $\overline R:=\max\{1,R\}$ is a static coefficient
assignment with floor one in Proposition~\ref{prop:static_profile_coefficients}.
\end{lemma}
\begin{proof}
Take $k:=d-1$ in Eq.~\eqref{eq:block_polynomial_testing}. Symmetrization
is a contraction, and finite output amplification preserves
$\|\mathcal T\|_{\mathrm{op}}=\sqrt R$. With the factorial convention
for $c_d$, the identity gives $c_d\le\sqrt R \cdot c_{d-1}$ for $d\ge2$.
Covariance normalization gives $c_1\le1$, so finite induction proves
the claim, also when $R<1$. This argument uses the same centered
polynomial tests and ordered output spaces as the identity cited above.
\end{proof}

The restricted radius supplies a static coefficient bound. Taking the
supremum of the normalized coefficients gives a common radius, and the
following comparison relates it to each finite block.

For a centered regular measure $\mu$ of covariance at most $I$, we use
$R:=R_{\mathrm{res}}=\|\mathcal T\|_{\mathrm{op}}^2$ and define
\begin{equation}\label{eq:coefficient_radius_definition}
 \mathcal A(\mu):=\max\{1,\sup_{d\ge2}c_d(\mu)^{2/(d-1)}\}.
\end{equation}
This parameter is independent of the chosen block length and order. It
will allow us to transfer coefficients between repetitions without
returning to the restricted radius at every step.

\begin{lemma}[Comparison of the coefficient radius]
\label{lem:coefficient_radius_comparison}
The parameter in Eq.~\eqref{eq:coefficient_radius_definition} is finite and
satisfies
\begin{equation}\label{eq:coefficient_radius_static_bounds}
 1\le\mathcal A(\mu)\le\max\{1,R\},\qquad
 c_d(\mu)\le\mathcal A(\mu)^{(d-1)/2},\ \forall d\ge1.
\end{equation}
Let $\mathcal G:\mathbb N\to[1,\infty)$ be nondecreasing and satisfy
$c_k(\mu)\le\mathcal G(k)^{k-1}$ for every $k\ge1$. For every integer
$m\ge1$, we define $z_m:=\|\mathcal T^m\|_{\mathrm{op}}^{2/m}$. Then
\begin{equation}\label{eq:coefficient_radius_block_maximum}
 \mathcal A(\mu)\le\max\{\mathcal G(m)^2,z_m\}.
\end{equation}
If $D^2W\succeq aI$ for $a>0$, then
\begin{equation}\label{eq:coefficient_radius_poincare}
 C_P(\mu)\le2^{85} \cdot \mathcal A(\mu).
\end{equation}
All constants are independent of the measure, dimension, curvature,
and block length.
\end{lemma}
\begin{proof}
Lemma~\ref{lem:restricted_radius_static} gives
$c_d\le R^{(d-1)/2}$ for $d\ge2$, also when $R<1$, and
covariance normalization gives $c_1\le1$. This proves finiteness and
Eq.~\eqref{eq:coefficient_radius_static_bounds}.

For a fixed $m$, finite iteration of
Eq.~\eqref{eq:block_polynomial_testing} gives
\begin{equation*}
 \frac{Q_d}{d!}
 =\mathsf P_d\frac{Q_{d-m}}{(d-m)!}\mathcal T^m,
 \qquad d\ge m+1.
\end{equation*}
Here $Q_{d-m}$ acts componentwise on the ordered output indices, and
$\mathsf P_d$ symmetrizes all $d$ resulting slots. Indeed,
$\mathcal T$ commutes with permutations of existing output slots,
and full symmetrization absorbs each intermediate symmetrization.
Taking norms, using contraction of $\mathsf P_d$ and the factorial
convention for $c_d$, gives
$c_d\le\|\mathcal T^m\|_{\mathrm{op}} \cdot c_{d-m}
=z_m^{m/2} \cdot c_{d-m}$. Write $d=jm+r$ with
$j\ge0$ and $1\le r\le m$. Iterating this inequality finitely many
times and using the remainder of degree one when $r=1$ gives
\begin{equation*}
 c_d\le z_m^{jm/2} \cdot c_r
     \le z_m^{jm/2} \cdot \mathcal G(m)^{r-1}
     \le\max\{\mathcal G(m)^2,z_m\}^{(d-1)/2}.
\end{equation*}
The maximum is at least one, so taking the supremum proves
Eq.~\eqref{eq:coefficient_radius_block_maximum}. In particular this
argument does not require $\mathcal G(m)^2\le z_m$.

To prove Eq.~\eqref{eq:coefficient_radius_poincare}, apply
Theorem~\ref{thm:dyadic_curvature} with $\epsilon:=1$, $\ell(k):=1$,
and its coefficient parameter equal to
$\mathcal R:=2^{40} \cdot \sqrt{\mathcal A(\mu)}$. This is above the
required absolute threshold. For every $k\ge1$,
\begin{equation*}
 (k+1)^2\le4^k\le2^{40k} \cdot \sqrt{\mathcal A(\mu)},
 \qquad
 c_k(\mu)\le\mathcal A(\mu)^{(k-1)/2}
             \le\frac{\mathcal R^k}{(k+1)^2}.
\end{equation*}
Thus every finite dyadic $d\ge2$ satisfies
\begin{equation*}
 C_P(\mu)\le32 \cdot \mathcal R^2 \cdot \max\{1,a^{-1/(d+1)}\}
            =2^{85} \cdot \mathcal A(\mu) \cdot \max\{1,a^{-1/(d+1)}\}.
\end{equation*}
Since $a>0$ is fixed, letting $d$ tend to infinity through the dyadic
integers proves the claim. This is a limit of scalar inequalities; each
application of the curvature theorem uses only finite degrees and
finite families.
\end{proof}

Equation~\eqref{eq:coefficient_radius_static_bounds} makes
$\mathcal A$ a common static coefficient assignment with floor one in
Proposition~\ref{prop:static_profile_coefficients}. Moreover, every
static coefficient assignment $Z\ge1$ satisfies $\mathcal A\le Z$.
The transfer needs no continuity of $\mathcal A$: as in that proposition,
regular approximation passes the required coefficient
inequalities at each finite degree with all profile parameters fixed. We use
Eq.~\eqref{eq:coefficient_radius_poincare} only after the profile for the coefficient radius
has been established, so its universal factor is paid once.

\subsection{Realizing the common radius by controlled blocks}
\label{sec:second_logstar_geometry}

To apply the curvature comparison at higher orders, we need blocks whose
initial losses are bounded by a prescribed inverse power of their radius.
We control a norm-attaining orbit and average its windows to construct
these starting families, then extend the blocks through successive odd
orders. The resulting radius also controls every polynomial degree,
so the block estimates remain compatible with the coefficient transfer.

\begin{lemma}[Control of a norm-attaining orbit]
\label{lem:second_logstar_green_orbit}
There are universal constants $C_0,c_0,C>0$ with the following property.
Let $0<\epsilon\le1/16$, let $d\ge2$ be dyadic, let
$m\ge2d$ be an integer, and suppose
\begin{equation}\label{eq:second_logstar_green_hypotheses}
 \begin{gathered}
 z:=\|\mathcal T^m\|_{\mathrm{op}}^{2/m},\qquad
 0\le R-z\le2,\qquad
 z\ge\max\{(1+\epsilon) \cdot \mathcal G(d)^2,C_0 \cdot \epsilon^{-1}\},\\
 \gamma:=C_F \cdot d \cdot c_d^2 \cdot z^{-(d-1)},\qquad
 \frac{\gamma \cdot m^2}{z}\le c_0.
 \end{gathered}
\end{equation}
Let $f$ be a centered unit top right singular vector of $\mathcal T^m$,
and define
\begin{equation*}
 w_j:=z^{-j/2}\mathcal T^jf,\qquad b_j:=\|w_j\|_2^2,\qquad
 \ell_j:=\|Lw_j\|^2,\qquad
 D_j:=\|\nabla[(B-zI)w_j]\|_2^2,\ \forall j\ge1.
\end{equation*}
Then
\begin{equation}\label{eq:second_logstar_orbit_bounds}
 \frac12\le b_j\le2,\ \forall 0\le j\le m,\qquad
 b_j\le2,\ \forall j\ge0.
\end{equation}
Moreover, for every integer $1\le J\le m/4$,
\begin{equation}\label{eq:second_logstar_window_loss}
 \sum_{k=2}^{J+1}\sum_{i=k}^{m+k-1}\ell_i
 \le\frac{C}z+C \cdot \gamma \cdot Jm.
\end{equation}
\end{lemma}
\begin{proof}
Compactness supplies the singular vector, as in
Proposition~\ref{prop:block_averaged_start}. We define $t:=R/z$. Maximality
and the one-step operator norm give
\begin{equation}\label{eq:second_logstar_orbit_shift}
 b_0=b_m=1,\qquad b_{j+m}\le b_j,\qquad
 b_{j+1}\le tb_j,\qquad b_j\le t^j,
 \qquad t\le1+2/z.
\end{equation}
Write $S:=\sum_{j=0}^{m-1}b_j$. In particular,
$\sum_{j=0}^{2m-1}b_j\le2S$.

We first use only local propagation in the raw joint frame. A swap at
degree $k$ is propagated through at most $3k-2$ normalized copies of
$\mathcal T$, whose norm is at most $t^{(3k-2)/2}$. The proof of
Lemma~\ref{lem:quarter_order_raw_joint}, with this bound in place of its
global factor, gives, for dyadic $e\ge2$ and $j\ge2e-1$,
\begin{align}
 \ell_j&\le C_F \cdot e \cdot c_e^2 \cdot z^{-(e-1)} \cdot b_{j-e+1}\notag\\
 &\quad+108 \cdot C_F \cdot \sum_{\substack{k<e\\k\ \mathrm{dyadic}}}
 k^5 \cdot c_k^2 \cdot z^{-k} \cdot t^{3k-2}
       \sum_{h=1}^{3k-1}D_{j-k+1-h}.
 \label{eq:second_logstar_local_raw_frame}
\end{align}
Indeed the exact projection remains
$\mathsf P_kQ_1w_j=Q_kw_{j-k+1}/(k! \cdot z^{(k-1)/2})$.
The $h$th old adjacent swap uses $D_{j-k+1-h}$ and is propagated
through $h-1$ copies of $\mathcal T/\sqrt z$. Squaring its sorting
bound and applying Theorem~\ref{thm:joint_uniform_frame} proves
Eq.~\eqref{eq:second_logstar_local_raw_frame}. Every fresh defect has
index at least one, so no gradient symmetry is required of $w_0$.

For $j\ge3$, use the largest available dyadic degree at most $d$, and
denote its coherent term by $h_j$. Set $h_2:=2 \cdot K_2 \cdot b_1/z$ and $D_i:=0$
for $i\le0$. The second-moment estimate
Eq.~\eqref{eq:quarter_order_second_moment} gives
\begin{equation}\label{eq:second_logstar_local_convolution}
 \ell_j\le h_j+\sum_{r\ge1}\alpha_r \cdot D_{j-r},\ \forall j\ge2,
\end{equation}
where the nonnegative finite kernel is defined as follows. Add $8/z$
at delay one and, for every dyadic $k<d$, add
$108 \cdot C_F \cdot k^5 \cdot c_k^2 \cdot z^{-k} \cdot t^{3k-2}$ at each delay
$r=k,\ldots,4k-2$. Thus
\begin{align}
 \eta:=\sum_r\alpha_r
 &\le\frac8z+324 \cdot C_F \cdot \sum_{\substack{k<d\\k\ \mathrm{dyadic}}}
       k^6 \cdot c_k^2 \cdot z^{-k} \cdot t^{3k-2}\le\frac{C}z,\notag\\
 \eta_1:=\sum_r r\alpha_r
 &\le\frac8z+1296 \cdot C_F \cdot \sum_{\substack{k<d\\k\ \mathrm{dyadic}}}
       k^7 \cdot c_k^2 \cdot z^{-k} \cdot t^{3k-2}\le\frac{C}z.
 \label{eq:second_logstar_local_kernel}
\end{align}
To verify these bounds, we define $F^2:=z/(1+\epsilon)$ and
$\rho_k:=\mathcal G(k)^2/F^2$. The hypotheses give
$\epsilon \cdot F^2\ge C$ and $t^3\le1+\epsilon/4$. The original cap
and the seed bound at $d$ imply
\begin{equation*}
 \rho_k\le\min\{1,C_* \cdot k^{1/3}/F^2\},\qquad
 t^3 \cdot \mathcal G(k)^2/z\le\rho_k \cdot e^{-\epsilon/2}
 ,\ \forall k\le d.
\end{equation*}
Lemma~\ref{lem:nearunit_seed_moments} bounds both sums uniformly;
the fixed term $k=1$ is treated separately. All constants are
independent of $\epsilon$.

These bounds also control the early coherent terms before any orbit
bound is known. Degree $k<d$ is used only at
$j=2k-1,\ldots,4k-2$, where $b_{j-k+1}\le t^{3k-1}$.
For $k=2$ use $c_2^2\le2$. For $k\ge4$, its contribution to $z$
times the weighted early sum
is bounded by a constant times
\begin{align*}
 k^3 \cdot t^{3k-1} \cdot \mathcal G(k)^{2k-2} \cdot z^{-(k-2)}
 &=t^5 \cdot k^3 \cdot \mathcal G(k)^2 \cdot (t^3 \cdot \mathcal G(k)^2/z)^{k-2}\\
 &\le C \cdot k^4 \cdot \rho_k^{k-2} \cdot e^{-\epsilon(k-2)/2}.
\end{align*}
Here $\mathcal G(k)^2\le C \cdot k$ follows from the original cap.
Lemma~\ref{lem:nearunit_seed_moments} therefore gives
\begin{equation}\label{eq:second_logstar_early_terms}
 \sum_{j=2}^{2d-2}h_j\le C/z,\qquad
 \sum_{j=2}^{2d-2}j h_j\le C/z.
\end{equation}
At every later index $h_j=\gamma \cdot b_{j-d+1}$. The absolute
threshold in the lemma makes both kernel moments small enough for the
following absorptions, with constants independent of $\epsilon$.

We now sum the defects on two blocks. Eq.~\eqref{eq:block_fixed_defect}
states
\begin{equation}\label{eq:second_logstar_second_difference}
 D_j\le z(b_{j+1}-2b_j+b_{j-1})+\ell_j.
\end{equation}
Writing $\delta_j:=b_j-b_{j-1}$, telescoping over the first block gives
\begin{equation}\label{eq:second_logstar_first_telescope}
 D_{\mathrm{int}}:=\sum_{j=1}^mD_j
 \le\sum_{j=1}^m\ell_j-z(b_1-b_{m+1})
 \le\sum_{j=1}^m\ell_j.
\end{equation}
Here $\ell_1\le b_1\le t\le2$. The terminal coherent sum on this
block is at most $\gamma \cdot S$. Equations
\eqref{eq:second_logstar_local_convolution}--\eqref{eq:second_logstar_early_terms}
and absorption at the stated absolute threshold therefore imply
\begin{equation}\label{eq:second_logstar_first_sums}
 D_{\mathrm{int}}\le C+C \cdot \gamma \cdot S,\qquad
 \sum_{j=1}^m\ell_j\le C+C \cdot \gamma \cdot S,\qquad
 \sum_{j=2}^m\ell_j\le C/z+C \cdot \gamma \cdot S.
\end{equation}
Some $\delta_j$ with $1\le j\le m$ is nonnegative, since their sum
is zero. Eq.~\eqref{eq:second_logstar_second_difference} gives
$\delta_{j+1}\ge\delta_j-\ell_j/z$, whence
$\delta_{m+1}\ge-z^{-1}\sum_{j=1}^m\ell_j$.
Also $b_{2m}\le1$ and
$\delta_{2m+1}\le(t-1)b_{2m}\le2/z$.
The telescoping sum over the second block is thus bounded by its moments plus
$2+\sum_{j=1}^m\ell_j$. The coherent sum through $2m$ is at most
$2\gamma \cdot S$, and all source defects have index at most $2m-1$.
A second absorption gives
\begin{equation}\label{eq:second_logstar_two_block_sums}
 D_{\mathrm{tot}}:=\sum_{j=1}^{2m}D_j\le C+C \cdot \gamma \cdot S,\qquad
 \sum_{j=2}^{2m}\ell_j\le C/z+C \cdot \gamma \cdot S.
\end{equation}

For the upper orbit bound, fix $1\le j<m$ and let
\begin{equation*}
 g_j(i):=\frac{\min\{i,j\}(m-\max\{i,j\})}{m}
 ,\ \forall 1\le i<m,
\end{equation*}
extended by zero outside this interval. This is a nonnegative
$1$-Lipschitz weight with $g_j(1)\le1$ and $g_j(i)\le i$.
Discrete summation by parts gives
\begin{equation*}
 \sum_i g_j(i)(b_{i+1}-2b_i+b_{i-1})=1-b_j.
\end{equation*}
We define $E_g:=\sum_i g_j(i)\ell_i$ and $F_g:=\sum_i g_j(i)D_i$.
Then $F_g\le E_g+z(1-b_j)\le E_g+z$.
The term at index one is bounded by two. For the remaining terms,
use $g_j(i+r)\le g_j(i)+r$, the two kernel moments, and
Eq.~\eqref{eq:second_logstar_early_terms}. The terminal coherent
contribution is at most $\gamma \cdot mS$, so
\begin{equation*}
 E_g\le C+C \cdot \gamma \cdot mS+\eta \cdot F_g+\eta_1 \cdot D_{\mathrm{int}}
 \le C+C \cdot \gamma \cdot mS+(C/z)E_g.
\end{equation*}
Since $D_i\ge0$, absorption and the Green identity imply
\begin{equation}\label{eq:second_logstar_upper_green}
 b_j\le1+E_g/z\le1+C/z+C \cdot \gamma \cdot mS/z.
\end{equation}
Let $M:=\max_{0\le j\le m}b_j$, so $S\le mM$.
Choosing $C_0$ large and $c_0$ small in
Eq.~\eqref{eq:second_logstar_green_hypotheses} makes
$M\le1+C/z+C \cdot c_0 \cdot M$ imply $M\le2$.
The shift in Eq.~\eqref{eq:second_logstar_orbit_shift} then gives
$b_j\le2$ for every future index.

For the lower bound, define the multiplicity of $m$ sliding windows by
\begin{equation*}
 n_m(i):=\#\{k\in\{2,\ldots,m+1\}:k\le i\le m+k-1\},\qquad
 E_m:=\sum_{i\ge1}n_m(i)\ell_i,\qquad
 F_m:=\sum_{i\ge1}n_m(i)D_i.
\end{equation*}
This weight is $1$-Lipschitz, supported on $2,\ldots,2m$, and satisfies
$n_m(i)\le\min\{i-1,m\}$ for $i\ge1$.
Summing Eq.~\eqref{eq:second_logstar_second_difference} over these
windows gives
\begin{align*}
 F_m&\le E_m+z(b_{2m+1}-2b_{m+1}+b_1)\\
 &\le E_m+z(b_1-b_{m+1})\le E_m+C+C \cdot \gamma \cdot S.
\end{align*}
Here $b_{2m+1}\le b_{m+1}$, and the last inequality follows from
Eq.~\eqref{eq:second_logstar_first_telescope}.
The weighted early coherent sum is at most $C/z$; the terminal part
is at most $2\gamma \cdot mS$. Applying the convolution estimate with
$n_m(i+r)\le n_m(i)+r$ and
Eq.~\eqref{eq:second_logstar_two_block_sums} yields
\begin{equation}\label{eq:second_logstar_long_window}
 E_m\le C/z+C \cdot \gamma \cdot mS+\eta \cdot F_m+\eta_1 \cdot D_{\mathrm{tot}}
 \le C/z+C \cdot \gamma \cdot mS+(C/z)E_m.
\end{equation}
Thus $E_m\le C/z+C \cdot \gamma \cdot mS$ after absorption.

Fix $1\le j<m$ and apply the Dirichlet Green formula on $[j,j+m]$
at the point $m$. The endpoints satisfy $b_{j+m}\le b_j$, and
Eq.~\eqref{eq:second_logstar_second_difference} gives
\begin{equation*}
 1=b_m\le b_j+\frac1z
 \sum_{i=j+1}^{j+m-1}
 \frac{(\min\{m,i\}-j)(j+m-\max\{m,i\})}{m}\ell_i.
\end{equation*}
The Green weight is at most $n_m(i)$. For $i\le m$ it is at most
$i-j\le i-1=n_m(i)$; for $i>m$ it is at most
$j+m-i\le2m-i\le n_m(i)$. Consequently
\begin{equation}\label{eq:second_logstar_lower_green}
 b_j\ge1-E_m/z\ge1-C/z^2-C \cdot \gamma \cdot mS/z\ge1/2,
\end{equation}
using $S\le2m$ and the same choices of $C_0,c_0$.
The endpoints are already equal to one. This proves
Eq.~\eqref{eq:second_logstar_orbit_bounds}.

Finally repeat the weighted calculation with $J$ windows rather than
$m$ windows. Its boundary is at most $z(b_1-b_{m+1})$ by the shift.
The early weighted coherent sum remains $C/z$, and the terminal
weighted sum is at most $2\gamma \cdot Jm$, now using $b_i\le2$.
The same kernel moments and defect bound over two blocks therefore give
Eq.~\eqref{eq:second_logstar_window_loss}. All indices lie within the
two blocks already used. The finite Bochner domain argument from
Lemma~\ref{lem:quantitative_raw_frame} justifies the pairings throughout.
\end{proof}

The orbit and window estimates allow us to average a finite segment into
one starting family. Its restart preserves the comparison between the
losses and the block's own energy budget.

For odd $p\ge3$, let $d_p^\circ$ be the least dyadic integer at least
$4p$, and define
\begin{equation}\label{eq:second_logstar_coherent_constants}
 4p\le d_p^\circ<8p,\qquad
 \mathsf A_p^\circ:=C_F \cdot d_p^\circ \cdot
                 \mathcal G(d_p^\circ)^{2d_p^\circ-2}.
\end{equation}

\begin{lemma}[Initialization from a local orbit]
\label{lem:second_logstar_initialization}
There are universal constants $C,c_0>0$ with the following property.
Suppose $0<\epsilon\le1/16$, $p\ge3$ is odd, $\kappa\ge1$,
$m,J_0$ are positive integers, $d:=d_p^\circ$, and
\begin{gather}
 z:=\|\mathcal T^m\|_{\mathrm{op}}^{2/m},\qquad
 0\le R-z\le2,\qquad
 z\ge\max\{(1+\epsilon)\mathcal G(d)^2,C\epsilon^{-1}\},\quad
 \gamma:=C_F \cdot d \cdot c_d^2 \cdot z^{-(d-1)},\notag\\
 \frac{z^{p-2}}{4\kappa}\le m\le\frac{4z^{p-2}}\kappa,
 \qquad m\ge2d,\qquad J_0\le m/4,\notag\\
 \frac{\gamma \cdot m^2}{z}\le c_0,\qquad
 J_0 \cdot \gamma \cdot zm\le1,\qquad C \cdot \kappa \cdot z^{-p}\le1/8.
 \label{eq:second_logstar_coherent_smallness}
\end{gather}
There is one centered unit finite starting family satisfying, for $u:=1/z$,
\begin{equation}\label{eq:second_logstar_initialization}
 \nu\le\frac{u}{1-C \cdot \kappa \cdot u^p},\qquad
 P_J\le C \cdot \kappa \cdot u^p,\ \forall 1\le J\le J_0,\qquad
 a \cdot V_N+X_N\le\nu-u+u \cdot P_{N+1},\ \forall N\ge0.
\end{equation}
\end{lemma}
\begin{proof}
The stated conditions give every hypothesis of
Lemma~\ref{lem:second_logstar_green_orbit} with $d:=d_p^\circ$.
Use its actual orbit and the same averaged family as before:
\begin{equation*}
 \mathcal W:=\bigoplus_{j=0}^{m-1}w_j,\qquad
 S_k:=\sum_{j=k}^{m+k-1}b_j,\qquad
 C_k:=\sum_{j=k}^{m+k-1}\ell_j,\qquad
 U:=\frac{\mathcal T\mathcal W}{\sqrt{z \cdot S_0}}.
\end{equation*}
Since $S_0=S_1$, the family $U$ is centered, unit, and a genuine
gradient, with energy at most $1/z$. Quadratic testing of the first
summand gives $\|\operatorname{Sym}(Lw_1)\|^2\le K_2/(4z)$.
Eq.~\eqref{eq:second_logstar_first_sums} controls the remaining
summands, and $S_0\ge m/2$. Hence
\begin{equation*}
 \|\operatorname{Sym}(LU)\|_{\mathrm{HS}}^2
 \le\frac{C}{zm}+C \cdot \gamma.
\end{equation*}
Start at $F_0:=B^{1/2}\mathcal W/\|B^{1/2}\mathcal W\|_2$ and take
$G_0:=F_2/\|F_2\|_2$. The exact skew credit in
Lemma~\ref{lem:skew_normalization_credit} gives
\begin{equation*}
 \nu\le\frac{u}{1-C/(z^2 \cdot m)-C \cdot \gamma/z}
 \le\frac{u}{1-C \cdot \kappa \cdot u^p},
\end{equation*}
using Eq.~\eqref{eq:second_logstar_coherent_smallness}.
The stated smallness condition makes the denominator positive.

For $J\le J_0$, the actual loss fractions are
$C_k/(z \cdot S_{k+1}+C_k)$ at $k=2,\ldots,J+1$.
Every denominator window contains at least $m-J-1\ge m/2$ interior
indices, so $S_{k+1}\ge m/4$. The unit restarted masses are at most
one. Eq.~\eqref{eq:second_logstar_window_loss} therefore gives
\begin{equation*}
 P_J\le\frac C{zm}\sum_{k=2}^{J+1}C_k
 \le\frac{C}{z^2 \cdot m}+C \cdot \gamma \cdot J/z
 \le C \cdot \kappa \cdot u^p.
\end{equation*}
The same positive interior terms justify every normalization used.
Finally $S_{k+1}\le S_k$ for all $k$ by the norm-attaining shift.
Writing the hierarchy as positive scalar multiples of
$B^{1/2}\mathcal T^j\mathcal W$ gives $v_{j+1}\le z \cdot e_j$, exactly as
in Eq.~\eqref{eq:quarter_order_matched_budget}. Combining this with
$a \cdot V_N+X_N\le\nu-e_N$ proves the last assertion of
Eq.~\eqref{eq:second_logstar_initialization}. The absorbing zero
convention handles any later vanishing family.
\end{proof}

We next extend these initialized blocks to the required order while
keeping the coefficient floor close to the squared majorant. The common
radius identifies the branch in which each actual family is available.

\begin{proposition}[A block floor with a small multiplicative loss]
\label{prop:nearunit_logstar_blocks}
There are universal constants $C_X,C_0,C>0$ with the following property.
Let $\mathcal G$ satisfy Eq.~\eqref{eq:second_logstar_majorant} for all
centered regular laws of covariance at most $I$. Fix
$0<\epsilon\le1/16$ and an odd integer $Q\ge3$, and define
\begin{equation}\label{eq:nearunit_geometry_parameters}
 X:=C_X \cdot Q,\qquad
 F^2:=\max\{\mathcal G(X^X)^2,C_0/\epsilon\},\qquad
 J_Q:=32Q.
\end{equation}
There are constants $H_Q,a_Q,s_Q>0$ and an assignment $Z_Q(\mu)$
satisfying
\begin{align}
 H_Q&:=(1+2\epsilon) \cdot F^2+1,\notag\\
 a_Q+s_Q+1&\le C \cdot F^{2Q-6}\le F^{2Q},\notag\\
 c_d(\mu)&\le Z_Q(\mu)^{(d-1)/2},\ \forall d\ge1,\notag\\
 H_Q&\le Z_Q(\mu)\le\max\{H_Q,C_P(\mu)\},\qquad
 C_P(\mu)\le Z_Q(\mu)+3.
 \label{eq:nearunit_geometry_static}
\end{align}
If $Z_Q(\mu)>H_Q$, we define $z:=Z_Q(\mu)$ and $u:=1/z$. One centered unit
finite starting family satisfies
\begin{equation}\label{eq:nearunit_geometry_startup}
 \nu\le\frac{u}{1-a_Q \cdot u^Q},\qquad
 P_J\le s_Q \cdot u^Q,\ \forall 1\le J\le J_Q,\qquad
 a \cdot V_N+X_N\le\nu-u+u \cdot P_{N+1},\ \forall N\ge0.
\end{equation}
Moreover,
\begin{equation}\label{eq:nearunit_geometry_propagation}
 \|\mathcal T^h\|_{\mathrm{op}}\le4(h+1) \cdot z^{h/2},\ \forall h\ge0,\qquad
 0\le R-z\le1,\qquad a_Q \cdot u^Q\le1/8.
\end{equation}
All constants are independent of $\epsilon,Q$ and $\mathcal G$.
\end{proposition}
\begin{proof}
We define $R_0:=(1+2\epsilon) \cdot F^2$ and $H_Q:=R_0+1$.
If $R<R_0$, we define $Z_Q:=H_Q$. Then $C_P\le R+1\le H_Q$,
and ordinary Poincar\'e iteration gives the coefficient bound.
The assertions requiring $Z_Q>H_Q$ do not apply in this case.
The constants $a_Q,s_Q$ used in both cases will be fixed below
independently of the law. We henceforth assume $R\ge R_0$.

Write $D:=1+\epsilon$ and define
\begin{equation}\label{eq:second_logstar_divisors}
 \kappa_q:=F^{2q-2},\ \forall q\ge1\text{ odd},\qquad J_p:=32p.
\end{equation}
In particular $\kappa_1=1$. The initialization degree
$d_p^\circ$ lies in $[4p,8p)$, and the extension degree $e_q$ is
the least dyadic integer at least $8q$, so $e_q<16q$.
Choose $C_X\ge16$ and above the fixed static threshold used below.
Every tested degree is then at most $X$, hence has seed at most $F$.
There is no implicit dependence of $X$ on $F$.

We check all initialization conditions at the single absolute floor
$\epsilon \cdot F^2\ge C_0$. Let $p\ge5$. We define $d:=d_p^\circ$,
$G_p:=\mathcal G(d)$, and $\rho_p:=G_p^2/F^2$. The original cap gives
\begin{equation*}
 \rho_p\le\min\{1,C \cdot p^{1/3}/F^2\}.
\end{equation*}
Suppose $z\ge DF^2$ and the radius and length comparisons of
Lemma~\ref{lem:second_logstar_initialization} hold with
$\kappa:=\kappa_{p-2}$. Its coherent constant obeys
\begin{align}
 \frac{\gamma \cdot m^2}{z}
 &\le C \cdot d \cdot G_p^2 \cdot \rho_p^{d-2} \cdot D^{-(d-2p+4)}
 \le C \cdot p^{4/3} \cdot \rho_p^{4p-2} \cdot e^{-\epsilon p},
 \label{eq:second_logstar_green_order_smallness}\\
 J_p \cdot \gamma \cdot zm
 &\le C \cdot J_p \cdot d \cdot G_p^6 \cdot F^{-2} \cdot \rho_p^{d-4} \cdot D^{p-d}
 \le C \cdot p^3 \cdot F^{-2} \cdot \rho_p^{4p-4} \cdot e^{-\epsilon p}.
 \label{eq:second_logstar_prefix_order_smallness}
\end{align}
These quantities are uniformly small when $C_0$ is large.
For $p\le F^3$, one has $\rho_p\le C/F\le F^{-1/2}$, and the
fixed polynomial factors times these powers tend uniformly to zero.
For $p>F^3$, discard $\rho_p\le1$ and use
$\epsilon \cdot F^3\ge C_0 \cdot F$. The exponential tail then dominates every
fixed polynomial factor, exactly as in
Lemma~\ref{lem:nearunit_seed_moments}. Thus the two displays are at
most $c_0$ and one, respectively.

The lower length bound gives $m\ge D^{p-2} \cdot F^2/4$, so
\begin{equation*}
 J_p/m\le C/(\epsilon \cdot F^2),\qquad
 \kappa_{p-2} \cdot z^{-p}\le F^{-6} \cdot D^{-p}.
\end{equation*}
One fixed choice of $C_0$ therefore gives $J_p\le m/4$, $m\ge2d$,
and every initialization smallness condition. The raw seed threshold
follows from $G_p\le F$ and $z\ge C_0/\epsilon$.
The initialization lemma supplies a family through $J_p$, with
\begin{equation}\label{eq:second_logstar_uniform_startup_cost}
 a_p+s_p+1\le C\kappa_{p-2}.
\end{equation}

The base is $m_3:=\lfloor R\rfloor$ and
$z_3:=\|\mathcal T^{m_3}\|_{\mathrm{op}}^{2/m_3}$.
Lemma~\ref{lem:mesoscopic_block_norm} gives $R-z_3\le C_*/R$;
Eq.~\eqref{eq:block_global_propagation} gives a global factor at most
two at a fixed large $R$. The actual divisor is one, consistently
with $\kappa_1=1$. Here $d_3^\circ=16$ has coherent constant at most
$16 \cdot C_F \cdot G_*(16)^{30}$. Since $m_3$ is comparable to $z_3$, the coherent
quantities are $O(z_3^{-14})$ and $O(J_3z_3^{-13})$, with $J_3=96$.
They are small at the chosen floor, which also gives $m_3\ge4J_3$.
Lemma~\ref{lem:second_logstar_initialization} now proves the base case of
Eq.~\eqref{eq:second_logstar_uniform_startup_cost}, together with the
actual family and matched budget. The required envelope follows from
Eq.~\eqref{eq:block_global_propagation}.

Apply Lemma~\ref{lem:nearunit_block_extension} with the common
$\epsilon,F$, the majorant $\mathcal G$, and envelope $4(h+1)$.
Its degree satisfies $2e_q-1<32q=J_q$. Its kernel constant is universal,
and Eq.~\eqref{eq:second_logstar_uniform_startup_cost} gives a bridge
constant at most $C_b \cdot \kappa_{q-2}$, for one universal $C_b$.
All required normalizer thresholds hold because its exit mass is at most
$C \cdot F^{-6} \cdot D^{-q}\le\epsilon/1024$. Fix $C_b$ before choosing the final
absolute floor constant.

At each intermediate order, we define
\begin{equation}\label{eq:second_logstar_block_ladder}
 m_{q+2}:=\lfloor z_q^q/\kappa_q\rfloor,\qquad
 z_{q+2}:=\|\mathcal T^{m_{q+2}}\|_{\mathrm{op}}^{2/m_{q+2}},\qquad
 \Delta_q:=C_b \cdot \kappa_{q-2} \cdot z_q^{1-q}.
\end{equation}
This length belongs to the exact horizon of
Lemma~\ref{lem:nearunit_block_extension}; hence
$z_{q+2}\ge z_q-\Delta_q$. At $z_q\ge DF^2$,
\begin{equation}\label{eq:second_logstar_gap_bounds}
 \frac{\Delta_q}{z_q}\le C_b \cdot F^{-6} \cdot D^{-q},\qquad
 \frac{q\Delta_q}{z_q}\le\frac C{\epsilon \cdot F^6}
 \le\epsilon^2/100.
\end{equation}
The last inequality follows from $\epsilon \cdot F^2\ge C_0$ after fixing
$C_0$ large. Provided the new radius is at least $DF^2$, its next
possible decrement satisfies
\begin{equation}\label{eq:second_logstar_gap_ratio}
 \frac{\Delta_{q+2}}{\Delta_q}
 =\frac{F^4}{z_{q+2}^2}(z_q/z_{q+2})^{q-1}
 \le D^{-2} \cdot e^{\epsilon^2/50}\le1-\epsilon.
\end{equation}
Increasing radii improve this estimate. Every constructed finite prefix
ending at $q'$ consequently satisfies
\begin{equation}\label{eq:second_logstar_gap_tail}
 (z_q-z_{q'})_+\le\Delta_q/\epsilon,\qquad
 R-z_{q'}\le C/R+C/(\epsilon \cdot F^4)<\min\{1,\log2/2\}.
\end{equation}
Every radius is at least $DF^2$, since $R\ge R_0=(1+2\epsilon) \cdot F^2$
and $\epsilon \cdot F^2\ge C_0$. In the induction, the already known decrements
establish this lower bound at the new radius before
Eq.~\eqref{eq:second_logstar_gap_ratio} is used for its next decrement.
No monotonicity of the actual radii is assumed.

All defining real block lengths are at least $D^q \cdot F^2$, so their
relative floor losses are at most $\epsilon/100$. For $q\ge5$,
Eq.~\eqref{eq:second_logstar_gap_bounds} gives
\begin{equation}\label{eq:second_logstar_length_lower}
 \frac{m_{q+2}}{m_q}
 \ge(1-\epsilon/100)\frac{z_q^2}{F^4}
                    (z_q/z_{q-2})^{q-2}\ge1+\epsilon.
\end{equation}
At the separate first step,
\begin{equation*}
 \frac{m_5}{m_3}\ge(1-\epsilon/100)
                 \frac{z_3^2}{F^4}\frac{z_3}{R}\ge1+\epsilon,
\end{equation*}
using the mesoscopic relative gap $O(R^{-2})$. The lower length
comparison is thus valid at every stage.

For any finite prefix ending at $q'$, greedily decompose $h\ge0$ as
\begin{equation*}
 h=\sum_{q=3,5,\ldots,q'}n_q \cdot m_q+r,\qquad0\le r<m_3.
\end{equation*}
Every nonterminal level satisfies $n_q \cdot m_q<m_{q+2}$. Exact block norms
and submultiplicativity bound $\|\mathcal T^h\|_{\mathrm{op}}$ by
$R^{r/2}\prod_qz_q^{n_q \cdot m_q/2}$. By
Eq.~\eqref{eq:second_logstar_gap_tail}, each nonterminal level has
logarithmic excess over $z_{q'}^{h/2}$ at most
\begin{equation*}
 \frac{C \cdot m_{q+2} \cdot \Delta_q}{\epsilon \cdot z_q}
 \le\frac C{\epsilon \cdot F^4}.
\end{equation*}
There are at most $1+\log(h+1)/\log(1+\epsilon)$ such levels.
Their polynomial exponent is at most $C/(\epsilon^2 \cdot F^4)\le1/2$.
The factor for one level and the base remainder are each at most two,
by the same gap estimate and $m_3\le R$. Therefore
\begin{equation*}
 \|\mathcal T^h\|_{\mathrm{op}}\le4(h+1) \cdot z_{q'}^{h/2}.
\end{equation*}
This bound holds at each already constructed prefix before it is
extended, with the same fixed factor four.

For the remaining radius--length comparison we keep the finite number
of levels, rather than the preceding bound in terms of $\log(h+1)$.
A prefix ending at $q$ has at most $q/2$ nonterminal levels. The same
calculation gives, for every $h\ge0$,
\begin{equation}\label{eq:nearunit_finite_prefix_factor}
 \|\mathcal T^h\|_{\mathrm{op}}\le E_q \cdot z_q^{h/2},\qquad
 \log E_q\le C/R+C \cdot q/(\epsilon \cdot F^4).
\end{equation}
Use this bound only at $m:=m_{q+2}\ge D^q \cdot F^2/2$. It gives
\begin{align*}
 q(\log(z_{q+2}/z_q))_+
 &\le\frac{C \cdot q \cdot D^{-q}}{F^2 \cdot R}
       +\frac{C \cdot q^2 \cdot D^{-q}}{\epsilon \cdot F^6}\\
 &\le\frac C{\epsilon \cdot F^2 \cdot R}+\frac C{\epsilon^3 \cdot F^6}.
\end{align*}
Both terms are uniformly small at $\epsilon \cdot F^2\ge C_0$.
Together with Eq.~\eqref{eq:second_logstar_gap_bounds} and the floor
loss, this proves
\begin{equation}\label{eq:second_logstar_radius_length_comparison}
 \frac{z_{q+2}^q}{4\kappa_q}\le m_{q+2}
 \le\frac{4z_{q+2}^q}{\kappa_q}.
\end{equation}
The prefix range follows from $J_{q+2}/m\le C/(\epsilon \cdot F^2)$.
The coherent and normalization conditions were checked above, so initialization at this radius proves
Eq.~\eqref{eq:second_logstar_uniform_startup_cost} for order $q+2$.
This closes the finite induction, including the actual family, budget,
radius, length, and envelope hypotheses. The factor $E_q$ is used only
for this comparison; the uniform polynomial envelope remains the one
used in every kernel estimate.

At the final order,
\begin{equation*}
 a_Q+s_Q+1\le C \cdot F^{2Q-6}\le F^{2Q},
\end{equation*}
since $F^6\ge C$. This includes $Q=3$ and its actual base divisor one.
The same floor gives $a_Qz_Q^{-Q}\le1/8$.

It remains to check all static coefficient degrees. We have $m_Q\le R^Q$.
If $R\le X$, then $Q\le X$ and $m_Q\le X^X$, so
$\mathcal G(m_Q)^2\le F^2<z_Q$ exactly. If $R>X$, use the original cap:
\begin{equation}\label{eq:nearunit_geometry_static_seed}
 \mathcal G(m_Q)^2\le G_*(R^R)^2
 \le C \cdot G_*(R)^2\le R/2\le z_Q.
\end{equation}
Choose $C_X$ above the fixed threshold for the last inequality.
Monotonicity controls all degrees at most $m_Q$;
Eq.~\eqref{eq:block_coefficient_recursion} controls every higher degree
with the same radius. No fixed multiple of $F^2$ enters this step.

We define $Z_Q:=\max\{H_Q,z_Q\}$ in the present case $R\ge R_0$.
We have $z_Q\le R\le C_P$ and $C_P\le R+1\le z_Q+2$.
Thus every static assertion holds. Whenever $Z_Q>H_Q$, it is the
actual constructed radius, so its starting family and budget apply.
The formula for $F^2$ gives the exact floor
\begin{equation*}
 H_Q=(1+2\epsilon)
       \max\{\mathcal G((C_X \cdot Q)^{C_X \cdot Q})^2,C_0/\epsilon\}+1.
\end{equation*}
All choices are finite, and all constants use only the original cap,
not the amplitude or shape of the smaller coefficient seed.
\end{proof}

Fixing the loss parameter gives the block estimates used in repeated
height reduction. The next proposition records the radius, initialization,
propagation, and coherent error bounds for any coefficient majorant
in the stated class.

\begin{proposition}[Blocks from a coefficient majorant]
\label{prop:second_logstar_blocks}
Fix a universal $C_{\mathrm{arg}}$ large enough that
$64Q \cdot \log(e+G_*(8Q))\le C_{\mathrm{arg}}(8Q)^2$ for $Q\ge3$.
Suppose $\mathcal G$ satisfies Eq.~\eqref{eq:second_logstar_majorant}
for all centered regular laws of covariance at most $I$ and, for one
fixed $C_{\mathrm{seed}}\ge1$,
\begin{align}
 \mathcal G(C_{\mathrm{arg}} \cdot x^2)
 &\le C_{\mathrm{seed}} \cdot \mathcal G(x),\ \forall x\ge3,\notag\\
 \mathcal G(R^Q)^2
 &\le C_{\mathrm{seed}} \cdot \mathcal G(8Q)^2+
       C_{\mathrm{seed}} \cdot G_*(R)^2,\ \forall R,Q\ge3.
 \label{eq:second_logstar_seed_comparisons}
\end{align}
There is a constant $K$, depending only on these fixed comparison
constants and the previously fixed universal constants, such that every odd
$Q\ge3$ has an assignment $Z_Q$ and constants $H_Q,a_Q,s_Q>0$ with
\begin{align}
 H_Q&\le K \cdot \mathcal G(8Q)^2,\qquad
 a_Q+s_Q+1\le(K \cdot \mathcal G(8Q)^2)^Q,\notag\\
 c_d(\mu)&\le Z_Q(\mu)^{(d-1)/2},\ \forall d\ge1,\notag\\
 H_Q&\le Z_Q(\mu)\le\max\{H_Q,C_P(\mu)\},\qquad
 C_P(\mu)\le Z_Q(\mu)+3.
 \label{eq:second_logstar_static}
\end{align}
If $Z_Q>H_Q$, we define $z:=Z_Q$, $u:=z^{-1}$, and
$J_Q:=\lceil128Q \cdot \log(e+\mathcal G(8Q))\rceil$. One centered unit
finite starting family satisfies
\begin{equation}\label{eq:second_logstar_startup}
 \nu\le\frac{u}{1-a_Q \cdot u^Q},\qquad
 P_J\le s_Q \cdot u^Q,\ \forall 1\le J\le J_Q,\qquad
 a \cdot V_N+X_N\le\nu-u+u \cdot P_{N+1},\ \forall N\ge0,
\end{equation}
and
\begin{equation}\label{eq:second_logstar_propagation}
 \|\mathcal T^h\|_{\mathrm{op}}\le4(h+1) \cdot z^{h/2},\ \forall h\ge0,\qquad
 0\le R-z\le1,\qquad a_Q \cdot u^Q\le1/8.
\end{equation}
The constant $K$ is uniform over all majorants with the fixed stated
comparison constants.
\end{proposition}
\begin{proof}
Fix $\epsilon_0:=1/16$ and repeat the construction of
Proposition~\ref{prop:nearunit_logstar_blocks} with
\begin{equation*}
 F:=C_2 \cdot (\mathcal G(8Q)+1),\qquad
 \widetilde J_p:=\lceil128p \cdot \log(e+F)\rceil,
\end{equation*}
where $C_2$ depends only on the fixed constants for the class of coefficient seeds.
At every initialization we construct the family through this longer
prefix, verifying the additional terms directly.

Every initialization degree is below $8Q$, and every extension degree
is below $16Q\le C_{\mathrm{arg}}(8Q)^2$. The first comparison in
Eq.~\eqref{eq:second_logstar_seed_comparisons} bounds all these seeds
by $C_{\mathrm{seed}} \cdot \mathcal G(8Q)\le F$, after choosing $C_2$.
The initialization lemma permits any prefix satisfying its stated
conditions. At fixed $D:=1+\epsilon_0$, the length bound gives
\begin{equation*}
 \widetilde J_p/m\le C \cdot \log(e+F)/F^2,
\end{equation*}
uniformly in $p$. The right side is small once $C_2$ is fixed large.
The first coherent smallness condition is unchanged. The second is
now bounded by
\begin{equation*}
 C \cdot p^3 \cdot \log(e+F) \cdot F^{-2} \cdot \rho_p^{4p-4} \cdot e^{-\epsilon_0p}.
\end{equation*}
For $p\le F^3$, use $\rho_p\le F^{-1/2}$; the fixed powers dominate
both the polynomial factor and $\log(e+F)$. For $p>F^3$, the fixed
exponential tail does so. This quantity is uniformly small at large
$C_2$. The base $p=3$ separately has
$\widetilde J_3=O(\log(e+F))$, $m_3$ comparable to $z_3\ge DF^2$,
and second coherent term $O(\widetilde J_3 \cdot z_3^{-13})$.
Thus its longer prefix also satisfies the initialization hypotheses.
All energy denominators and normalization conditions are unchanged.
Initialization of this same family consequently gives
$a_p+s_p+1\le C\kappa_{p-2}$ through $\widetilde J_p$.
This prefix contains $32p$, so the block extension and the entire
radius, length, gap, and propagation induction remain valid.

For the static remainder $m_Q\le R^Q$, the second comparison gives
\begin{equation*}
 \mathcal G(m_Q)^2
 \le C_{\mathrm{seed}} \cdot \mathcal G(8Q)^2+C_{\mathrm{seed}} \cdot G_*(R)^2
 \le R/2\le z_Q.
\end{equation*}
At $R\ge(1+2\epsilon_0)F^2$, the first term is at most $R/4$ by
$F\ge C_2 \cdot \mathcal G(8Q)$; the second is at most $R/4$ by the fixed
sublinear growth of $G_*(R)^2$, after one fixed enlargement of $C_2$.
The exact coefficient recurrence gives every higher degree as before.
We also choose $C_2$ above all fixed absolute thresholds.

Finally $F\le2 \cdot C_2 \cdot \mathcal G(8Q)$, and
$\widetilde J_Q\ge\lceil128Q \cdot \log(e+\mathcal G(8Q))\rceil$.
The floor and startup bounds yield
Eq.~\eqref{eq:second_logstar_static} with one constant $K$, while
this prefix gives Eq.~\eqref{eq:second_logstar_startup}.
The budget and propagation bounds were preserved throughout the
construction. All parameters in these block estimates are fixed independently of
any subsequent curvature iteration.
\end{proof}

\subsection{From a coefficient seed to a full curvature profile}
\label{sec:profile_improvement}
\label{sec:second_logstar_profiles}

We now convert a coefficient seed and one initial curvature bound into
a profile valid at every subsequent inner depth. The block order reduces
the depth cost to $(r+1)^{1/Q}$, with the dependence on the seed kept
explicit. Controlling the logarithmic distortion closes this depth
induction, and the original seed supplies the first such profile.

The curvature profiles below concern centered regular
laws of covariance at most $I$ and curvature at least $aI$, where $a>0$.
We use $g(x):=\log(e+x)$ and $\ell_r:=g^{\circ r}$, together with
$\mathfrak t(x):=1+\log^*(x+2)$, its compositions $\mathfrak t_j$, and
$t_Q:=\mathfrak t(Q)$. The coefficient seed $G_*$ is defined in
Eq.~\eqref{eq:quarter_order_seed}. The block radius $Z_Q$ and its starting
family are supplied by Proposition~\ref{prop:second_logstar_blocks};
Proposition~\ref{prop:static_profile_coefficients} transfers a profile
for this radius to polynomial coefficients.

\begin{lemma}[Finite burn-in for the logarithmic distortion]
\label{lem:repeated_height_distortion}
We define $\eta:=(e+1)^{-1}$. For $x,M\ge1$ and $r\ge1$, define
$F(v):=\log(1+\eta v)$ and $b(M):=1+\log^*M$. Then
\begin{equation}\label{eq:repeated_height_exact_distortion}
 0\le\ell_r(Mx)-\ell_r(x)\le F^{\circ(r-1)}(\log M).
\end{equation}
For $r\ge b(M)$, the right side is at most $2\eta^{r-b(M)}$.
Consequently, for each fixed universal $C_T\ge1$, there is a universal
$C_b$ such that, for
\begin{equation}\label{eq:second_logstar_profile_depth}
 r_Q:=\max\{r_*,\lceil C_b \cdot t_Q\rceil\},
\end{equation}
the terminal degree
$d<2 \cdot C_T \cdot Q^2 \cdot (r+1)^2 \cdot \log(e+a^{-1})$, $0<a\le1$, satisfies
\begin{equation}\label{eq:repeated_height_terminal_distortion}
 \ell_r(d)\le\zeta_r \cdot \ell_{r+1}(a^{-1}),\qquad
 \zeta_r:=1+2\eta^{r/2},\ \forall r\ge r_Q.
\end{equation}
The products $\prod_{r=r_Q}^{R}(1+r^{-2})^2 \cdot \zeta_r^2$ are uniformly
bounded, independently of $Q$ and the finite terminal depth $R$.
\end{lemma}
\begin{proof}
We have $g(Mx)-g(x)\le\log M$. For $y\ge1$ and $v\ge0$,
$g(y+v)-g(y)\le F(v)$. Induction gives
Eq.~\eqref{eq:repeated_height_exact_distortion}.
Moreover $F(v)\le\eta v$ for all $v\ge0$, and $F(v)\le\log v$
for $v\ge2$. Until the value is at most two, the orbit therefore
falls at least as fast as ordinary logarithm; it then contracts by
$\eta$ at every step. This proves the stated bound after $b(M)$ steps,
including $1\le M\le e$ directly.

For $M:=2 \cdot C_T \cdot Q^2 \cdot (r+1)^2$, there is a fixed $c_T$ such that
\begin{equation*}
 b(M)\le c_T+\max\{t_Q,\log^*(r+1)\}.
\end{equation*}
Indeed, with $N:=\max\{16,2 \cdot C_T,Q,r+1\}$, we have
$M\le N^5\le e^N$. For $r\ge16$,
$\log^*(r+1)\le1+\log(r+1)\le r/4$.
Choose $C_b$ so that $r\ge r_Q$ implies
$r\ge\max\{16,4t_Q,4c_T\}$. Then $b(M)\le r/2$.
Use $x:=g(a^{-1})\ge1$ in the first assertion and $\ell_r(x)\ge1$
to obtain Eq.~\eqref{eq:repeated_height_terminal_distortion}.
Finally the logarithm of the full product is at most
\begin{equation*}
 \frac2{r_Q-1}+\frac{4\eta^{r_Q/2}}{1-\sqrt\eta},
\end{equation*}
which is universally bounded. The integer $r_*$ is chosen above the
thresholds for coefficient transfer and the initial profile.
\end{proof}

With the terminal distortion controlled after a uniform initial depth,
the block comparison advances any admissible initial profile through all
subsequent depths. We keep its amplitude explicit for later repetition.

\begin{proposition}[A profile from a coefficient seed and an initial bound]
\label{prop:second_logstar_profile}
Let $\mathcal G$ satisfy the seed assumptions of
Proposition~\ref{prop:second_logstar_blocks}, and suppose
Eqs.~\eqref{eq:second_logstar_seed_cutoff} and
\eqref{eq:second_logstar_coherent_root} hold for
$\mathfrak g:=\mathcal G^2$, with uniform constants.
For each odd $Q\ge3$, suppose a number $M_Q\ge\mathfrak g(Q)$ satisfies,
uniformly over centered regular laws of covariance at most $I$,
\begin{equation}\label{eq:repeated_height_profile_initial_input}
 D^2W\succeq aI\quad\Longrightarrow\quad
 C_P\le M_Q \cdot \ell_{r_Q}(a^{-1})^2,\ \forall a>0.
\end{equation}
Then there is a universal $C$ such that, for all integers $r\ge r_Q$,
\begin{equation}\label{eq:second_logstar_radius_profile}
 Z_Q\le C \cdot M_Q \cdot (r+1)^{1/Q} \cdot \ell_r(a^{-1})^2,
\end{equation}
and
\begin{equation}\label{eq:second_logstar_profile_bound}
 C_P\le C \cdot M_Q \cdot (r+1)^{1/Q} \cdot \ell_r(a^{-1})^2.
\end{equation}
The constant depends only on the fixed constants for the coefficient seeds and the analytic estimates,
not on the seed, $M_Q$, $Q$, $r$, or the dimension.
\end{proposition}
\begin{proof}
Use the assignment and constants from
Proposition~\ref{prop:second_logstar_blocks}. Set
\begin{equation}\label{eq:second_logstar_profile_amplitude}
 A_Q:=A \cdot M_Q,
\end{equation}
where $A$ is a universal constant to be chosen. The cutoff estimate at
$s=1$ and a fixed argument constant at least eight gives
$\mathfrak g(8Q)\le C \cdot \mathfrak g(Q)$. Therefore
\begin{equation*}
 H_Q\le C \cdot M_Q,\qquad
 a_Q+s_Q+1\le(C \cdot M_Q)^Q.
\end{equation*}
The polynomial propagation bound has $E=4$ and $\alpha=1$.
Suppose a uniform profile for $Z_Q$ is established at depth $r$, with
squared amplitude $\Gamma_{Q,r}^2\ge A_Q \cdot (r+1)^{1/Q}$.
Choose $r_*$ above the thresholds of
Proposition~\ref{prop:static_profile_coefficients} and
Proposition~\ref{prop:sharp_iterated_curvature}. Write
\begin{align}
 s&:=r+1,&\widehat\Gamma&:=(1+r^{-2}) \cdot \Gamma_{Q,r},&
 \delta&:=\frac1{64Qs},&p&:=\frac\delta{128},\notag\\
 b&:=(1-2p)^{-1},&M&:=2^{64} \cdot (1+C_F \cdot E^2),&
 k_0&:=\lceil100\delta^{-1} \cdot \log(M/\delta)\rceil.
 \label{eq:second_logstar_profile_parameters}
\end{align}
The static coefficient transfer and the fixed seed give
\begin{equation}\label{eq:second_logstar_profile_coefficients}
 c_k^2\le[\widehat\Gamma^2 \cdot \ell_r(k)^2]^{k-1},\qquad
 c_k^2\le\mathfrak g(k)^{k-1}.
\end{equation}
Let $d_0$ be the least dyadic integer at least $Q$, and let $d\ge d_0$
be dyadic. Assume for a contradiction that
\begin{equation}\label{eq:second_logstar_large_radius}
 Z_Q(\mu)>b^4 \cdot e^\delta \cdot \widehat\Gamma^2 \cdot \ell_r(d)^2.
\end{equation}
Choose $A$ large enough that $A_Q>2H_Q$. The assumption excludes the
floor, so $z:=Z_Q(\mu)$ is an actual radius and the starting family of
Proposition~\ref{prop:second_logstar_blocks} applies. With $u:=z^{-1}$,
\begin{equation}\label{eq:second_logstar_inverse_bounds}
 b^4 \cdot u \cdot \widehat\Gamma^2 \cdot \ell_r(d)^2<e^{-\delta},\qquad
 u\le A_Q^{-1} \cdot s^{-1/Q},\qquad u^Q\le A_Q^{-Q} \cdot s^{-1}.
\end{equation}
Its energy and retained losses obey
\begin{equation}\label{eq:second_logstar_profile_startup}
 \nu\le\frac{u}{1-a_Q \cdot u^Q},\qquad
 P_{2d_0-1}\le s_Q \cdot u^Q,\qquad
 D_0:=\max\{\nu-u,0\}\le2a_Q \cdot u^{Q+1},
\end{equation}
provided $a_Q \cdot u^Q\le1/2$. The retained range is valid because
$2d_0-1<4Q$. Also $C_P\le z+3$ gives $\lambda\ge u/2$ when $u\le1/3$.

We will ensure $a_Q \cdot u^Q\le p$. Before a first exit of mass $p$, all
required normalizers are then at most
\begin{equation}\label{eq:second_logstar_profile_normalizers}
 B_*:=\frac\nu{1-p}\le\frac{u}{(1-p)^2}\le bu,
 \qquad \max\{1,B_* \cdot z\}\le b.
\end{equation}
Lemma~\ref{lem:second_logstar_polynomial_loss} consequently applies with
\begin{align}
 \theta_d&:=5184 \cdot C_F \cdot E^2 \cdot \sum_{\substack{k<d\\k\ \mathrm{dyadic}}}
             k^8 \cdot (b^4 \cdot u)^k \cdot c_k^2,\notag\\
 \Delta_d&:=s_Q \cdot u^Q+2 \cdot C_F \cdot \sum_{\substack{d_0\le k<d\\k\ \mathrm{dyadic}}}
             k^2 \cdot (bu)^k \cdot c_k^2,\notag\\
 \tau_d&:=C_F \cdot d \cdot (bu)^d \cdot c_d^2.
 \label{eq:second_logstar_profile_kernel}
\end{align}

Since $k_0\le C \cdot (Qs)^2$, the assumed cutoff estimate
shows that $\mathfrak g(k_0)\le C \cdot \mathfrak g(Q) \cdot s^{1/Q}$.
Increasing $A$ once ensures $\Gamma_{Q,r}^2\ge4\mathfrak g(k_0)$.
For $k\le k_0$, Eqs.~\eqref{eq:second_logstar_profile_coefficients}
and~\eqref{eq:second_logstar_inverse_bounds} therefore give
\begin{align*}
 (b^4 \cdot u)^k \cdot c_k^2&\le b^4 \cdot u \cdot 4^{-(k-1)},\\
 (bu)^k \cdot c_k^2&\le b^Q \cdot u^Q \cdot \mathfrak g(k)^{Q-1} \cdot 4^{-(k-Q)}
 ,\ \forall k\ge Q.
\end{align*}
The low-degree contribution to $\theta_d$ is at most a fixed constant
times $u$. The low-degree coherent contribution is at most $D_Q \cdot u^Q$,
where $D_Q^{1/Q}\le C \cdot \mathfrak g(Q)$ by
Eq.~\eqref{eq:second_logstar_coherent_root}; here $b\le2$.

For $k_0<k<d$, the profile bound gives
\begin{equation*}
 (b^4 \cdot u)^k \cdot c_k^2\le b^4 \cdot u \cdot e^{-\delta(k-1)},\qquad
 (bu)^k \cdot c_k^2\le e^{-\delta k}.
\end{equation*}
For $0<\delta\le1$,
\begin{align*}
 \sum_{k>k_0}k^8 \cdot e^{-\delta k}
 &\le2^{26} \cdot \delta^{-9} \cdot e^{-\delta k_0/2},\\
 \sum_{k>k_0}k^2 \cdot e^{-\delta k}
 &\le2^{10} \cdot \delta^{-3} \cdot e^{-\delta k_0/2},\qquad
 e^{-\delta k_0/2}\le(\delta/M)^{50}.
\end{align*}
For the first two bounds, take out $e^{-\delta k_0/2}$ and bound the
remaining degree-$m$ moment by the integral of
$(x+1)^m e^{-\delta x/2}$, which is at most $2m! \cdot (2/\delta)^{m+1}$.
The fixed $M$ in Eq.~\eqref{eq:second_logstar_profile_parameters} makes
the high-degree contributions at most $u$ and $p/512$, respectively.
Thus there is one universal $C_\theta$ such that
\begin{equation}\label{eq:second_logstar_profile_split}
 \theta_d\le C_\theta \cdot u,\qquad
 \Delta_d\le(s_Q+D_Q) \cdot u^Q+p/512.
\end{equation}

Choose $A$ universally so that, for every odd $Q\ge3$,
\begin{equation}\label{eq:second_logstar_amplitude_conditions}
 A_Q\ge\max\{2H_Q,16,16 \cdot C_\theta,C \cdot \mathfrak g(Q)\},\qquad
 A_Q^Q\ge2^{24} \cdot Q \cdot (a_Q+s_Q+D_Q+4 \cdot C_\theta \cdot a_Q+1).
\end{equation}
These inequalities are compatible with
Eq.~\eqref{eq:second_logstar_profile_amplitude}, since
$Q^{1/Q}$ is bounded and $M_Q\ge\mathfrak g(Q)$.
They imply
\begin{equation*}
 \nu\le1/4,\quad D_0\le u,\quad
 \frac{\theta_d \cdot u}{\lambda}\le1/8,\quad
 \Delta_d+\frac{\theta_d \cdot D_0}{\lambda}\le p/64,\quad
 u+D_0/p\le2u.
\end{equation*}
Indeed $(\theta_d/\lambda) \cdot D_0\le4 \cdot C_\theta \cdot a_Q \cdot u^{Q+1}$, and
$D_0/p\le16384Q \cdot a_Q \cdot A_Q^{-Q} \cdot u\le u$.
The same conditions give $a_Q \cdot u^Q\le p$, as required before exit.
The actual matched budget and Lemma~\ref{lem:matched_scale_comparison}
now yield
\begin{equation}\label{eq:second_logstar_terminal_comparison}
 a<12(u+D_0/p) \cdot \tau_d\le24 \cdot C_F \cdot d \cdot (bu)^d \cdot c_d^2.
\end{equation}

For $0<a\le1$, choose a dyadic $d$ with
\begin{equation}\label{eq:second_logstar_terminal_degree}
 C_T \cdot Q^2 \cdot s^2 \cdot \log(e+a^{-1})\le d
 <2 \cdot C_T \cdot Q^2 \cdot s^2 \cdot \log(e+a^{-1}),
\end{equation}
where $C_T$ is a sufficiently large universal constant. Then $d\ge d_0$
and $\log(24 \cdot C_F \cdot d/a)\le\delta d$. Indeed the left side is at most
$\log(48 \cdot C_F \cdot C_T)+2\log Q+2\log s+2\log(e+a^{-1})$, whereas
$\delta d\ge(C_T/64)Qs \cdot \log(e+a^{-1})$.
We define $L:=\widehat\Gamma^2 \cdot \ell_r(d)^2\ge1$.
The strict assumption gives $b^4 \cdot uL<e^{-\delta}$, while
$c_d^2\le L^{d-1}$. Equation~\eqref{eq:second_logstar_terminal_comparison}
would imply
\begin{equation*}
 a<24 \cdot C_F \cdot d \cdot L^{-1} \cdot (buL)^d\le24 \cdot C_F \cdot d \cdot e^{-\delta d}\le a,
\end{equation*}
a contradiction. Therefore
\begin{equation}\label{eq:second_logstar_profile_step}
 Z_Q(\mu)\le b^4 \cdot e^\delta \cdot \widehat\Gamma^2 \cdot \ell_r(d)^2.
\end{equation}
When $a\ge1$, the curvature bound gives $C_P\le1$, so $Z_Q=H_Q$;
the chosen amplitude already covers this case.

We close the depth induction using
Lemma~\ref{lem:repeated_height_distortion}. The scalar increment is unchanged:
\begin{equation*}
 \log(b^4 \cdot e^\delta)\le16p+\delta
 =\frac9{512Q(r+1)}\le\frac1Q \cdot \log(1+1/(r+1)).
\end{equation*}
We define
\begin{equation}\label{eq:second_logstar_distortion_product}
 J_{Q,r_Q}:=1,\qquad
 J_{Q,r+1}:=J_{Q,r} \cdot (1+r^{-2})^2 \cdot \zeta_r^2,\qquad
 \Gamma_{Q,r}^2:=A_Q \cdot (r+1)^{1/Q} \cdot J_{Q,r}.
\end{equation}
All $J_{Q,r}$ have one universal upper bound. Equations
\eqref{eq:second_logstar_profile_step} and
\eqref{eq:repeated_height_terminal_distortion} consequently give the
claimed depth step with these amplitudes.

To initialize at $r_Q$, the static comparison and
Eq.~\eqref{eq:repeated_height_profile_initial_input} give
\begin{equation*}
 Z_Q\le\max\{H_Q,C_P\}\le C \cdot M_Q \cdot \ell_{r_Q}(a^{-1})^2.
\end{equation*}
Increasing the same universal $A$ covers this at
$\Gamma_{Q,r_Q}^2=A_Q \cdot (r_Q+1)^{1/Q}$. Finite induction proves
Eq.~\eqref{eq:second_logstar_radius_profile}. Applying $C_P\le Z_Q+3$
once and using $M_Q\ge1$ and $\ell_r\ge1$ proves
Eq.~\eqref{eq:second_logstar_profile_bound}.
\end{proof}

The curvature bound from Section~\ref{sec:iterated_curvature} supplies the
initial profile for the original seed. This gives the first round
of height reduction.

\begin{proposition}[The original seed at the sharper depth]
\label{prop:repeated_height_original_seed_profile}
For $\mathcal G=G_*$ there is a universal $C$ such that
\begin{equation}\label{eq:repeated_height_original_seed_profile}
 C_P\le C \cdot t_Q^{1/3} \cdot (r+1)^{1/Q} \cdot \ell_r(a^{-1})^2
 ,\ \forall Q\ge3\text{ odd},\ r\ge r_Q.
\end{equation}
\end{proposition}
\begin{proof}
The seed satisfies $\mathfrak g(Q)\le C \cdot t_Q^{1/3}$.
At $r_Q=O(t_Q)$, the curvature profile in
Proposition~\ref{prop:sharp_iterated_curvature} has squared amplitude at most
$C \cdot (r_Q+1)^{1/3}\le C' \cdot t_Q^{1/3}$. Thus the hypothesis on the initial curvature bound
holds with $M_Q:=C' \cdot t_Q^{1/3}$, after increasing $C'$ to cover the seed.
Apply Proposition~\ref{prop:second_logstar_profile}.
\end{proof}

\subsection{Proof of repeated height reduction}
\label{sec:repeated_height_outer_stages}

We now close the outer induction in
Proposition~\ref{prop:repeated_height_outer_profiles}: choosing the depth
and block order separately for each degree gives improved coefficients,
which in turn give the next full curvature profile. Each round introduces
another log-star with one universal multiplicative cost and the same
starting depths. We also show that the scalar height iteration terminates,
identifying the accumulated cost that the next section must remove.

\begin{lemma}[Termination of the height iteration]
\label{lem:repeated_height_termination}
Set $\mathfrak t_0(x):=x$. For every $x\ge1$, the integer
\begin{equation}\label{eq:repeated_height_count}
 \kappa(x):=\min\{j\ge0:\mathfrak t_j(x)\le3\}
\end{equation}
is finite. The value three is fixed by $\mathfrak t$, and
$\mathfrak t(x)<x$ for every $x>3$. The function $\kappa$ is
unbounded.
\end{lemma}
\begin{proof}
We have $\mathfrak t(3)=3$ and $\mathfrak t(x)=3$ for $3<x<4$.
For $x\ge4$, the log-star recursion and $\log^*v\le v$ give
\begin{equation*}
 \mathfrak t(x)\le2+\log(x+2)<x.
\end{equation*}
After one composition the values are integers, so they strictly
decrease until they reach three, which is then fixed. This proves
finiteness. Every finite composition $\mathfrak t_j(x)$ tends to
infinity as $x$ tends to infinity. Therefore no fixed number of
compositions brings every argument to three, and $\kappa$ is unbounded.
\end{proof}

The coefficient seed used in the next reduction comes from the full
curvature profile of the preceding reduction. Capping it by the original
seed keeps the constants in the block construction uniform. The
remaining cost is one fixed factor for each reduction.

\begin{proof}[Proof of Proposition~\ref{prop:repeated_height_outer_profiles}]
The case $j=1$ of Eq.~\eqref{eq:repeated_height_outer_profile} is
Proposition~\ref{prop:repeated_height_original_seed_profile}. Increase
$A_1$ once so that $A_1\ge1$.

We first extract coefficients from any established stage, keeping its
constant explicit. Suppose that
Eq.~\eqref{eq:repeated_height_outer_profile} holds at stage $j$.
For $d\ge2$, we define $h:=\mathfrak t(d)$ and choose
$r:=\lceil C_r \cdot h\rceil$, where $C_r$ is one fixed sufficiently large
constant. The argument in Lemma~\ref{lem:iterated_coefficient_seed}
gives $\ell_r(d)\le5$, and $r$ exceeds the threshold in
Proposition~\ref{prop:constant_log_specialization}.
Let $k$ be the least odd integer at least
$\max\{3,\log(r+1)\}$. Above a fixed threshold for $h$,
\begin{equation*}
 k\le h,\qquad r_k\le r,\qquad (r+1)^{1/k}\le e.
\end{equation*}
Indeed $k\le\log(C_r \cdot h+2)+3$, whereas
$r_k\le C \cdot \mathfrak t(k)\le C \cdot k$; choosing $C_r$ and then the fixed
threshold proves both comparisons. Neither choice depends on $j$.
The established profile at order $k$ therefore implies, uniformly over
its class,
\begin{equation*}
 C_P\le e \cdot A_j \cdot \mathfrak t_j(k)^{1/3} \cdot \ell_r(a^{-1})^2
 \le e \cdot A_j \cdot \mathfrak t_{j+1}(d)^{1/3} \cdot \ell_r(a^{-1})^2.
\end{equation*}
Proposition~\ref{prop:constant_log_specialization} now yields
Eq.~\eqref{eq:repeated_height_outer_coefficients}, with a universal
$C_e$. When $h$ is below the fixed threshold, the original seed
$G_*(d)$ is bounded by a fixed constant. Increasing $C_e$ handles
these degrees because $A_j\ge1$ and
$\mathfrak t_{j+1}(d)\ge3$. The case $d=1$ is the covariance bound.
Thus the coefficient estimate is uniform even in the exceptional
range of logarithm heights.
The specialization first gives the estimate for isotropic laws.
Whitening and contraction of each coefficient tensor extend it to
covariance at most $I$, as in
Lemma~\ref{lem:iterated_coefficient_seed}; a singular covariance is
treated on its supporting subspace.

Assume the profile has been proved through stage $j-1$, where $j\ge2$.
We define the next majorant by
\begin{equation}\label{eq:repeated_height_outer_capped_seed}
 \mathcal G_j(x):=
 \min\{G_*(x),C_e \cdot \sqrt{A_{j-1}} \cdot \mathfrak t_j(x)^{1/6}\},
 \qquad \mathfrak g_j(x):=\mathcal G_j(x)^2.
\end{equation}
The coefficient extraction just proved and the original seed both
bound the same coefficients, so their minimum does too.
Lemma~\ref{lem:second_logstar_seed_estimates} supplies every seed
comparison, cutoff estimate, and coherent-root estimate required by
Propositions~\ref{prop:second_logstar_blocks} and
\ref{prop:second_logstar_profile}. Their constants are independent of
$j$ and $A_{j-1}$. In particular, the cap keeps
$\mathcal G_j(k)^6\le C_0 \cdot k$ with the original universal $C_0$.
This is why enlarging the amplitude at an earlier stage does not
change a later block or kernel constant.

It remains to provide the initial curvature bound at the common
$r_Q$. Let $k$ be the least odd integer at least
$\max\{3,\log(r_Q+1)\}$. For all $Q$ above one fixed threshold,
\begin{equation*}
 k\le Q,\qquad r_k\le r_Q,\qquad
 k\le C \cdot \mathfrak t(Q),\qquad (r_Q+1)^{1/k}\le e.
\end{equation*}
Here $r_Q\le C \cdot \mathfrak t(Q)$, and the monotonicity of the common
depth gives $r_k\le r_Q$. The fixed-power argument bound in
Eq.~\eqref{eq:repeated_height_power_argument} gives, uniformly in $j$,
\begin{equation*}
 \mathfrak t_{j-1}(k)
 \le\mathfrak t_{j-1}(C \cdot \mathfrak t(Q))
 \le\mathfrak t_j(Q)+b_{C,1}
 \le C' \cdot \mathfrak t_j(Q).
\end{equation*}
Apply the preceding stage at order $k$ and depth $r_Q$ to obtain
\begin{equation}\label{eq:repeated_height_outer_initialization}
 C_P\le C'' \cdot A_{j-1} \cdot \mathfrak t_j(Q)^{1/3} \cdot
                  \ell_{r_Q}(a^{-1})^2.
\end{equation}
For the fixed remaining range of $Q$, the earlier profile in
Proposition~\ref{prop:sharp_iterated_curvature} at $r_Q$ has a fixed bounded
amplitude. The same inequality follows after a universal enlargement
of $C''$, since $A_{j-1}\ge1$ and $\mathfrak t_j(Q)\ge3$.
No threshold in this initialization depends on $j$.

Take
\begin{equation*}
 M_Q:=C_M \cdot A_{j-1} \cdot \mathfrak t_j(Q)^{1/3},
\end{equation*}
where the universal $C_M$ covers Eq.~\eqref{eq:repeated_height_outer_initialization}
and $C_e^2$. Then $M_Q\ge\mathfrak g_j(Q)$, and the hypothesis on the initial
curvature bound of Proposition~\ref{prop:second_logstar_profile} holds.
That proposition gives Eq.~\eqref{eq:repeated_height_outer_profile}
with $A_j:=C_* \cdot A_{j-1}$ for one universal $C_*$. Its constant depends
only on the fixed estimates for the coefficient seeds already verified, so the same
$C_*$ works for every stage. The coefficient extraction above then
proves Eq.~\eqref{eq:repeated_height_outer_coefficients} at stage $j$.
Induction proves both assertions for every finite $j$.
\end{proof}

\Needspace{22\baselineskip}
\section{Controlling Repetition Costs through Geometric Decay}
\label{sec:repetition_costs}
\label{sec:controlling_repetition_costs}
\label{sec:cutoff_blocks}
\label{sec:small_loss_refinement}

The repeated height reduction of the preceding section still incurs a
fixed multiplicative cost at each round. We remove this accumulation by
choosing losses that decrease geometrically and controlling the resulting
degree cutoffs and starting depths. Their total cost stays bounded,
yielding a Poincar\'e inequality with a constant independent of the
dimension and the measure.

\begin{theorem}[Universal Poincar\'e and KLS bounds, formal version of Theorem~\ref{thm:main_informal}]
\label{tg:thm:main}
There is a universal constant $C>0$ such that, for every $n\ge1$ and every isotropic log-concave probability measure $\mu$ on $\mathbb R^n$, each locally Lipschitz function $f$ with $\int|\nabla f|^2\d\mu<\infty$ belongs to $L^2(\mu)$ and satisfies
\begin{equation*}
 \Var_\mu[f]\le C \cdot \int|\nabla f|^2\d\mu.
\end{equation*}
In particular, after increasing the universal constant if necessary,
\begin{equation}\label{eq:constant_chain_dimension_bound}
 C_P(\mu)\le C,\qquad \psi_n\le C.
\end{equation}
The constant is independent of the dimension, the measure, and the test function.
\end{theorem}

It suffices to bound the common coefficient radius $\mathcal A$ from
Section~\ref{sec:coefficient_radius_static}. We construct curvature
estimates whose multiplicative and additive losses remain bounded
through every finite repetition. The precise inductive estimate is
Proposition~\ref{prop:constant_chain_profiles}; a finite choice of
the number of repetitions and the depth then gives the universal
bound. Regular approximation completes the proof of
Theorem~\ref{tg:thm:main} at the end of the section.

We use the common coefficient radius $\mathcal A$, the operator
$\mathcal T$, and the class of centered regular measures of covariance
at most $I$ and strictly positive curvature. The nondecreasing majorant
$\mathcal G$ satisfies Eq.~\eqref{eq:second_logstar_majorant}.
Fix also a nondecreasing coefficient majorant $\mathcal H$ such that,
uniformly over the same laws,
\begin{equation}\label{eq:secondary_geometry_cap}
 \begin{gathered}
 1\le\mathcal H(x)\le G_*(x),\ \forall x\ge1,\qquad
 \mathcal G\le\mathcal H,\\
 c_k\le\mathcal H(k)^{k-1},\ \forall k\ge1.
 \end{gathered}
\end{equation}
Any valid additional majorant may first be capped by $G_*$, and the
current $\mathcal G$ may be replaced by its minimum with $\mathcal H$.
The choice $\mathcal H:=G_*$ gives the cutoff floor associated with $G_*$.
All constants below are universal and independent of $\mathcal H$.
Fix sufficiently large numerical
constants $C_{\mathrm{cut}},C_{\mathrm{deg}}$, and define
\begin{equation}\label{eq:slow_offset_cutoff}
 K_\epsilon:=\lceil C_{\mathrm{cut}} \cdot \epsilon^{-2}\rceil
 ,\ \forall 0<\epsilon\le1/16.
\end{equation}
The coefficient seed below this cutoff will control the small degrees;
the factor $e^{-\epsilon k}$ controls the remaining fixed moments.

\subsection{Constructing blocks above a coefficient cutoff}
\label{sec:slow_offset_geometry}

A small multiplicative loss is useful only if the block constants remain
uniform as that loss decreases. We split the coefficient moments at a
cutoff, use those bounds to control a norm-attaining orbit, and initialize
and extend its family. The resulting block floor records exactly where
the smaller loss is available.

\begin{lemma}[Moments with a fixed-ratio initial range]
\label{lem:slow_offset_seed_moments}
Fix nonnegative integers $M,r$ and $c>0$. Suppose $0\le\rho_k\le1$
for all $k$, and $\rho_k\le C_0^{-1}$ for $k\le16K_\epsilon$.
For sufficiently large fixed $C_{\mathrm{cut}},C_0$,
\begin{equation}\label{eq:slow_offset_seed_moments}
 \sum_{k\ge r+1}k^M \cdot \rho_k^{k-r} \cdot e^{-c\epsilon(k-r)}\le C.
\end{equation}
The bound is uniform in $\epsilon$ and the sequence. Its left side
tends uniformly to zero when these two fixed thresholds increase.
\end{lemma}
\begin{proof}
Below $K_\epsilon$ the sum is bounded by
$\sum_{k\ge r+1}k^M \cdot C_0^{-(k-r)}$. Choose the cutoff constant so
that $K_\epsilon\ge2r$. Above $K_\epsilon$, discard $\rho_k\le1$
and use $k-r\ge k/2$. This part is at most
\begin{equation*}
 C \cdot \epsilon^{-M-1} \cdot e^{-c\epsilon K_\epsilon/4}
 \le C \cdot \epsilon^{-M-1} \cdot e^{-cC_{\mathrm{cut}}/(4\epsilon)}.
\end{equation*}
Both bounds have the stated uniform smallness. We fix constants large
enough for the finitely many moments used below.
\end{proof}

The cutoff moments control the convolution terms along an actual orbit.
The next estimate converts them into uniform orbit and window bounds
under the stated polynomial propagation envelope.

\begingroup
\emergencystretch=2em
\begin{lemma}[Control of a norm-attaining orbit]
\label{lem:coefficient_radius_green}
There are universal constants $C_0,c_0,C>0$ with the following property.
Let $0<\epsilon\le1/16$, let $d\ge2$ be dyadic, and suppose
$F^2\ge\max\{C_0 \cdot \mathcal H(16K_\epsilon)^2,C_0\}$,
$\mathcal G(d)\le F$, $m\ge2d$ is an integer, and
\begin{equation}\label{eq:coefficient_radius_green_hypotheses}
 \begin{gathered}
 z=\|\mathcal T^m\|_{\mathrm{op}}^{2/m},\qquad
 0\le R-z\le2,\qquad
 z\ge(1+\epsilon) \cdot F^2,\\
 \|\mathcal T^h\|_{\mathrm{op}}\le4(h+1)z^{h/2},\ \forall h\ge0,\qquad
 \gamma:=C_F \cdot d \cdot c_d^2 \cdot z^{-(d-1)},\qquad
 \frac{\gamma m^2}{z}\le c_0.
 \end{gathered}
\end{equation}
Let $f$ be a centered unit top right singular vector of $\mathcal T^m$,
and define
\begin{equation*}
 w_j:=z^{-j/2} \cdot \mathcal T^jf,\qquad b_j:=\|w_j\|_2^2,\qquad
 \ell_j:=\|Lw_j\|^2,\qquad
 D_j:=\|\nabla[(B-zI)w_j]\|_2^2,\ \forall j\ge1.
\end{equation*}
Then
\begin{equation}\label{eq:coefficient_radius_orbit_bounds}
 \frac12\le b_j\le2,\ \forall 0\le j\le m,\qquad
 b_j\le2,\ \forall j\ge0.
\end{equation}
Moreover, for every integer $1\le J\le m/4$,
\begin{equation}\label{eq:coefficient_radius_window_loss}
 \sum_{k=2}^{J+1}\sum_{i=k}^{m+k-1}\ell_i
 \le\frac{C}z+C \cdot \gamma Jm.
\end{equation}
\end{lemma}
\begin{proof}
Compactness supplies the singular vector, as in
Proposition~\ref{prop:block_averaged_start}. We define $t:=R/z$. Maximality
and the one-step operator norm give
\begin{equation}\label{eq:coefficient_radius_orbit_shift}
 b_0=b_m=1,\qquad b_{j+m}\le b_j,\qquad
 b_{j+1}\le tb_j,\qquad b_j\le t^j,
 \qquad t\le1+2/z.
\end{equation}
Write $S:=\sum_{j=0}^{m-1}b_j$. In particular,
$\sum_{j=0}^{2m-1}b_j\le2S$.

Write $E:=4$. A swap at degree $k$ is propagated through at most
$3k-2$ normalized copies of $\mathcal T$. The stated power bound gives
a squared propagation factor at most $E^2(3k)^2$. The proof of
Lemma~\ref{lem:quarter_order_raw_joint} therefore gives, for dyadic
$e\ge2$ and $j\ge2e-1$,
\begin{align}
 \ell_j&\le C_F \cdot e \cdot c_e^2 \cdot z^{-(e-1)} \cdot b_{j-e+1}\notag\\
 &\quad+972C_F \cdot E^2 \cdot \sum_{\substack{k<e\\k\ \mathrm{dyadic}}}
 k^7 \cdot c_k^2 \cdot z^{-k} \cdot \sum_{h=1}^{3k-1}D_{j-k+1-h}.
 \label{eq:coefficient_radius_local_raw_frame}
\end{align}
Indeed the exact projection is
$\mathsf P_kQ_1w_j=Q_kw_{j-k+1}/(k! \cdot z^{(k-1)/2})$.
The $h$th preceding adjacent swap uses $D_{j-k+1-h}$ and propagates through
$h-1$ copies of $\mathcal T/\sqrt z$. Applying the polynomial power
bound inside the same sorting sum proves the display. Every defect
index is at least one; no gradient symmetry of $w_0$ is needed.

For $j\ge3$, use the largest available dyadic degree at most $d$, and
denote its coherent term by $h_j$. Set $h_2:=2K_2 \cdot b_1/z$ and $D_i:=0$
for $i\le0$. The second-moment estimate
Eq.~\eqref{eq:quarter_order_second_moment} gives
\begin{equation}\label{eq:coefficient_radius_local_convolution}
 \ell_j\le h_j+\sum_{r\ge1}\alpha_r \cdot D_{j-r},\ \forall j\ge2,
\end{equation}
where the nonnegative finite kernel puts $8/z$ at delay one and,
for every dyadic $k<d$, puts $972C_F \cdot E^2 \cdot k^7 \cdot c_k^2 \cdot z^{-k}$ at each
delay $r=k,\ldots,4k-2$. In particular
\begin{align}
 \eta:=\sum_r\alpha_r
 &\le\frac8z+2916C_F \cdot E^2 \cdot \sum_{\substack{k<d\\k\ \mathrm{dyadic}}}
                    k^8 \cdot c_k^2 \cdot z^{-k}\le\frac Cz,\notag\\
 \eta_1:=\sum_r r\alpha_r
 &\le\frac8z+11664C_F \cdot E^2 \cdot \sum_{\substack{k<d\\k\ \mathrm{dyadic}}}
                    k^9 \cdot c_k^2 \cdot z^{-k}\le\frac Cz.
 \label{eq:coefficient_radius_local_kernel}
\end{align}
For every $k$ define $\rho_k:=\min\{1,\mathcal H(k)^2/F^2\}$.
Then $\rho_k\le C_0^{-1}$ for $k\le16K_\epsilon$ and, at every
used degree $k\le d$,
$\mathcal G(k)^2/z\le\rho_k \cdot e^{-\epsilon/2}$.
Lemma~\ref{lem:slow_offset_seed_moments} bounds both sums.
The degree-one term is fixed.

Before any uniform orbit bound is known, the stated propagation gives
$b_j\le E^2(j+1)^2$. Degree $k<d$ is used only at
$j=2k-1,\ldots,4k-2$, where $b_{j-k+1}\le9E^2 \cdot k^2$.
Its contribution to $z$ times the weighted early coherent sum is at most
\begin{equation*}
 C \cdot E^2 \cdot k^5 \cdot \mathcal G(k)^2
       (\mathcal G(k)^2/z)^{k-2}
 \le C \cdot k^6 \cdot \rho_k^{k-2} \cdot e^{-\epsilon(k-2)/2}.
\end{equation*}
Here $\mathcal G(k)^2\le Ck$ follows from the original cap.
Degree two and $h_2$ are bounded directly, using $c_2^2\le2$ and
$b_j\le E^2(j+1)^2$ at their finitely many indices.
The same moment lemma consequently gives
\begin{equation}\label{eq:coefficient_radius_early_terms}
 \sum_{j=2}^{2d-2}h_j\le C/z,\qquad
 \sum_{j=2}^{2d-2}j h_j\le C/z.
\end{equation}
At every later index $h_j=\gamma b_{j-d+1}$. The absolute
threshold in the lemma makes both kernel moments small enough for the
following absorptions, with constants independent of $\epsilon$.

The remainder is the Green calculation in
Lemma~\ref{lem:second_logstar_green_orbit}. We verify its required
estimates. After its kernel and early coherent estimates, that proof
uses only $\eta,\eta_1\le C/z$, the two early sums just proved,
$D_j\le z(b_{j+1}-2b_j+b_{j-1})+\ell_j$, the norm-attaining shift,
$R-z\le2$, and $\gamma m^2/z\le c_0$. In particular the boundary of its second block uses $b_{2m}\le b_0=1$ and
\begin{equation*}
 z(b_{2m+1}-b_{2m})\le(R-z)b_{2m}\le2.
\end{equation*}
Also $b_1\le R/z\le2$. All remaining absorptions require only a fixed
large lower bound on $z$. Starting from the assumed polynomial envelope,
these absorptions yield the uniform forward orbit bound.

Explicitly, with $D_{\mathrm{int}}:=\sum_{j=1}^mD_j$ and
$D_{\mathrm{tot}}:=\sum_{j=1}^{2m}D_j$, the two telescopes in
Eqs.~\eqref{eq:second_logstar_first_telescope}--\eqref{eq:second_logstar_two_block_sums}
and the same absorptions give
\begin{equation*}
 D_{\mathrm{int}},D_{\mathrm{tot}}\le C+C \cdot \gamma S,\qquad
 \sum_{j=2}^{2m}\ell_j\le C/z+C \cdot \gamma S.
\end{equation*}
The upper Dirichlet weight $g_j$ and the sliding-window weight $n_m$
from that proof are nonnegative and $1$-Lipschitz. Their weighted losses
satisfy $E_g\le C+C \cdot \gamma mS$ and $E_m\le C/z+C \cdot \gamma mS$.
Consequently the same upper and shifted lower Green formulas give
\begin{equation*}
 b_j\le1+C/z+C \cdot \gamma mS/z,\qquad
 b_j\ge1-C/z^2-C \cdot \gamma mS/z,\ \forall 1\le j<m.
\end{equation*}
Since $S\le m\max_{j\le m}b_j$, a fixed large $C_0$ and small $c_0$
first give $b_j\le2$, and then $b_j\ge1/2$ on the first block.
The shift gives the upper bound at every forward index. Repeating the
same weight calculation with $J\le m/4$ windows gives
Eq.~\eqref{eq:coefficient_radius_window_loss}; its terminal coherent
term is at most $2\gamma Jm$. Every source index remains within the
two blocks used above. All form pairings and finite Bochner domains are
those in the proof of Lemma~\ref{lem:second_logstar_green_orbit}. This establishes the claimed orbit and
window bounds under the stated polynomial propagation hypothesis.
\end{proof}
\endgroup

The orbit bounds now supply one starting family. Keeping the dyadic
degree free allows the subsequent extension to choose its degree according
to the order and the loss parameter.

\begin{lemma}[Initialization with a chosen dyadic degree]
\label{lem:coefficient_radius_initialization}
\label{lem:slow_offset_initialization}
There is a universal constant $C>0$ with the following property.
Suppose $p\ge3$ is odd, $\kappa\ge1$, and the hypotheses of
Lemma~\ref{lem:coefficient_radius_green} hold for a dyadic $d\ge2$.
Assume also
\begin{equation}\label{eq:coefficient_radius_initialization_hypotheses}
 \frac{z^{p-2}}{4\kappa}\le m\le\frac{4z^{p-2}}\kappa,
 \qquad1\le J_0\le m/4,\qquad
 J_0 \cdot \gamma zm\le1,\qquad C \cdot \kappa z^{-p}\le1/8.
\end{equation}
There is one centered unit finite starting family satisfying, for $u:=1/z$,
\begin{equation}\label{eq:coefficient_radius_initialization}
 \nu\le\frac{u}{1-C \cdot \kappa u^p},\qquad
 P_J\le C \cdot \kappa u^p,\ \forall 1\le J\le J_0,\qquad
 aV_N+X_N\le\nu-u+uP_{N+1},\ \forall N\ge0.
\end{equation}
\end{lemma}
\begin{proof}
Use the actual orbit in Lemma~\ref{lem:coefficient_radius_green} and
exactly the family in Lemma~\ref{lem:second_logstar_initialization}:
$\mathcal W:=\bigoplus_{j=0}^{m-1}w_j$, followed by its compensated
$F_2$ restart. Its proof uses the chosen degree only through the orbit
bounds, $m\ge2d$, and the displayed coherent conditions; it does not
require $d$ to be a particular function of $p$.
For completeness, $S_0=S_1\ge m/2$ makes
$U:=\mathcal T\mathcal W/\sqrt{zS_0}$ a centered unit gradient of
energy at most $1/z$. Quadratic testing of the first summand and the
moment sum over two blocks give
$\|\operatorname{Sym}(LU)\|^2\le C/(zm)+C \cdot \gamma$.
The exact skew credit and the actual restarted loss fractions therefore
give
\begin{equation*}
 \nu\le\frac{u}{1-C/(z^2 \cdot m)-C \cdot \gamma/z},\qquad
 P_J\le C/(z^2 \cdot m)+C \cdot \gamma J/z.
\end{equation*}
Every denominator window has $S_{k+1}\ge m/4$ for $k\le J_0+1$,
using the same interior indices as in that proof. The hypotheses bound
both errors by $C \cdot \kappa u^p$ and keep the normalizers positive.
Finally the norm-attaining shift gives $S_{k+1}\le S_k$ for every
$k$, hence $v_{j+1}\le ze_j$. Combining this exact inequality with
$aV_N+X_N\le\nu-e_N$ gives the asserted budget. Every output slot
is included, and the construction uses the absorbing convention at zero
of Lemma~\ref{lem:second_logstar_initialization}.
\end{proof}

The degree choice changes when the order crosses the cutoff. We treat
the two ranges in one extension, with the normalizer threshold matched
to the estimate available in each range.

For odd $p\ge3$, choose dyadic degrees and a retained prefix by
\begin{equation}\label{eq:slow_offset_degrees}
 \begin{gathered}
 d_p:=\min\{2^j:2^j\ge4p\},\quad
 e_p:=\min\{2^j:2^j\ge8p\}\quad\text{if }p\le K_\epsilon,\\
 d_p:=e_p:=\min\{2^j:2^j\ge C_{\mathrm{deg}} \cdot p/\epsilon\}
                  \quad\text{if }p>K_\epsilon,\qquad J_p:=4e_p.
 \end{gathered}
\end{equation}
In particular $J_p\ge32p$.

\begin{lemma}[Extension with two normalizer thresholds]
\label{lem:slow_offset_extension}
Let $q\ge3$ be odd. Suppose
$F^2\ge\max\{C_0 \cdot \mathcal H(16K_\epsilon)^2,C_0\}$ at a sufficiently large
universal $C_0$, and define
\begin{equation*}
 \kappa_q:=(F^2/4)^{q-1},\qquad
 \kappa_{q-2}:=(F^2/4)^{q-3},\qquad u:=1/z,\qquad
 \omega_q:=\begin{cases}1,&q\le K_\epsilon,\\
                       \epsilon,&q>K_\epsilon.\end{cases}
\end{equation*}
Assume $\mathcal G(e_q)\le F$, $0\le R-z\le2$, and
$\|\mathcal T^h\|_{\mathrm{op}}\le4(h+1)z^{h/2}$ for all $h\ge0$.
Suppose one centered unit family has energy $\nu$ and actual losses
through $J_q$ satisfying, with constants $a_q,s_q$,
\begin{equation*}
 \nu\le\frac{u}{1-a_q \cdot u^q},\qquad
 P_J\le s_q \cdot u^q,\ \forall 1\le J\le J_q,\qquad
 aV_N+X_N\le\nu-u+uP_{N+1},\ \forall N\ge0.
\end{equation*}
There is a universal $C_\theta$ such that, on defining
\begin{equation*}
 \mathsf K_q:=4(s_q+C_F \cdot \kappa_{q-2}+4C_\theta \cdot a_q+1),\qquad
 C_q:=a_q+4\mathsf K_q,
\end{equation*}
the conditions
\begin{equation}\label{eq:slow_offset_extension_thresholds}
 z\ge\max\{(1+\epsilon) \cdot F^2,8,8C_\theta\},\qquad
 a_q \cdot u^q,\mathsf K_q \cdot u^q\le\omega_q/1024
\end{equation}
imply
\begin{equation}\label{eq:slow_offset_extension}
 \|\mathcal T^m\|_{\mathrm{op}}^{2/m}
 \ge z-C_q \cdot z^{1-q}
 ,\ \forall 1\le m\le\lfloor z^q/\kappa_q\rfloor.
\end{equation}
Moreover $C_q\le C \cdot (a_q+s_q+\kappa_{q-2}+1)$ with a universal constant.
\end{lemma}
\begin{proof}
Up to the first exit of $p_*:=\mathsf K_q \cdot u^q$, the actual normalizers obey
\begin{equation*}
 B_*\le\frac{u}{(1-a_q \cdot u^q)(1-p_*)}
       \le\beta u,\qquad\beta:=1+\omega_q/128.
\end{equation*}
Thus $t_*:=\max\{1,B_* \cdot z\}\le\beta$, and the absolute gap gives
$\lambda\ge u/2$. Apply Lemma~\ref{lem:second_logstar_polynomial_loss}
with $d_0:=d:=e_q$. Its terminal term is
$\tau:=C_F \cdot d \cdot B_*^d \cdot c_d^2$, and its loss kernel is at most
\begin{equation*}
 Cu\sum_{\substack{k<e_q\\k\ \mathrm{dyadic}}}
              k^8(\beta^4 \cdot \mathcal G(k)^2/z)^{k-1}.
\end{equation*}
If $q\le K_\epsilon$, every tested degree is below $16K_\epsilon$.
Consequently $\mathcal G(k)^2/F^2\le1/C_0$, and the fixed value
$\beta=1+1/128$ gives a geometric ratio at most $1/4$ after increasing
$C_0$. The kernel is at most $C_\theta \cdot u$.
If $q>K_\epsilon$, use $\beta=1+\epsilon/128$ and the sequence
$\rho_k:=\min\{1,\mathcal H(k)^2/F^2\}$. At tested degrees,
$\beta^4 \cdot \mathcal G(k)^2/z\le\rho_k \cdot e^{-\epsilon/2}$.
The moment lemma gives the same uniform kernel bound.

We define $d:=e_q$, $D:=1+\epsilon$, $g:=\mathcal G(d)^2$, and $\rho:=g/F^2$.
The exact terminal ratio can be written without a positive power of $F$:
\begin{align}\label{eq:slow_offset_terminal_ratio}
 \frac{\tau}{\kappa_q \cdot \kappa_{q-2} \cdot u^{2q}}
 &\le C_F \cdot d \cdot F^6 \cdot 4^{2q-4} \cdot \rho^{d-1} \cdot \beta^d \cdot D^{2q-d}\notag\\
 &=C_F \cdot d \cdot g^3 \cdot 4^{2q-4} \cdot \rho^{d-4} \cdot \beta^d \cdot D^{2q-d}.
\end{align}
The original cap $\mathcal G\le G_*$ gives $g\le Cd$. In the low
branch $8q\le d<16q$,
$\rho\le C_0^{-1}$, and $D^{2q-d}\le1$. The right side is at most
\begin{equation*}
 C \cdot q^4 \cdot C_0^4 \cdot (16\beta^{16}/C_0^8)^q,
\end{equation*}
which is uniformly small for $q\ge3$ at a fixed large $C_0$.
In the high branch $d\ge C_{\mathrm{deg}} \cdot q/\epsilon$,
$\rho\le1$, and the $D$ damping absorbs $\beta^d \cdot 4^{2q}$ for a
fixed sufficiently large $C_{\mathrm{deg}}$. Since
$\epsilon^{-1}<\sqrt{q/C_{\mathrm{cut}}}\le\sqrt q$,
the remaining bound is
\begin{equation*}
 C \cdot (q/\epsilon)^4 \cdot e^{-cC_{\mathrm{deg}}q}
       \le C \cdot q^6 \cdot e^{-cC_{\mathrm{deg}}q}.
\end{equation*}
All constants are fixed, and $q>K_\epsilon\ge256C_{\mathrm{cut}}$.
Thus both branches give
$\tau\le C_F \cdot \kappa_q \cdot \kappa_{q-2} \cdot u^{2q}$.

The prefix contains $2d-1$. On
$N\le\lfloor z^q/\kappa_q\rfloor$, the delayed loss estimate and
actual matched budget give
\begin{equation*}
 (1-2C_\theta \cdot u)P_N
 \le s_q \cdot u^q+C_F \cdot \kappa_{q-2} \cdot u^q+4C_\theta \cdot a_q \cdot u^{q+1}
 \le\mathsf K_q \cdot u^q/4.
\end{equation*}
This excludes first exit, including its justified delayed prefix.
Successor masses are positive because $p_*\le1/1024$ and $\nu<1/4$.
The raw quotient identities in Lemma~\ref{lem:nearunit_block_extension}
then give
\begin{equation*}
 \sum_{j<m}\frac{\|LY_j\|^2}{\|Y_j\|^2}
 \le\frac{P_m}{\lambda(1-p_*)}\le4\mathsf K_q \cdot u^{q-1},\qquad
 \nu^{-1}\ge z-a_q \cdot u^{q-1}.
\end{equation*}
Every successive squared norm ratio is at least
$z-(a_q+4\mathsf K_q)u^{q-1}$. Multiplication proves the assertion
on precisely the stated finite horizon. The low branch needs only its
fixed small exit for these norm identities and positivity; its stronger
coefficient ratio already controlled the kernel.
\end{proof}

We assemble the orbit, initialization, and extension estimates into a
block of the prescribed order. The proposition records the coefficient
floor and the uniform bounds needed to keep the loss per repetition small.

\begin{proposition}[Blocks above the cutoff floor]
\label{prop:slow_offset_blocks}
There are universal constants $C_X,C_0,C>0$ with the following property.
Let $\mathcal G,\mathcal H$ satisfy
Eqs.~\eqref{eq:second_logstar_majorant} and
\eqref{eq:secondary_geometry_cap} uniformly over the stated regular laws.
For $0<\epsilon\le1/16$ and
odd $Q\ge3$, we define
\begin{equation}\label{eq:slow_offset_geometry_parameters}
 \begin{gathered}
 X:=C_X \cdot Q^2,\qquad
 F^2:=\max\{\mathcal G(X^X)^2,C_0 \cdot \mathcal H(16K_\epsilon)^2,C_0\},\\
 H_Q:=(1+2\epsilon) \cdot F^2.
 \end{gathered}
\end{equation}
There is an assignment $Z_Q$ and constants $a_Q,s_Q>0$ satisfying
\begin{align}
 H_Q&\le Z_Q\le\max\{H_Q,R\},\qquad\mathcal A\le Z_Q,\notag\\
 c_d&\le Z_Q^{(d-1)/2},\ \forall d\ge1,\qquad
 a_Q+s_Q+1\le C \cdot (F^2/4)^{Q-3}\le F^{2Q}.
 \label{eq:slow_offset_geometry_static}
\end{align}
If $Z_Q>H_Q$, let $z:=Z_Q$ and $u:=1/z$. One centered unit finite
starting family satisfies
\begin{equation}\label{eq:slow_offset_geometry_startup}
 \nu\le\frac{u}{1-a_Q \cdot u^Q},\qquad
 P_J\le s_Q \cdot u^Q,\ \forall 1\le J\le J_Q,\qquad
 aV_N+X_N\le\nu-u+uP_{N+1},\ \forall N\ge0,
\end{equation}
where $J_Q$ is defined in Eq.~\eqref{eq:slow_offset_degrees}, and
\begin{equation}\label{eq:slow_offset_geometry_propagation}
 \|\mathcal T^h\|_{\mathrm{op}}\le4(h+1)z^{h/2},\ \forall h\ge0,\qquad
 0\le R-z\le\min\{1,C/R\},\qquad a_Q \cdot u^Q\le1/8.
\end{equation}
All constants are independent of the order, both majorants, and $\epsilon$.
\end{proposition}
\begin{proof}
Construct first in the branch $\mathcal A>H_Q$, so $R\ge\mathcal A$.
Use $\kappa_q:=(F^2/4)^{q-1}$ and Eq.~\eqref{eq:slow_offset_degrees}.
In the high branch $p>K_\epsilon$ implies $\epsilon^{-1}<\sqrt p$,
so all degrees are at most $C \cdot p^{3/2}\le C \cdot p^2$. One fixed $C_X$ puts
every tested degree below $X$, with its seed at most $F$.

For odd $p\ge5$, suppose $z\ge DF^2$, $D:=1+\epsilon$, and $m$
is within a factor four of $z^{p-2}/\kappa_{p-2}$. Write $d:=d_p$,
$g:=\mathcal G(d)^2$ and $\rho:=g/F^2$. The two coherent conditions for initialization are bounded by
\begin{align*}
 \gamma m^2/z
 &\le Cd \cdot F^2 \cdot \rho^{d-1} \cdot 4^{2p-6} \cdot D^{2p-d-4}
   =C \cdot d \cdot g \cdot \rho^{d-2} \cdot 4^{2p-6} \cdot D^{2p-d-4},\\
 J_p \cdot \gamma zm
 &\le C \cdot J_p \cdot d \cdot F^4 \cdot \rho^{d-1} \cdot 4^{p-3} \cdot D^{p-d}
   =C \cdot J_p \cdot d \cdot g^2 \cdot \rho^{d-3} \cdot 4^{p-3} \cdot D^{p-d}.
\end{align*}
The original cap $\mathcal G\le G_*$ gives $g\le Cd$. In the low
branch, $d\ge4p$,
$J_p\le64p$, $\rho\le C_0^{-1}$, and both $D$ exponents are
negative. The two bounds are at most
\begin{equation*}
 C \cdot p^2 \cdot C_0^2 \cdot (16/C_0^4)^p,\qquad
 C \cdot p^4 \cdot C_0^3 \cdot (4/C_0^4)^p.
\end{equation*}
For $p\ge5$ these are uniformly small at fixed large $C_0$.
In the high branch, the $D$ decay absorbs the powers of four, and
$\rho\le1$. The bounds are at most
\begin{equation*}
 C \cdot (p/\epsilon)^2 \cdot e^{-cC_{\mathrm{deg}}p}\le C \cdot p^3 \cdot e^{-cC_{\mathrm{deg}}p},
 \qquad
 C \cdot (p/\epsilon)^4 \cdot e^{-cC_{\mathrm{deg}}p}\le C \cdot p^6 \cdot e^{-cC_{\mathrm{deg}}p}.
\end{equation*}
The chosen constants make them small uniformly for $p>K_\epsilon$.
The length lower bound is $m\ge cF^2 \cdot 4^{p-3} \cdot D^{p-2}$. Thus
$J_p/m\le Cp \cdot 4^{-p}/F^2$ in the low branch, and
$J_p/m\le C \cdot p^{3/2} \cdot 4^{-p}/F^2$ in the high branch. These give
$m\ge2d_p$ and $J_p\le m/4$. Also
$\kappa_{p-2} \cdot z^{-p}\le F^{-6} \cdot 4^{3-p} \cdot D^{-p}$ is uniformly small.
With the polynomial envelope established below, the initialization lemma
therefore gives $a_p+s_p+1\le C \cdot \kappa_{p-2}$ through its full prefix.

At the base $m_3:=\lfloor R\rfloor$, the low degree rule always applies:
$d_3=16,e_3=32,J_3=128$, and $\kappa_1=1$.
The mesoscopic estimate gives $R-z_3\le C_*/R$ and global remainder
factor at most two. Its two coherent terms are $O(z_3^{-14})$ and
$O(J_3 \cdot z_3^{-13})$. The lower radius is justified below, after which
Lemmas~\ref{lem:coefficient_radius_green} and
\ref{lem:coefficient_radius_initialization} apply at this fixed base.

For each constructed order define
\begin{equation*}
 m_{q+2}:=\lfloor z_q^q/\kappa_q\rfloor,\qquad
 z_{q+2}:=\|\mathcal T^{m_{q+2}}\|_{\mathrm{op}}^{2/m_{q+2}},\qquad
 \Delta_q:=C_b \cdot \kappa_{q-2} \cdot z_q^{1-q}.
\end{equation*}
The exit thresholds follow from
$(a_q+\mathsf K_q)u^q\le C \cdot F^{-6} \cdot 4^{3-q} \cdot D^{-q}$.
This is uniformly small in the low branch. In the high branch, multiplying
by $\epsilon^{-1}<\sqrt q$ still leaves a uniformly small quantity.
Lemma~\ref{lem:slow_offset_extension} therefore supplies every bridge.
In addition,
\begin{equation*}
 \Delta_q/z_q\le C_b \cdot F^{-6} \cdot 4^{3-q} \cdot D^{-q},\qquad
 q\Delta_q/z_q\le C \cdot F^{-6}.
\end{equation*}

Each proposed length is at most $R^Q$. If $R\le X$, its first $m$
coefficients have seed $\mathcal G(m)^2\le F^2<\mathcal A$.
If $X<R<2\mathcal A$, then $Q<R$ and the original cap gives
$\mathcal G(m)^2\le G_*(R^R)^2\le R/2<\mathcal A$ after a fixed
enlargement of $C_X$. The exact block maximum in
Eq.~\eqref{eq:coefficient_radius_block_maximum} forces $z_m\ge\mathcal A$
in both cases. If $R\ge2\mathcal A$, the mesoscopic gap gives this
at the base. At each new step the old decrement tail and current bridge
give $z_{q+2}\ge R-C_*/R-2\Delta_3\ge\mathcal A$, since
$\Delta_3=C_b/z_3^2=O(R^{-2})$. Only old ratios enter that tail.
Thus every new actual radius is at least $\mathcal A>H_Q$ before
its new ratio is evaluated.

The exact successive ratio, including increasing actual radii, is
\begin{equation}\label{eq:coefficient_radius_gap_ratio}
 \frac{\Delta_{q+2}}{\Delta_q}
 =\frac{F^4}{16z_{q+2}^2}(z_q/z_{q+2})^{q-1}
 \le\frac2{16}<\frac12.
\end{equation}
For a downward step use $z_{q+2}\ge z_q-\Delta_q$ and
$q\Delta_q/z_q\le C \cdot F^{-6}$; for an upward step the second factor
is at most one. Hence all finite constructed prefixes satisfy
\begin{equation}\label{eq:coefficient_radius_gap}
 0\le R-z_q\le C_*/R+C'/R^2\le C_g/R\le1.
\end{equation}
This provides the absolute gap required by the Green lemma. The lower
radius comes from the common coefficient radius, so no reserve of size
$\epsilon F^2$ is used to pay this gap. Fix the degree and moment constants, then the constants for initialization, extension,
and the total gap, and finally $C_0,C_X$ above all required
fixed thresholds.

The defining real lengths are at least $F^2 \cdot 4^{q-1} \cdot D^q$. Their integer
floor losses are uniformly small. For $q\ge5$, the exact lower length
ratio has factor $16z_q^2/F^4\ge16D^2$, while the correction for a decrease in radius is $e^{-O(F^{-6})}$. Thus $m_{q+2}/m_q\ge4$.
At the distinct first step the ratio is at least
$(1-o(1))16(z_3/F^2)^2(z_3/R)\ge4$, using the mesoscopic relative
gap $O(R^{-2})$.

In the greedy decomposition into exact block lengths, each nonterminal
level has logarithmic excess at most
\begin{equation*}
 C \cdot m_{q+2}\Delta_q/z_q
 \le C \cdot \kappa_{q-2}/\kappa_q=16C/F^4.
\end{equation*}
The base remainder costs at most $C/R$ by
Eq.~\eqref{eq:coefficient_radius_gap}. Counting levels by
$C \cdot \log(h+1)/\log4$ gives the uniform envelope
$\|\mathcal T^h\|_{\mathrm{op}}\le4(h+1)z_q^{h/2}$.
Counting instead their actual finite number, at most $q/2$, gives
$\log E_q\le C/R+Cq/F^4$. At
$m=m_{q+2}\ge cF^2 \cdot 4^{q-1} \cdot D^q$, this implies
\begin{equation*}
 q(\log(z_{q+2}/z_q))_+
 \le\frac{Cq \cdot 4^{-q}}{F^2 \cdot R}+\frac{C \cdot q^2 \cdot 4^{-q}}{F^6}
 \le\frac C{F^2 \cdot R}+\frac C{F^6}.
\end{equation*}
Together with the downward gap and floor loss, it proves the required
comparison within a factor of four of $m_{q+2}$ with $z_{q+2}^q/\kappa_q$.
The coherent, prefix, and normalization conditions were checked above.
At this point both the new actual radius and its polynomial envelope are
known. Apply Lemma~\ref{lem:coefficient_radius_green} and then
Lemma~\ref{lem:coefficient_radius_initialization} to its own norm-attaining
orbit. This closes the finite induction with one actual family and its
matched budget. Thus the polynomial envelope is established before
the Green estimate is applied.

At the final order $a_Q+s_Q+1\le C \cdot (F^2/4)^{Q-3}\le F^{2Q}$,
including $Q=3$. In the nontrivial branch set $Z_Q:=z_Q\ge\mathcal A$.
If instead $\mathcal A\le H_Q$, set $Z_Q:=H_Q$ and use the same
upper bounds for $a_Q,s_Q$. The coefficient inequalities follow from
the definition of $\mathcal A$. The assertions about an actual family
and a comparison of $R$ with $H_Q$ apply only to the nontrivial branch. This proves the proposition.
\end{proof}

\subsection{Repeated height reduction with a small multiplicative loss}
\label{sec:slow_offset_profiles}

We use the preceding block estimates to repeat height reduction with a
factor $e^{24\epsilon}$ per round instead of the fixed factor $C_*$.
Keeping the dependence on $\epsilon$ explicit will let us choose the
loss according to the number of rounds, while tracking the required
coefficient cutoff and starting depth.

The curvature profiles below concern centered regular
laws of covariance at most $I$ and curvature at least $aI$, where $a>0$.
We use the common coefficient radius $\mathcal A$ from
Eq.~\eqref{eq:coefficient_radius_definition} and the seed $G_*$ from
Eq.~\eqref{eq:quarter_order_seed}. We use the iterated logarithms
$\ell_r:=g^{\circ r}$, where $g(x):=\log(e+x)$, and the height functions
$\mathfrak t(x):=1+\log^*(x+2)$ and $\mathfrak t_j:=\mathfrak t^{\circ j}$.
Proposition~\ref{prop:repeated_height_original_seed_profile}
provides the initial profile for $G_*$. The block estimates in
Proposition~\ref{prop:slow_offset_blocks} reduce the cost of each repetition.

We first normalize the iterated logarithm at its fixed point. This lets
the coefficient transfer preserve a multiplicative factor close to one
throughout the repetition.

We apply the block estimates of Proposition~\ref{prop:slow_offset_blocks}
with $\mathcal H:=G_*$.
Its constants are universal, and its parameters are
\begin{equation}\label{eq:slow_offset_geometry_input}
 X:=C_X \cdot Q^2,\qquad \Xi_\epsilon:=16K_\epsilon,\qquad
 F^2:=\max\{\mathcal G(X^X)^2,C_0 \cdot G_*(\Xi_\epsilon)^2,C_0\},
 \qquad H_Q:=(1+2\epsilon) \cdot F^2.
\end{equation}
When $\mathcal A>H_Q$, we use its actual radius $z\ge\mathcal A$,
startup bound $(a_Q+s_Q+1)^{1/Q}\le F^2$, actual prefix through at
least $4Q$, matched budget, propagation envelope $4(h+1)z^{h/2}$,
and absolute gap $R_{\mathrm{res}}-z\le2$.
Here $K_\epsilon$ is defined in Eq.~\eqref{eq:slow_offset_cutoff}.
These estimates concern the radius of the constructed block.

Let $\varrho$ be the unique positive solution of
$\varrho=\log(e+\varrho)$, and define
\begin{equation}\label{eq:nearunit_normalized_log}
 \mathcal L_r(x):=\ell_r(x)/\varrho,\ \forall x\ge0.
\end{equation}
We first record the precise normalization allowed in the static
coefficient transfer.

\begin{lemma}[Normalized static transfer]
\label{lem:nearunit_static_transfer}
The conclusion of Proposition~\ref{prop:static_profile_coefficients}
remains valid if its hypothesis $\Gamma\ge\sqrt H$ is replaced by
\begin{equation}\label{eq:nearunit_static_floor}
 \Gamma>0,\qquad \Gamma \cdot \ell_r(0)\ge\sqrt H.
\end{equation}
\end{lemma}
\begin{proof}
With $\Psi(x):=\Gamma^2 \cdot \ell_r(x)^2$, the stated assumption gives
$\Psi(x)\ge H$ for every $x>0$. The posterior Appell estimate and
regular approximation in that proof therefore remain valid. In the
simultaneous degree induction of
Proposition~\ref{prop:block_profile_coefficients}, the auxiliary
radius is $A:=(1+r^{-2})\Gamma \cdot \ell_r(d)\ge\sqrt H\ge1$.
This is the lower bound used for its weights and for the second term
of Eq.~\eqref{eq:block_finite_transfer}. The estimates for $J/A$ and
the logarithmic elasticity are independent of $\Gamma$. Thus every
step of the proof applies under Eq.~\eqref{eq:nearunit_static_floor}.

We will apply this conclusion to $\mathcal A$. Its defining coefficient
inequalities and $\mathcal A\ge1$ give the static hypothesis with $H:=1$.
The proof passes only coefficient bounds at finitely many degrees through each
approximation, so it requires no continuity of the supremum defining
$\mathcal A$.
\end{proof}

The extracted coefficient bound must preserve the degree at which the
preceding profile was tested. We choose the order and depth together,
then estimate the resulting majorant at every degree used by the block
construction.

We define
\begin{align}
 U_\epsilon&:=1+\mathfrak t(\epsilon^{-1}),\qquad
 S_\epsilon:=C_S \cdot U_\epsilon,\qquad
 W_{j,\epsilon}(x):=[\mathfrak t_j(x)+S_\epsilon]^{1/3},
                  \label{eq:slow_offset_height}\\
 Q_\epsilon&:=\text{the least odd integer at least }C_Q \cdot \epsilon^{-4},\notag\\
 r_{Q,\epsilon}&:=\max\{r_*,\lceil C_b \cdot \mathfrak t(Q)\rceil
                                      +\lceil C_R/\epsilon\rceil\}.
                  \label{eq:slow_offset_depth}
\end{align}
Here $0<\epsilon\le1/16$, and $C_b$ is from
Lemma~\ref{lem:repeated_height_distortion}. Choose $r_*$ above the
thresholds for coefficient transfer in
Proposition~\ref{prop:static_profile_coefficients}. Write $\eta:=(e+1)^{-1}$. Since
$g(x)=\log(e+x)$ fixes $\varrho$ and is $\eta$-Lipschitz on
$[1,\infty)$, we have $\mathcal L_r(0)\ge1-C \cdot \eta^{r-1}$.
Also $\ell_m(x)\le5$ implies
$\mathcal L_r(x)\le1+C \cdot \eta^{r-m}$ for $r\ge m$.
Choose $C_R$ sufficiently large that
\begin{equation}\label{eq:slow_offset_depth_losses}
 \mathcal L_r(0)^2\ge e^{-\epsilon},\qquad
 \prod_{v=r_{Q,\epsilon}}^r(1+v^{-2})^2\zeta_v^2\le e^\epsilon,
 \qquad \zeta_v:=1+2(e+1)^{-v/2}.
\end{equation}
Indeed the logarithm of the product is at most
$2/(r_{Q,\epsilon}-1)+4\eta^{r_{Q,\epsilon}/2}/(1-\sqrt\eta)$.
For every integer $d\ge3$, we also have $\mathcal L_r(d)\ge1$
because $d>\varrho$. The depths are integers and nondecreasing in $Q$.
All constants below are fixed before the finite repetition count.

\begin{lemma}[Retaining the extracted degree]
\label{lem:slow_offset_extracted_seed}
Suppose that, for some $j\ge1$ and $A\ge A_0$, uniformly over the
centered regular laws of covariance at most $I$,
\begin{equation}\label{eq:slow_offset_profile_input}
 \mathcal A\le A \cdot W_{j,\epsilon}(Q) \cdot (r+1)^{1/Q}
                              \mathcal L_r(a^{-1})^2
 ,\ \forall Q\ge3\text{ odd},\ r\ge r_{Q,\epsilon}.
\end{equation}
There is a nondecreasing majorant $\mathcal G$ satisfying
Eq.~\eqref{eq:second_logstar_majorant} for which, on writing
$\mathfrak g:=\mathcal G^2$ and $M_Q:=A \cdot W_{j+1,\epsilon}(Q)$,
\begin{align}
 \mathfrak g(X^X)&\le e^{5\epsilon} \cdot M_Q,
                          \label{eq:slow_offset_seed_lookahead}\\
 \sup_{s\ge1}\frac{\mathfrak g(C_k \cdot (Qs)^2)}{s^{1/Q}}
       &\le e^{5\epsilon} \cdot M_Q.
                          \label{eq:slow_offset_seed_cutoff}
\end{align}
The second assertion holds for each fixed universal $C_k$ needed in
the proof, after one universal choice of $C_S$.
Moreover, if $Q\ge Q_\epsilon$, then
\begin{equation}\label{eq:slow_offset_coherent_root}
 E_{Q,\epsilon}:=2C_F \cdot \sum_{k\ge Q}k^2 \cdot \mathfrak g(k)^{Q-1}
                         e^{-\epsilon(k-Q)}
 \quad\text{satisfies}\quad
 E_{Q,\epsilon}^{1/Q}\le e^\epsilon \cdot \mathfrak g(N),
 \qquad N:=\lceil D\epsilon^{-2} \cdot Q^2\rceil\le X^X,
\end{equation}
where $D$ is a universal constant.
\end{lemma}
\begin{proof}
For $x\ge1$, define
\begin{equation}\label{eq:slow_offset_extraction_order}
 \begin{split}
 r_\epsilon(x)&:=\max\{r_*,\lceil C_r \cdot \mathfrak t(x)\rceil
                                      +\lceil D_R/\epsilon\rceil\},\\
 k_\epsilon(x)&:=\text{the least odd integer at least }
       \max\{3,\epsilon^{-1} \cdot \log(r_\epsilon(x)+1)\}.
 \end{split}
\end{equation}
These are nondecreasing functions. Choose $D_R\ge2C_R+2$ and enlarge
$r_*$ and $C_r$ by universal amounts. Then
$r_\epsilon(x)\ge r_{k_\epsilon(x),\epsilon}$ for every $x$.
To see this, we write $r:=r_\epsilon(x)$. Since $r\ge D_R/\epsilon$,
$k_\epsilon(x)\le r\log(r+1)/D_R+3$. The elementary estimate
$\mathfrak t(k)\le C \cdot \log(e+k)$ gives
$C_b \cdot \mathfrak t(k_\epsilon(x))+2\le r/2$ above a fixed threshold.
The remaining $\lceil C_R/\epsilon\rceil$ is at most $r/2$ after
increasing $D_R$ if necessary. The maximum with $r_*$ causes no issue.

Lemma~\ref{lem:iterated_coefficient_seed} brings $\ell_m(x)$ below
five for some $m\le\mathfrak t(x)$. Contraction toward $\varrho$,
together with the additional $D_R/\epsilon$ depths, gives
\begin{equation*}
 (r_\epsilon(x)+1)^{1/k_\epsilon(x)}\le e^\epsilon,\qquad
 r_\epsilon(x)^{-2}\le\epsilon,\qquad
 \mathcal L_{r_\epsilon(x)}(x)\le1+\epsilon.
\end{equation*}
Apply Eq.~\eqref{eq:slow_offset_profile_input} at this order and depth,
and then Lemma~\ref{lem:nearunit_static_transfer} with the common
assignment $\mathcal A$ and floor one. Its floor condition follows
from Eq.~\eqref{eq:slow_offset_depth_losses} and a fixed $A_0$.
The three factors just displayed give, for every $d\ge2$,
\begin{equation}\label{eq:slow_offset_actual_degree_seed}
 c_d^2\le[e^{5\epsilon} \cdot A
                 W_{j,\epsilon}(k_\epsilon(d))]^{d-1}.
\end{equation}
This is a simultaneous coefficient estimate on the full law class.
Degree one follows from covariance normalization.

The fixed first-stage coefficient estimate in
Eq.~\eqref{eq:repeated_height_outer_coefficients} also gives
$c_d^2\le[C_2 \cdot \mathfrak t_2(d)^{1/3}]^{d-1}$ for a universal
$C_2\ge1$. Hence we may take
\begin{equation}\label{eq:slow_offset_seed_definition}
 \mathfrak g(x):=\min\{G_*(x)^2,C_2 \cdot \mathfrak t_2(x)^{1/3},
                    e^{5\epsilon} \cdot A \cdot W_{j,\epsilon}(k_\epsilon(x))\}.
\end{equation}
Every branch is nondecreasing and at least one. Thus its square root
is a valid monotone coefficient majorant, also across the degree ranges
used next.

Write $h:=\mathfrak t(Q)$. Fixed-power and self-power bounds give
$\mathfrak t(X^X)\le h+B_1$, with universal $B_1$.
For $Y:=C_k \cdot (Qs)^2$ and $y:=\log(s)/Q$, the product and fixed-argument height estimates give
\begin{equation}\label{eq:slow_offset_cutoff_height}
 \mathfrak t(Y)\le h+B_2+y.
\end{equation}
Indeed $\log(Y+2)\le KQ(1+y)$; apply the product and fixed-argument
bounds to the remaining iterated logarithms. Choose $B_h$ above both
$B_1$ and $B_2$.

Choose a universal $C_h$ large enough that, whenever
\begin{equation}\label{eq:slow_offset_height_threshold}
 h\ge h_\epsilon:=C_h \cdot \epsilon^{-1} \cdot \log(e/\epsilon),
\end{equation}
we have $k_\epsilon(x)\le h/2$ for all $x$ with
$\mathfrak t(x)\le h+B_h+1$. This follows from
\begin{equation*}
 k_\epsilon(x)\le3+\epsilon^{-1}
            \log[C \cdot (\mathfrak t(x)+\epsilon^{-1})]
\end{equation*}
at $h=h_\epsilon$, and then from the decrease of $\log(h)/h$.
In particular $k_\epsilon(X^X)\le h$, which proves
Eq.~\eqref{eq:slow_offset_seed_lookahead} in this range.
For the cutoff, $0\le y\le1$ likewise gives $k_\epsilon(Y)\le h$.
For $y\ge1$, the same logarithmic bound, allowing an increase of at most two when rounding to an odd integer, gives
\begin{equation*}
 k_\epsilon(Y)\le h/2+2+
          \epsilon^{-1} \cdot \log(1+y/(h+B_h+\epsilon^{-1}))
 \le h/2+2+y/(\epsilon h)\le h+y\le h(1+y).
\end{equation*}
The fixed choice of $C_h$ ensures all these inequalities.
Consequently the uniform product estimate gives
$\mathfrak t_j(k_\epsilon(Y))\le\mathfrak t_j(h)+5y$ for $y\ge1$,
using $\mathfrak t_j(1+y)\le y+4$. For
$D:=\mathfrak t_j(h)+S_\epsilon\ge3$,
\begin{equation*}
 [(D+5y)/D]^{1/3} \cdot e^{-y}\le1,
\end{equation*}
since its logarithm is at most $5y/(3D)-y\le0$.
On $0\le y\le1$ there is no height increase. The last branch in
Eq.~\eqref{eq:slow_offset_seed_definition} therefore proves
Eq.~\eqref{eq:slow_offset_seed_cutoff} at every high height.

If $h<h_\epsilon$, use its fixed $\mathfrak t_2$ branch instead.
By Eq.~\eqref{eq:slow_offset_cutoff_height} and the product bound,
\begin{equation*}
 \mathfrak t_2(Y)\le\mathfrak t(h+B_h+1)+\mathfrak t(1+y)
                  \le\mathfrak t(h+B_h+1)+y+4.
\end{equation*}
Taking its cube root and multiplying by $e^{-y}$ shows that the
cutoff supremum is at most $C[\mathfrak t(h+B_h+1)+4]^{1/3}$.
The same bound controls $\mathfrak g(X^X)$.
Since $h+B_h+1\le C \cdot \epsilon^{-2}$, a fixed-power height estimate gives
\begin{equation*}
 [\mathfrak t(h+B_h+1)+4]^{1/3}\le C \cdot U_\epsilon^{1/3}.
\end{equation*}
One universal choice of $C_S,A_0$ makes these bounds at most
$M_Q\ge A_0 \cdot C_S^{1/3} \cdot U_\epsilon^{1/3}$. No fixed factor is added to
$A$. This proves both estimates in the remaining range. Keeping the exact value of
$k_\epsilon$ is essential here: a fixed additive shift in
$\mathfrak t_{j+1}(Q)$ could not be paid by $\epsilon S_\epsilon$.

For the coherent sum, use monotonicity, $\mathfrak g\ge1$, and
the cap $\mathfrak g(k)\le C \cdot k^{1/3}$. The terms through $N$ are at most
$C \cdot \epsilon^{-3}(Q+1)^2 \cdot \mathfrak g(N)^{Q-1}$. For $k>N$, choose
$D$ so $C \cdot Q\log(e+k)\le\epsilon k/4$; the remaining sum is at most
$C \cdot \epsilon^{-1} \cdot e^{\epsilon Q-3\epsilon N/4}\le1$.
Since $\log[C \cdot \epsilon^{-3}(Q+1)^2+1]\le\epsilon Q$ for
$Q\ge C_Q \cdot \epsilon^{-4}$, taking roots gives
Eq.~\eqref{eq:slow_offset_coherent_root}. In this range
$N\le C' \cdot Q^{5/2}\le(C_X \cdot Q^2)^{C_X \cdot Q^2}=X^X$ after one fixed choice
of $C_X$. This proves the lemma.
\end{proof}

The coefficient majorant and cutoff blocks now close a full profile
induction. Its cost per repetition is small, but the required depth and
the additive offset depend on the loss parameter. The following estimate records both dependencies explicitly.

\begin{proposition}[A slowly growing offset]
\label{prop:slow_offset_outer_profiles}
There is a universal $A_0\ge1$ such that, with $A_{j,\epsilon}:=A_0 \cdot e^{24\epsilon j}$,
\begin{equation}\label{eq:slow_offset_outer_profile}
 \mathcal A\le A_{j,\epsilon} \cdot W_{j,\epsilon}(Q) \cdot (r+1)^{1/Q}
                                  \mathcal L_r(a^{-1})^2
\end{equation}
for every $j\ge1$, odd $Q\ge3$, $0<\epsilon\le1/16$, and
$r\ge r_{Q,\epsilon}$, uniformly over centered regular measures of
covariance at most $I$ and curvature at least $aI$.
Every constant is independent of these parameters and the dimension.
\end{proposition}
\begin{proof}
Assume that the preceding full stage holds with amplitude $A\ge A_0$.
Use the majorant in Lemma~\ref{lem:slow_offset_extracted_seed}, and set
$M_Q:=A \cdot W_{j+1,\epsilon}(Q)$ and $L_Q:=e^{20\epsilon} \cdot M_Q$.
We first initialize its next stage at $r=r_{Q,\epsilon}$.
When $h:=\mathfrak t(Q)\ge h_\epsilon$, take the least odd
$k\ge\max\{3,\epsilon^{-1} \cdot \log(r+1)\}$. The same threshold choice
as above gives $k\le h\le Q$, $r\ge r_{k,\epsilon}$ and
$(r+1)^{1/k}\le e^\epsilon$. Thus
\begin{equation}\label{eq:slow_offset_initial_bound}
 \mathcal A\le e^\epsilon \cdot M_Q \cdot \mathcal L_r(a^{-1})^2.
\end{equation}
For $h<h_\epsilon$, use the profile already proved for the initial coefficient seed
in Proposition~\ref{prop:repeated_height_original_seed_profile} at
$k_0$, the least odd integer at least $\max\{3,\log(r+1)\}$.
Increasing $r_*$ ensures that the initial depth required by that proposition,
$\max\{r_*,\lceil C_b \cdot \mathfrak t(k_0)\rceil\}$, is at most $r$;
indeed $\mathfrak t(k_0)=O(\log\log(r+1)+1)$.
The factor $(r+1)^{1/k_0}$ is at most $e$.
In this range $r\le C \cdot \epsilon^{-1} \cdot \log(e/\epsilon)$, so
$k_0\le C \cdot \epsilon^{-1}$ and
$\mathfrak t(k_0)\le\mathfrak t(\epsilon^{-1})+C$.
Using $\mathcal A\le\max\{1,C_P\}$, its normalized amplitude is
therefore at most $C \cdot U_\epsilon^{1/3}$, covered by a fixed
$A_0 \cdot C_S^{1/3} \cdot U_\epsilon^{1/3}$. This proves
Eq.~\eqref{eq:slow_offset_initial_bound} at every height without a
fixed multiplicative cost in the previous amplitude.

For the moment suppose $Q\ge Q_\epsilon$. The seed lemma and
$G_*(\Xi_\epsilon)^2\le C \cdot U_\epsilon^{1/3}$ give
\begin{equation}\label{eq:slow_offset_geometry_costs}
 F^2\le e^{6\epsilon} \cdot M_Q,\qquad H_Q\le e^{9\epsilon} \cdot M_Q,
 \qquad(a_Q+s_Q+1)^{1/Q}\le e^{6\epsilon} \cdot M_Q.
\end{equation}
For the second member of the maximum in
Eq.~\eqref{eq:slow_offset_geometry_input}, use
$\mathfrak t(\Xi_\epsilon)\le\mathfrak t(\epsilon^{-1})+C$
and choose $C_S$ once so $C_0 \cdot G_*(\Xi_\epsilon)^2,C_0\le M_Q$.
This is an exact maximum estimate, with no added radius constant.

Induct on the depth with the original profile parameter
\begin{equation}\label{eq:slow_offset_profile_induction}
 \Gamma_{Q,r}^2:=\frac{L_Q}{\varrho^2}(r+1)^{1/Q} \cdot J_{Q,r},\qquad
 J_{Q,r_{Q,\epsilon}}:=1,\qquad
 J_{Q,r+1}:=J_{Q,r}(1+r^{-2})^2\zeta_r^2.
\end{equation}
Eq.~\eqref{eq:slow_offset_depth_losses} gives $J_{Q,r}\le e^\epsilon$.
Let a uniform profile $\mathcal A\le\Gamma_{Q,r}^2 \cdot \ell_r(a^{-1})^2$
be established at depth $r$. We use the same parameters as in the
previous depth argument:
\begin{align}
 s&:=r+1,& \widehat\Gamma&:=(1+r^{-2})\Gamma_{Q,r},&
 \delta&:=\frac1{64Qs},& p&:=\frac\delta{128},\notag\\
 b&:=(1-2p)^{-1},& E&:=4,& M&:=2^{64}(1+C_F \cdot E^2),&
 k_0&:=\lceil100\delta^{-1} \cdot \log(M/\delta)\rceil.
 \label{eq:slow_offset_profile_parameters}
\end{align}
In particular $k_0\le C_k \cdot (Qs)^2$ with a universal $C_k$ independent
of $\epsilon$.
Lemma~\ref{lem:nearunit_static_transfer} applies: by
Eqs.~\eqref{eq:slow_offset_depth_losses} and
\eqref{eq:slow_offset_geometry_costs},
\begin{equation*}
 \frac{\Gamma_{Q,r}^2 \cdot \ell_r(0)^2}{H_Q}\ge e^{10\epsilon}.
\end{equation*}
Consequently
\begin{equation*}
 c_k^2\le[\widehat\Gamma^2 \cdot \ell_r(k)^2]^{k-1},\qquad
 c_k^2\le\mathfrak g(k)^{k-1},\qquad \mathfrak g:=\mathcal G^2.
\end{equation*}

Let $d_0$ be the least dyadic integer at least $Q$ and let $d\ge d_0$
be dyadic. Suppose for a contradiction that
\begin{equation}\label{eq:slow_offset_large_radius}
 \mathcal A>b^4 \cdot e^\delta \cdot \widehat\Gamma^2 \cdot \ell_r(d)^2.
\end{equation}
The right side exceeds $H_Q$ by
Eq.~\eqref{eq:slow_offset_geometry_costs}. Proposition~\ref{prop:slow_offset_blocks}
therefore supplies an actual radius $z\ge\mathcal A$ and its starting
family. Since $d\ge3$ and
$\mathcal L_r(d)\ge1$, with $u:=z^{-1}$ we have
\begin{equation}\label{eq:slow_offset_inverse_bounds}
 u\le L_Q^{-1} \cdot s^{-1/Q},\qquad u^Q\le L_Q^{-Q} \cdot s^{-1},\qquad
 b^4 \cdot u \cdot \widehat\Gamma^2 \cdot \ell_r(d)^2<e^{-\delta}.
\end{equation}
The retained prefix covers $2d_0-1<4Q$. Its energy excess obeys
$D_0:=\max\{\nu-u,0\}\le2a_Q \cdot u^{Q+1}$ once $a_Q \cdot u^Q\le1/2$.
The actual gap $R_{\mathrm{res}}-z\le2$ gives $C_P\le z+3$,
hence $\lambda\ge u/2$.

We verify $a_Q \cdot u^Q\le p$ below. Before a first exit of mass $p$, every
required normalizer is then at most
$B_*\le u/(1-p)^2\le bu$, and $\max\{1,B_* \cdot z\}\le b$.
Lemma~\ref{lem:second_logstar_polynomial_loss} gives
\begin{align*}
 \theta_d&:=5184C_F \cdot E^2 \cdot \sum_{\substack{k<d\\k\ \mathrm{dyadic}}}
                      k^8(b^4 \cdot u)^k \cdot c_k^2,\\
 \Delta_d&:=s_Q \cdot u^Q+2C_F \cdot \sum_{\substack{d_0\le k<d\\k\ \mathrm{dyadic}}}
                      k^2(bu)^k \cdot c_k^2,\qquad
 \tau_d:=C_F \cdot d(bu)^d \cdot c_d^2.
\end{align*}
For $k\le k_0$, Eq.~\eqref{eq:slow_offset_seed_cutoff} gives
\begin{equation*}
 \mathfrak g(k)\le e^{6\epsilon} \cdot M_Q \cdot s^{1/Q}.
\end{equation*}
Together with Eq.~\eqref{eq:slow_offset_large_radius} and
$\mathcal L_r(d)\ge1$, this implies
$b^4 \cdot u \cdot \mathfrak g(k)\le e^{-\epsilon}$.
We split this kernel once more using the original cap. We define
$m_\epsilon:=\lceil C_m \cdot \epsilon^{-1} \cdot \log(e/\epsilon)\rceil$, with
one sufficiently large universal $C_m$. For
$k\le\min\{m_\epsilon,k_0\}$, the cap gives
$\mathfrak g(k)\le G_*(m_\epsilon)^2\le C \cdot U_\epsilon^{1/3}$.
Here $m_\epsilon\le C \cdot \epsilon^{-2}$, so the fixed-power height
bound applies. Since $u\le L_Q^{-1}\le
(A_0 \cdot C_S^{1/3} \cdot U_\epsilon^{1/3})^{-1}$, a fixed choice of $C_S$
makes $b^4 \cdot u \cdot \mathfrak g(k)\le1/4$. This part of the kernel
is at most $Cu$. For $m_\epsilon<k\le k_0$, the near-unit estimate
above gives the bound
\begin{equation*}
 Cu\sum_{k>m_\epsilon}k^8 \cdot e^{-\epsilon(k-1)}
 \le Cu\epsilon^{-9} \cdot e^{-\epsilon m_\epsilon/2}\le C'u.
\end{equation*}
The last constant is universal after fixing $C_m$ large enough.
If $k_0<m_\epsilon$, the second range is empty. Thus the entire
low-degree kernel is at most $Cu$, uniformly in $\epsilon$.
Its coherent part is at most $D_Q \cdot u^Q$, where
\begin{equation*}
 D_Q:=b^Q \cdot E_{Q,\epsilon},\qquad
 D_Q^{1/Q}\le e^{8\epsilon} \cdot M_Q.
\end{equation*}
Here Eq.~\eqref{eq:slow_offset_coherent_root} gives
$E_{Q,\epsilon}^{1/Q}\le e^\epsilon \cdot \mathfrak g(N)$,
$N\le X^X$, and $b\le e^\epsilon$ for $Q\ge Q_\epsilon$.
Keeping the factor $b^Q$ explicit preserves the small relative cost.

For $k>k_0$, the current profile gives the same bounds as before:
$(b^4 \cdot u)^k \cdot c_k^2\le b^4 \cdot ue^{-\delta(k-1)}$ and
$(bu)^k \cdot c_k^2\le e^{-\delta k}$.
The two integer tail estimates following
Eq.~\eqref{eq:second_logstar_profile_kernel} therefore bound their
contributions by $u$ and $p/512$, respectively.
Their constants depend on $C_F,E$ and the universal cutoff, and do not
depend on $\epsilon$. We conclude that
\begin{equation}\label{eq:slow_offset_kernel_split}
 \theta_d\le C_\theta \cdot u,\qquad
 \Delta_d\le(s_Q+D_Q)u^Q+p/512,
\end{equation}
with a universal $C_\theta$ independent of $\epsilon$.

The parameters satisfy the exact amplitude conditions
\begin{equation}\label{eq:slow_offset_amplitude_conditions}
 L_Q\ge\max\{16,16C_\theta\},\qquad
 L_Q^Q\ge2^{24} \cdot Q(a_Q+s_Q+D_Q+4C_\theta \cdot a_Q+1).
\end{equation}
The first uses the fixed large minimum $M_Q\ge A_0 \cdot C_S^{1/3} \cdot U_\epsilon^{1/3}$.
For the second, Eqs.~\eqref{eq:slow_offset_geometry_costs} and
\eqref{eq:slow_offset_kernel_split} leave a relative margin of at least
$e^{12\epsilon Q}$, whereas the remaining factor is at most
$C \cdot Q$. The inequality
$\log(C \cdot Q)\le\epsilon Q$ holds for every
$Q\ge C_Q \cdot \epsilon^{-4}$ after a universal choice of $C_Q$.
This proves Eq.~\eqref{eq:slow_offset_amplitude_conditions} without
multiplying $L_Q$ by a fixed constant.

The verification following
Eq.~\eqref{eq:second_logstar_amplitude_conditions} now applies with
$L_Q$ and $C_\theta$: it gives
\begin{equation*}
 \nu\le1/4,\quad D_0\le u,\quad
 \theta_d \cdot u/\lambda\le1/8,\quad
 \Delta_d+\theta_d \cdot D_0/\lambda\le p/64,\quad
 D_0/p\le u,\quad a_Q \cdot u^Q\le p.
\end{equation*}
The actual matched budget and
Lemma~\ref{lem:matched_scale_comparison} consequently yield
\begin{equation}\label{eq:slow_offset_terminal_comparison}
 a<12(u+D_0/p)\tau_d\le24C_F \cdot d(bu)^d \cdot c_d^2.
\end{equation}
Choose the dyadic terminal degree as in
Eq.~\eqref{eq:second_logstar_terminal_degree}. Its constant $C_T$
remains universal: the inequality
$\log(24C_F \cdot d/a)\le\delta d$ and the factor $24$ in
Eq.~\eqref{eq:slow_offset_terminal_comparison} have the same form as in
Proposition~\ref{prop:second_logstar_profile}.
All dependence on $\epsilon$ has already been absorbed in
Eq.~\eqref{eq:slow_offset_amplitude_conditions}.
The contradiction following
Eq.~\eqref{eq:second_logstar_terminal_comparison} proves
\begin{equation*}
 \mathcal A\le b^4 \cdot e^\delta \cdot \widehat\Gamma^2 \cdot \ell_r(d)^2.
\end{equation*}
The same universal $C_b$ in
Lemma~\ref{lem:repeated_height_distortion} therefore supplies the
bound at the next depth with the amplitudes in
Eq.~\eqref{eq:slow_offset_profile_induction}.

Initialize the induction directly with
Eq.~\eqref{eq:slow_offset_initial_bound}. The induction variable
is $\mathcal A$, so no initial bound on an actual block radius is
required. The floor estimates ensure that every contradictory threshold
is above $H_Q$. They also cover $a\ge1$, where
$\mathcal A\le\max\{1,C_P\}=1$.
Finite induction gives
\begin{equation*}
 \mathcal A\le e^{21\epsilon} \cdot M_Q \cdot (r+1)^{1/Q} \cdot \mathcal L_r(a^{-1})^2.
\end{equation*}
This proves the next full profile for the coefficient radius with amplitude
$e^{21\epsilon} \cdot A$, which is covered by $e^{24\epsilon} \cdot A$.
Each stage propagates the profile for the coefficient radius directly.

\emph{The low orders and the first stage.}
For every odd $Q<Q_\epsilon$,
Proposition~\ref{prop:repeated_height_original_seed_profile} gives
\begin{equation*}
 \mathcal A\le\max\{1,C_P\}
 \le C \cdot \mathfrak t(Q)^{1/3}(r+1)^{1/Q} \cdot \ell_r(a^{-1})^2
 ,\ \forall r\ge r_{Q,\epsilon}.
\end{equation*}
The initial depth required by that proposition is below the present one. Since
$Q_\epsilon\le C \cdot \epsilon^{-4}$,
$\mathfrak t(Q_\epsilon)\le\mathfrak t(\epsilon^{-1})+C$.
Its normalized amplitude is thus at most $C \cdot U_\epsilon^{1/3}$,
covered by one fixed $A_0,C_S$ at every stage. In particular all smaller
orders used in extraction and initialization are already available.

The first full stage follows directly from
Proposition~\ref{prop:repeated_height_original_seed_profile} at the same
order $Q$. Its initial depth is at most $r_{Q,\epsilon}$, and
\begin{equation*}
 C \cdot \varrho^2 \cdot \mathfrak t(Q)^{1/3}
      \le A_0 \cdot [\mathfrak t(Q)+S_\epsilon]^{1/3}
\end{equation*}
after one universal choice of $A_0$. Thus it supplies the required
bound at every permitted depth without a separate initialization or
an additional multiplicative loss.

All constants, including those of the geometry, were fixed
before the stage count. Finite induction proves
Eq.~\eqref{eq:slow_offset_outer_profile}. The only conversion to the
ordinary Poincar\'e constant will be made after these repetitions.
\end{proof}

The smaller loss per repetition requires larger cutoffs and initial depths.
We control their cost using inherited coefficient bounds and height functions
that contract changes in their arguments.

\subsection{Returning to a common domain of prescribed depths}
\label{sec:inherited_heights}
\label{sec:secondary_feedback_heights}
\label{sec:secondary_feedback_profiles}

Completing a finite number of outer rounds reduces the log-star
factor, while the required precision contributes another scalar cost.
We encode that cost by a height function and return the profile to every
depth above one fixed threshold. The same lower threshold can therefore
be used in the next profile refinement. The index $m$ below counts
these refinements of the profile; $j$ still counts the outer rounds of
repeated height reduction, and $r$ the inner depth of iterative curvature
refinement.

We use $\mathfrak t$ and $\kappa$ from
Eq.~\eqref{eq:second_logstar_seed_notation} and
Lemma~\ref{lem:repeated_height_termination}. For $x\ge1$, we define
\begin{equation}\label{eq:secondary_feedback_heights}
 \begin{gathered}
 \chi(x):=\max\{3,1+\kappa(x)\},\qquad
 \chi_0(x):=x,\qquad \chi_{m+1}(x):=\chi(\chi_m(x)),\\
 V_m(x):=\mathfrak t(\chi_m(x)),\ \forall m\ge0,\qquad
 N_\chi(x):=\min\{m\ge0:\chi_m(x)\le3\}.
 \end{gathered}
\end{equation}
The lower cutoff fixes the terminal value at three. In particular,
$V_0=\mathfrak t$ and
$V_1(x)=\mathfrak t(1+\kappa(x))$.

\begin{lemma}[Uniform estimates for the repetition count]
\label{lem:secondary_feedback_heights}
The functions $\kappa$ and $\chi$ are nondecreasing, and each has
Lipschitz constant one on the positive integers. The map $\chi$ fixes
three and satisfies $\chi(x)<x$ for $x>3$. For every $m\ge0$,
\begin{equation}\label{eq:secondary_feedback_height_bounds}
 3\le V_m(x)\le\mathfrak t(x),\qquad
 V_{m+1}(x)\le V_m(x),\ \forall x\ge1.
\end{equation}
For fixed $C,p\ge1$ there is an integer $b_{C,p}$, independent of $m$,
such that
\begin{equation}\label{eq:secondary_feedback_fixed_power}
 V_m(C \cdot x^p)\le V_m(x)+b_{C,p},\ \forall x\ge1,m\ge0.
\end{equation}
There is a universal $b$ such that, for $x,y\ge1$ and $m\ge0$,
\begin{equation}\label{eq:secondary_feedback_product}
 \max\{V_m(x+y),V_m(xy)\}
       \le\max\{V_m(x),V_m(y)\}+b.
\end{equation}
For each fixed $B\ge0$ there is a constant $C_B$, independent of $m$,
satisfying
\begin{equation}\label{eq:secondary_feedback_weighted_cutoff}
 \sup_{y\ge0}e^{-y} \cdot V_m(h+B+y)^{1/3}
       \le C_B \cdot V_m(h)^{1/3},\ \forall h\ge1.
\end{equation}
Finally, $N_\chi(x)$ is finite for every $x$ and is unbounded as
$x\to\infty$.
\end{lemma}
\begin{proof}
Monotonicity of every $\mathfrak t_j$ implies monotonicity of $\kappa$.
The relevant successive exponential thresholds for $\log^*$ are more
than one apart. Thus
$\mathfrak t(n+1)-\mathfrak t(n)\in\{0,1\}$ for positive integers
$n$. Since $\mathfrak t(4)=3$, it also follows that
$\mathfrak t(n+1)\le n$ for $n\ge3$.
The values of $\kappa$ are zero through three and $\kappa(4)=1$.
For $n\ge4$, use
\begin{equation*}
 \kappa(n)=1+\kappa(\mathfrak t(n)),\qquad
 \kappa(n+1)=1+\kappa(\mathfrak t(n+1)).
\end{equation*}
The two arguments are at most $n$ and differ by at most one. Strong
induction gives $0\le\kappa(n+1)-\kappa(n)\le1$. Taking the maximum
with three proves the same assertion for $\chi$.
These Lipschitz assertions concern integer arguments; no continuity
across a real height threshold is required.

The strict descent of $\mathfrak t$ above three gives
$\kappa(n)\le n-3$ for integers $n\ge3$. Consequently, for real $x>3$,
\begin{equation*}
 \kappa(x)=1+\kappa(\mathfrak t(x))\le\mathfrak t(x)-2.
\end{equation*}
It follows that
$\chi(x)\le\max\{3,\mathfrak t(x)-1\}<x$ for $x>3$.
On $[1,3]$, $\chi(x)=3$ and $\mathfrak t(x)=3$.
Monotonicity now proves Eq.~\eqref{eq:secondary_feedback_height_bounds}.

The case $m=0$ of Eq.~\eqref{eq:secondary_feedback_fixed_power} is
Eq.~\eqref{eq:repeated_height_power_argument}. Enlarge its integer
constant if necessary. For $x>3$, integer nonexpansiveness gives
\begin{align*}
 \kappa(C \cdot x^p)
 &=1+\kappa(\mathfrak t(C \cdot x^p))\\
 &\le1+\kappa(\mathfrak t(x)+b_{C,p})
 \le\kappa(x)+b_{C,p}.
\end{align*}
Increasing $b_{C,p}$ once covers $1\le x\le3$ as well. Hence
$\chi(C \cdot x^p)\le\chi(x)+b_{C,p}$.
Both values of $\chi$ are integers, so applying its integer
nonexpansiveness repeatedly gives
\begin{equation*}
 \chi_m(C \cdot x^p)\le\chi_m(x)+b_{C,p},\ \forall m\ge1.
\end{equation*}
A final application of integer nonexpansiveness of $\mathfrak t$
proves Eq.~\eqref{eq:secondary_feedback_fixed_power} for every $m$,
with the same constant. If $M:=\max\{x,y\}$, then $xy\le M^2$ and
$x+y\le2M$. These inequalities and the fixed-power estimate prove
Eq.~\eqref{eq:secondary_feedback_product}.

For the weighted estimate, the product and sum bounds, monotonicity,
and $V_m\le\mathfrak t$ imply
\begin{equation*}
 V_m(h+B+y)\le V_m(h)+y+C_B,\ \forall h\ge1,y\ge0,
\end{equation*}
where we used $\mathfrak t(1+y)\le y+C$.
Since $V_m(h)\ge3$, taking cube roots and multiplying by $e^{-y}$
gives Eq.~\eqref{eq:secondary_feedback_weighted_cutoff} with a
constant independent of $m$.

After one application of $\chi$, its argument is an integer at least
three. Its strict descent therefore proves termination. Each fixed
composition of $\mathfrak t$ is unbounded, so $\kappa$ is unbounded:
for any fixed $j$, sufficiently large $x$ has $\mathfrak t_j(x)>3$.
Thus $\chi$, and every fixed finite composition of $\chi$, is
unbounded. If $\lambda$ were bounded by a fixed $m$, then
$\chi_m(x)\le3$ for all $x$, a contradiction.
\end{proof}

We use the functions $\chi$, $\chi_m$, and $V_m$ from
Eq.~\eqref{eq:secondary_feedback_heights}. The additional repetition
will carry a profile at every prescribed depth. Keeping the same lower
depth for every repetition avoids a growing threshold. Within each such
repetition we use only sufficiently large block orders.

Choose the universal constants $C_b,C_R,C_Q$ of the preceding profile
argument first. We then choose one sufficiently large integer $R_0$,
above its transfer thresholds. Constants depending on this fixed $R_0$
will be absorbed in amplitudes and offsets, without changing $C_b,C_R$
or $C_Q$.

\begin{proposition}[Profiles at large block orders from an initial bound]
\label{prop:secondary_feedback_restricted_profile}
Suppose that for some $m\ge0$ and $C_m\ge1$, uniformly over centered
regular measures of covariance at most $I$ and curvature at least $aI$,
\begin{equation}\label{eq:secondary_feedback_prescribed_input}
 \mathcal A\le C_m \cdot V_m(r)^{1/3} \cdot \mathcal L_r(a^{-1})^2
 ,\ \forall r\ge R_0,\ r\text{ an integer}.
\end{equation}
There are universal constants $D_{\mathrm{prof}},C_S\ge1$ such that the following holds.
For $0<\epsilon\le1/16$, we define
\begin{align}
 S_{m,\epsilon}&:=C_S \cdot V_m(\epsilon^{-1}),\qquad
 W_{m,j,\epsilon}(x):=[\mathfrak t_j(x)+S_{m,\epsilon}]^{1/3},
                        \label{eq:secondary_feedback_offset}\\
 Q_\epsilon&:=\text{the least odd integer at least }C_Q \cdot \epsilon^{-4},
 \notag\\
 r_{Q,\epsilon}&:=\max\{R_0,\lceil C_b \cdot \mathfrak t(Q)\rceil
                                  +\lceil C_R/\epsilon\rceil\}.
                        \label{eq:secondary_feedback_restricted_depth}
\end{align}
Then, for every integer $j\ge1$, odd $Q\ge Q_\epsilon$, and integer
$r\ge r_{Q,\epsilon}$,
\begin{equation}\label{eq:secondary_feedback_restricted_profile}
 \mathcal A\le D_{\mathrm{prof}} \cdot C_m \cdot e^{24\epsilon j} \cdot W_{m,j,\epsilon}(Q) \cdot
                       (r+1)^{1/Q} \cdot \mathcal L_r(a^{-1})^2.
\end{equation}
The constants are independent of $m,C_m,j,Q,r,\epsilon$, and the law.
No bound at orders below $Q_\epsilon$ is asserted here.
\end{proposition}
\begin{proof}
We first extract a coefficient majorant that will remain fixed throughout
this repetition. Let
$r(x):=\max\{R_0,\lceil C_r \cdot \mathfrak t(x)\rceil\}$, where $C_r$
is universal and sufficiently large. The scalar estimate in the proof of
Lemma~\ref{lem:iterated_coefficient_seed} gives
$\mathcal L_{r(d)}(d)\le C$. Lemma~\ref{lem:nearunit_static_transfer},
applied to Eq.~\eqref{eq:secondary_feedback_prescribed_input} at the
fixed depth $r(d)$, yields
\begin{equation*}
 c_d^2\le[C_e \cdot C_m \cdot V_m(r(d))^{1/3}]^{d-1},\ \forall d\ge2.
\end{equation*}
The condition that the floor is one follows from $V_m\ge3$ and
$\mathcal L_r(0)^2\ge e^{-1/16}$ at the fixed large depth.
The transfer uses only finite degrees and gives the estimate over the
full log-concave class of covariance at most $I$. The uniform
fixed-argument estimate of Lemma~\ref{lem:secondary_feedback_heights}
therefore permits the monotone fixed majorant
\begin{equation}\label{eq:secondary_feedback_fixed_seed}
 \mathcal H_m(x)^2:=\min\{G_*(x)^2,
                       C_e' \cdot C_m \cdot V_m(\mathfrak t(x))^{1/3}\}.
\end{equation}
Increase $C_e'$ once to cover degree one and the fixed depth $R_0$.
Both branches are at least one and nondecreasing. This majorant depends
on the fixed amplitude $C_m$, and not on the subsequent index $j$.

Assume the preceding stage at the allowed block orders has amplitude $A$,
where $A\ge D_{\mathrm{prof}} \cdot C_m$. Define
\begin{equation}\label{eq:secondary_feedback_selected_order}
 \begin{split}
 r_\epsilon(x)&:=\max\{R_0,\lceil C_r \cdot \mathfrak t(x)\rceil
                                    +\lceil D_R/\epsilon\rceil\},\\
 k_\epsilon(x)&:=\text{the least odd integer at least }
       \max\{Q_\epsilon,\epsilon^{-1} \cdot \log(r_\epsilon(x)+1)\}.
 \end{split}
\end{equation}
Choose $D_R\ge2C_R+2$ before $R_0$. Then
$r_\epsilon(x)\ge r_{k_\epsilon(x),\epsilon}$ for every $x$.
Indeed, writing $r:=r_\epsilon(x)$, both entries in the maximum
for $k_\epsilon(x)$ are at most $C \cdot r^4$, since
$\epsilon^{-1}\le r/D_R$. Thus
$C_b \cdot \mathfrak t(k_\epsilon(x))+2\le r/2$ once $R_0$ is fixed
large enough. The remaining $\lceil C_R/\epsilon\rceil$ uses the
other half. The same choice leaves
\begin{equation*}
 (r_\epsilon(x)+1)^{1/k_\epsilon(x)}\le e^\epsilon,\qquad
 r_\epsilon(x)^{-2}\le\epsilon,\qquad
 \mathcal L_{r_\epsilon(x)}(x)\le1+\epsilon.
\end{equation*}
The normalized static transfer consequently gives, for every degree,
\begin{equation*}
 c_d^2\le[e^{5\epsilon} \cdot A
                   W_{m,j,\epsilon}(k_\epsilon(d))]^{d-1}.
\end{equation*}
Every order used here is at least $Q_\epsilon$. Take the working squared
majorant
\begin{equation}\label{eq:secondary_feedback_working_seed}
 \mathfrak g(x):=\min\{\mathcal H_m(x)^2,
                  e^{5\epsilon} \cdot A \cdot W_{m,j,\epsilon}(k_\epsilon(x))\}.
\end{equation}
It is valid, monotone, and capped by the same original $G_*^2$.

Write $h:=\mathfrak t(Q)$ and $M_Q:=A \cdot W_{m,j+1,\epsilon}(Q)$.
Choose the universal constant $C_h$ so that the high-height range is
\begin{equation}\label{eq:secondary_feedback_height_threshold}
 h\ge h_\epsilon:=C_h \cdot \epsilon^{-5}.
\end{equation}
In this range $k_\epsilon(x)\le h/2$ whenever
$\mathfrak t(x)\le h+B+1$, for each fixed $B$ needed below.
The entry $Q_\epsilon\le C \cdot \epsilon^{-4}$ is covered by this
choice. The other entry is bounded by
$3+\epsilon^{-1} \cdot \log[C \cdot (h+B+1+\epsilon^{-1})]$.
As in Eq.~\eqref{eq:slow_offset_cutoff_height}, for
$X:=C_X \cdot Q^2$, $Y:=C_k \cdot (Qs)^2$, and $y:=\log(s)/Q$,
\begin{equation*}
 \mathfrak t(X^X)\le h+B,\qquad \mathfrak t(Y)\le h+B+y.
\end{equation*}
For $y\le1$, we have $k_\epsilon(Y)\le h$. For $y\ge1$, the
logarithmic entry is at most
\begin{equation*}
 h/2+2+\epsilon^{-1} \cdot \log(1+y/(h+B+\epsilon^{-1}))
       \le h+y\le h(1+y),
\end{equation*}
and the $Q_\epsilon$ entry is already at most $h/2$.
The uniform product bound for $\mathfrak t_j$ gives
$\mathfrak t_j(k_\epsilon(Y))\le\mathfrak t_j(h)+5y$.
Taking the cube root and dividing by $e^y$ proves the same high-height
estimates as before, with no fixed additive height error.

In the low-height range, use the fixed branch of
Eq.~\eqref{eq:secondary_feedback_working_seed}. By
Lemma~\ref{lem:secondary_feedback_heights}, uniformly in $m$,
\begin{equation*}
 V_m(\mathfrak t(Y))\le C[V_m(h+B+1)+1+y].
\end{equation*}
Its cube root divided by $e^y$ has supremum at most
$C \cdot V_m(h+B+1)^{1/3}$. Since $h+B+1\le C \cdot \epsilon^{-5}$, the
uniform fixed-power estimate bounds this by
$C \cdot V_m(\epsilon^{-1})^{1/3}$. The same estimate controls $X^X$.
Choose $D_{\mathrm{prof}},C_S$ once so that this fixed cost, including $C_e' \cdot C_m$,
is at most $M_Q$. We obtain, in both height ranges,
\begin{equation}\label{eq:secondary_feedback_seed_bounds}
 \mathfrak g(X^X)\le e^{5\epsilon} \cdot M_Q,\qquad
 \sup_{s\ge1}\frac{\mathfrak g(C_k \cdot (Qs)^2)}{s^{1/Q}}
                                      \le e^{5\epsilon} \cdot M_Q.
\end{equation}
The coherent-root estimate in
Eq.~\eqref{eq:slow_offset_coherent_root} uses only monotonicity and
the original polynomial cap. It therefore applies with
$N:=\lceil D_1 \cdot \epsilon^{-2} \cdot Q^2\rceil\le X^X$, since
$Q\ge Q_\epsilon$.

For the initial depth of the next stage $r:=r_{Q,\epsilon}$, the high-height
range uses the preceding stage at the least odd
$k\ge\max\{Q_\epsilon,\epsilon^{-1} \cdot \log(r+1)\}$.
The same threshold gives $k\le h\le Q$ and
$r\ge r_{k,\epsilon}$, so
\begin{equation}\label{eq:secondary_feedback_initialization}
 \mathcal A\le e^\epsilon \cdot M_Q \cdot \mathcal L_r(a^{-1})^2.
\end{equation}
In the low-height range, $r\le R_0+C \cdot \epsilon^{-5}$. Apply the
assumed bound in Eq.~\eqref{eq:secondary_feedback_prescribed_input} at
this very depth. The uniform fixed-power estimate gives
\begin{equation*}
 \mathcal A\le C \cdot C_m \cdot V_m(\epsilon^{-1})^{1/3} \cdot \mathcal L_r(a^{-1})^2
              \le M_Q \cdot \mathcal L_r(a^{-1})^2.
\end{equation*}
Thus Eq.~\eqref{eq:secondary_feedback_initialization} holds everywhere.
No profile at a smaller block order has been used in this low range.

Apply Proposition~\ref{prop:slow_offset_blocks} with the fixed cap
$\mathcal H:=\mathcal H_m$ and the working majorant
$\mathcal G:=\sqrt{\mathfrak g}$. Since $\Xi_\epsilon=16K_\epsilon$
is at most $C \cdot \epsilon^{-2}$, Eq.~\eqref{eq:secondary_feedback_fixed_seed}
and the uniform scalar estimates give
\begin{equation*}
 C_0 \cdot \mathcal H_m(\Xi_\epsilon)^2,C_0
          \le C \cdot C_m \cdot V_m(\epsilon^{-1})^{1/3}\le M_Q.
\end{equation*}
The exact maximum in its floor therefore satisfies
\begin{equation}\label{eq:secondary_feedback_geometry_costs}
 F^2\le e^{6\epsilon} \cdot M_Q,\qquad H_Q\le e^{9\epsilon} \cdot M_Q,
 \qquad(a_Q+s_Q+1)^{1/Q}\le e^{6\epsilon} \cdot M_Q.
\end{equation}
Only the coefficient bound from the start of the repetition occurs in
this absolute floor.

We verify the other change needed in the depth induction. The cutoff for the small degrees in its kernel is $l_\epsilon:=\lceil C_l \cdot \epsilon^{-1} \cdot \log(e/\epsilon)\rceil$.
For $k\le l_\epsilon$, the fixed cap gives
\begin{equation*}
 \mathfrak g(k)\le C \cdot C_m \cdot V_m(\epsilon^{-1})^{1/3}.
\end{equation*}
Under the assumption on the radius made for contradiction with
$L_Q:=e^{20\epsilon} \cdot M_Q$, its inverse radius obeys
$u\le L_Q^{-1} \cdot s^{-1/Q}$. Thus $D_{\mathrm{prof}},C_S$ make the geometric ratio
at most $1/4$ on this range, independently of $m,C_m,j,\epsilon$.
Above $l_\epsilon$, Eq.~\eqref{eq:secondary_feedback_seed_bounds}
provides the same $e^{-\epsilon k}$ summation. Hence the full
$C_\theta$ in Eq.~\eqref{eq:slow_offset_kernel_split} remains universal.
The coherent-root estimate and Eq.~\eqref{eq:secondary_feedback_geometry_costs}
then give the exact amplitude conditions in
Eq.~\eqref{eq:slow_offset_amplitude_conditions}, with the same
$Q_\epsilon$ and $L_Q$.

All remaining parameters and inequalities of the proof of
Proposition~\ref{prop:slow_offset_outer_profiles} are now available:
the normalized floor, contradiction for the constructed family, available prefix,
matched budget, and the universal terminal degree. Its inner induction
therefore applies with
\begin{equation*}
 \Gamma_{Q,r}^2:=\frac{L_Q}{\varrho^2}(r+1)^{1/Q} \cdot J_{Q,r},\qquad
 J_{Q,r_{Q,\epsilon}}:=1,\qquad
 J_{Q,r+1}:=J_{Q,r}(1+r^{-2})^2\zeta_r^2.
\end{equation*}
The product is at most $e^\epsilon$. Initialize by
Eq.~\eqref{eq:secondary_feedback_initialization}. The same comparison
proves the full profile at the next depth with amplitude at most
$e^{21\epsilon} \cdot A$, covered by $e^{24\epsilon} \cdot A$.
Throughout this induction we estimate $\mathcal A$ directly, without
converting to $C_P$ or to the restricted radius. In particular its terminal constant
$C_T$ is independent of $m$ and $\epsilon$.

The first full stage follows directly from
Proposition~\ref{prop:repeated_height_original_seed_profile} at the
same $Q\ge Q_\epsilon$: its normalized amplitude
$C \cdot \mathfrak t(Q)^{1/3}$ is at most
$D_{\mathrm{prof}} \cdot C_m \cdot [\mathfrak t(Q)+S_{m,\epsilon}]^{1/3}$.
Its initial depth is smaller than Eq.~\eqref{eq:secondary_feedback_restricted_depth}.
Finite induction proves Eq.~\eqref{eq:secondary_feedback_restricted_profile}.
The constants $C_b,C_R,C_Q$ were fixed before $R_0$. All remaining
scalar costs depending on the fixed $R_0$ were absorbed in $D_{\mathrm{prof}},C_S$, which only
increase the amplitude and do not change the order or depth thresholds.
\end{proof}

We now choose the number of outer rounds at each prescribed
depth. This improves the height profile and preserves the same
range of depths, at the cost of one fixed factor in the amplitude.

\begin{proposition}[Uniform profiles at prescribed depths]
\label{prop:secondary_feedback_prescribed_profile}
There are universal constants $D_0,D\ge1$ and an integer $R_0$ such
that, with $C_m:=D_0 \cdot D^m$, every centered regular measure of covariance
at most $I$ and curvature at least $aI$ satisfies
\begin{equation}\label{eq:secondary_feedback_prescribed_profile}
 \mathcal A\le C_m \cdot V_m(r)^{1/3} \cdot \mathcal L_r(a^{-1})^2
 ,\ \forall m\ge0,\ r\ge R_0,\ r\text{ an integer}.
\end{equation}
All constants, including the lower depth $R_0$, are independent of $m$.
\end{proposition}
\begin{proof}
For $m=0$, choose the least odd
$q\ge\max\{3,\log(r+1)\}$ in
Proposition~\ref{prop:repeated_height_original_seed_profile}.
For every $r$ above one universal threshold, $q\le r$, its initial
depth is at most $r$, and $(r+1)^{1/q}\le e$.
Use $\mathcal A\le\max\{1,C_P\}$. The fixed lower bound
$\mathcal L_r(0)^2\ge e^{-1/16}$ allows the constant one to be
absorbed by a universal $D_0$. Since
$\mathfrak t(q)\le\mathfrak t(r)=V_0(r)$, this proves the initial
assertion, including the case $C_P<1$.

Choose $R_0$ once, also large enough for the following finite choices
at every $r\ge R_0$:
\begin{equation}\label{eq:secondary_feedback_invariant_parameters}
 j:=\kappa(r),\qquad \epsilon:=\frac1{\max\{16,j+1\}},\qquad
 Q:=\text{the least odd integer at least }
       \max\{Q_\epsilon,\log(r+1)\}.
\end{equation}
We require $j\ge1$, $Q\le r$, and $r\ge r_{Q,\epsilon}$.
These are universal scalar conditions. Indeed
$\kappa(r)\le\mathfrak t(r)-2$ gives
$Q_\epsilon\le C \cdot \mathfrak t(r)^4\le r$ beyond a fixed threshold,
with enough room for the odd rounding. Further,
\begin{equation*}
 r_{Q,\epsilon}\le\max\{R_0,
          C_b \cdot \mathfrak t(r)+C_R\max\{16,1+\kappa(r)\}+2\}\le r
\end{equation*}
for every such $r$. This choice depends only on the already fixed
$C_b,C_R,C_Q$, and not on the profile index or amplitude.

Suppose Eq.~\eqref{eq:secondary_feedback_prescribed_profile} holds
at $m$. Apply Proposition~\ref{prop:secondary_feedback_restricted_profile}
at the finite choices in Eq.~\eqref{eq:secondary_feedback_invariant_parameters}.
Then $\mathfrak t_j(Q)\le\mathfrak t_j(r)=3$,
$(r+1)^{1/Q}\le e$, and $\epsilon j\le1$. Also, by the uniform
fixed-argument estimate,
\begin{equation*}
 V_m(\epsilon^{-1})\le V_m(16\chi(r))
                 \le V_m(\chi(r))+C\le C' \cdot V_{m+1}(r).
\end{equation*}
It follows that
\begin{equation*}
 \mathcal A\le D' \cdot C_m \cdot V_{m+1}(r)^{1/3} \cdot \mathcal L_r(a^{-1})^2,
\end{equation*}
with one universal $D'$. Choose $D\ge D'$.
This proves the next prescribed-depth profile on the same domain
$r\ge R_0$, and finite induction proves the proposition.
\end{proof}

\subsection{Using an inherited bound to lower the block floor}
\label{sec:nearunit_secondary_geometry}

We lower the block floor required for a small loss by retaining the
coefficient bound available at the start of the profile refinement.
Different degree ranges use the original cap, the inherited margin, and
the current margin. This lets the smallest loss govern the later degree
range while the inherited bound controls the earlier degrees, and yields
the orbit, restart, and extension estimates needed for the next refinement.

We use the common coefficient radius $\mathcal A$, the operator
$\mathcal T$, and the class of centered regular measures of covariance
at most $I$ and strictly positive curvature. The nondecreasing majorant
$\mathcal G$ satisfies Eq.~\eqref{eq:second_logstar_majorant}.
Fix also a nondecreasing coefficient majorant $\mathcal H$ such that,
uniformly over the same laws,
\begin{equation}\label{eq:nearunit_secondary_cap}
 \begin{gathered}
 1\le\mathcal H(x)\le G_*(x),\ \forall x\ge1,\qquad
 \mathcal G\le\mathcal H,\\
 c_k\le\mathcal H(k)^{k-1},\ \forall k\ge1.
 \end{gathered}
\end{equation}
Any valid additional majorant may first be capped by $G_*$, and the
current $\mathcal G$ may be replaced by its minimum with $\mathcal H$.
The original seed is allowed as the choice $\mathcal H:=G_*$.
All constants below are universal and independent of $\mathcal H$.
Fix sufficiently large universal constants $C_{\mathrm{cut}},C_{\mathrm{deg}}$.
For the construction with an inherited floor, we assume
$0<\epsilon\le\delta\le1/16$ and define
\begin{equation}\label{eq:nearunit_secondary_cutoff}
 K_\alpha:=\lceil C_{\mathrm{cut}} \cdot \alpha^{-2}\rceil,\qquad
 \Xi:=\lceil4C_{\mathrm{deg}} \cdot \delta^{-1} \cdot K_\epsilon\rceil,\qquad
 \mathsf A_\delta:=\max\{C_0 \cdot G_*(16K_\delta)^2,C_0\}.
\end{equation}
Thus $K_\delta\le K_\epsilon$ and $\Xi\ge K_\epsilon$.
The small degrees use the original coefficient cap, the next range
uses the inherited $\delta$ margin, and the final range uses $\epsilon$.

\begin{lemma}[Three-range moments]
\label{lem:nearunit_secondary_seed_moments}
Fix nonnegative integers $M,r$ and $c>0$. Suppose $0\le\rho_k\le1$,
$\rho_k\le C_0^{-1}$ for $k\le16K_\delta$, and
$\rho_k\le(1+\delta)^{-1}$ for $k\le\Xi$. Then
\begin{equation}\label{eq:nearunit_secondary_seed_moments}
 \sum_{k\ge r+1}k^M \cdot \rho_k^{k-r} \cdot e^{-c\epsilon(k-r)}\le C.
\end{equation}
The bound is uniform in both parameters and the sequence. It tends
uniformly to zero when the fixed thresholds $C_{\mathrm{cut}},C_0$
increase.
\end{lemma}
\begin{proof}
Up to $K_\delta$, use the geometric sum with ratio $C_0^{-1}$.
Choose $C_{\mathrm{cut}}$ so $K_\delta\ge2r$. Between $K_\delta$
and $\Xi$, discard the $\epsilon$ factor and use
$(1+\delta)^{-(k-r)}\le e^{-\delta(k-r)/2}$. This contribution is at
most $C \cdot \delta^{-M-1} \cdot e^{-c' C_{\mathrm{cut}}/\delta}$.
Above $\Xi\ge K_\epsilon$, use $\rho_k\le1$ and the $\epsilon$
factor. The bound is $C \cdot \epsilon^{-M-1} \cdot e^{-c''C_{\mathrm{cut}}/\epsilon}$.
Both tails are uniformly small as the fixed cutoff increases. Only
finitely many moments will be used.
\end{proof}

The moment bounds over the three ranges give the kernel and early coherent
estimates for an actual orbit. The orbit and restart then supply the
same family and its matched budget under the inherited floor.

\begin{lemma}[The Green and restart estimates]
\label{lem:nearunit_secondary_green}
\label{lem:nearunit_secondary_initialization}
Suppose that the hypotheses of Lemma~\ref{lem:coefficient_radius_green}
hold, except that its condition on the coefficient floor is replaced by
\begin{equation}\label{eq:nearunit_secondary_green_hypotheses}
 F^2\ge\max\{(1+\delta)\mathcal H(\Xi)^2,\mathsf A_\delta\}.
\end{equation}
Its orbit and window conclusions,
Eqs.~\eqref{eq:coefficient_radius_orbit_bounds} and
\eqref{eq:coefficient_radius_window_loss}, still hold with universal
constants. Under the additional hypotheses in
Eq.~\eqref{eq:coefficient_radius_initialization_hypotheses}, the same
actual family has all conclusions of
Eq.~\eqref{eq:coefficient_radius_initialization}.
\end{lemma}
\begin{proof}
Keep the raw projection, convolution, and delays from
Eqs.~\eqref{eq:coefficient_radius_local_raw_frame}--\eqref{eq:coefficient_radius_local_kernel}.
In particular a degree-$k$ swap has squared propagation factor at most
$16(3k)^2$ under the assumed envelope $4(h+1)$; its delays are
$k,\ldots,4k-2$. Set $\rho_k:=\min\{1,\mathcal H(k)^2/F^2\}$.
Eq.~\eqref{eq:nearunit_secondary_green_hypotheses} gives the two
ratio hypotheses of Lemma~\ref{lem:nearunit_secondary_seed_moments}.
At every used degree $k\le d$, $\mathcal G(k)\le F$ and
\begin{equation*}
 \mathcal G(k)^2/z\le\rho_k \cdot e^{-\epsilon/2}.
\end{equation*}
The moments with weights $k^8$ and $k^9$ therefore give the same
kernel bounds $\eta,\eta_1\le C/z$, uniformly in both parameters.

Before the orbit is bounded, the assumed envelope gives
$b_j\le16(j+1)^2$. Degree $k<d$ is used only at
$j=2k-1,\ldots,4k-2$. As in the proof of
Lemma~\ref{lem:coefficient_radius_green}, its contribution
to $z$ times the weighted early coherent sum is at most
\begin{equation*}
 C \cdot k^5 \cdot \mathcal G(k)^2(\mathcal G(k)^2/z)^{k-2}
 \le C \cdot k^6 \cdot \rho_k^{k-2} \cdot e^{-\epsilon(k-2)/2}.
\end{equation*}
Here the prefactor uses only the original cap $\mathcal G\le G_*$.
Degree two is treated directly by $c_2^2\le2$. The moment lemma
gives both early bounds in Eq.~\eqref{eq:coefficient_radius_early_terms}.
These are the only uses of the original condition on the coefficient floor.

The remaining Green calculation uses these early
bounds, $\eta,\eta_1\le C/z$,
$D_j\le z(b_{j+1}-2b_j+b_{j-1})+\ell_j$, the norm-attaining shift,
$R-z\le2$, and $\gamma m^2/z\le c_0$. The boundary of its second block
is $b_{2m}\le1$ and
$z(b_{2m+1}-b_{2m})\le(R-z)b_{2m}\le2$.
Thus the same weighted absorptions give the asserted orbit and actual
window bounds at a fixed large lower threshold for $z$.

The proof of Lemma~\ref{lem:coefficient_radius_initialization} now
applies to its same direct sum and compensated $F_2$ restart. It uses
only these bounds, $m\ge2d$, and its displayed coherent, length, and
prefix conditions. In particular it gives
\begin{equation*}
 \nu\le\frac{u}{1-C/(z^2 \cdot m)-C \cdot \gamma/z},\qquad
 P_J\le C/(z^2 \cdot m)+C \cdot \gamma J/z.
\end{equation*}
All denominator windows contain at least $m/4$ mass. The exact shift
$S_{k+1}\le S_k$ gives $v_{j+1}\le ze_j$, hence the same matched
budget. The construction uses the same actual family and output slots,
with the same subsequent normalizations.
\end{proof}

We extend the family by assigning a normalizer threshold to each degree
range. This keeps the propagation and terminal coherent term within the
margin supplied by the corresponding cap.

For odd $p\ge3$, choose dyadic degrees, a retained prefix, and an exit
parameter by
\begin{equation}\label{eq:nearunit_secondary_degrees}
 \begin{gathered}
 d_p:=\min\{2^j:2^j\ge4p\},\quad
 e_p:=\min\{2^j:2^j\ge8p\},\quad\omega_p:=1
                         ,\ \forall p\le K_\delta,\\
 d_p:=e_p:=\min\{2^j:2^j\ge C_{\mathrm{deg}} \cdot p/\delta\},\quad
 \omega_p:=\delta,\ \forall K_\delta<p\le K_\epsilon,\\
 d_p:=e_p:=\min\{2^j:2^j\ge C_{\mathrm{deg}} \cdot p/\epsilon\},\quad
 \omega_p:=\epsilon,\ \forall p>K_\epsilon,\qquad J_p:=4e_p.
 \end{gathered}
\end{equation}
In particular $J_p\ge32p$. In the middle range $e_p<\Xi$.
In both nonsmall ranges $d_p,e_p,J_p\le C \cdot p^{3/2}$, since
$\delta^{-1}<\sqrt p$ or $\epsilon^{-1}<\sqrt p$, respectively.

For an odd order $Q\ge3$, $z>0$, $u:=1/z$, and constants
$a_Q,s_Q>0$, the startup conditions for a centered unit family are
\begin{equation}\label{eq:nearunit_secondary_geometry_startup}
 \nu\le\frac{u}{1-a_Q \cdot u^Q},\qquad
 P_J\le s_Q \cdot u^Q,\ \forall 1\le J\le J_Q,\qquad
 aV_N+X_N\le\nu-u+uP_{N+1},\ \forall N\ge0,
\end{equation}
where $J_Q$ is defined in Eq.~\eqref{eq:nearunit_secondary_degrees}.

\begin{lemma}[Extension with three normalizer thresholds]
\label{lem:nearunit_secondary_extension}
Let $q\ge3$ be odd. Suppose
$F^2\ge\max\{(1+\delta)\mathcal H(\Xi)^2,\mathsf A_\delta\}$
at a sufficiently large universal $C_0$, and define
\begin{equation*}
 \kappa_q:=(F^2/4)^{q-1},\qquad
 \kappa_{q-2}:=(F^2/4)^{q-3},\qquad u:=1/z.
\end{equation*}
Assume $\mathcal G(e_q)\le F$, $0\le R-z\le2$, and
$\|\mathcal T^h\|_{\mathrm{op}}\le4(h+1)z^{h/2}$ for all $h\ge0$.
Suppose one centered unit family has the energy bound, actual losses
through $J_q$, and matched budget of
Eq.~\eqref{eq:nearunit_secondary_geometry_startup}, with constants $a_q,s_q$.
There is a universal $C_\theta$ such that, on defining
\begin{equation*}
 \mathsf K_q:=4(s_q+C_F \cdot \kappa_{q-2}+4C_\theta \cdot a_q+1),\qquad
 C_q:=a_q+4\mathsf K_q,
\end{equation*}
the conditions
\begin{equation}\label{eq:nearunit_secondary_extension_thresholds}
 z\ge\max\{(1+\epsilon) \cdot F^2,8,8C_\theta\},\qquad
 a_q \cdot u^q,\mathsf K_q \cdot u^q\le\omega_q/1024
\end{equation}
imply
\begin{equation}\label{eq:nearunit_secondary_extension}
 \|\mathcal T^m\|_{\mathrm{op}}^{2/m}
 \ge z-C_q \cdot z^{1-q}
 ,\ \forall 1\le m\le\lfloor z^q/\kappa_q\rfloor.
\end{equation}
Moreover $C_q\le C \cdot (a_q+s_q+\kappa_{q-2}+1)$ with a universal constant.
\end{lemma}
\begin{proof}
Up to the first exit of $p_*:=\mathsf K_q \cdot u^q$, the actual normalizers obey
\begin{equation*}
 B_*\le\frac{u}{(1-a_q \cdot u^q)(1-p_*)}
       \le\beta u,\qquad\beta:=1+\omega_q/128.
\end{equation*}
Thus $t_*:=\max\{1,B_* \cdot z\}\le\beta$, and the absolute gap gives
$\lambda\ge u/2$. Apply Lemma~\ref{lem:second_logstar_polynomial_loss}
with $d_0:=d:=e_q$. Its terminal term is
$\tau:=C_F \cdot d \cdot B_*^d \cdot c_d^2$, and its loss kernel is at most
\begin{equation*}
 Cu\sum_{\substack{k<e_q\\k\ \mathrm{dyadic}}}
              k^8(\beta^4 \cdot \mathcal G(k)^2/z)^{k-1}.
\end{equation*}
If $q\le K_\delta$, every tested degree is below $16K_\delta$,
so the floor supplied by the original cap gives ratio at most $1/C_0$. The fixed
$\beta=1+1/128$ is harmless. If $K_\delta<q\le K_\epsilon$, every
tested degree is below $\Xi$. Below $K_\delta$ use the original fixed
ratio; above it use
$\beta^4/(1+\delta)\le e^{-\delta/2}$ and the tail after $K_\delta$.
This gives a universal bound $C_\theta \cdot u$ in both cases.
If $q>K_\epsilon$, then $\beta=1+\epsilon/128$ and
$\beta^4/(1+\epsilon)\le e^{-\epsilon/2}$. Apply
Lemma~\ref{lem:nearunit_secondary_seed_moments} to
$\rho_k:=\min\{1,\mathcal H(k)^2/F^2\}$. The same universal
$C_\theta$ bounds this third kernel.

We define $d:=e_q$, $D:=1+\epsilon$, $g:=\mathcal G(d)^2$, and $\rho:=g/F^2$.
The exact terminal ratio can be written without a positive power of $F$:
\begin{align}\label{eq:nearunit_secondary_terminal_ratio}
 \frac{\tau}{\kappa_q \cdot \kappa_{q-2} \cdot u^{2q}}
 &\le C_F \cdot d \cdot F^6 \cdot 4^{2q-4} \cdot \rho^{d-1} \cdot \beta^d \cdot D^{2q-d}\notag\\
 &=C_F \cdot d \cdot g^3 \cdot 4^{2q-4} \cdot \rho^{d-4} \cdot \beta^d \cdot D^{2q-d}.
\end{align}
The original cap $\mathcal G\le G_*$ gives $g\le Cd$. In the small
range $8q\le d<16q$,
$\rho\le C_0^{-1}$, and $D^{2q-d}\le1$. The right side is at most
\begin{equation*}
 C \cdot q^4 \cdot C_0^4 \cdot (16\beta^{16}/C_0^8)^q,
\end{equation*}
which is uniformly small for $q\ge3$ at a fixed large $C_0$.
In the middle range $d\ge C_{\mathrm{deg}} \cdot q/\delta$,
$\rho\le(1+\delta)^{-1}$, and $\beta=1+\delta/128$.
Discard the negative power of $D$. The product
$\rho^{d-4} \cdot \beta^d$ is at most $C \cdot e^{-c\delta d}$.
The divisor powers are absorbed for fixed large $C_{\mathrm{deg}}$,
and $\delta^{-1}<\sqrt q$ gives a remaining bound
$C \cdot q^6 \cdot e^{-cC_{\mathrm{deg}}q}$. In the large range, use
$\rho\le1$, $\beta=1+\epsilon/128$, and the $D$ damping instead;
$\epsilon^{-1}<\sqrt q$ gives the identical bound. Both nonsmall
ranges begin above a fixed large cutoff. Thus all three cases give
$\tau\le C_F \cdot \kappa_q \cdot \kappa_{q-2} \cdot u^{2q}$.

The prefix contains $2d-1$. Apply the first-exit calculation of
Lemma~\ref{lem:slow_offset_extension} with $\omega_q$.
The coefficient bounds it uses are
$\theta\le C_\theta \cdot u$ and the terminal bound just established.
On $N\le\lfloor z^q/\kappa_q\rfloor$ it gives
\begin{equation*}
 (1-2C_\theta \cdot u)P_N
 \le s_q \cdot u^q+C_F \cdot \kappa_{q-2} \cdot u^q+4C_\theta \cdot a_q \cdot u^{q+1}
 \le\mathsf K_q \cdot u^q/4.
\end{equation*}
This excludes first exit on the actual finite horizon. The same
positive successor masses and raw quotient identities give successive
squared norm ratios at least $z-(a_q+4\mathsf K_q)u^{q-1}$.
Their multiplication proves the assertion. The fixed exit in the small
range suffices for positivity, while its stronger coefficient ratio
already controlled the kernel.
\end{proof}

The extension can now be iterated to the prescribed block order. The
following estimate preserves the inherited coefficient floor and makes
its dependence on the two loss parameters explicit.

\begin{proposition}[Blocks with an inherited floor]
\label{prop:nearunit_secondary_blocks}
There are universal constants $C_X,C_0,C>0$ with the following property.
Let $\mathcal G,\mathcal H$ satisfy
Eqs.~\eqref{eq:second_logstar_majorant} and
\eqref{eq:nearunit_secondary_cap} uniformly over the stated regular laws.
For $0<\epsilon\le\delta\le1/16$ and odd $Q\ge3$, we define
\begin{equation}\label{eq:nearunit_secondary_geometry_parameters}
 \begin{gathered}
 X:=C_X \cdot Q^2,\qquad
 F^2:=\max\{\mathcal G(X^X)^2,(1+\delta)\mathcal H(\Xi)^2,\mathsf A_\delta\},\\
 H_Q:=(1+2\epsilon) \cdot F^2.
 \end{gathered}
\end{equation}
There is an assignment $Z_Q$ and constants $a_Q,s_Q>0$ satisfying
\begin{align}
 H_Q&\le Z_Q\le\max\{H_Q,R\},\qquad\mathcal A\le Z_Q,\notag\\
 c_d&\le Z_Q^{(d-1)/2},\ \forall d\ge1,\qquad
 a_Q+s_Q+1\le C \cdot (F^2/4)^{Q-3}\le F^{2Q}.
 \label{eq:nearunit_secondary_geometry_static}
\end{align}
If $Z_Q>H_Q$, let $z:=Z_Q$ and $u:=1/z$. One centered unit finite
starting family satisfies the startup conditions
in Eq.~\eqref{eq:nearunit_secondary_geometry_startup}, and
\begin{equation}\label{eq:nearunit_secondary_geometry_propagation}
 \|\mathcal T^h\|_{\mathrm{op}}\le4(h+1)z^{h/2},\ \forall h\ge0,\qquad
 0\le R-z\le\min\{1,C/R\},\qquad a_Q \cdot u^Q\le1/8.
\end{equation}
Moreover,
\begin{equation}\label{eq:nearunit_secondary_absolute_floor}
 \mathsf A_\delta\le C[1+\mathfrak t(\delta^{-1})]^{1/3}
                         \le C \cdot \delta^{-1/3},\qquad
 \Xi\le C \cdot \delta^{-1} \cdot \epsilon^{-2}.
\end{equation}
All constants are independent of the order, both majorants, $\delta$,
and $\epsilon$.
\end{proposition}
\begin{proof}
Construct first in the branch $\mathcal A>H_Q$, so $R\ge\mathcal A$.
Use $\kappa_q:=(F^2/4)^{q-1}$ and Eq.~\eqref{eq:nearunit_secondary_degrees}.
The degree bounds following Eq.~\eqref{eq:nearunit_secondary_degrees} place
every tested degree below one universal $X:=C_X \cdot Q^2$, with seed at most
$F$. This bound on the tested degrees is independent of $\delta$ and $\epsilon$.

For odd $p\ge5$, suppose $z\ge DF^2$, $D:=1+\epsilon$, and $m$
is within a factor four of $z^{p-2}/\kappa_{p-2}$. Write $d:=d_p$,
$g:=\mathcal G(d)^2$ and $\rho:=g/F^2$. The two coherent conditions for initialization are bounded by
\begin{align*}
 \gamma m^2/z
 &\le Cd \cdot F^2 \cdot \rho^{d-1} \cdot 4^{2p-6} \cdot D^{2p-d-4}
   =C \cdot d \cdot g \cdot \rho^{d-2} \cdot 4^{2p-6} \cdot D^{2p-d-4},\\
 J_p \cdot \gamma zm
 &\le C \cdot J_p \cdot d \cdot F^4 \cdot \rho^{d-1} \cdot 4^{p-3} \cdot D^{p-d}
   =C \cdot J_p \cdot d \cdot g^2 \cdot \rho^{d-3} \cdot 4^{p-3} \cdot D^{p-d}.
\end{align*}
The original cap $\mathcal G\le G_*$ gives $g\le Cd$. In the small
range, $d\ge4p$,
$J_p\le64p$, $\rho\le C_0^{-1}$, and both $D$ exponents are
negative. The two bounds are at most
\begin{equation*}
 C \cdot p^2 \cdot C_0^2 \cdot (16/C_0^4)^p,\qquad
 C \cdot p^4 \cdot C_0^3 \cdot (4/C_0^4)^p.
\end{equation*}
For $p\ge5$ these are uniformly small at fixed large $C_0$.
In the middle range, $\rho\le(1+\delta)^{-1}$ and
$d,J_p=O(p/\delta)$. Discard the negative $D$ powers.
The $\rho$ decay absorbs the divisor powers, leaving
$C \cdot p^3 \cdot e^{-cC_{\mathrm{deg}}p}$ and $C \cdot p^6 \cdot e^{-cC_{\mathrm{deg}}p}$,
because $\delta^{-1}<\sqrt p$. In the large range use the $D$ decay
and $\epsilon^{-1}<\sqrt p$, obtaining the same bounds.
All three cases therefore satisfy both coherent conditions universally.
The length lower bound is $m\ge cF^2 \cdot 4^{p-3} \cdot D^{p-2}$. Thus
$J_p/m\le Cp \cdot 4^{-p}/F^2$ in the small range, and
$J_p/m\le C \cdot p^{3/2} \cdot 4^{-p}/F^2$ otherwise. These give
$m\ge2d_p$ and $J_p\le m/4$. Also
$\kappa_{p-2} \cdot z^{-p}\le F^{-6} \cdot 4^{3-p} \cdot D^{-p}$ is uniformly small.
With the envelope established below, the initialization lemma gives
$a_p+s_p+1\le C \cdot \kappa_{p-2}$ through the full prefix.

At the base $m_3:=\lfloor R\rfloor$, the small degree rule always applies:
$d_3=16,e_3=32,J_3=128$, and $\kappa_1=1$.
The mesoscopic estimate gives $R-z_3\le C_*/R$ and global remainder
factor at most two. Its two coherent terms are $O(z_3^{-14})$ and
$O(J_3 \cdot z_3^{-13})$. The lower radius is justified below, after which
Lemma~\ref{lem:nearunit_secondary_green} supplies the orbit estimate
and initialization at this fixed base.

For each constructed order define
\begin{equation*}
 m_{q+2}:=\lfloor z_q^q/\kappa_q\rfloor,\qquad
 z_{q+2}:=\|\mathcal T^{m_{q+2}}\|_{\mathrm{op}}^{2/m_{q+2}},\qquad
 \Delta_q:=C_b \cdot \kappa_{q-2} \cdot z_q^{1-q}.
\end{equation*}
The exit thresholds follow from
$(a_q+\mathsf K_q)u^q\le C \cdot F^{-6} \cdot 4^{3-q} \cdot D^{-q}$.
This is uniformly small in the small range. In the middle and large ranges,
multiply respectively by $\delta^{-1}<\sqrt q$ and
$\epsilon^{-1}<\sqrt q$. The result is still uniformly small.
Lemma~\ref{lem:nearunit_secondary_extension} therefore supplies every bridge.
In addition,
\begin{equation*}
 \Delta_q/z_q\le C_b \cdot F^{-6} \cdot 4^{3-q} \cdot D^{-q},\qquad
 q\Delta_q/z_q\le C \cdot F^{-6}.
\end{equation*}

Each proposed length is at most $R^Q$. If $R\le X$, its first $m$
coefficients have seed $\mathcal G(m)^2\le F^2<\mathcal A$.
If $X<R<2\mathcal A$, then $Q<R$ and the original cap gives
$\mathcal G(m)^2\le G_*(R^R)^2\le R/2<\mathcal A$ after a fixed
enlargement of $C_X$. The exact block maximum in
Eq.~\eqref{eq:coefficient_radius_block_maximum} forces $z_m\ge\mathcal A$
in both cases. If $R\ge2\mathcal A$, the mesoscopic gap gives this
at the base. At each new step the old decrement tail and current bridge
give $z_{q+2}\ge R-C_*/R-2\Delta_3\ge\mathcal A$, since
$\Delta_3=C_b/z_3^2=O(R^{-2})$. Only old ratios enter that tail.
Thus every new actual radius is at least $\mathcal A>H_Q$ before
its new ratio is evaluated.

The exact successive ratio, including increasing actual radii, is
\begin{equation}\label{eq:nearunit_secondary_gap_ratio}
 \frac{\Delta_{q+2}}{\Delta_q}
 =\frac{F^4}{16z_{q+2}^2}(z_q/z_{q+2})^{q-1}
 \le\frac2{16}<\frac12.
\end{equation}
For a downward step use $z_{q+2}\ge z_q-\Delta_q$ and
$q\Delta_q/z_q\le C \cdot F^{-6}$; for an upward step the second factor
is at most one. Hence all finite constructed prefixes satisfy
\begin{equation}\label{eq:nearunit_secondary_gap}
 0\le R-z_q\le C_*/R+C'/R^2\le C_g/R\le1.
\end{equation}
This provides the absolute gap required by the Green lemma. The lower
radius comes from the common coefficient radius, so no reserve of size
$\epsilon F^2$ is used to pay this gap. Fix the degree and moment constants, then the constants for initialization, extension,
and the total gap, and finally $C_0,C_X$ above all required
fixed thresholds.

The defining real lengths are at least $F^2 \cdot 4^{q-1} \cdot D^q$. Their integer
floor losses are uniformly small. For $q\ge5$, the exact lower length
ratio has factor $16z_q^2/F^4\ge16D^2$, while the correction for a decrease in radius is $e^{-O(F^{-6})}$. Thus $m_{q+2}/m_q\ge4$.
At the distinct first step the ratio is at least
$(1-o(1))16(z_3/F^2)^2(z_3/R)\ge4$, using the mesoscopic relative
gap $O(R^{-2})$.

The mixed-radix proof of Proposition~\ref{prop:slow_offset_blocks}
now applies with the same bounds: block growth at least four,
the geometric decrement tail, and per-level logarithmic excess
\begin{equation*}
 C \cdot m_{q+2}\Delta_q/z_q
 \le C \cdot \kappa_{q-2}/\kappa_q=16C/F^4.
\end{equation*}
The base remainder costs $C/R$. Counting at most
$C \cdot \log(h+1)/\log4$ levels gives the envelope $4(h+1)$.
Counting the actual finite prefix, at most $q/2$ levels, instead gives
$\log E_q\le C/R+Cq/F^4$. Therefore, at the proposed next length,
\begin{equation*}
 q(\log(z_{q+2}/z_q))_+
 \le\frac{Cq \cdot 4^{-q}}{F^2 \cdot R}+\frac{C \cdot q^2 \cdot 4^{-q}}{F^6}
 \le\frac C{F^2 \cdot R}+\frac C{F^6}.
\end{equation*}
Together with the downward gap and integer floor, this gives the
comparison within a factor of four with $z_{q+2}^q/\kappa_q$. At this point the
new actual radius and its envelope are both known; only then apply
Lemma~\ref{lem:nearunit_secondary_green} with the coherent and prefix
conditions checked above. Its actual restart closes the finite
induction. Thus the polynomial envelope is established before
the Green estimate is applied.

At the final order $a_Q+s_Q+1\le C \cdot (F^2/4)^{Q-3}\le F^{2Q}$,
including $Q=3$. In the nontrivial branch set $Z_Q:=z_Q\ge\mathcal A$.
If instead $\mathcal A\le H_Q$, set $Z_Q:=H_Q$ and use the same
upper bounds for $a_Q,s_Q$. The coefficient inequalities follow from
the definition of $\mathcal A$. The assertions about an actual family
and a comparison of $R$ with $H_Q$ apply only to the nontrivial branch.
Finally the original cap and the fixed-power height estimate give
$\mathsf A_\delta\le C[1+\mathfrak t(\delta^{-1})]^{1/3}$.
The elementary bound on $\log^*$ then gives
Eq.~\eqref{eq:nearunit_secondary_absolute_floor}. This proves the proposition.
\end{proof}

\subsection{Making changes in height arguments summable}
\label{sec:nearunit_secondary_profiles}

The scalar arguments change when a profile is extracted and returned,
so their contribution must also remain bounded over repeated refinements.
We replace the discrete height functions by averaged majorants whose
contraction makes the costs of fixed powers, sums, and products decrease
geometrically with the refinement index. These estimates will control
the additive errors in the profile induction.

We use $\mathfrak t$, $\kappa$, and $\chi$ from
Eq.~\eqref{eq:secondary_feedback_heights}. We introduce a second loss
parameter and include its height in a maximum. Smoothing the scalar
heights makes their argument costs decrease with the repetition count.
The additive offsets will be bounded uniformly. The coefficient seed
is fixed throughout these outer rounds, so its normalization is paid
only once per profile refinement.

Fix $b_*:=L_*:=4$ and write $\vartheta:=b_*^{-1}$. For $f=\mathfrak t,\chi$, set
\begin{equation}\label{eq:smoothed_secondary_heights}
 \overline f(x):=\max\{L_*,\frac1{b_*}
                            \int_0^{b_*}f(x+t)dt\},\ \forall x\ge1,
 \qquad
 \widetilde V_m(x):=\overline{\mathfrak t}
                         (\overline\chi^{\circ m}(x)).
\end{equation}
The bars below denote these averages of the scalar height functions.

\begin{lemma}[Contractive majorants]
\label{lem:smoothed_secondary_maps}
Both averaged maps are
nondecreasing and globally $\vartheta$-Lipschitz on $[1,\infty)$,
\begin{equation}\label{eq:smoothed_secondary_map_bounds}
 \max\{L_*,f(x)\}\le\overline f(x)\le f(x)+1,
 \forall f\in\{\mathfrak t,\chi\},
\end{equation}
and
\begin{equation}\label{eq:smoothed_secondary_fixed_point}
 \overline{\mathfrak t}(L_*)=L_*,\qquad
 \overline\chi(x)=L_*,\ \forall L_*\le x\le L_*+2.
\end{equation}
Consequently $\widetilde V_m\ge L_*$, $\widetilde V_{m+1}\le\widetilde V_m
\le\overline{\mathfrak t}$, and $\widetilde V_m$ is globally
$\vartheta^{m+1}$-Lipschitz.
\end{lemma}
\begin{proof}
For $x\ge1$, the first jump of $\mathfrak t$ occurs at
$e^e-2>13$, and consecutive later jumps are more than four apart.
Each jump has size one. For $x>3$,
$\kappa(x)=1+\kappa(\mathfrak t(x))$ and integer nonexpansiveness
of $\kappa$ imply that every jump of $\chi$ is a jump of
$\mathfrak t$ and has size at most one. The possible jump of
$\kappa$ at three is hidden by $\chi=\max\{3,1+\kappa\}$.
Thus for either function $f$ one has
$0\le f(x+4)-f(x)\le1$ on $[1,\infty)$.
Its forward average is absolutely continuous, with derivative almost everywhere equal to $[f(x+4)-f(x)]/4$ in $[0,1/4]$.
Taking the maximum with four preserves this Lipschitz bound.
Both original functions are at least three. Forward averaging lies
between $f(x)$ and $f(x)+1$, and the floor four is also at most
$f(x)+1$, giving
Eq.~\eqref{eq:smoothed_secondary_map_bounds}.

Since $\mathfrak t(8)=3$ and $\chi(10)=3$, the fixed-point assertions
follow. For $x>4$, contraction at four gives
$\overline\chi(x)-4\le\vartheta(x-4)<x-4$.
For $x\le4$ both averaged maps equal four. Monotone iteration and
composition of the global Lipschitz bounds prove the last assertions.
\end{proof}

The contraction turns fixed-power, sum, and product changes in the
argument into errors that decrease with the profile index. We record
the estimates used when choosing the next profile's parameters.

\begin{lemma}[Decreasing argument costs]
\label{lem:smoothed_secondary_costs}
For all $x,y\ge1$ and integers $m\ge0$, the following estimates hold.
For each fixed $C,p\ge1$ there is $b_{C,p}$, independent of $m$, with
\begin{equation}\label{eq:smoothed_secondary_power_cost}
 \widetilde V_m(C \cdot x^p)\le\widetilde V_m(x)+b_{C,p} \cdot \vartheta^m.
\end{equation}
There is a universal $b$ such that
\begin{equation}\label{eq:smoothed_secondary_product_cost}
 \max\{\widetilde V_m(xy),\widetilde V_m(x+y)\}
 \le\max\{\widetilde V_m(x),\widetilde V_m(y)\}+b\vartheta^m.
\end{equation}
For every fixed $C,p,a\ge1$ and $0<\delta\le1/16$,
\begin{equation}\label{eq:smoothed_secondary_parameter_cost}
 \widetilde V_m(C \cdot \delta^{-p} \cdot x^a)
 \le\max\{\widetilde V_m(x),\widetilde V_m(\delta^{-1})\}
                                      +b_{C,p,a} \cdot \vartheta^m.
\end{equation}
Adding a fixed polynomial in $\delta^{-1}$ to the argument has the
same effect, with a different fixed constant. All assertions remain
valid after replacing $\widetilde V_m$ by $\max\{\widetilde V_m,h\}$
for any fixed $h\ge L_*$.
\end{lemma}
\begin{proof}
The proof of Lemma~\ref{lem:secondary_feedback_heights} gives the
fixed-power estimate for both $\mathfrak t$ and $\chi$. The forward
average and Eq.~\eqref{eq:smoothed_secondary_map_bounds} therefore
give $\overline f(C \cdot x^p)\le\overline f(x)+b_{C,p}$ as well.
For $m\ge1$, apply this estimate to the first $\overline\chi$.
The remaining $m-1$ copies of $\overline\chi$ and the final
$\overline{\mathfrak t}$ contract its additive error by
$\vartheta^m$. The case $m=0$ is the averaged fixed-power estimate.
This proves Eq.~\eqref{eq:smoothed_secondary_power_cost}.
For the product and sum, use $xy\le\max\{x,y\}^2$ and
$x+y\le2\max\{x,y\}$ before the same first-map argument.
Applying these estimates to the fixed powers of $x$ and
$\delta^{-1}$ proves the parameter assertion. Maxima with a constant
preserve all the stated inequalities.
\end{proof}

Smoothing must preserve the finite stopping property of the discrete
heights. The next lemma shows that a finite number of profile refinements
still brings the scalar majorant to its constant floor.

\begin{lemma}[Finite termination]
\label{lem:smoothed_secondary_termination}
\label{lem:smoothed_secondary_heights}
For every $x\ge1$,
\begin{equation}\label{eq:smoothed_secondary_termination}
 \overline\chi^{\circ(N_\chi(x)+1)}(x)=L_*,\qquad
 \widetilde V_{N_\chi(x)+1}(x)=L_*.
\end{equation}
\end{lemma}
\begin{proof}
We define $\chi_L:=\max\{L_*,\chi\}$. Since $\chi(L_*)\le L_*$,
$\chi_L^{\circ j}(x)=\max\{L_*,\chi^{\circ j}(x)\}$ for $j\ge1$.
Let $x_j:=\overline\chi^{\circ j}(x)$ and
$y_j:=\chi_L^{\circ j}(x)$. Monotonicity and
Eq.~\eqref{eq:smoothed_secondary_map_bounds} give $x_j\ge y_j$.
Moreover,
\begin{equation*}
 x_{j+1}-y_{j+1}
 \le\vartheta(x_j-y_j)+1,
 \qquad0\le x_j-y_j\le\frac1{1-\vartheta}<2.
\end{equation*}
At $j=N_\chi(x)>0$, one has $y_j=L_*$. One further application
of Eq.~\eqref{eq:smoothed_secondary_fixed_point} gives $x_{j+1}=L_*$.
If $N_\chi(x)=0$, then $x\le3$ and the first application already
equals $L_*$. Apply $\overline{\mathfrak t}(L_*)=L_*$.
\end{proof}

We now allocate an additive reserve for the argument changes. The
following parameter choices make the reserve increments summable and
prepare the common domain used by the profile induction.

The constants $C_b,C_R,C_Q$ for repeated height reduction and the terminal-degree constant
$C_T$ are fixed first. For $0<\delta\le1/16$, we define
\begin{equation}\label{eq:nearunit_secondary_profile_parameters}
 \begin{gathered}
 \widehat V_{m,\delta}(x):=\max\{\widetilde V_m(x),\overline{\mathfrak t}(\delta^{-1})\},\\
 R_\delta:=\lceil C_{R,0} \cdot \delta^{-12}\rceil,\qquad
 S_m:=1+2D_{\mathrm{add}} \cdot \sum_{i=0}^{m-1}4^{-i}.
 \end{gathered}
\end{equation}
We choose $C_{R,0}$ and then $D_{\mathrm{add}}$ below, independently
of the repetition count. The sum is empty at $m=0$.
All constants in this subsection are universal.

\begin{lemma}[Additive parameter costs]
\label{lem:nearunit_secondary_scalar_costs}
For each fixed $C,p,a\ge1$, there is a constant $b_{C,p,a}$,
independent of $m$ and $\delta$, such that
\begin{equation}\label{eq:nearunit_secondary_scalar_costs}
 \widehat V_{m,\delta}(C \cdot \delta^{-p} \cdot x^a)
       \le\widehat V_{m,\delta}(x)+b_{C,p,a} \cdot 4^{-m}
 ,\ \forall x\ge1,m\ge0.
\end{equation}
The same estimate holds after adding a fixed polynomial in
$\delta^{-1}$ to the argument, with a different fixed constant.
\end{lemma}
\begin{proof}
Apply the product and fixed-power estimates of
Lemma~\ref{lem:smoothed_secondary_costs} separately to
$C$, $\delta^{-p}$, and $x^a$. The result is bounded by
$\max\{\widetilde V_m(x),\widetilde V_m(\delta^{-1})\}+b_{C,p,a} \cdot 4^{-m}$.
Use $\widetilde V_m(\delta^{-1})\le\overline{\mathfrak t}(\delta^{-1})$ and take the
maximum with $\overline{\mathfrak t}(\delta^{-1})$. The sum estimate in that
lemma proves the second assertion. Thus the parameter contributes only
an additive error of order $4^{-m}$ to the enlarged height.
\end{proof}

Only finitely many powers occur below. Choose the fixed constant
$D_{\mathrm{add}}$ above all their combined constants. Then
\begin{equation}\label{eq:nearunit_secondary_offset_reserve}
 S_{m+1}-S_m=2D_{\mathrm{add}} \cdot 4^{-m},\qquad
 1\le S_m\le1+\frac83D_{\mathrm{add}}.
\end{equation}
One increment of $D_{\mathrm{add}} \cdot 4^{-m}$ covers changes in the
argument of the fixed coefficient bound; the other covers the final
passage to the next height profile.

\subsection{Refining the profile with a small multiplicative loss}

We combine the inherited block floor with the contractive height
estimates to advance the prescribed-depth profile by one refinement,
with a multiplicative factor arbitrarily close to one.
The inherited coefficient bound remains fixed throughout the outer
rounds, and choosing their number returns the improved profile to the
same range of depths for each fixed loss parameter. Keeping both loss
parameters explicit records the cost of this advance.

\begin{proposition}[Repeated height reduction with a fixed offset]
\label{prop:nearunit_secondary_restricted_profile}
Suppose, for some $m\ge0$ and $A\ge A_0$, that uniformly over centered
regular laws of covariance at most $I$ and curvature at least $aI$,
\begin{equation}\label{eq:nearunit_secondary_prescribed_input}
 \mathcal A\le A \cdot [\widehat V_{m,\delta}(r)+S_m]^{1/3} \cdot \mathcal L_r(a^{-1})^2
 ,\ \forall r\ge R_\delta,\ r\text{ an integer}.
\end{equation}
For $0<\epsilon\le\delta$, define
\begin{align}
 D&:=\widehat V_{m,\delta}(\epsilon^{-1})+S_m+D_{\mathrm{add}} \cdot 4^{-m},\notag\\
 W_j(x)&:=\max\{\mathfrak t_j(x),D\}^{1/3},
                      \label{eq:nearunit_secondary_primary_weight}\\
 Q_\epsilon&:=\text{the least odd integer at least }C_Q \cdot \epsilon^{-4},
 \notag\\
 r_{Q,\epsilon,\delta}&:=\max\{R_\delta,
              \lceil C_b \cdot \mathfrak t(Q)\rceil+\lceil C_R/\epsilon\rceil\}.
                      \label{eq:nearunit_secondary_primary_depth}
\end{align}
Then every integer $j\ge1$, odd $Q\ge Q_\epsilon$, and integer
$r\ge r_{Q,\epsilon,\delta}$ satisfy
\begin{equation}\label{eq:nearunit_secondary_primary_profile}
 \mathcal A\le e^{16\delta+24\epsilon j} \cdot A \cdot W_j(Q) \cdot
                           (r+1)^{1/Q} \cdot \mathcal L_r(a^{-1})^2.
\end{equation}
The constants $A_0,C_{R,0},D_{\mathrm{add}}$ can be fixed independently of
$m,A,\delta,\epsilon$, and the law. The constants $C_b,C_R,C_Q,C_T$
are those fixed before Eq.~\eqref{eq:nearunit_secondary_profile_parameters}.
\end{proposition}
\begin{proof}
We first extract the fixed majorant. We define
\begin{equation}\label{eq:nearunit_secondary_extraction_depth}
 r_\delta(x):=\max\{R_\delta,
                         \lceil C_r \cdot \mathfrak t(x)+D_r/\delta\rceil\}.
\end{equation}
Choose $C_r>1$ and $D_r$ universally so that scalar burn-in and
contraction give
\begin{equation*}
 \mathcal L_{r_\delta(d)}(d)^2\le e^\delta,\qquad
 (1+r_\delta(d)^{-2})^2\le e^\delta.
\end{equation*}
Lemma~\ref{lem:nearunit_static_transfer} with the common assignment
$\mathcal A$ and floor one therefore yields
\begin{equation*}
 c_d^2\le[e^{3\delta} \cdot A
                  (\widehat V_{m,\delta}(r_\delta(d))+S_m)^{1/3}]^{d-1}.
\end{equation*}
The normalized floor follows from the fixed large $A_0$ and
$\mathcal L_r(0)^2\ge e^{-\delta}$ for $r\ge R_\delta$.
Every depth used here lies in the range of the assumed bound. The transfer gives the
finite-degree estimate on the full coefficient class without any
continuity assumption on the supremum defining $\mathcal A$.
Thus the valid nondecreasing fixed cap is
\begin{equation}\label{eq:nearunit_secondary_fixed_seed}
 \mathcal H(x)^2:=\min\{G_*(x)^2,
             e^{3\delta} \cdot A \cdot [\widehat V_{m,\delta}(r_\delta(x))+S_m]^{1/3}\}.
\end{equation}
It depends on the profile amplitude at the start of this step, and not on the new
outer index $j$.

Set $A_{\mathrm{base}}:=e^{16\delta} \cdot A$. Suppose the preceding outer
round has amplitude $A_j:=A_{\mathrm{base}} \cdot e^{24\epsilon j}$.
We extract coefficients at the chosen order, now with
\begin{equation}\label{eq:nearunit_secondary_selected_order}
 \begin{split}
 r_\epsilon(x)&:=\max\{R_\delta,
               \lceil C_r \cdot \mathfrak t(x)\rceil+\lceil D_R/\epsilon\rceil\},\\
 k_\epsilon(x)&:=\text{the least odd integer at least }
       \max\{Q_\epsilon,\epsilon^{-1} \cdot \log(r_\epsilon(x)+1)\}.
 \end{split}
\end{equation}
The depth verification in
the proof of Proposition~\ref{prop:secondary_feedback_restricted_profile}
applies:
both entries defining $k_\epsilon$ are at most $C \cdot r_\epsilon^4$,
so $r_\epsilon(x)\ge r_{k_\epsilon(x),\epsilon,\delta}$ once
$C_{R,0}$ is large. The normalized transfer gives the
$e^{5\epsilon}$ extracted seed, and the working squared majorant is
\begin{equation}\label{eq:nearunit_secondary_working_seed}
 \mathfrak g(x):=\min\{\mathcal H(x)^2,
                             e^{5\epsilon} \cdot A_j \cdot W_j(k_\epsilon(x))\}.
\end{equation}
Every order used in this extraction is at least $Q_\epsilon$.

Write $h:=\mathfrak t(Q)$ and $M_Q:=A_j \cdot W_{j+1}(Q)$.
Use the high-height threshold
\begin{equation}\label{eq:nearunit_secondary_height_threshold}
 h\ge C_h \cdot R_\delta \cdot \epsilon^{-5}.
\end{equation}
The argument at the chosen order in
Eq.~\eqref{eq:secondary_feedback_seed_bounds} then applies verbatim:
for $X:=C_X \cdot Q^2$, $Y:=C_k \cdot (Qs)^2$, and $y:=\log(s)/Q$, one has
$k_\epsilon(X^X)\le h$, $k_\epsilon(Y)\le h$ for $y\le1$,
and $k_\epsilon(Y)\le h(1+y)$ for $y\ge1$.
For $y\ge1$, the maximum satisfies
\begin{equation*}
 W_j(k_\epsilon(Y))^3
 \le\max\{\mathfrak t_j(h)+5y,D\}
 \le\max\{\mathfrak t_j(h),D\}+5y.
\end{equation*}
The last maximum is at least five, so its cube root divided by $e^y$
is at most $W_{j+1}(Q)$. For $y\le1$ there is no height increase.
In particular this leaves no fixed additive error in $\mathfrak t_{j+1}(Q)$.

If Eq.~\eqref{eq:nearunit_secondary_height_threshold} fails, then
$\mathfrak t(Y)\le h+B+y$ and
$r_\delta(Y)\le C[R_\delta \cdot \epsilon^{-5}+y]$.
Remove the fixed outside dilation by
Lemma~\ref{lem:smoothed_secondary_costs}, before treating the sum.
The remaining argument estimates and the global $4^{-m-1}$ Lipschitz
bound give
\begin{equation}\label{eq:nearunit_secondary_low_argument}
 \widehat V_{m,\delta}(r_\delta(Y))\le \widehat V_{m,\delta}(\epsilon^{-1})+D_{\mathrm{add}} \cdot 4^{-m}+4^{-m-1} \cdot y.
\end{equation}
The analogous estimate for $X^X$ has no $y$ term.
The definition of $D$ prepays this additive error. Since $D\ge5$,
\begin{equation*}
 \sup_{y\ge0}e^{-y} \cdot [D+4^{-m-1} \cdot y]^{1/3}\le D^{1/3}.
\end{equation*}
Thus the fixed cap
contributes at most $e^{4\delta} \cdot A \cdot D^{1/3}$, which is below
$A_{\mathrm{base}} \cdot D^{1/3}$. We have proved at all heights
\begin{equation}\label{eq:nearunit_secondary_seed_bounds}
 \mathfrak g(X^X)\le e^{5\epsilon} \cdot M_Q,\qquad
 \sup_{s\ge1}\frac{\mathfrak g(C_k \cdot (Qs)^2)}{s^{1/Q}}
                                  \le e^{5\epsilon} \cdot M_Q.
\end{equation}
The coherent-root estimate in
Eq.~\eqref{eq:slow_offset_coherent_root} uses only monotonicity and
the original cap; it remains valid with the same universal constants
and $N:=\lceil D_1 \cdot \epsilon^{-2} \cdot Q^2\rceil\le X^X$.

Initialization at $r:=r_{Q,\epsilon,\delta}$ also follows the same
high-height argument, using the preceding stage at
$k\ge\max\{Q_\epsilon,\epsilon^{-1} \cdot \log(r+1)\}$ with $k\le h\le Q$.
At low height, $r\le C \cdot R_\delta \cdot \epsilon^{-5}$. Apply the assumed bound in
Eq.~\eqref{eq:nearunit_secondary_prescribed_input} at this depth.
By the additive parameter estimate it gives
\begin{equation}\label{eq:nearunit_secondary_initial_bound}
 \mathcal A\le A \cdot D^{1/3} \cdot \mathcal L_r(a^{-1})^2
       \le M_Q \cdot \mathcal L_r(a^{-1})^2.
\end{equation}
The high-height bound is at most $e^\epsilon \cdot M_Q \cdot \mathcal L_r(a^{-1})^2$.
No profile below $Q_\epsilon$ is used in the low range.

Apply Proposition~\ref{prop:nearunit_secondary_blocks} with the fixed
cap in Eq.~\eqref{eq:nearunit_secondary_fixed_seed} and
$\mathcal G:=\sqrt{\mathfrak g}$. Its argument satisfies
$\Xi\le C \cdot \delta^{-1} \cdot \epsilon^{-2}$, and $\epsilon\le\delta$.
The same additive estimate gives
\begin{equation*}
 (1+\delta)\mathcal H(\Xi)^2
       \le e^{5\delta} \cdot A \cdot D^{1/3}
       \le M_Q.
\end{equation*}
The separate absolute floor satisfies
$\mathsf A_\delta\le C \cdot \mathfrak t(\delta^{-1})^{1/3}$, which is
at most $A_0 \cdot D^{1/3}$
after the fixed choice of $A_0$, since
$D\ge\overline{\mathfrak t}(\delta^{-1})\ge\mathfrak t(\delta^{-1})$.
The exact maximum in that proposition therefore yields
\begin{equation}\label{eq:nearunit_secondary_profile_floor}
 F^2\le e^{6\epsilon} \cdot M_Q,\qquad H_Q\le e^{9\epsilon} \cdot M_Q,
 \qquad(a_Q+s_Q+1)^{1/Q}\le e^{6\epsilon} \cdot M_Q.
\end{equation}
The $\delta$ cost is paid only in $A_{\mathrm{base}}$, independently
of the number of outer rounds.

We verify the kernel estimate without a fixed multiplicative increase
in that amplitude. Use the parameters for the inner depth induction in
Eq.~\eqref{eq:slow_offset_profile_parameters}, renaming its local
depth parameter $\delta=1/(64Qs)$ as $\xi:=1/(64Qs)$ to distinguish
it from the fixed loss parameter $\delta$. Suppose the contradictory actual
radius has inverse $u$. Set
$L_Q:=e^{20\epsilon} \cdot M_Q$ as before. Then
$u\le L_Q^{-1} \cdot s^{-1/Q}$ and $b^4\le e^\epsilon$.
Let $l_\epsilon:=\lceil C_l \cdot \epsilon^{-1} \cdot \log(e/\epsilon)\rceil$.
Partition the kernel below its finite cutoff $k_0$ into three ranges.
For $k\le16K_\delta$, the original cap and the absolute geometric
floor give $b^4 \cdot u \cdot \mathfrak g(k)\le b^4/C_0\le1/4$.
The degree-one term is bounded directly using $c_1\le1$.
For $16K_\delta<k\le l_\epsilon$, if this range is nonempty,
Eq.~\eqref{eq:nearunit_secondary_fixed_seed} and the additive estimate
bound the fixed cap by
$e^{4\delta} \cdot A \cdot D^{1/3}$.
The baseline $e^{16\delta} \cdot A$ consequently gives
$b^4 \cdot u \cdot \mathfrak g(k)\le e^{-\delta}$. Its moment is at most
\begin{equation*}
 Cu\sum_{k>16K_\delta}k^8 \cdot e^{-\delta(k-1)}
       \le Cu\delta^{-9} \cdot e^{-c\delta K_\delta}\le C'u.
\end{equation*}
For $k>\max\{16K_\delta,l_\epsilon\}$, still below $k_0$,
Eq.~\eqref{eq:nearunit_secondary_seed_bounds} gives the usual
$e^{-\epsilon}$ ratio. Its moment is at most
$Cu\epsilon^{-9} \cdot e^{-\epsilon l_\epsilon/2}\le C'u$.
If a range is empty it contributes zero. The profile tail above $k_0$
is bounded by the argument preceding
Eq.~\eqref{eq:slow_offset_kernel_split}. Hence
\begin{equation}\label{eq:nearunit_secondary_profile_kernel}
 \theta_d\le C_\theta \cdot u
\end{equation}
with a universal $C_\theta$, independent of both parameters.

The coherent-root estimate, Eq.~\eqref{eq:nearunit_secondary_profile_floor},
and this universal kernel establish exactly the amplitude conditions in
Eq.~\eqref{eq:slow_offset_amplitude_conditions}, with the same $C_Q$.
The normalized floor, actual retained prefix, matched budget, finite
first-exit comparison, and terminal-degree proof in
Proposition~\ref{prop:slow_offset_outer_profiles} therefore apply
without further change. Its distortion product is at most $e^\epsilon$
from the initial depth in Eq.~\eqref{eq:nearunit_secondary_primary_depth}.
Initialize by Eq.~\eqref{eq:nearunit_secondary_initial_bound} and the
high-height initialization. The next full depth profile has amplitude
at most $e^{21\epsilon} \cdot A_j$, covered by $e^{24\epsilon} \cdot A_j$.
In particular $C_T,C_b,C_R,C_Q$ acquire no dependence on $\delta$.
The outer induction propagates the profile for the coefficient radius directly.

For the first outer round,
Proposition~\ref{prop:repeated_height_original_seed_profile} at the same
$Q$ gives normalized amplitude
$C \cdot \mathfrak t(Q)^{1/3}$. A fixed sufficiently large $A_0$ makes this
at most $A \cdot W_1(Q)$ for every $\delta,m$; the initial depth required
by that proposition is
below Eq.~\eqref{eq:nearunit_secondary_primary_depth}.
Finite outer induction proves Eq.~\eqref{eq:nearunit_secondary_primary_profile}.
The choices of $C_{R,0}$ and $D_{\mathrm{add}}$ only increase the permitted initial
depth and amplitude. They never change the universal constants for repeated height reduction.
\end{proof}

Choosing the number of outer rounds gives a profile at every
prescribed depth in the same range. The height in
$\widehat V_{m,\delta}$ still depends on the loss parameter. We will
control this dependence by using a finite list of inherited bounds
simultaneously.

\begin{proposition}[A profile recurrence with a small multiplicative loss]
\label{prop:nearunit_secondary_prescribed_profile}
There are universal constants $A_0,C\ge1$ such that, for every
$0<\delta\le1/16$, every integer $m\ge0$, and every integer
$r\ge R_\delta$, each centered regular law of covariance at most $I$
and curvature at least $aI$ satisfies
\begin{equation}\label{eq:nearunit_secondary_prescribed_profile}
 \mathcal A\le A_0 \cdot e^{C \cdot \delta m}
        [\widehat V_{m,\delta}(r)+S_m]^{1/3} \cdot \mathcal L_r(a^{-1})^2.
\end{equation}
The constants and the domain $r\ge R_\delta$ are independent of $m$.
\end{proposition}
\begin{proof}
The $m=0$ profile follows from the initial prescribed-depth profile in
Proposition~\ref{prop:secondary_feedback_prescribed_profile}, increasing
$A_0$ once. That initial amplitude is universal, and
$\widehat V_{0,\delta}\ge\overline{\mathfrak t}\ge V_0$ and $S_m>0$ only weaken the bound. This also covers
$\mathcal A=1$ when $C_P<1$.

Fix a permitted integer $r$ and suppose the assertion holds at $m$
with amplitude $A$. Choose the finite parameters
\begin{equation}\label{eq:nearunit_secondary_return_parameters}
 j:=\kappa(r),\qquad \epsilon:=\frac{\delta}{j+1},\qquad
 Q:=\text{the least odd integer at least }
       \max\{Q_\epsilon,\delta^{-1} \cdot \log(r+1)\}.
\end{equation}
A universal choice of $C_{R,0}$ ensures $j\ge1$, $Q\le r$, and
$r\ge r_{Q,\epsilon,\delta}$ for every $r\ge R_\delta$.
Indeed $\kappa(r)\le\mathfrak t(r)\le C \cdot r^{1/8}$ gives
$Q_\epsilon\le C \cdot \delta^{-4} \cdot r^{1/2}\le r/2$ once
$r\ge C \cdot \delta^{-8}$, with room for rounding.
The other order term is at most $r/2$ at a smaller polynomial threshold.
Also
\begin{equation*}
 r_{Q,\epsilon,\delta}
       \le\max\{R_\delta,C \cdot \delta^{-1} \cdot \mathfrak t(r)\}\le r.
\end{equation*}
The exponent twelve in
Eq.~\eqref{eq:nearunit_secondary_profile_parameters} covers all these
conditions. None depends on $m$ or its amplitude.

Apply Proposition~\ref{prop:nearunit_secondary_restricted_profile}
with these parameters. Then $\mathfrak t_j(Q)\le\mathfrak t_j(r)=3<D$, so the
maximum defining the outer profile equals $D^{1/3}$ exactly. Also
$\epsilon j\le\delta$, and $(r+1)^{1/Q}\le e^\delta$.
Since $\epsilon^{-1}=(j+1)/\delta$ and $j+1\le\chi(r)$,
Lemma~\ref{lem:smoothed_secondary_costs} gives
\begin{align*}
 \widehat V_{m,\delta}(\epsilon^{-1})
 &\le\max\{\widetilde V_m(\chi(r)),\widetilde V_m(\delta^{-1}),
                   \overline{\mathfrak t}(\delta^{-1})\}
                                  +D_{\mathrm{add}} \cdot 4^{-m}\\
 &\le\widehat V_{m+1,\delta}(r)+D_{\mathrm{add}} \cdot 4^{-m}.
\end{align*}
The last inequality uses $\overline\chi\ge\chi$ and
$\widetilde V_m\circ\overline\chi=\widetilde V_{m+1}$.
Including the increment already allowed in $D$ gives the following
exact bound for the next height profile:
\begin{equation*}
 D\le\widehat V_{m+1,\delta}(r)+S_m+2D_{\mathrm{add}} \cdot 4^{-m}
   =\widehat V_{m+1,\delta}(r)+S_{m+1}.
\end{equation*}
There is no terminal additive three and no further multiplicative loss.
Together with $e^{16\delta+24\epsilon j}$ and the order factor,
this gives the next prescribed-depth profile with amplitude at most
$e^{42\delta} \cdot A$. Fix a universal $C\ge42$.
Finite induction proves Eq.~\eqref{eq:nearunit_secondary_prescribed_profile}
on the same domain $r\ge R_\delta$. Every parameter used in each
application is finite and is fixed before invoking coefficient transfer.
\end{proof}

Retaining a finite list of coefficient bounds allows different bounds to control disjoint ranges of degrees. This gives uniform errors across repetitions and closes the induction with bounded total cost.

\subsection{Retaining coefficient bounds across a finite chain}
\label{sec:finite_chain}
\label{sec:constant_chain_geometry}

The prescribed-depth profile carries a parameter-dependent height floor. To
control that cost as the losses decrease, we keep every coefficient
bound already obtained. We assign these caps to consecutive degree ranges
and prove a block estimate whose constants are independent of the
finite list's length.

We use the common coefficient radius $\mathcal A$, the operator
$\mathcal T$, and the regular class used above. Let $N\ge1$ be an integer,
and let $\mathcal H_0,\ldots,\mathcal H_{N-1}$ be nondecreasing
coefficient majorants on this class, with
\begin{equation}\label{eq:constant_chain_caps}
 1\le\mathcal H_i\le G_*,\qquad
 c_k\le\mathcal H_i(k)^{k-1},\qquad
 \mathcal G\le\min_{0\le i<N}\mathcal H_i.
\end{equation}
Here $\mathcal G$ satisfies Eq.~\eqref{eq:second_logstar_majorant}.
The caps need not be ordered; the working majorant may be replaced by
its minimum with them. Fix a finite ordered chain and write
\begin{equation}\label{eq:constant_chain_parameters}
 \begin{gathered}
 0<\epsilon\le\alpha_{N-1}\le\cdots\le\alpha_0\le1/16,
 \qquad\alpha_N:=\epsilon,\\
 K_i:=\lceil C_{\mathrm{cut}} \cdot \alpha_i^{-2}\rceil,\ \forall 0\le i\le N,\qquad
 \Xi_i:=\lceil4C_{\mathrm{deg}} \cdot \alpha_i^{-1} \cdot K_{i+1}\rceil
                                      ,\ \forall 0\le i<N,\\
 \mathsf A_{\alpha_0}:=\max\{C_0 \cdot G_*(16K_0)^2,C_0\}.
 \end{gathered}
\end{equation}
The constants will be universal. The $\Xi_i$ are nondecreasing and
$\Xi_i\ge K_{i+1}$. Equal adjacent margins are allowed.
For odd $p\ge3$, use the degrees and exits
\begin{equation}\label{eq:constant_chain_degrees}
 \begin{gathered}
 d_p:=\min\{2^j:2^j\ge4p\},\quad
 e_p:=\min\{2^j:2^j\ge8p\},\quad\omega_p:=1,\ \forall p\le K_0,\\
 d_p:=e_p:=\min\{2^j:2^j\ge C_{\mathrm{deg}} \cdot p/\alpha_i\},\quad
 \omega_p:=\alpha_i,\ \forall K_i<p\le K_{i+1},\ 0\le i<N,\\
 d_p:=e_p:=\min\{2^j:2^j\ge C_{\mathrm{deg}} \cdot p/\epsilon\},\quad
 \omega_p:=\epsilon,\ \forall p>K_N,\qquad J_p:=4e_p.
 \end{gathered}
\end{equation}
Empty order intervals are discarded. In band $i<N$, $e_p<\Xi_i$.
Every nonsmall band has $\omega_p^{-1}<\sqrt p$, so
$d_p,e_p,J_p\le C \cdot p^{3/2}$, and always $J_p\ge32p$.

The original cap controls the initial degree range, and each inherited
bound supplies a margin on its assigned interval. Summing the errors over
these disjoint intervals gives one convergent degree sum, which keeps the
block estimates uniform in the number of retained bounds.

\begin{proposition}[Blocks with a finite chain of floors]
\label{prop:constant_chain_blocks}
There are universal constants $C_X,C_0,C>0$ with the following property.
For every finite chain above and odd $Q\ge3$, we define
\begin{equation}\label{eq:constant_chain_block_parameters}
 \begin{gathered}
 X:=C_X \cdot Q^2,\qquad
 F^2:=\max\{\mathcal G(X^X)^2,\mathsf A_{\alpha_0},
                 (1+\alpha_i)\mathcal H_i(\Xi_i)^2:0\le i<N\},\\
 H_Q:=(1+2\epsilon) \cdot F^2.
 \end{gathered}
\end{equation}
There is an assignment $Z_Q$ and constants $a_Q,s_Q>0$ satisfying
\begin{align}
 H_Q&\le Z_Q\le\max\{H_Q,R\},\qquad\mathcal A\le Z_Q,\notag\\
 c_d&\le Z_Q^{(d-1)/2},\ \forall d\ge1,\qquad
 a_Q+s_Q+1\le C \cdot (F^2/4)^{Q-3}\le F^{2Q}.
                         \label{eq:constant_chain_static}
\end{align}
If $Z_Q>H_Q$, let $z:=Z_Q$ and $u:=1/z$. One centered unit finite
starting family satisfies
\begin{equation}\label{eq:constant_chain_startup}
 \nu\le\frac{u}{1-a_Q \cdot u^Q},\qquad
 P_J\le s_Q \cdot u^Q,\ \forall 1\le J\le J_Q,\qquad
 aV_M+X_M\le\nu-u+uP_{M+1},\ \forall M\ge0,
\end{equation}
and
\begin{equation}\label{eq:constant_chain_propagation}
 \|\mathcal T^h\|_{\mathrm{op}}\le4(h+1)z^{h/2},\ \forall h\ge0,\qquad
 0\le R-z\le\min\{1,C/R\},\qquad a_Q \cdot u^Q\le1/8.
\end{equation}
Every constant is independent of $N$, the margins, the caps, and $Q$.
For $N=1$ this is the construction with one inherited cap. When
$\mathcal A$ does not exceed the floor, the bound concerns only this
radius; no actual family or comparison of $R$ with
$H_Q$ is asserted there.
\end{proposition}
\begin{proof}
We check exactly the coefficient estimates needed in the proof of
Proposition~\ref{prop:nearunit_secondary_blocks}. Set
\begin{equation*}
 \rho_k:=\min\{1,G_*(k)^2/F^2,
                      \mathcal H_i(k)^2/F^2:0\le i<N\}.
\end{equation*}
Then $\rho_k\le C_0^{-1}$ for $k\le16K_0$, and
$\rho_k\le(1+\alpha_i)^{-1}$ for $k\le\Xi_i$.
For each of the finitely many fixed $M,r,c$ needed below,
\begin{equation}\label{eq:constant_chain_moments}
 \sum_{k\ge r+1}k^M \cdot \rho_k^{k-r} \cdot e^{-c\epsilon(k-r)}\le C
\end{equation}
with no dependence on $N$. To see this, use the geometric sum through
$K_0$, and partition the remaining integers into
$(K_0,\Xi_0]$, $(\Xi_{i-1},\Xi_i]$ for $1\le i<N$, and
$(\Xi_{N-1},\infty)$. A finite interval has damping
$e^{-\alpha_i(k-r)/2}$ and begins beyond $K_i$; the final interval
has the $\epsilon$ damping and begins beyond $K_N$.
Thus $\alpha_i k\ge\sqrt{C_{\mathrm{cut}} \cdot k}$ in the corresponding
band. Taking $K_0\ge2r$, the whole disjoint tail is bounded by
\begin{equation*}
 \sum_{k>K_0}k^M \cdot e^{-c'\sqrt{C_{\mathrm{cut}} \cdot k}}.
\end{equation*}
We have one sum over degrees, not a sum of constants over the bands.
It and the low geometric sum become uniformly small as the fixed
thresholds increase.

The raw kernel and early estimates in the proof of
Lemma~\ref{lem:nearunit_secondary_green} now hold without change.
Indeed, for every used $k\le d$, $\mathcal G(k)\le F$ and
$\mathcal G(k)^2/z\le\rho_k \cdot e^{-\epsilon/2}$ when
$z\ge(1+\epsilon) \cdot F^2$. The squared propagation factor remains
$16(3k)^2$ and the delays remain $k,\ldots,4k-2$, so the $k^8,k^9$
moments give $\eta,\eta_1\le C/z$. Before flattening, $b_j\le16(j+1)^2$;
the weighted early term is bounded by
$C \cdot k^6 \cdot \rho_k^{k-2} \cdot e^{-\epsilon(k-2)/2}$, using only
$\mathcal G(k)^2\le G_*(k)^2\le Ck$ in the prefactor.
Degree two still uses $c_2^2\le2$.
These are the only uses of the coefficient floor in that Green proof.
Its remaining hypotheses are the defect inequality, norm-attaining
shift, $R-z\le2$, and small $\gamma m^2/z$.
The boundary of the second block is still $b_{2m}\le1$ and
$z(b_{2m+1}-b_{2m})\le2$. The subsequent restart uses only those
bounds, $m\ge2d$, and the retained coherent, length, and prefix
conditions. Thus it applies to the same direct sum and compensated
$F_2$ family, with the same denominator windows and matched budget.

For extension at odd $q$, we use
$\kappa_q:=(F^2/4)^{q-1}$ and $\kappa_{q-2}:=(F^2/4)^{q-3}$.
The first-exit conditions
$a_q \cdot u^q,\mathsf K_q \cdot u^q\le\omega_q/1024$ give
$B_*\le\beta u$, where $\beta:=1+\omega_q/128$.
In the small band every tested degree is below $16K_0$, so the fixed
ratio absorbs $\beta$. In band $i<N$, every tested degree is below
$\Xi_i$; an earlier coefficient band has margin
$\alpha_l\ge\alpha_i=\omega_q$. The factor $\beta^4$ consumes
only a fixed small fraction of its damping. In the last band use
$\omega_q=\epsilon$ and the $z/F^2$ margin after the last cutoff.
The same disjoint bound gives $\theta\le C_\theta \cdot u$ universally.
For $d:=e_q$, $g:=\mathcal G(d)^2$, $\rho:=g/F^2$, and $D:=1+\epsilon$,
\begin{equation}\label{eq:constant_chain_terminal_ratio}
 \frac{\tau}{\kappa_q \cdot \kappa_{q-2} \cdot u^{2q}}
 \le C_F \cdot d \cdot g^3 \cdot 4^{2q-4} \cdot \rho^{d-4} \cdot \beta^d \cdot D^{2q-d}.
\end{equation}
This is the exact cancellation of $F^6 \cdot \rho^{d-1}$ used in
Eq.~\eqref{eq:nearunit_secondary_terminal_ratio}.
In the small range its bound is
$C \cdot q^4 \cdot C_0^4 \cdot (16\beta^{16}/C_0^8)^q$.
In band $i<N$, use $\rho\le(1+\alpha_i)^{-1}$,
$d\ge C_{\mathrm{deg}} \cdot q/\alpha_i$, and $\alpha_i^{-1}<\sqrt q$.
With $g\le Cd$, the bound is $C \cdot q^6 \cdot e^{-cC_{\mathrm{deg}}q}$.
The last band has the same bound from the $D$ damping.
There is no summation over bands in these orderwise estimates.
Thus $\tau\le C_F \cdot \kappa_q \cdot \kappa_{q-2} \cdot u^{2q}$.

Apply the first-exit calculation of
Lemma~\ref{lem:nearunit_secondary_extension}, with
$\mathsf K_q:=4(s_q+C_F \cdot \kappa_{q-2}+4C_\theta \cdot a_q+1)$.
It uses just the kernel and terminal bounds above, the actual prefix
$J_q=4e_q$, and the same normalizers and matched budget. It excludes
exit on $1\le m\le\lfloor z^q/\kappa_q\rfloor$ and supplies the
bridge loss $C_q \cdot z^{1-q}$ with
$C_q\le C \cdot (a_q+s_q+\kappa_{q-2}+1)$.

For a length $m$ within a factor four of $z^{p-2}/\kappa_{p-2}$,
the initializer's two coefficient conditions are bounded by
\begin{equation*}
 C \cdot d \cdot g \cdot \rho^{d-2} \cdot 4^{2p-6} \cdot D^{2p-d-4},\qquad
 C \cdot J_p \cdot d \cdot g^2 \cdot \rho^{d-3} \cdot 4^{p-3} \cdot D^{p-d},\qquad d:=d_p.
\end{equation*}
In the small band these are at most
$C \cdot p^2 \cdot C_0^2 \cdot (16/C_0^4)^p$ and $C \cdot p^4 \cdot C_0^3 \cdot (4/C_0^4)^p$.
In every other band its own margin gives respectively
$C \cdot p^3 \cdot e^{-cC_{\mathrm{deg}}p}$ and $C \cdot p^6 \cdot e^{-cC_{\mathrm{deg}}p}$.
Also $m\ge cF^2 \cdot 4^{p-3} \cdot D^{p-2}$ and
$J_p/m\le C \cdot p^{3/2} \cdot 4^{-p}/F^2$, so $m\ge2d_p$ and $J_p\le m/4$.
Once the propagation envelope has been established below, the initializer
therefore gives $a_p+s_p+1\le C \cdot \kappa_{p-2}$.
Its extension thresholds follow from
$(a_q+\mathsf K_q)u^q\le C \cdot F^{-6} \cdot 4^{3-q} \cdot D^{-q}$:
it is uniformly small in the first band, and multiplication by
$\omega_q^{-1}<\sqrt q$ handles every other band.
At the base $m_3:=\lfloor R\rfloor$, the degrees are still
$d_3=16,e_3=32,J_3=128$; the same fixed absolute floor and mesoscopic
estimate establish the base case; the lower radius is justified below.

We now verify the remaining hypotheses of
Proposition~\ref{prop:nearunit_secondary_blocks} in the order in which
its finite induction uses them. The uniform bound
$d_p,e_p,J_p\le C \cdot p^{3/2}$ puts every tested degree below
$X=C_X \cdot Q^2$. In the branch $\mathcal A>H_Q$, use the same actual
blocks, with
\begin{equation*}
 \begin{gathered}
 m_3:=\lfloor R\rfloor,\qquad
 z_q:=\|\mathcal T^{m_q}\|_{\mathrm{op}}^{2/m_q},\\
 m_{q+2}:=\lfloor z_q^q/\kappa_q\rfloor,\qquad
 \Delta_q:=C_b \cdot \kappa_{q-2} \cdot z_q^{1-q}.
 \end{gathered}
\end{equation*}
The block maximum in Eq.~\eqref{eq:coefficient_radius_block_maximum}
uses $\mathcal G(X^X)^2<\mathcal A$ when $R\le X$, and the unchanged
original bound $G_*(R^R)^2\le R/2<\mathcal A$ when
$X<R<2\mathcal A$. When $R\ge2\mathcal A$, the mesoscopic estimate
establishes the base lower radius, and the old decrement tail together
with the current bridge gives
$z_{q+2}\ge R-C_*/R-2\Delta_3\ge\mathcal A$.
Only ratios from previously constructed orders enter that tail.
Thus each actual radius is at least $\mathcal A>H_Q$ before its new
ratio is evaluated. Since $q\Delta_q/z_q\le C \cdot F^{-6}$,
\begin{equation*}
 \frac{\Delta_{q+2}}{\Delta_q}
 =\frac{F^4}{16z_{q+2}^2}(z_q/z_{q+2})^{q-1}\le\frac2{16},\qquad
 0\le R-z_q\le C/R\le1.
\end{equation*}
For an upward step the second factor is at most one; for a downward
step use $z_{q+2}\ge z_q-\Delta_q$. These are the same gap hypotheses
as in Eqs.~\eqref{eq:nearunit_secondary_gap_ratio} and
\eqref{eq:nearunit_secondary_gap}.

The real lengths are at least $F^2 \cdot 4^{q-1}(1+\epsilon)^q$.
The relative gap above, including the mesoscopic base gap, controls
rounding and gives $m_{q+2}/m_q\ge4$. The mixed-radix estimate has
per-level logarithmic excess
$C \cdot m_{q+2}\Delta_q/z_q\le16C/F^4$ and base cost $C/R$.
Counting levels gives the propagation envelope $4(h+1)$; counting only
the constructed prefix gives $\log E_q\le C/R+Cq/F^4$.
At the proposed next length this implies
\begin{equation*}
 q(\log(z_{q+2}/z_q))_+
 \le\frac{Cq \cdot 4^{-q}}{F^2 \cdot R}
      +\frac{C \cdot q^2 \cdot 4^{-q}}{F^6}
 \le\frac C{F^2 \cdot R}+\frac C{F^6}.
\end{equation*}
Together with the downward gap and rounding, this compares $m_{q+2}$
with $z_{q+2}^q/\kappa_q$ within a factor of four. At this point the
new radius and its propagation envelope are known. The coherent bounds
checked above, $m\ge2d_p$, $J_p\le m/4$, and the absolute gap therefore
supply all the hypotheses of Lemma~\ref{lem:nearunit_secondary_green}
and its restart. They give the same actual prefix and matched budget
on one starting family. This closes the finite induction without using
the new Green estimate to prove its own propagation hypothesis.
Every estimate uses a fixed degree moment, a geometric decrement sum,
or a maximum over the caps; none adds a cost for the chain length.

Set $Z_Q:=z_Q$ in this branch; otherwise set $Z_Q:=H_Q$.
The coefficient bounds follow from $\mathcal A\le Z_Q$.
The startup estimate includes $Q=3$, and its last constant is absorbed
in the fixed absolute floor to give $a_Q+s_Q+1\le F^{2Q}$.
Choose the moment and degree constants first, the initializer and
bridge constants next, and $C_0,C_X$ last. Every estimate above is
uniform in the finite chain length, so these choices are universal.
\end{proof}

\subsection{Closing the profile induction with bounded total cost}
\label{sec:constant_chain_profiles}

We now construct curvature profiles whose amplitude and additive reserve
stay bounded independently of the number of refinements. The finite-chain
block estimate retains all earlier coefficient bounds, while a geometric
choice of losses makes both the logarithms of the multiplicative costs
and the additive increments summable. The argument estimate below
controls the growing depth thresholds before the profile induction is
carried out.

We use the common coefficient radius $\mathcal A$, normalized logarithms
$\mathcal L_r$, and averaged maps from
Eq.~\eqref{eq:smoothed_secondary_heights}. Write
\begin{equation}\label{eq:constant_chain_carrier}
 W_i(x):=\overline{\mathfrak t}(\overline\chi^{\circ i}(x)),
 \ \forall i\ge0,\ x\ge1.
\end{equation}
Thus $W_i\ge4$, $W_i$ is $4^{-i-1}$-Lipschitz, and the fixed-power,
sum, and product argument costs in
Lemma~\ref{lem:smoothed_secondary_costs} are $O(4^{-i})$.
We keep all earlier coefficient caps during each finite construction.

We use Proposition~\ref{prop:constant_chain_blocks}. Its constants are
independent of the finite number of coefficient caps. The maximum in
Eq.~\eqref{eq:constant_chain_block_parameters}, the initialization of a common family in Eq.~\eqref{eq:constant_chain_startup}, and the actual
propagation bounds in Eq.~\eqref{eq:constant_chain_propagation}
will be used without changing the caps during an outer induction.

Fix $\alpha_0:=1/16$ and define
\begin{equation}\label{eq:constant_chain_profile_parameters}
 \begin{gathered}
 \alpha_i:=\alpha_0 2^{-i},\qquad
 R_i:=\lceil C_{R,0} \cdot \alpha_i^{-12}\rceil,\\
 S_0:=1,\qquad S_{i+1}:=S_i+B \cdot 2^{-i},\qquad
 A_{i+1}:=e^{C_A \cdot \alpha_i} \cdot A_i.
 \end{gathered}
\end{equation}
The fixed constant $A_0$ will be chosen below. The constants in the
proof of the inner depth induction and Proposition~\ref{prop:constant_chain_blocks}
are fixed first.
We then choose $C_{R,0}$, $C_A$, $B$, and $A_0$, in that order.
All choices will be independent of the finite chain length.

The initial depth $R_i$ grows with the chain index, so its effect on the height function must also be controlled. The following lemma first uses the
fixed-power estimates to reduce the polynomial argument and then applies
contraction; applying the Lipschitz bound directly at $R_i$ would not give
the required decreasing error.

\begin{lemma}[Polynomial arguments along the chain]
\label{lem:constant_chain_arguments}
For each fixed $C,p\ge1$, there is $b_{C,p}$ independent of $i$ such that
\begin{equation}\label{eq:constant_chain_polynomial_argument}
 W_i(C \cdot \alpha_i^{-p})\le4+b_{C,p} \cdot 2^{-i}.
\end{equation}
For a fixed sufficiently large $C$, there is also a universal $b$ with
\begin{equation}\label{eq:constant_chain_endpoint_argument}
 W_i(R_i+C \cdot \mathfrak t(x))\le W_i(x)+b \cdot 2^{-i}
 ,\ \forall x\ge1\text{ real},\ i\ge0.
\end{equation}
\end{lemma}
\begin{proof}
Apply the fixed-power estimate before using the Lipschitz bound:
\begin{equation*}
 W_i(C \cdot \alpha_i^{-p})
 \le W_i(2^i)+b_{C\alpha_0^{-p},p} \cdot 4^{-i}
 \le4+4^{-i-1}(2^i-4)_++b_{C\alpha_0^{-p},p} \cdot 4^{-i}.
\end{equation*}
This proves Eq.~\eqref{eq:constant_chain_polynomial_argument}, including
$i=0$ after increasing its fixed constant. The estimate for a fixed power of the argument is what makes the
subsequent contraction give a decreasing error at $R_i$.
For the second assertion, use the sum estimate of
Lemma~\ref{lem:smoothed_secondary_costs}, the first assertion with
$p:=12$, and $\mathfrak t(x)\le C' \cdot (x+1)$ for real $x\ge1$.
The fixed-power and sum estimates give
$W_i(C \cdot \mathfrak t(x))\le W_i(x)+b' \cdot 4^{-i}$.
Since $W_i(x)\ge4$, taking the maximum proves the claim.
\end{proof}

With the threshold costs controlled, the induction retains every bound
already proved and holds the resulting coefficient caps fixed during each
outer round. The multiplicative increments in $A_i$ have summable
logarithms, and the additive increments in $S_i$ are summable; the proof
must realize both prescribed increments with constants independent of the
number of retained caps.

\begin{proposition}[A bounded sequence of prescribed-depth profiles]
\label{prop:constant_chain_profiles}
The constants in Eq.~\eqref{eq:constant_chain_profile_parameters} can be chosen
universally so that every integer $i\ge0$, every integer $r\ge R_i$,
and every centered regular law of covariance at most $I$ and curvature
at least $aI$, where $a>0$, satisfy
\begin{equation}\label{eq:constant_chain_profile}
 \mathcal A\le A_i \cdot [W_i(r)+S_i]^{1/3} \cdot \mathcal L_r(a^{-1})^2.
\end{equation}
Moreover,
\begin{equation}\label{eq:constant_chain_bounded_constants}
 S_i\le1+2B,\qquad A_i\le A_0 \cdot e^{2C_A \cdot \alpha_0}.
\end{equation}
\end{proposition}
\begin{proof}
The $i=0$ case follows from Proposition~\ref{prop:secondary_feedback_prescribed_profile}
at its zeroth index, since $W_0=\overline{\mathfrak t}\ge\mathfrak t$.
Choose $C_{R,0}$ above its depth threshold and $A_0$ above its fixed
amplitude. The floor $\mathcal A\ge1$ is covered as well.
We prove Eq.~\eqref{eq:constant_chain_profile} by finite induction,
constructing one fixed coefficient cap immediately after each profile.

Suppose the profile has been proved through index $i$. For $0\le l\le i$,
set
\begin{equation}\label{eq:constant_chain_extraction_depth}
 r_l(x):=\max\{R_l,\lceil C_r \cdot \mathfrak t(x)+D_r/\alpha_l\rceil\}.
\end{equation}
The initial contraction and normalized transfer of
Lemma~\ref{lem:nearunit_static_transfer} give the valid nondecreasing cap
\begin{equation}\label{eq:constant_chain_fixed_caps}
 \mathcal H_l(x)^2:=\min\{G_*(x)^2,
 e^{3\alpha_l} \cdot A_l \cdot [W_l(r_l(x))+S_l]^{1/3}\}.
\end{equation}
Indeed $r_l(d)\ge R_l$ is in the previously proved profile domain,
$\mathcal L_{r_l(d)}(d)^2\le e^{\alpha_l}$, and
$(1+r_l(d)^{-2})^2\le e^{\alpha_l}$.
The normalized floor follows by a fixed choice of $A_0$.
The transfer passes finite-degree coefficient estimates on the whole
law class. It requires no continuity of the supremum defining
$\mathcal A$. Once constructed, $\mathcal H_l$ is fixed throughout all
later outer and inner depth inductions.

We first check every previously established coefficient floor. At the present stage, its
margin chain will be
\begin{equation}\label{eq:constant_chain_used_margins}
 \alpha_0,\alpha_1,\ldots,\alpha_i,\epsilon,
 \qquad0<\epsilon\le\alpha_i.
\end{equation}
In this application $N:=i+1$ and $K_{i+1}=\lceil C_{\mathrm{cut}} \cdot \epsilon^{-2}\rceil$.
For $l<i$, the endpoint in Eq.~\eqref{eq:constant_chain_parameters}
obeys $\Xi_l\le C \cdot \alpha_l^{-3}$ because $\alpha_{l+1}=\alpha_l/2$.
Consequently $r_l(\Xi_l)\le C' \cdot \alpha_l^{-12}$.
Lemma~\ref{lem:constant_chain_arguments} implies
\begin{equation*}
 W_l(r_l(\Xi_l))\le4+b_{\mathrm{old}} \cdot 2^{-l}.
\end{equation*}
Choose $B\ge b_{\mathrm{old}}$. Then for every $l<i$,
\begin{equation*}
 W_l(r_l(\Xi_l))+S_l\le4+S_i,
 \qquad A_i/A_l\ge e^{C_A \cdot \alpha_l}.
\end{equation*}
For $C_A\ge4$, Eq.~\eqref{eq:constant_chain_fixed_caps} gives
\begin{equation}\label{eq:constant_chain_old_floors}
 (1+\alpha_l)\mathcal H_l(\Xi_l)^2
       \le A_i(4+S_i)^{1/3},\ \forall 0\le l<i.
\end{equation}
This controls each member of the maximum without adding their costs.
It includes the possibly large but fixed $l=0$ argument cost.
The absolute floor from the original cap
$C_0 \cdot G_*(16K_0)^2$ is also fixed. Increase $A_0$ once so that
it and $C_0$ are below $A_0(4+S_0)^{1/3}$.

We now establish a complete family for repeated height reduction at index $i$. For
$0<\epsilon\le\alpha_i$, we define
\begin{align}
 D&:=W_i(\epsilon^{-1})+S_i+b_{\mathrm{pre}} \cdot 4^{-i},\notag\\
 P_j(x)&:=\max\{\mathfrak t_j(x),D\}^{1/3},
                  \label{eq:constant_chain_primary_weight}\\
 Q_\epsilon&:=\text{the least odd integer at least }C_Q \cdot \epsilon^{-4},\notag\\
 r_{Q,\epsilon,i}&:=\max\{R_i,
       \lceil C_b \cdot \mathfrak t(Q)\rceil+\lceil C_R/\epsilon\rceil\}.
                  \label{eq:constant_chain_primary_depth}
\end{align}
The constant $b_{\mathrm{pre}}$ will cover finitely many fixed-power
argument costs below. We claim that every $j\ge1$, odd $Q\ge Q_\epsilon$,
and integer $r\ge r_{Q,\epsilon,i}$ satisfy
\begin{equation}\label{eq:constant_chain_primary_profile}
 \mathcal A\le e^{16\alpha_i+24\epsilon j} \cdot A_i \cdot P_j(Q) \cdot
                           (r+1)^{1/Q} \cdot \mathcal L_r(a^{-1})^2.
\end{equation}
The constants $C_Q,C_b,C_R,C_T$ are the previously fixed universal constants for the outer iteration, independent of $i$ and the number of caps.

For the outer induction, write
$A_j':=e^{16\alpha_i+24\epsilon j} \cdot A_i$ and
$M_Q:=A_j' \cdot P_{j+1}(Q)$. We extract coefficients at the chosen order
\begin{equation}\label{eq:constant_chain_selected_order}
 \begin{split}
 r_\epsilon(x)&:=\max\{R_i,
           \lceil C_r \cdot \mathfrak t(x)\rceil+\lceil D_R/\epsilon\rceil\},\\
 k_\epsilon(x)&:=\text{the least odd integer at least }
     \max\{Q_\epsilon,\epsilon^{-1} \cdot \log(r_\epsilon(x)+1)\}.
 \end{split}
\end{equation}
Since $R_i\le C \cdot \epsilon^{-12}$, the polynomial depth comparison of
Proposition~\ref{prop:nearunit_secondary_restricted_profile} applies
with the same constants. In particular
$r_\epsilon(x)\ge r_{k_\epsilon(x),\epsilon,i}$.
The normalized transfer gives a valid working squared majorant
\begin{equation}\label{eq:constant_chain_working_seed}
 \mathfrak g(x):=\min\{\mathcal H_0(x)^2,\ldots,\mathcal H_i(x)^2,
                         e^{5\epsilon} \cdot A_j' \cdot P_j(k_\epsilon(x))\}.
\end{equation}
Taking the additional minima preserves validity and monotonicity.
No coefficient cap is extracted from a dimensional conclusion.

Set $h:=\mathfrak t(Q)$, $Y:=C_k \cdot (Qs)^2$, and $y:=\log(s)/Q$.
At $h\ge C_h \cdot R_i \cdot \epsilon^{-5}$, the bounds for the chosen order in the proof
of Proposition~\ref{prop:nearunit_secondary_restricted_profile} apply
with $R_\delta:=R_i$. They give $k_\epsilon(X^X)\le h$, no height
increase for $y\le1$, and, for $y\ge1$,
\begin{equation*}
 P_j(k_\epsilon(Y))^3
 \le\max\{\mathfrak t_j(h),D\}+5y.
\end{equation*}
Since $D\ge5$, division of the cube root by $e^y$ gives at most
$P_{j+1}(Q)$, exactly as in that proof.

When this condition on the height fails, the same argument gives
\begin{equation*}
 r_i(Y)\le C[R_i \cdot \epsilon^{-5}+y]
          \le C'[\epsilon^{-17}+y].
\end{equation*}
Remove the fixed outside dilation by the fixed-power estimate of
Lemma~\ref{lem:smoothed_secondary_costs}, and only then use the global
$4^{-i-1}$ Lipschitz bound for the additive $y$.
For one universal $b_{\mathrm{pre}}$ this gives
\begin{equation}\label{eq:constant_chain_low_argument}
 W_i(r_i(Y))\le W_i(\epsilon^{-1})
                   +b_{\mathrm{pre}} \cdot 4^{-i}+4^{-i-1} \cdot y.
\end{equation}
The analogous inequality at $X^X$ has no $y$ term.
Because $D\ge5$,
\begin{equation*}
 \sup_{y\ge0}e^{-y} \cdot [D+4^{-i-1} \cdot y]^{1/3}\le D^{1/3}.
\end{equation*}
The fixed cap $\mathcal H_i$ therefore contributes at most
$e^{4\alpha_i} \cdot A_i \cdot D^{1/3}$, below the initial amplitude for the outer iteration.
At all heights we obtain
\begin{equation}\label{eq:constant_chain_seed_bounds}
 \mathfrak g(X^X)\le e^{5\epsilon} \cdot M_Q,\qquad
 \sup_{s\ge1}\frac{\mathfrak g(C_k \cdot (Qs)^2)}{s^{1/Q}}
                                      \le e^{5\epsilon} \cdot M_Q.
\end{equation}
Only $\mathcal H_i$ was needed for the low range; retaining the earlier
caps can only decrease the working majorant.

The initial depth $r=r_{Q,\epsilon,i}$ is also valid.
At high height, use the preceding outer round at the chosen smaller order $k\le h\le Q$, obtaining at most
$e^\epsilon \cdot M_Q \cdot \mathcal L_r(a^{-1})^2$.
At low height, $r\le C \cdot R_i \cdot \epsilon^{-5}\le C' \cdot \epsilon^{-17}$,
and the profile already proved at index $i$ gives directly
\begin{equation}\label{eq:constant_chain_initial_depth}
 \mathcal A\le A_i \cdot D^{1/3} \cdot \mathcal L_r(a^{-1})^2
                    \le M_Q \cdot \mathcal L_r(a^{-1})^2.
\end{equation}
Here $r\ge R_i$; no assertion below the range of the preceding estimate is used.

Apply Proposition~\ref{prop:constant_chain_blocks} with
Eq.~\eqref{eq:constant_chain_used_margins}, fixed caps
$\mathcal H_0,\ldots,\mathcal H_i$, and
$\mathcal G:=\sqrt{\mathfrak g}$.
Every previously established coefficient floor is bounded by
Eq.~\eqref{eq:constant_chain_old_floors}, hence by $A_i \cdot D^{1/3}$.
For the newest cap, $\Xi_i\le C \cdot \alpha_i^{-1} \cdot \epsilon^{-2}$ and
$r_i(\Xi_i)\le C' \cdot \epsilon^{-12}$.
The same fixed-power estimate and the additive allowance already included give
\begin{equation}\label{eq:constant_chain_newest_floor}
 (1+\alpha_i)\mathcal H_i(\Xi_i)^2
          \le e^{5\alpha_i} \cdot A_i \cdot D^{1/3}\le M_Q.
\end{equation}
The fixed absolute floor is already below $A_i \cdot D^{1/3}$.
Together with Eq.~\eqref{eq:constant_chain_seed_bounds}, the exact
maximum defining $F$ yields
\begin{equation}\label{eq:constant_chain_profile_floor}
 F^2\le e^{6\epsilon} \cdot M_Q,\qquad H_Q\le e^{9\epsilon} \cdot M_Q,
 \qquad(a_Q+s_Q+1)^{1/Q}\le e^{6\epsilon} \cdot M_Q.
\end{equation}
The $\alpha_i$ normalization is paid only in the initial factor of $A_j'$.
There is no factor for the number of earlier caps.

We check the analytic kernel separately from the Green moments in
Proposition~\ref{prop:constant_chain_blocks}. Use the parameters for the inner depth induction in
Eq.~\eqref{eq:slow_offset_profile_parameters}, denoting its local depth
parameter by $\xi:=1/(64Qs)$. Thus $s=r+1$,
$p=\xi/128=1/(8192Qs)$, and $b=(1-2p)^{-1}$.
At a contradictory actual radius $z:=Z_Q$, write $u:=1/z$ and
$L_Q:=e^{20\epsilon} \cdot M_Q$. The comparison above gives
\begin{equation*}
 u\le L_Q^{-1} \cdot s^{-1/Q},\qquad b^4\le e^\epsilon.
\end{equation*}
The contradiction implies $Z_Q>H_Q$, so the actual family is available.
For $k\le16K_0$, the original geometric floor and
$z>H_Q\ge F^2$ give
$b^4 \cdot u \cdot \mathfrak g(k)\le b^4/C_0\le1/4$.
The degree-one term uses $c_1\le1$ directly.

For the remaining degrees below the finite profile cutoff $k_0$,
partition into the disjoint intervals
\begin{equation*}
 (16K_0,\Xi_0],\qquad
 (\Xi_{l-1},\Xi_l],\ \forall 1\le l\le i,\qquad
 (\Xi_i,\infty),
\end{equation*}
intersecting each with $k\le k_0$ and discarding empty intervals.
On the interval ending at $\Xi_l$, the corresponding old or newest
floor is at most $M_Q$. Thus
\begin{equation*}
 b^4 \cdot u \cdot \mathfrak g(k)
 \le e^{-19\epsilon}(1+\alpha_l)^{-1}
 \le e^{-c\alpha_l}.
\end{equation*}
Its left endpoint is at least $K_l$.
Consequently $\alpha_l k\ge\sqrt{C_{\mathrm{cut}} \cdot k}$.
Beyond $\Xi_i\ge K_{i+1}$, use $k\le k_0\le C_k \cdot (Qs)^2$ and
the dynamic estimate in Eq.~\eqref{eq:constant_chain_seed_bounds}.
The factor $s^{1/Q}$ cancels against the inverse-radius bound, giving
the same $e^{-c\epsilon}$ ratio and the same stretched-exponential bound.
Therefore all these disjoint intervals together contribute at most
\begin{equation*}
 Cu\sum_{k>16K_0}k^8
               e^{-c'\sqrt{C_{\mathrm{cut}} \cdot k}}\le C'u.
\end{equation*}
We have bounded one sum over degrees, not a sum of constants indexed
by the cap chain. The high-$k_0$ profile tail is unchanged, so
\begin{equation}\label{eq:constant_chain_kernel}
 \theta_d\le C_\theta \cdot u
\end{equation}
with $C_\theta$ independent of $i$, $\epsilon$, and all cap amplitudes.

The coherent-root estimate in
Eq.~\eqref{eq:slow_offset_coherent_root} uses only the original cap and
monotonicity. Its cutoff is
$\lceil D_1 \cdot \epsilon^{-2} \cdot Q^2\rceil\le X^X$ for $Q\ge Q_\epsilon$.
As in Eq.~\eqref{eq:slow_offset_kernel_split}, write
$D_Q:=b^Q \cdot E_{Q,\epsilon}$ for the coherent contribution. Then
\begin{equation*}
 D_Q^{1/Q}\le e^{8\epsilon} \cdot M_Q,\qquad
 \Delta_d\le(s_Q+D_Q)u^Q+p/512.
\end{equation*}
In particular, Eq.~\eqref{eq:constant_chain_profile_floor} gives
\begin{equation*}
 a_Q+s_Q+D_Q+4C_\theta \cdot a_Q+1
 \le(2+4C_\theta)e^{8\epsilon Q} \cdot M_Q^Q.
\end{equation*}
Choose $C_Q$ once so that
$\log(2^{24}(2+4C_\theta)Q)\le12\epsilon Q$ whenever
$0<\epsilon\le1/16$ and $Q\ge C_Q \cdot \epsilon^{-4}$.
This is possible uniformly in $\epsilon$: $\log(CQ)/Q$ decreases
for large $Q$, and at $Q=C_Q \cdot \epsilon^{-4}$ we use
$\epsilon^3\log(\epsilon^{-1})\le1/(3e)$. Since $M_Q\ge A_0$,
a fixed choice $A_0\ge\max\{16,16C_\theta\}$ now gives
\begin{equation}\label{eq:constant_chain_amplitude_conditions}
 L_Q\ge\max\{16,16C_\theta\},\qquad
 L_Q^Q\ge2^{24} \cdot Q(a_Q+s_Q+D_Q+4C_\theta \cdot a_Q+1).
\end{equation}
Thus the amplitude conditions hold without a fixed multiplicative
increase of $L_Q$ at each stage.

We also record the resulting first-exit inequalities. Set
$D_0:=\max\{\nu-u,0\}$. The actual startup and gap estimates give
$D_0\le2a_Q \cdot u^{Q+1}$ and $\lambda\ge u/2$.
Using $u^Q\le L_Q^{-Q}/s$ in
Eq.~\eqref{eq:constant_chain_amplitude_conditions}, we obtain
\begin{equation*}
 \begin{gathered}
 a_Q \cdot u^Q\le p/2048,\qquad
 \nu\le1/4,\qquad D_0\le u,\qquad D_0/p\le u,\\
 \theta_d \cdot u/\lambda\le2C_\theta \cdot u\le1/8,
 \end{gathered}
\end{equation*}
and
\begin{equation*}
 \begin{aligned}
 \Delta_d+\theta_d \cdot D_0/\lambda
 &\le\frac{a_Q+s_Q+D_Q+4C_\theta \cdot a_Q+1}{L_Q^Q \cdot s}
       +\frac p{512}\\
 &\le\frac p{2048}+\frac p{512}<\frac p{64}.
 \end{aligned}
\end{equation*}
Before a first exit, every required normalizer is at most
$\nu/(1-p)\le u/(1-p)^2\le bu$. The actual prefix and matched
budget in Eq.~\eqref{eq:constant_chain_startup} therefore allow
Lemma~\ref{lem:matched_scale_comparison} to give
$a<12(u+D_0/p)\tau_d\le24C_F \cdot d(bu)^d \cdot c_d^2$.
The dyadic terminal degree in
Eq.~\eqref{eq:second_logstar_terminal_degree} uses only
$\log(24C_F \cdot d/a)\le\xi d$, so its constant $C_T$ has no
dependence on the caps, their number, or their amplitudes.
The terminal contradiction and depth induction of
Proposition~\ref{prop:slow_offset_outer_profiles} now apply.
Initialize by Eq.~\eqref{eq:constant_chain_initial_depth} or the
high-height estimate. The distortion product is at most $e^\epsilon$,
so the next full depth profile has amplitude at most
$e^{21\epsilon} \cdot A_j'$, covered by $e^{24\epsilon} \cdot A_j'$.

For the first outer round, Proposition~\ref{prop:repeated_height_original_seed_profile}
at the same $Q$ applies because $P_1(Q)\ge\mathfrak t(Q)^{1/3}$ and
$r_{Q,\epsilon,i}$ is above its initial depth. One fixed large $A_0$ covers its amplitude.
This proves Eq.~\eqref{eq:constant_chain_primary_profile} by finite
outer induction. The induction estimates $\mathcal A$ directly, without converting
to $C_P$ or to the restricted radius.

It remains to prove the profile at index $i+1$ for every prescribed depth.
Fix any integer $r\ge R_{i+1}$ and choose
\begin{equation}\label{eq:constant_chain_return}
 j:=\kappa(r),\qquad\epsilon:=\frac{\alpha_i}{j+1},\qquad
 Q:=\text{the least odd integer at least }
       \max\{Q_\epsilon,\alpha_i^{-1} \cdot \log(r+1)\}.
\end{equation}
The polynomial threshold calculation in the proof of
Proposition~\ref{prop:nearunit_secondary_prescribed_profile} gives
$j\ge1$, $Q\le r$, and $r\ge r_{Q,\epsilon,i}$ whenever
$r\ge C \cdot \alpha_i^{-12}$. A fixed universal $C_{R,0}$ makes this true
for every $r\ge R_{i+1}$.
Explicitly, $\kappa(r)\le\mathfrak t(r)\le C \cdot r^{1/8}$ gives
$Q_\epsilon\le C \cdot \alpha_i^{-4} \cdot r^{1/2}\le r/2$ beyond a polynomial
threshold; the other order term and rounding fit in the remaining half.
Also $C_b \cdot \mathfrak t(Q)+C_R/\epsilon\le C' \cdot \alpha_i^{-1} \cdot \mathfrak t(r)
\le r$, and $R_i\le R_{i+1}\le r$.
These estimates do not depend on the number or amplitudes of the caps.

Now $\mathfrak t_j(Q)\le3<D$, so $P_j(Q)=D^{1/3}$.
Moreover $\epsilon j\le\alpha_i$ and $(r+1)^{1/Q}\le e^{\alpha_i}$.
The product estimate, monotonicity, and $\overline\chi\ge\chi$ give
\begin{align*}
 W_i(\epsilon^{-1})
 &\le\max\{W_i(\chi(r)),W_i(\alpha_i^{-1})\}+b \cdot 4^{-i}\\
 &\le W_{i+1}(r)+b_{\mathrm{return}} \cdot 2^{-i}.
\end{align*}
The second line uses Lemma~\ref{lem:constant_chain_arguments} and
$W_{i+1}(r)\ge4$.
Choose $B\ge b_{\mathrm{return}}+b_{\mathrm{pre}}$ as well as the requirement for the previously established floors. Then
\begin{equation*}
 D\le W_{i+1}(r)+S_i+B \cdot 2^{-i}=W_{i+1}(r)+S_{i+1}.
\end{equation*}
The next amplitude is at most $e^{42\alpha_i} \cdot A_i$.
Taking $C_A\ge42$ proves Eq.~\eqref{eq:constant_chain_profile} at $i+1$.
This closes the finite induction and makes its new cap available for the
next stage. Summing the two geometric series proves
Eq.~\eqref{eq:constant_chain_bounded_constants}.

To make the dependence of the constants explicit, first fix the degree,
moment, Green, and bridge constants in
Proposition~\ref{prop:constant_chain_blocks}. These fix $C_\theta$.
Next fix $C_Q$ as in Eq.~\eqref{eq:constant_chain_amplitude_conditions},
and the transfer, terminal-degree, and distortion constants
$C_r,D_r,D_R,C_T,C_b,C_R$. The polynomial depth comparisons then fix
$C_{R,0}$. After that, take $C_A\ge42$ and
$B\ge\max\{b_{\mathrm{old}},b_{\mathrm{return}}+b_{\mathrm{pre}}\}$;
these argument constants depend only on the already fixed choices.
Finally take $A_0$ above the initial profile, normalized floor, and
amplitude thresholds. None of these choices depends on $i$, $j$,
$\epsilon$, $Q$, or the finite cap count.
The choices of $B$ and $A_0$ only increase offsets and amplitudes. They do not change the prescribed range of depths, the threshold for choosing the order,
the terminal degree constant, or the constants in the block estimates.
The posterior and regular approximation arguments use inequalities at
finitely many degrees, with all their parameters fixed.
\end{proof}

\Needspace{14\baselineskip}
\subsection{Obtaining the universal bound and the main theorem}
\label{sec:proof_main}

We finish by obtaining a bound for the common coefficient radius, and
hence the Poincar\'e constant, that is independent of the curvature and
the dimension. For each regular law, finite termination of the height
function selects a finite chain length; the threshold estimate and
bounded amplitudes and offsets then give the uniform bound. Regular
approximation extends it to every isotropic log-concave measure and
completes the main theorem.

\begin{proposition}[A universal Poincar\'e bound]
\label{prop:constant_chain_endpoint}
There is a universal constant $C$ such that every centered regular law
of covariance at most $I$ satisfies $\mathcal A\le C$ and $C_P\le C$.
Consequently every isotropic log-concave measure $\mu$ on $\mathbb R^n$
satisfies $C_P(\mu)\le C$ and $\psi_n\le C$.
\end{proposition}
\begin{proof}
For a fixed regular law with curvature at least $aI$, we define
\begin{equation*}
 x:=\max\{1,a^{-1}\},\qquad i:=N_\chi(x)+1,\qquad
 r:=R_i+\lceil C \cdot \mathfrak t(x)\rceil.
\end{equation*}
Choose $C$ to be a fixed sufficiently large integer. Every parameter
is finite, and this choice makes
$\mathcal L_r(a^{-1})\le\mathcal L_r(x)$ bounded universally.
Lemma~\ref{lem:constant_chain_arguments}, valid for every real $x\ge1$,
and Lemma~\ref{lem:smoothed_secondary_heights} give
\begin{equation*}
 W_i(r)\le W_i(x)+b \cdot 2^{-i}=4+b \cdot 2^{-i}.
\end{equation*}
Apply Eq.~\eqref{eq:constant_chain_profile} at this permitted depth and
use Eq.~\eqref{eq:constant_chain_bounded_constants}. It follows that
$\mathcal A\le C$ with a constant independent of $a$, the dimension,
and the law. Lemma~\ref{lem:coefficient_radius_comparison} then gives
$C_P\le2^{85} \cdot \mathcal A\le C$.

For an arbitrary isotropic log-concave measure, use the regular
isotropic approximations of Lemma~\ref{lem:regular_isotropic_approximation}.
Each approximant satisfies this same universal Poincar\'e inequality.
Its finite internal parameters may be chosen separately; no parameter
or coefficient-radius supremum is passed through a limit.
Lemma~\ref{lem:poincare_limit} passes the uniform inequality to the
limiting measure for every locally Lipschitz function of finite energy.
The Cheeger--Poincar\'e comparison gives the assertion for $\psi_n$.
Every cap list, number of outer rounds, degree, and depth used to establish the
regular-law bound is finite. No infinite operator construction is required.
\end{proof}

\Needspace{8\baselineskip}

The preceding proposition gives the universal bound for regular laws.
Regular approximation and the Cheeger--Poincar\'e comparison now yield
the main theorem for arbitrary isotropic log-concave measures.

\begin{proof}[Proof of Theorem~\ref{tg:thm:main}]
Proposition~\ref{prop:constant_chain_endpoint} gives a universal bound on the common coefficient radius for every centered regular measure of covariance at most $I$. Lemma~\ref{lem:coefficient_radius_comparison} converts this into a universal Poincar\'e bound.

Given an arbitrary isotropic log-concave measure $\mu$, choose the regular isotropic approximations from Lemma~\ref{lem:regular_isotropic_approximation}. The preceding bound applies to every approximant with the same constant. Lemma~\ref{lem:poincare_limit} passes the inequality to $\mu$ and gives the stated $L^2$ conclusion for locally Lipschitz functions of finite Dirichlet energy. Finally, Lemma~\ref{lem:cheeger_poincare_comparison} yields $\psi_\mu^2\le\pi \cdot C_P(\mu)$. Taking the supremum over isotropic log-concave measures and enlarging the universal constant proves the assertion for $\psi_n$.
\end{proof}